\documentclass[twoside]{book}
\usepackage{geometry}
\usepackage{lectures}
\hypersetup{
breaklinks=true,
pdftitle={Euclidean plane and its relatives; a minimalist introduction.}
}

\newcommand{\arxiv}[2]{#1} 
\arxiv{}{\pdfminorversion=3}

\newcommand{\spell}[2]{#1} 

\makeindex

\begin{document}
\spell{}{\pagestyle{empty}\renewcommand\includegraphics[2][{}]{}\def\emph{\textit}\renewcommand\footnote[1]{\ (#1)}\renewcommand\z{}\renewcommand\texorpdfstring[2]{}\renewcommand\section[1]{SECTION. {#1} SECTION.}}

\title{\tt Euclidean plane and its relatives}
\subtitle{\tt A minimalist introduction}
\author{\tt Anton Petrunin}
\date{}
\maketitle


\null\vfill\noindent{\includegraphics[scale=0.5]{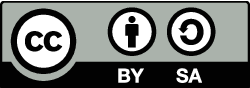}
\vspace*{1mm}
\\
\hbox{\parbox{.7\textwidth}
{This work is licensed under the Creative Commons Attribution-ShareAlike
4.0 International License.
To view a copy of this license, visit\\
\texttt{http://creativecommons.org/licenses/by-sa/4.0/}}}



\thispagestyle{empty}
\newpage
{\sloppy
\tableofcontents
}

\chapter*{Preface}
\addcontentsline{toc}{chapter}{Preface}

This book is meant to be 
rigorous, 
conservative, 
elementary, 
and minimalist.
At the same time, it includes about the maximum that students can absorb in one semester.

The present book is based 
on courses given by the author 
at Pennsylvania State University.
These lectures were oriented toward sophomore and senior university students who are familiar with real numbers and continuity.
It makes it possible to cover the material faster 
and  in a more rigorous way
than it could be done in high school.
Approximately one-third of the material used to be covered in high school, but not anymore.

\section{Prerequisite}

The students should be familiar 
with the following topics:
\begin{itemize}
\item Elementary set theory: 
$\in$,
$\cup$, 
$\cap$,
$\backslash$,
$\subset$,~$\times$.
\item Real numbers: intervals, inequalities, algebraic identities.
\item Limits, continuous functions, and the intermediate value theorem.
\item Standard functions: 
absolute value, 
natural logarithm,
exponential function. 
Occasionally, trigonometric functions  are used.
\item  Chapter~\ref{chap:trans} uses basic vector algebra.
\item To read Chapter~\ref{chap:sphere}, it is better to have some previous experience with the \textit{scalar product}, also known as the \textit{dot product}.
\item To read Chapter~\ref{chap:complex}, it is better to have some previous experience with complex numbers.
\end{itemize} 

\section{Overview}

We use the so-called \textit{metric approach} introduced by Birkhoff.
It means that we define the Euclidean plane as a \textit{metric space} that satisfies a list of properties (\textit{axioms}).
This way we minimize the tedious parts
which are unavoidable in the more classical Hilbert's approach.
At the same time, the students have a chance to learn the geometry of metric spaces.

Here is a dependency graph of the chapters.

\begin{figure}[!ht]
\centering

\begin{tikzpicture}[->,>=stealth',shorten >=1pt,auto,scale=1.4,
  thick,main node/.style={circle,draw,font=\sffamily\bfseries,minimum size=8mm}]

  \node[main node] (1) at (-1.5,10/6) {\ref{chap:metr}};
  \node[main node] (2) at (-.5,10/6){\ref{chap:axioms}};
  \node[main node] (3) at (.5,10/6) {\ref{chap:half-planes}};
  \node[main node] (4) at (1.5,10/6) {\ref{chap:cong}};
  \node[main node] (5) at (2.5,10/6) {\ref{chap:perp}};
  \node[main node] (6) at (3.5,10/6) {\ref{chap:parallel}};
  \node[main node] (7) at (4.5,10/6) {\ref{chap:angle-sum}};
  \node[main node] (8) at (5.5,10/6) {\ref{chap:triangle}};
  \node[main node] (9) at (5,5/6){\ref{chap:inscribed-angle}};
  \node[main node] (10) at (4.5,0) {\ref{chap:inversion}};
  \node[main node] (11) at (3,5/6) {\ref{chap:non-euclid}};
  \node[main node] (12) at (3.6,0){\ref{chap:poincare}};
  \node[main node] (13) at (3,-5/6) {\ref{chap:h-plane}};
  \node[main node] (14) at (5,-5/6) {\ref{chap:trans}};
  \node[main node] (15) at (4.5,-10/6) {\ref{chap:proj}};
  \node[main node] (16) at (4,-5/6) {\ref{chap:sphere}};
  \node[main node] (17) at (3.5,-10/6) {\ref{chap:klein}};
  \node[main node] (18) at (5.5,0) {\ref{chap:complex}};
  \node[main node] (19) at (4,15/6) {\ref{chap:car}};
  \node[main node] (20) at (5,15/6) {\ref{chap:area}};

  \path[every node/.style={font=\sffamily\small}]
   (1) edge node{}(2)
   (2) edge node{}(3)
   (3) edge node{}(4)
   (4) edge node{}(5)
   (5) edge node{}(6)
   (6) edge node{}(7)
   (7) edge node{}(9)
   (7) edge node{}(8)
   (7) edge node{}(20)
   (10) edge node{}(14)
   (14) edge node{}(15)
   (7) edge node{}(19)
   (9) edge node{}(10)
   (10) edge node{}(12)
   (10) edge node{}(16)
   (10) edge node{}(18)
   (5) edge node{}(11)
   (11) edge node{}(12)
   (12) edge node{}(13)
   (15) edge node{}(17)
   (16) edge[dashed] node{}(17)
   (13) edge node{}(17);
\end{tikzpicture}
\end{figure}

In (\ref{chap:metr}) we introduce the concepts necessary to formulate the axioms;
it includes metric spaces, 
lines, 
angle measure, 
continuous maps,
and congruent triangles.

Next, we dive into Euclidean geometry:
(\ref{chap:axioms}) Axioms and immediate corollaries;
(\ref{chap:half-planes}) Half-planes and continuity;
(\ref{chap:cong}) Congruent triangles;
(\ref{chap:perp}) Circles, motions, and perpendicular lines;
(\ref{chap:parallel}) Similar triangles and (\ref{chap:angle-sum}) Parallel lines  
--- these are the first two chapters where we use Axiom~\ref{def:birkhoff-axioms:4}, an equivalent of Euclid's parallel postulate.
In (\ref{chap:triangle}) we give the most classical theorems of triangle geometry;
this chapter is included mainly as an illustration.

In the following two chapters, we discuss the geometry of circles on the Euclidean plane:
(\ref{chap:inscribed-angle}) Inscribed angles; (\ref{chap:inversion}) Inversion.
It  will be used to construct the model of the hyperbolic plane.

Next, we discuss non-Euclidean geometry:
(\ref{chap:non-euclid})
Neutral geometry --- geometry without the parallel postulate;
(\ref{chap:poincare})
Conformal disc model ---
this is a construction of the hyperbolic plane,
an example of a neutral plane that is not Euclidean.
In (\ref{chap:h-plane}) we discuss geometry of the constructed hyperbolic plane --- this is the highest point in the book.

In the remaining chapters, we discuss additional topics:
(\ref{chap:trans}) Affine geometry;
(\ref{chap:proj}) Projective geometry;
(\ref{chap:sphere}) Spherical geometry;
(\ref{chap:klein}) Projective model of the hyperbolic plane;
(\ref{chap:complex}) Complex coordinates;
(\ref{chap:car}) Geometric constructions;
(\ref{chap:area}) Area.
The proofs in these chapters are not completely rigorous.

We encourage the use of visual assignments on the author's website.

\section{Disclaimer}
 
Most of the proofs in the book already appeared in Euclid's Elements.
I did not try to find original references; it is an impossible task.

\section{Recommended resources}

\begin{itemize}
\item {}\emph{Byrne's Euclid} \cite{byrne} --- a colored version of the first six books of Euclid's Elements edited by Oliver Byrne. 

\item {}\emph{Kiselyov's geometry} \cite{kiselev} ---
a classical textbook for school students written by Andrey Kiselyov; it should help if you have trouble following this book.

\item {}\emph{Lessons in Geometry} by Jacques Hadamard \cite{hadamard} --- an encyclopedia of elementary geometry originally written for school teachers.



\item {}\emph{Problems in geometry} by Victor Prasolov\cite{prasolov}  is perfect for mastering your problem-solving skills.

\item {}\emph{Geometry in figures} by Arseniy Akopyan \cite{akopyan} --- an encyclopedia of Euclidean geometry with barely any words.

\item {}\emph{Euclidea} \cite{euclidea} --- a fun and challenging way to learn geometric constructions.
 
\item {}\emph{Geometry} by Igor Sharygin \cite{sharygin} --- the greatest textbook in geometry for school students, I recommend it to anyone who can read Russian.

\item {}\emph{problems.ru} \cite{zadachi} --- an on-line collection of problems in Russian.

\end{itemize}

\section{Acknowledgments}

{\sloppy

Let me thank 
Stephanie Alexander,
Thomas Barthelme,
Berk Ceylan,
Matthew Chao, 
Quinn Culver,
Svetlana Katok, 
Nina Lebedeva,
Alexander Lytchak,
Alexei Novikov,
and Lukeria Petrunina
for useful suggestions and for correcting misprints.

This work was partially supported by
NSF grants
DMS-0103957,
DMS-0406482,
DMS-0905138,
DMS-1309340,
DMS-2005279,
and Simons Foundation grants 
245094, 584781.

}

\chapter{Preliminaries}\label{chap:metr}

\section{What is the axiomatic approach?}
\label{preaxioms}

In the axiomatic approach, one defines the plane as anything that satisfies a specified set of properties known as {}\emph{axioms}.
These axioms essentially serve as the rules for a mathematical game.
Once the axiom system is fixed, a statement is true if it logically follows from the axioms, and nothing else is.

The formulations of the first axioms were not rigorous at all.
For example, Euclid described a {}\emph{line} as \textit{breadthless length}
and a {}\emph{straight line} as a line that \textit{lies evenly with the points on itself}.
However, these formulations were clear enuf for one mathematician to understand another.

The best way to understand an axiomatic system
is to make one by yourself.
For now, let's assume an intuitive understanding of concepts like lines and points.
Imagine an infinite and perfect surface of a chalkboard. 
Let us try to collect the key observations about this model.

\begin{enumerate}[(i)]
 \item\label{preaxiomI} We can measure distances between points.
 \item\label{preaxiomII} We can draw a unique line thru two given points.
 \item\label{preaxiomIII} We can measure angles.
 \item\label{preaxiomIV} If we rotate or shift we will not see the difference.
 \item\label{preaxiomV} If we change the scale we will not see the difference.
\end{enumerate}
These observations are good to start with.
Later, we will develop a language to reformulate them rigorously.

\section{What is a model?}
\label{page:model}

The Euclidean plane can be defined rigorously the following way:

\textit{Define a {}\emph{point} in the Euclidean plane as a pair of real numbers $(x,y)$ and define the {}\emph{distance} between two points $(x_1,y_1)$ and $(x_2,y_2)$ by the following formula:}
\[\sqrt{(x_1-x_2)^2+(y_1-y_2)^2}.\]

That is it!
We gave a {}\emph{numerical model} of the Euclidean plane;
it builds the Euclidean plane from real numbers
while the latter is assumed to be known.

Shortness is the main advantage of the model approach,
but it is not intuitively clear why we define points and distances this way.

On the other hand, the observations made in the previous section are intuitively obvious ---
this is the primary advantage of the axiomatic approach.
Another advantage --- the axiomatic approach is easily adjustable. 
For example, we can remove or replace an axiom from the list. 
We will do it in Chapter~\ref{chap:non-euclid}.

\section{Metric spaces}

The concept of a metric space provides 
a rigorous way to say: \textit{``we can measure distances between points''}.
That is, instead of (\ref{preaxiomI}) in Section~\ref{preaxioms},
we can say \textit{``Euclidean plane is a metric space''}.

\begin{thm}{Definition}\label{def:metric-space}
Let $\mathcal X$ be a nonempty set and 
let $d$ be a function
that assigns a real number $d(A,B)$
to each pair $A,B\in\mathcal X$.
Then $d$ is called a \index{metric}\emph{metric} on 
$\mathcal X$ if, for any
$A,B,C\in \mathcal X$, the following conditions are satisfied:
\begin{enumerate}[(a)]
\item\label{def:metric-space:a} Positiveness: 
$d(A,B)\ge 0.$
\item\label{def:metric-space:b} $A=B$ if and only if 
$d(A,B)=0.$
\item\label{def:metric-space:c} Symmetry: $d(A, B) = d(B, A).$
\item\label{def:metric-space:d} \index{triangle!inequality}Triangle inequality: 
$d(A, C) \le d(A, B) + d(B, C).$
\end{enumerate}
A \index{metric!space}\emph{metric space} is defined as a set equipped with a metric. 
More formally, a metric space is a pair $(\mathcal X, d)$ where $\mathcal X$ is a set and $d$ is a metric on~$\mathcal X$.

The elements of $\mathcal X$ are referred to as \index{point}\emph{points} of the metric space.
The value $d(A, B)$ is called the \index{distance}\emph{distance} from $A$ to~$B$.
\end{thm}

\pagebreak

\subsection*{Examples}

\begin{itemize}
\item {}\emph{Discrete metric.} Let $\mathcal X$ be an arbitrary set. 
For any $A,B\z\in\mathcal X$ set $d(A,B)\z=0$ if $A=B$ and $d(A,B)=1$ otherwise.
The metric $d$ is called the \index{discrete metric}\emph{discrete metric} on~$\mathcal X$.
\item\index{real!line}\emph{Real line.}
Set of all real numbers ($\mathbb{R}$) with the metric $d$ defined as 
$$d(A,B)\df|A-B|.$$
\end{itemize}

\begin{thm}{Exercise}\label{ex:dist-square}
Show that $d(A,B)=|A-B|^2$ is \textit{not} a metric on $\mathbb{R}$.
\end{thm}

\begin{itemize}
\item \textit{Metrics on the plane.}
Suppose that $\mathbb{R}^2$ denotes the set of all pairs $(x,y)$ of real numbers.
Assume $A=(x_A,y_A)$ and $B=(x_B,y_B)$.
Consider the following metrics on $\mathbb{R}^2$:
\begin{itemize}
\item\index{Euclidean!metric}\emph{Euclidean metric,} denoted by \index{59@$d_1$, $d_2$, $d_\infty$}$d_2$, and defined as \label{def:d_2}
$$d_2(A,B)=\sqrt{(x_A-x_B)^2+(y_A-y_B)^2}.$$
\item\label{Manhattan plane}\index{Manhattan plane}\emph{Manhattan metric,} denoted by $d_1$ and defined as 
$$d_1(A,B)=|x_A-x_B|+|y_A-y_B|.$$
\item{}\emph{Maximum metric,} denoted by $d_\infty$ and defined as 
$$d_\infty(A,B)=\max\{|x_A-x_B|,|y_A-y_B|\}.$$
\end{itemize}
\end{itemize}

\begin{thm}{Exercise}\label{ex:d_1+d_2+d_infty}
Prove that the following functions are metrics on $\mathbb{R}^2$:
(a)~$d_1$; (b)~$d_2$; (c)~$d_\infty$.
\end{thm}

\section{A shortcut for distance}

Most of the time, we need only one metric on a space.
Therefore, we will not need to name the metric each time.

Given a metric space $\mathcal X$,
the distance between points $A$ and $B$ will be further denoted by 
$$AB
\quad
\text{or}
\quad
d_{\mathcal X}(A,B);$$
the latter is used only if we need to emphasize that $A$ and $B$ are points of the metric space~$\mathcal X$.

For example, the triangle inequality can be written as 
$$AC\le AB+BC.$$

To avoid confusion, we will always use ``$\cdot$'' for multiplication,
so $AB$ is not mistaken for $A\cdot B$.

\begin{thm}{Exercise}\label{ex:4-triangle}
Show that the inequality
\[AB+PQ\le AP+AQ+BP+BQ\]
holds for any four points $A$, $B$, $P$, $Q$ in a metric space.
\end{thm}

\section{Isometries, motions, and lines}

In this section, we define lines in a metric space.
Once it is done the sentence \textit{``We can draw a unique line that passes thru two given points.''} becomes rigorous; see (\ref{preaxiomII}) in Section~\ref{preaxioms}. 

Recall that a map $f\:\mathcal{X}\to\mathcal{Y}$
is a \index{bijection}\emph{bijection}
if it gives an exact pairing of the elements of two sets.
Equivalently, $f\:\mathcal{X}\to\mathcal{Y}$ is a bijection if it has an \index{inverse}\emph{inverse};
that is, a map $g\:\mathcal{Y}\to\mathcal{X}$
such that 
$g(f(A))\z=A$ for any $A\in\mathcal{X}$
and
$f(g(B))\z=B$ for any $B\in\mathcal{Y}$. 

Let $\mathcal X$ and $\mathcal Y$ be two metric spaces and $d_{\mathcal X}$, $d_{\mathcal Y}$ be their metrics. 
A~map 
$$f\:\mathcal X \z\to \mathcal Y$$ 
is
called \index{distance-preserving map}\emph{distance-preserving} if 
$$d_{\mathcal Y}(f(A), f(B))
 = d_{\mathcal X}(A,B)$$
for all $A,B\in {\mathcal X}$.

A bijective distance-preserving map is called an \index{isometry}\emph{isometry}. 

Two metric spaces are called \index{isometric spaces}\emph{isometric} if there exists an isometry between them.

An isometry from a metric space to itself is also referred a \index{motion}\emph{motion}.

\begin{thm}{Exercise}\label{ex:dist-preserv=>injective}
Show that each distance-preserving map is \index{injective map}\emph{injective};
that is, if $f\:\mathcal X\to\mathcal Y$ is a distance-preserving map, 
then $f(A)\ne f(B)$ for every pair of distinct points $A, B\in \mathcal X$.
\end{thm}

\begin{thm}{Exercise}\label{ex:motion-of-R}
Show that if $f\:\mathbb{R}\to\mathbb{R}$ is a motion of the real line,
then either (a)
$f(x)=f(0)+x$ for any $x\in \mathbb{R}$, 
or (b)
$f(x)=f(0)-x$ for any $x\in \mathbb{R}$. 

\end{thm}

\begin{thm}{Exercise}\label{ex:d_1=d_infty}
Prove that $(\mathbb{R}^2,d_1)$ is isometric to $(\mathbb{R}^2,d_\infty)$.
\end{thm}

\begin{thm}{Advanced exercise}\label{ad-ex:motions of Manhattan plane}
Describe all the motions of the Manhattan plane (defined in \ref{ex:dist-square}).
\end{thm}

If $\mathcal X$ is a metric space and $\mathcal Y$ is a subset of $\mathcal X$,
then a metric on $\mathcal Y$ can be obtained by restricting the metric from~$\mathcal X$. 
In other words, 
the distance between two points in $\mathcal Y$ is defined to be the distance between these points in $\mathcal X$.
In this way any subset of a metric space can be also considered as a metric space.

\begin{thm}{Definition}\label{def:line}
A subset $\ell$ of a metric space is called a \index{line}\emph{line} if it is isometric to the real line.
\end{thm}

A triple of points that lie on one line is called \index{collinear points}\emph{collinear}.
Note that if $A$, $B$, and  $C$ are  collinear, $AC\ge AB$, and $AC\ge BC$, then $AC\z= AB+BC$.

Some metric spaces have no lines; for example, discrete metrics.
The picture shows examples of lines on the Manhattan plane $(\mathbb{R}^2,d_1)$. 
\begin{figure}[!ht]
\centering
\includegraphics{mppics/pic-2}
\end{figure}

\begin{thm}{Exercise}\label{ex:y=|x|}
Consider the graph $y=|x|$ in $\mathbb{R}^2$.
In which of the following spaces
(a) $(\mathbb{R}^2,d_1)$, 
(b) $(\mathbb{R}^2,d_2)$, 
(c) $(\mathbb{R}^2,d_\infty)$ 
does it form a line? 
Why?
\end{thm}

\begin{thm}{Exercise}\label{ex:line-motion}
Show that each motion maps a line to a line.
\end{thm}

\section{Half-lines and segments}

Assume there is a line $\ell$ passing thru
two distinct points $P$ and $Q$.
In this case, we might denote $\ell$ as $(PQ)$.
There might be more than one line thru $P$ and $Q$,
but if we write \index{60@$(PQ)$, $[PQ)$, $[PQ]$}$(PQ)$ we assume that we made a choice of such a line. 

We will denote by $[P Q)$ the \index{half-line}\emph{half-line}
that starts at $P$ and contains~$Q$. 
Formally speaking, $[P Q)$ is a subset of $(P Q)$ that corresponds to $[0,\infty)$ under an isometry $f\:(P Q)\to \mathbb{R}$ such that $f(P)=0$ and $f(Q)>0$.

The subset of line $(P Q)$ between $P$ and $Q$ is called the \index{segment}\emph{segment between} $P$ and $Q$; it is denoted by~$[P Q]$.
Formally, the segment can be defined as the intersection of two half-lines: $[P Q]=[P Q)\cap[Q P)$.

\begin{thm}{Exercise}\label{ex:trig==}
Show that given $r\ge 0$ and a half-line $[P Q)$ there is a unique $X\in [P Q)$  such that $P X=r$.

\end{thm}


\section{Angles}

Our next goal is to introduce angles and angle measures; 
after that, the statement \textit{``we can measure angles''} will become rigorous;
see (\ref{preaxiomIII}) in Section~\ref{preaxioms}.

An ordered pair of half-lines that start at the same point is called an \index{angle}\emph{angle}.
The angle $AOB$ (also denoted by \index{10@$\angle$, $\measuredangle$}$\angle AOB$) is the pair of half-lines $[OA)$ and $[OB)$; the point $O$ is called the \index{vertex!of angle}\emph{vertex} of the angle.

Intuitively, the angle measure tells how much one has to rotate the first half-line counterclockwise, so it gets the position of the second half-line of the angle. 
The full turn is assumed to be $2\cdot\pi$;
it corresponds to the angle measure in radians.%
\footnote{For now, you can think of $\pi$ as a the size of a half-turn in specific units. Its numerical value, $\pi\approx 3.14$, is not essential at this stage.}

The angle measure of $\angle AOB$ is denoted as  $\measuredangle AOB$;
it is a real number in the interval $(-\pi,\pi]$. 

\begin{wrapfigure}{o}{25mm}
\vskip-0mm
\centering
\includegraphics{mppics/pic-4}
\end{wrapfigure}

The notations $\angle AOB$ and $\measuredangle AOB$ look similar;
they also have close but different meanings which better not be confused.
For example, the equality 
$\angle AOB\z=\angle A'O'B'$
means that
$[OA)=[O'A')$ and $[OB)\z=[O'B')$;
in particular, $O=O'$.
On the other hand, the equality 
$\measuredangle AOB\z=\measuredangle A'O'B'$ 
means only equality of two real numbers;
in this case, $O$ may be distinct from~$O'$.

Here is the first property of angle measure which will become a part of the axiom.
\textit{Given a half-line $[O A)$ and $\alpha\in(-\pi,\pi]$ there is a unique half-line $[O B)$ such that $\measuredangle A O B= \alpha$.}

\section[\texorpdfstring{Reals modulo $2\cdot\pi$}{Reals modulo 2·π}]{Reals modulo $\bm{2\cdot\pi}$}

Consider three half-lines starting from the same point, $[O A)$, $[O B)$, and $[O C)$.
They make three angles $A O B$, $B O C$, and $A O C$,
so the value $\measuredangle A O C$ should coincide with
the sum $\measuredangle A O B+\measuredangle B O C$ up to full rotation.
This property will be expressed by the formula 
$$\measuredangle A O B+\measuredangle B O C\equiv \measuredangle A O C,$$
where \index{34@$\equiv$}``$\equiv$'' is a new notation that we are about to introduce.
The last identity will become a part of the axioms.

We will write $\alpha\equiv\beta\pmod{2\cdot\pi}$, or briefly
\begin{align*}
\alpha&\equiv\beta
\end{align*}
if $\alpha=\beta+2\cdot\pi\cdot n$
for an integer~$n$.
In this case, we say 
$$\textit{``$\alpha$ is equal to $\beta$ modulo $2\cdot\pi$''}.$$
For example, $-\pi
\equiv
\pi\equiv 3\cdot\pi$ and $\tfrac12\cdot\pi
\equiv
-\tfrac32\cdot\pi$.

The introduced relation ``$\equiv$'' behaves like an equality sign,
but
\[\dots\equiv\alpha-2\cdot\pi\equiv \alpha\equiv \alpha+2\cdot\pi\equiv \alpha+4\cdot\pi\equiv\dots;\] 
that is, if the angle measures differ by full turn,
then they are considered to be the same.

With ``$\equiv$'', we can do addition, subtraction, and multiplication with integer numbers without getting into trouble.
That is, if
$$\alpha\equiv\beta
\quad
\text{and}
\quad
\alpha'\equiv \beta',$$ 
then
$$\alpha+\alpha'\equiv\beta+\beta',
\quad
\alpha-\alpha'\equiv \beta-\beta'
\quad 
\text{and}
\quad
n\cdot\alpha\equiv n\cdot\beta$$
for every integer~$n$.
However, ``$\equiv$'' does not respect multiplication with non-integer numbers; for example, 
$$\pi
\equiv 
-\pi
\quad
\text{but}
\quad
\tfrac12\cdot\pi
\not\equiv
-\tfrac12\cdot\pi.$$ 

\begin{thm}{Exercise}\label{ex:2a=0}
Show that $2\cdot\alpha\equiv0$ if and only if $\alpha\equiv0$ or $\alpha\equiv\pi$.
\end{thm}

\begin{thm}{Exercise}\label{ex:a+b==c}
Suppose that $0\le \alpha\le \pi$, $0\le \beta\le \pi$, and $0<\gamma<2\cdot \pi$.
Show that
\[\alpha+\beta\equiv \gamma\quad\text{implies that}\quad\alpha+\beta=\gamma.\]
\end{thm}

\section{Continuity}

The angle measure is also assumed to be continuous.
Namely, the following property of angle measure will become a part of the axioms:

\textit{The function}
$$\measuredangle\:(A,O,B)\mapsto\measuredangle A O B$$
\textit{is continuous at every triple of points $(A,O,B)$
such that $O\ne A$ and $O\ne B$ and $\measuredangle A O B\ne\pi$.}

To clarify this property, we need to extend the concept of {}\emph{continuity} to functions between metric spaces.
The definition is a straightforward generalization of the standard definition for real-to-real functions.

Let $\mathcal X$ and $\mathcal Y$ be two metric spaces,
and $d_{\mathcal X}$, $d_{\mathcal Y}$ be their metrics.

A map $f\:\mathcal X\to\mathcal Y$ is called \index{continuous}\emph{continuous} at point $A\in \mathcal X$
if, for any $\epsilon>0$, there is $\delta>0$, such that 
\[d_{\mathcal X}(A,A')
<
\delta
\quad
\Rightarrow
\quad
d_{\mathcal Y}(f(A),f(A'))
<
\epsilon.\]
(Informally it means that sufficiently small changes of $A$ result in arbitrarily small changes of $f(A)$.)

A map $f\:\mathcal X\to\mathcal Y$ is called \index{continuous}\emph{continuous} if it is continuous at every point $A\in \mathcal X$.

One may define a continuous map of several variables in a similar manner.
Assume $f$ is a function that returns a point in $\mathcal Y$ for a triple of points $(A,B,C)$
in~$\mathcal X$.
The map $f$ might be defined only for some triples in~$\mathcal X$.

Assume $f(A,B,C)$ is defined.
Then, we say that $f$ is continuous at the triple $(A,B,C)$ 
if, for any $\epsilon>0$, there is $\delta>0$ such that 
\[d_{\mathcal Y}(f(A,B,C),f(A',B',C'))<\epsilon.\]
if $d_{\mathcal X}(A,A')<\delta$, $d_{\mathcal X}(B,B')<\delta$, and $d_{\mathcal X}(C,C')<\delta$.

\begin{thm}{Exercise}\label{ex:dist-cont}
Let $\mathcal{X}$ be a metric space.
\begin{enumerate}[(a)]
\item\label{ex:dist-cont:a} Let $A\in \mathcal{X}$ be a fixed point.
Show that the function 
$$f\:B\mapsto
d_{\mathcal{X}}(A,B)$$ 
is continuous at every point~$B$.
\item Show that the function $g\:(A,B)\mapsto d_{\mathcal{X}}(A,B)$ is continuous at every pair $A,B\in \mathcal{X}$.
\end{enumerate}

\end{thm}

\begin{thm}{Exercise}\label{ex:comp+cont}
Let $\mathcal{X}$, $\mathcal{Y}$, and $\mathcal{Z}$ be metric spaces.
Assume that the functions $f\:\mathcal{X}\to\mathcal{Y}$
and $g\:\mathcal{Y}\to\mathcal{Z}$ are continuous at every point,
and $h=g\circ f$ is their composition;
that is, $h(A)=g(f(A))$ for each $A\in \mathcal{X}$.
Show that $h\:\mathcal{X}\to\mathcal{Z}$ is continuous at every point.
\end{thm}

\begin{thm}{Exercise}\label{ex:isom-cont}
Show that every distance-preserving map is continuous.
\end{thm}

\section{Congruent triangles}
\label{sec:cong-triangles}

Our next goal is to find a rigorous interpretation for statement (\ref{preaxiomIV}) in Section~\ref{preaxioms}.
To do this, we will introduce the concept of congruent triangles
so instead of \textit{``if we rotate or shift we will not see the difference''} we say that for triangles, the side-angle-side congruence holds.

An \textit{ordered} triple of \textit{distinct} points in a metric space $\mathcal{X}$, 
say $A$, $B$, $C$,
is called \index{triangle}\emph{triangle $ABC$}\label{page:def:triangle} (briefly \index{20@$\triangle$}$\triangle A B C$).
The points $A$, $B$, and $C$ are called \index{vertex!of triangle}\emph{vertices} of $\triangle ABC$.
Note that the triangles $A B C$ and $A C B$ are different.

Two triangles $A' B' C'$ and $A B C$ are  
\index{triangle!congruent triangles}
\index{congruent!triangles}\emph{congruent}
(briefly  \index{32@$\cong$}$\triangle A' B' C'\z\cong\triangle A B C$) if there is a motion $f\:\mathcal{X}\to\mathcal{X}$ such that 
\[A'\z=f(A),
\quad
B'=f(B)
\quad
\text{and}
\quad
C'=f(C).\]

Let $\mathcal X$ be a metric space,
and $f,g\:\mathcal X\to\mathcal X$ be two motions.
Note that the inverse $f^{-1}:\mathcal X\to\mathcal X$,
as well as the composition $f\circ g:\mathcal X\to\mathcal X$,
are also motions.

It follows that ``$\cong$'' is an \index{equivalence relation}\emph{equivalence relation};
that is, every triangle is congruent to itself,
and the following two conditions hold:
\begin{itemize} 
\item If $\triangle A' B' C'\z\cong\triangle A B C$, then $\triangle A B C\z\cong\triangle A' B' C'$.
\item If $\triangle A'' B'' C''\z\cong\triangle A' B' C'$ and $\triangle A' B' C'\z\cong\triangle A B C$,
then 
$$\triangle A'' B'' C''\cong\triangle A B C.$$
\end{itemize}

If $\triangle A' B' C'\z\cong\triangle A B C$,
then $AB\z=A'B'$,
$BC=B'C'$ and $CA=C'A'$.

For a discrete metric, as well as some other metrics, 
the converse also holds.
The following example shows that it does not hold in the Manhattan plane:

\parbf{Example.}\label{example:isometric but not congruent}
Consider three points 
$A=(0,1)$, $B=(1,0)$, and $C\z=(-1,0)$ on the Manhattan plane $(\mathbb{R}^2,d_1)$.
Then
$$d_1(A,B)=d_1(A,C)=d_1(B,C)=2.$$

On one hand,
$$\triangle ABC\cong \triangle ACB.$$
Indeed, the map $(x,y)\z\mapsto (-x,y)$ is a motion of $(\mathbb{R}^2,d_1)$
that sends $A\z\mapsto A$, $B\mapsto C$, and $C\z\mapsto B$.

On the other hand,
$$\triangle ABC\ncong \triangle BCA.$$
Indeed, arguing by contradiction, assume that $\triangle ABC\cong \triangle BCA$; that is, there is a motion $f$ of $(\mathbb{R}^2,d_1)$ that sends $A\mapsto B$, $B\mapsto C$, and $C\mapsto A$.

\begin{wrapfigure}[6]{o}{33mm}
\vskip-4mm
\centering
\includegraphics{mppics/pic-6}
\end{wrapfigure}

We say that $M$ is a midpoint of $A$ and $B$ if 
\[d_1(A,M)=d_1(B,M)=\tfrac12\cdot d_1(A,B).\]
Note that a point $M$ is a midpoint of $A$ and $B$ if and only if $f(M)$ is a midpoint of $B$ and~$C$.

The set of midpoints for $A$ and $B$ is infinite, it contains all points $(t,t)$ for $t\in[0,1]$ (it is the gray segment in the picture above).
In contrast, the midpoint for $B$ and $C$ is unique.
Thus, the map $f$ cannot be bijective, leading to a contradiction.

\begin{thm}{Exercise}\label{ex:ncong}
Consider a metric space with four points $A$, $B$, $C$, $D$ and a metric defined by $AB=AC=AD=BC=BD=1$, and 
$CD=2$.
Show that (a) $\triangle ABC\cong \triangle BAC$, and (b) $\triangle ABC\ncong \triangle BCA$.
\end{thm}

\addtocontents{toc}{\protect\contentsline{part}{\protect\numberline{}Euclidean geometry}{}{}}
  
\chapter{Axioms}
\label{chap:axioms}

\vfill

A system of axioms appears already in Euclid's ``Elements'' --- the most successful and influential textbook ever written.

The systematic study of geometries as axiomatic systems
 was
triggered by the discovery of non-Euclidean geometry.
The branch of mathematics, emerging this way, is called ``Foundations of geometry''.

The most popular system of axioms
was proposed by David Hilbert in 1899.
This is also the first rigorous system by modern standards.
It consists of twenty axioms in five groups, six ``primitive notions'', and three ``primitive terms'';
these are not defined in terms of previously defined concepts.

Later many different systems were proposed.
It is worth mentioning
the system developed by Alexandr Alexandrov \cite{alexandrov} which is intuitive and elementary, 
the system by Friedrich Bachmann \cite{bachmann} based on the concept of symmetry,
and the system by Alfred Tarski \cite{tarski} --- a minimalist system designed for analysis using mathematical logic.

We will use another system close to the one proposed by George Birkhoff \cite{birkhoff}.
It is based on the key observations (\ref{preaxiomI})--(\ref{preaxiomV}) listed in Section~\ref{preaxioms}.
The axioms use the notions of 
metric space, 
lines, 
angles,
triangles,
equalities modulo $2\cdot\pi$ ($\equiv$), 
the continuity of maps between metric spaces,
and the congruence of triangles ($\cong$).
All of these concepts are discussed in the preliminaries.

Our system is built upon metric spaces.
In particular, we use real numbers as a building block. 
For this reason our approach is not purely axiomatic --- we build the theory upon something else;
it resembles a model-based introduction to Euclidean geometry discussed in Section~\ref{page:model}.
We used this approach to minimize the tedious parts which are unavoidable in purely axiomatic foundations.

\newpage

\section{The axioms}
\label{sec:axioms}

\begin{framed}
\begin{enumerate}[I.]
\item\label{def:birkhoff-axioms:0} The \index{plane!Euclidean plane}\index{Euclidean!plane}\emph{Euclidean plane} is a metric space with at least two points.

\item\label{def:birkhoff-axioms:1} 
There is one and only one line that contains any two given distinct points $P$ and $Q$ in the Euclidean plane.

\item\label{def:birkhoff-axioms:2} 
Any angle $\angle AOB$ in the Euclidean plane
defines a real number in the interval $(-\pi,\pi]$.
This number is called the \index{angle!measure}\emph{angle measure of $\angle AOB$}
and is denoted by $\measuredangle A O B$.
It satisfies the following conditions:
\begin{enumerate}[(a)]
\item\label{def:birkhoff-axioms:2a} 
Given a half-line $[O A)$ and $\alpha\in(-\pi,\pi]$, 
there is a unique  half-line $[O B)$, 
such that $\measuredangle A O B= \alpha$.
\item\label{def:birkhoff-axioms:2b} 
For any points $A$, $B$, and $C$, distinct from $O$ we have
$$\measuredangle A O B+\measuredangle B O C
\equiv\measuredangle A O C.$$
\item\label{def:birkhoff-axioms:2c} 
The function 
$$\measuredangle\:(A,O,B)\mapsto\measuredangle A O B$$
is continuous at any triple of points $(A,O,B)$,
such that $O\ne A$ and $O\ne B$ and $\measuredangle A O B\ne\pi$.

\end{enumerate}

\item\label{def:birkhoff-axioms:3}  
In the Euclidean plane, we have
$\triangle A B C\cong\triangle A' B' C'$
if and only if 
\begin{align*}
A' B'&=A B, & A' C'&= A C, &&\text{and}
&\measuredangle C' A' B'&=\pm\measuredangle C A B.
\end{align*}
\item\label{def:birkhoff-axioms:4}
If for two triangles $\triangle ABC$, $\triangle AB'C'$ in the Euclidean plane
and for $k>0$ we have
\begin{align*}
B'&\in [AB),
& C'&\in [AC),
\\
AB'&=k\cdot AB,&
AC'&=k\cdot AC,
\end{align*}
then
\begin{align*}
B'C'&=k\cdot BC,&
\measuredangle ABC&=\measuredangle AB'C',
&
\measuredangle ACB&=\measuredangle AC'B'.
\end{align*}
\end{enumerate}
\end{framed}

From now on,  
we can use no information about the Euclidean plane that does not follow from the five axioms above.

\begin{thm}{Classroom exercise}\label{ex:infinite}
Show that there are (a) an infinite set of points,
(b) an infinite set of lines on the plane.
\end{thm}

\section{Lines and half-lines}

\begin{thm}[\abs]{Proposition}\label{lem:line-line}
\let\thefootnote\relax\footnotetext{${}^\a$ The mark ``$\a$'' indicates that Axiom~\ref{def:birkhoff-axioms:4} is not used in the proof.
Please ignore this mark for now; it will become important in Chapter~\ref{chap:non-euclid}.}
Any two distinct lines intersect at most at one point.
\end{thm}

\parit{Proof.}
Assume that two lines $\ell$ and $m$ intersect at two distinct points $P$ and~$Q$.
Applying Axiom~\ref{def:birkhoff-axioms:1}, we get that $\ell=m$.
\qeds

\begin{thm}{Exercise}\label{ex:[OA)=[OA')}
Suppose $A'\in[OA)$ and $A'\not=O$. 
Show that 
\[[O A)\z=[O A').\]

\end{thm}

\section{A zero angle}

\begin{thm}[\abs]{Proposition}\label{lem:AOA=0}
$\measuredangle A O A= 0$ for every $A\not=O$.
\end{thm}

\parit{Proof.}
According to Axiom~\ref{def:birkhoff-axioms:2b},
$$\measuredangle A O A
+
\measuredangle A O A 
\equiv
\measuredangle A O A.$$
Subtract  $\measuredangle A O A$ from both sides, we get that
$\measuredangle A O A \equiv 0$.

By Axiom~\ref{def:birkhoff-axioms:2}, $-\pi<\measuredangle A O A\le \pi$;
therefore $\measuredangle A O A \z= 0$.
\qeds

\begin{thm}{Exercise}\label{ex:2.4} 
Assume $\measuredangle A O B= 0$.
Show that $[OA)=[OB)$.
\end{thm}

\begin{thm}[\abs]{Proposition}\label{lem:AOB+BOA=0}
For any $A$ and $B$ distinct from $O$,
we have 
$$\measuredangle A O B\equiv-\measuredangle B O A.$$

\end{thm}

\parit{Proof.}
According to Axiom~\ref{def:birkhoff-axioms:2b},
$$\measuredangle A O B+\measuredangle B O A \equiv\measuredangle A O A.$$
By Proposition~\ref{lem:AOA=0}, $\measuredangle A O A=0$.
Hence the result.
\qeds

\section{A straight angle}

If $\measuredangle A O B=\pi$,
we say that $\angle A O B$ is a 
\index{angle!straight angle}\emph{straight angle}.
By Proposition~\ref{lem:AOB+BOA=0}, 
if $\angle A O B$ is straight,
then so is $\angle B O A$.

We say that point $O$ \index{between}\emph{lies between} points $A$ and $B$, 
if $O\not= A$, $O\not= B$, and $O\in[A B]$.

\begin{thm}[\abs]{Theorem}\label{thm:straight-angle}
An angle $A O B$ is straight 
if and only if $O$ 
\index{between}\emph{lies between} $A$ and~$B$.
\end{thm}

\begin{wrapfigure}{r}{40mm}
\centering
\vskip-0mm
\includegraphics{mppics/pic-8}
\end{wrapfigure}

\parit{Proof.}
By Exercise~\ref{ex:trig==},  we may assume that
$O A = O B = 1$.

\parit{``If'' part.}
Assume $O$  
lies between $A$ and~$B$.
Set  $\alpha=\measuredangle A O B$.

Applying Axiom~\ref{def:birkhoff-axioms:2a},
we get a half-line $[OA')$ such that $\alpha\z=\measuredangle B O A'$.
By Exercise~\ref{ex:trig==}, we can assume that $OA'=1$.
According to Axiom~\ref{def:birkhoff-axioms:3},
\[\triangle AOB\z\cong\triangle BOA'.\]
Suppose that $f$ denotes the corresponding motion of the plane;
that is, $f$ is a motion such that $f(A)=B$, $f(O)=O$, and $f(B)=A'$. 

\begin{wrapfigure}{o}{40mm}
\centering
\vskip-2mm
\includegraphics{mppics/pic-10}
\end{wrapfigure}

Then 
\[(A'B)=f((AB))\ni f(O)=O.\]
Therefore, both lines $(AB)$ and $(A'B)$ contain $B$ and~$O$.
By Axiom~\ref{def:birkhoff-axioms:1}, $(AB)=(A'B)$.

By the definition of a line,
$(AB)$ contains exactly two points, $A$ and $B$, that are at distance $1$ from~$O$.
Since $OA' = 1$, it follows that $A' = A$ or $A' = B$.
The isometry $f$ must be injective; see Exercise~\ref{ex:dist-preserv=>injective}.
Therefore, $A' = f(B) \ne f(A) = B$.
That is, $A' \ne B$, and hence $A = A'$.

By Axiom~\ref{def:birkhoff-axioms:2b} and Proposition~\ref{lem:AOA=0}, we get that
\begin{align*}
2\cdot\alpha&=
\measuredangle AOB+\measuredangle BOA'=
\\
&=\measuredangle AOB+\measuredangle BOA\equiv
\\
&\equiv\measuredangle AOA=
\\
&= 0.
\end{align*}
Therefore, by Exercise~\ref{ex:2a=0}, $\alpha$ is either $0$ or~$\pi$.

Since $[OA)\ne [OB)$,  
we have that $\alpha\ne 0$, see Exercise~\ref{ex:2.4}.
Therefore, $\alpha=\pi$.

\parit{``Only if'' part.}
Suppose that $\measuredangle A O B= \pi$.
Consider the line $(OA)$ and choose a point $B'$ on $(OA)$ so that $O$ lies between $A$ and~$B'$.

From above, we have that $\measuredangle AOB'=\pi$.
Applying Axiom~\ref{def:birkhoff-axioms:2a}, 
we get that $[O B)\z=[O B')$.
In particular, $O$ lies between $A$ and~$B$.
\qeds 

A triangle $ABC$ is called 
\index{triangle!degenerate triangle}\index{degenerate! triangle}\emph{degenerate}
if $A$, $B$, and $C$ lie on one line.
The following corollary is just a reformulation of Theorem~\ref{thm:straight-angle}.

\begin{thm}[\abs]{Corollary}\label{cor:degenerate=pi}
A triangle is degenerate if and only if one of its angles is equal to $\pi$ or~$0$.
Moreover, in a degenerate triangle, the angle measures are $0$, $0$, and $\pi$.
\end{thm}

\begin{thm}{Exercise}\label{ex:lineAOB}
Show that three distinct points $A$, $O$, and $B$ lie on one line if and only if 
$$2\cdot \measuredangle AOB\equiv 0.$$ 

\end{thm}

\begin{thm}{Exercise}\label{ex:ABCO-line}
Let $A$, $B$, and $C$ be three points distinct from~$O$.
Show that $B$, $O$, and $C$ lie on one line if and only if
$$2\cdot \measuredangle AOB\equiv 2\cdot \measuredangle AOC.$$ 

\end{thm}

\begin{thm}{Exercise}\label{ex:infinite-number-of-lines} 
Show that there is a nondegenerate triangle.
\end{thm}

\section{Vertical angles}

A pair of angles $AOB$ and $A'OB'$ 
is called \index{angle!vertical angles}\index{vertical angles}\emph{vertical}
if point $O$ 
lies between $A$ and $A'$ 
and between $B$ and $B'$ at the same time.

\begin{thm}[\abs]{Proposition}\label{prop:vert}
Vertical angles have equal measures.
\end{thm}

\parit{Proof.}
Assume that the angles $AOB$ and $A'OB'$ are vertical.
Note that $\angle AOA'$ and $\angle BOB'$ are straight.
Therefore, $\measuredangle AOA'\z=\measuredangle BOB'=\pi$.

{

\begin{wrapfigure}[4]{o}{24mm}
\vskip-2mm
\centering
\includegraphics{mppics/pic-12}
\end{wrapfigure}

It follows that
\begin{align*}
0&=\measuredangle AOA'-\measuredangle BOB'\equiv
\\
&\equiv 
\measuredangle AOB+\measuredangle BOA'-\measuredangle BOA'-\measuredangle A'OB'
\equiv
\\
&\equiv\measuredangle AOB-\measuredangle A'OB'.
\end{align*}
Since $-\pi<\measuredangle AOB\le \pi$ and $-\pi<\measuredangle A'OB'\le \pi$, we get that $\measuredangle AOB\z=\measuredangle A'OB'$.
\qeds

}

\section{Reflection across a point}

\begin{thm}[\abs]{Lemma}\label{lem:point-reflection}
Given distinct points $X$ and $O$, there is a unique point $X'$ such that $O$ is the midpoint of a line segment $[XX']$.
\end{thm}

\begin{wrapfigure}{o}{33mm}
\vskip-0mm
\centering
\includegraphics{mppics/pic-915}
\end{wrapfigure}

\parit{Proof.}
The point $O$ is the midpoint of a line segment $[XX']$ if and only if $X'\in [XO)$ and $XX'=2\cdot XO$.

By Axiom~\ref{def:birkhoff-axioms:1}, the line $(XO)$ and therefore half-line $[XO)$ are uniquely defined.
Further, by Exercise~\ref{ex:trig==} there is a unique point $X'\in [XO)$ such that $XX'=2\cdot XO$.
\qeds

The point $X'$ provided by the lemma is called \index{reflection across a point}\emph{reflection} of point $X$ across $O$.
We also assume that $O'=O$; that is, $O$ is a reflection of itself across itself.

The point $O$ is called the \index{center!of reflection}\emph{center of reflection}.

\begin{thm}[\abs]{Proposition}\label{prop:point-reflection}
A reflection across a point is a motion of the plane.
\end{thm}

\parit{Proof.}
Let $O$ be the center of reflection.
Observe that if $X'$ is a reflection of $X$ across $O$,
then $X$ is a reflection of $X'$.
In other words, the composition of the reflection with itself is the identity map.
In particular, the reflection is a bijection.

\begin{wrapfigure}{o}{33mm}
\vskip-6mm
\centering
\includegraphics{mppics/pic-76}
\end{wrapfigure}

Now choose two points $X$ and $Y$;
let $X'$ and $Y'$ be their reflections across $O$.
To check that the reflection is distance preserving, we need to show that $X'Y'=XY$.

We may assume that $X$, $Y$, and $O$ are distinct; otherwise, the statement is trivial.
By the definition of the reflection, we have that $OX=OX'$, $OY=OY'$.
Note also that the angles $XOY$ and $X'OY'$ are vertical.
By \ref{prop:vert}, $\measuredangle XOY\z=\measuredangle X'OY'$.

Now Axiom~\ref{def:birkhoff-axioms:3} implies that $\triangle XOY\cong\triangle X'OY'$ and $X'Y'=XY$.
\qeds

\begin{thm}{Exercise}\label{ex:refelection-of-line}
Show that the reflection across a point maps any line to a line.
\end{thm}

\chapter{Half-planes}\label{chap:half-planes}

This chapter contains lengthy proofs of intuitively evident statements.
It is okay to skip it, but make sure you understand the definitions of positive/negative angles and that your intuition agrees with statements \ref{thm:signs-of-triug}, \ref{prop:half-plane}, \ref{cor:half-plane}, \ref{thm:pasch}, and \ref{thm:abc}.

\section{Sign of an angle}

Positive and negative angles can be visualized as {}\emph{counterclockwise} and {}\emph{clockwise} directions; they are defined the following way:
\begin{itemize}
\item The angle $A O B$ is called \index{angle!positive and negative angles}\emph{positive} 
if $0<\measuredangle A O B<\pi$;
\item The  angle $A O B$ is called {}\emph{negative} 
if $\measuredangle A O B<0$.
\end{itemize}

According to the above definitions the straight angle, as well as the zero angle,
are neither positive nor negative.

\begin{thm}{Classroom exercise}\label{ex:AOB+<=>BOA-}
Show that $\angle A O B$ is positive if and only if $\angle B O A$ is negative.
\end{thm}

\begin{thm}{Lemma}\label{lem:straight-sign}
Let $\angle AOB$ be straight.
Then $\angle AOX$ is positive 
if and only if $\angle BOX$ is negative.
\end{thm}

\parit{Proof.}
Set $\alpha=\measuredangle AOX$ 
and 
$\beta=\measuredangle BOX$.
Since $\angle AOB$ is straight,
$$\alpha-\beta\equiv \pi.\eqlbl{eq:alpha-beta}$$

It follows that $\alpha=\pi$ $\Leftrightarrow$ $\beta=0$
and $\alpha=0$ $\Leftrightarrow$ $\beta=\pi$.
In these two cases, the signs of $\angle AOX$ and $\angle BOX$ are undefined.

In the remaining cases we have that $|\alpha|<\pi$ and $|\beta|<\pi$.
If $\alpha$ and $\beta$ have the same sign, then $|\alpha-\beta|<\pi$;
the latter contradicts \ref{eq:alpha-beta}.
Hence the statement follows.
\qeds

\begin{thm}{Exercise}\label{ex:PP(PN)}
Assume that angles $ABC$ and $A'B'C'$ have the same sign
and 
\[2\cdot \measuredangle ABC\equiv 2\cdot \measuredangle A'B'C'.\]
Show that $\measuredangle ABC= \measuredangle A'B'C'$.
\end{thm}

\section{Intermediate value theorem}

\begin{thm}{Intermediate value theorem}\label{thm:intermidiate}
Let $f\:[a,b]\to \mathbb{R}$ be a continuous function.
Assume 
$f(a)$ and $f(b)$ have opposite signs.
Then $f(t_0)=0$ for some $t_0\in[a,b]$.
\end{thm}

\begin{wrapfigure}{r}{38mm}
\vskip-6mm
\centering
\includegraphics{mppics/pic-14}
\end{wrapfigure}

The intermediate value theorem is assumed to be known;
it should be covered in any calculus course.
We will use only the following corollary:

\begin{thm}[\abs]{Corollary}\label{cor:intermidiate}
Assume that for any $t\z\in [0,1]$ we have three points in the plane  $O_t$, $A_t$, and $B_t$, such that 
\begin{enumerate}[(a)]
\item Each function $t\mapsto O_t$, $t\mapsto A_t$, and $t\z\mapsto B_t$ is continuous.
\end{enumerate}

\begin{enumerate}[(a)]\addtocounter{enumi}{1}
\item For any $t\in [0,1]$, the points $O_t$, $A_t$, and $B_t$ do not lie on one line.  
\end{enumerate}
Then $\angle A_0O_0B_0$ and $\angle A_1O_1B_1$ have the same sign.
\end{thm}

\parit{Proof.}
Consider the function 
$f(t)=\measuredangle A_tO_tB_t$.

Since 
the points $O_t$, $A_t$, and $B_t$ do not lie on one line,
Theorem~\ref{thm:straight-angle} implies that $f(t)=\measuredangle A_tO_tB_t\ne 0$ nor $\pi$ for any $t\in[0,1]$.

Therefore, by Axiom~\ref{def:birkhoff-axioms:2c} and Exercise~\ref{ex:comp+cont},
$f$ is a continuous function.
By the intermediate value theorem, $f(0)$ and $f(1)$ have the same sign;
hence the result follows.
\qeds

\section{Same sign lemmas}

\begin{thm}[\abs]{Lemma}\label{lem:signs}
Assume $Q'\in [PQ)$ and $Q'\z\ne P$.
Then for every $X\z\notin (PQ)$ the angles $PQX$ and $PQ'X$ have the same sign.
\end{thm}

{

\begin{wrapfigure}{o}{33mm}
\centering
\vskip-5mm
\includegraphics{mppics/pic-16}
\end{wrapfigure}

\parit{Proof.}
By Exercise~\ref{ex:trig==},
for any $t\in [0,1]$ there is a unique point $Q_t\in[PQ)$ 
such that 
\[PQ_t=  (1-t)\cdot PQ+t\cdot PQ'.\]
The map $t\mapsto Q_t$ is continuous,
\begin{align*}
Q_0&=Q,
&
Q_1&=Q'
\end{align*}
and for any $t\in [0,1]$, 
we have that $P\z\ne Q_t$.

}

Applying Corollary~\ref{cor:intermidiate},
for $P_t=P$, $Q_t$, and $X_t=X$, we get that $\angle PQX$ has the same sign as $\angle PQ'X$.
\qeds

\begin{thm}[\abs]{Signs of angles of a triangle}\label{thm:signs-of-triug}
In arbitrary nondegenerate triangle $ABC$,
the angles $ABC$, $BCA$, and $CAB$ have the same sign. 
\end{thm}

{

\begin{wrapfigure}{o}{33mm}
\vskip-4mm
\centering
\includegraphics{mppics/pic-18}
\end{wrapfigure}

\parit{Proof.}
Choose a point $Z\in (AB)$ so that $A$ lies between $B$ and~$Z$.

According to Lemma~\ref{lem:signs},
the angles $ZBC$ and $ZAC$ have the same sign.

Note that $\measuredangle ABC=\measuredangle ZBC$
and 
$$\measuredangle ZAC+\measuredangle CAB\equiv \pi.$$
Therefore, $\angle CAB$ has the same sign as $\angle ZAC$
which in turn has the same sign as $\measuredangle ABC\z=\measuredangle ZBC$.

}

Repeating the same argument for $\angle BCA$ and $\angle CAB$,
we get the result.
\qeds

\begin{wrapfigure}{r}{26mm}
\vskip-2mm
\centering
\includegraphics{mppics/pic-838}
\end{wrapfigure}

\begin{thm}{Classroom exercise}\label{ex:between}
Let $A$, $B$, $C$, and $P$ be points such that $B$ lies between $A$ and $C$, and let $P \not\in (AC)$.
Show that the angles $\angle APB$, $\angle BPC$, and $\angle APC$ have the same sign.
Conclude that
\[|\measuredangle APB|+|\measuredangle BPC|=|\measuredangle APC|.\]
\end{thm}

\begin{thm}[\abs]{Lemma}\label{lem:signsXY}
Assume $[XY]$ does not intersect $(PQ)$,
then angles $PQX$ and $PQY$ 
have the same sign.
\end{thm}

\begin{wrapfigure}{o}{26mm}
\vskip-4mm
\centering
\includegraphics{mppics/pic-20}
\end{wrapfigure}

The proof is nearly identical to the one above.

\parit{Proof.}
According to Exercise~\ref{ex:trig==},
for any $t\z\in [0,1]$ there is a point  $X_t\in[XY]$, 
such that 
\[XX_t= t\cdot XY.\]
The map $t\mapsto X_t$ is continuous.
Moreover, $X_0=X$, $X_1=Y$, and $X_t\notin(QP)$ for any $t\in [0,1]$.

Applying Corollary~\ref{cor:intermidiate},
for $P_t\z=P$, $Q_t\z=Q$, and $X_t$, we get that
$\angle PQX$ has the same sign as $\angle PQY$.
\qeds

\section{Half-planes}

\begin{thm}{Proposition}\label{prop:half-plane}
Assume $X,Y\notin(PQ)$.
Then angles $PQX$ and $PQY$ have the same sign if and only if $[XY]$ does not intersect $(PQ)$.
\end{thm}

\begin{wrapfigure}{o}{30mm}
\includegraphics{mppics/pic-22}
\centering
\end{wrapfigure}

\parit{Proof.} The if-part follows from Lemma~\ref{lem:signsXY}. 

Assume $[XY]$ intersects $(PQ)$;
suppose that $Z$ denotes the point of intersection.
Without loss of generality, we can assume $Z\ne P$.

Note that $Z$ lies between $X$ and $Y$.
According to Lemma~\ref{lem:straight-sign}, $\angle PZX$ and $\angle PZY$ have opposite signs.
It proves the statement if $Z=Q$.

If $Z\ne Q$, then $\angle ZQX$ and $\angle QZX$ have opposite signs by \ref{thm:signs-of-triug}.
In the same way, we get that $\angle ZQY$ and $\angle QZY$ have opposite signs.

If $Q$ lies between $Z$ and $P$, then by Lemma~\ref{lem:straight-sign} two pairs of angles $\angle PQX$, $\angle ZQX$ and $\angle PQY$, $\angle ZQY$ have opposite signs. 
It follows that $\angle PQX$ and $\angle PQY$ have opposite signs as required.

In the remaining case $[QZ)=[QP)$ and therefore $\angle PQX=\angle ZQX$ and $\angle PQY=\angle ZQY$. 
Therefore again $\angle PQX$ and $\angle PQY$ have opposite signs as required.
\qeds

\begin{thm}[\abs]{Corollary}\label{cor:half-plane}
The complement of a line $(PQ)$ in the plane 
can be presented in a unique way as a union of two disjoint subsets 
called \index{half-plane}\emph{half-planes}
such that 
\begin{enumerate}[(a)]
\item\label{cor:half-plane:angle} Two points $X,Y\notin(PQ)$ lie in the same half-plane if and only if the angles $PQX$ and $PQY$ have the same sign.
\item\label{cor:half-plane:intersect} Two points $X,Y\notin(PQ)$ lie in the same half-plane if and only if $[XY]$ does not intersect~$(PQ)$.
\end{enumerate}

\end{thm}

{

\begin{wrapfigure}{r}{26mm}
\vskip-4mm
\centering
\includegraphics{mppics/pic-24}
\end{wrapfigure}

We say that $X$ and $Y$ lie on {}\emph{one side of} $(PQ)$ if they lie in one of the half-planes of $(PQ)$ and we say that  $P$ and $Q$ lie on {}\emph{opposite sides of} $\ell$ if they lie in different half-planes of~$\ell$.

\begin{thm}{Exercise}\label{ex:vert-intersect}
Suppose that angles $AOB$ and $A'OB'$ are vertical and $B\notin (OA)$.
Show that the line $(AB)$ does not intersect the segment~$[A'B']$.
\end{thm}

}

Consider triangle $ABC$.
Recall that its vertices are points $A$, $B$ and $C$;
segments $[AB]$, $[BC]$, and $[CA]$ will be called the \index{side!of triangle}\emph{sides of the triangle}.

{

\begin{wrapfigure}{r}{27mm}
\vskip-4mm
\centering
\includegraphics{mppics/pic-26}
\end{wrapfigure}

\begin{thm}[\abs]{Pasch's theorem}\label{thm:pasch}\index{Pasch's theorem}
Assume that a line $\ell$ does not pass thru any vertex of a triangle.
Then it intersects either two or zero sides of the triangle.
\end{thm}

\parit{Proof.}
Assume that line $\ell$ intersects side $[AB]$ of the triangle $ABC$ and does not pass thru $A$, $B$, and $C$.

By Corollary~\ref{cor:half-plane}, the vertices $A$ and $B$ lie on opposite sides of~$\ell$.

}

The vertex $C$ may lie on the same side with $A$ and on the opposite side with $B$ or the other way around.
By Corollary~\ref{cor:half-plane}, in the first case, $\ell$ intersects side $[BC]$ and does not intersect $[AC]$; in the second case, $\ell$ intersects side $[AC]$ and does not intersect $[BC]$.
Hence the statement follows.
\qeds

{

\begin{wrapfigure}[5]{r}{29mm}
\vskip-4mm
\centering
\includegraphics{mppics/pic-28}
\bigskip
\includegraphics{mppics/pic-30}
\end{wrapfigure}

\begin{thm}{Exercise}\label{ex:signs-PXQ-PYQ}
Show that two points $X,Y\z\notin(PQ)$ lie on the same side of $(PQ)$
if and only if angles $PXQ$ and $PYQ$ have the same sign.
\end{thm}

\begin{thm}{Exercise}\label{ex:chevinas}
Let $\triangle ABC$ be a nondegenerate triangle,
$A'\in[BC]$  and 
$B'\in [AC]$.
Show that $[AA']$ intersects $[BB']$.
\end{thm}

\begin{thm}{Exercise}\label{ex:Z}
Assume that points $X$ and $Y$ lie on opposite sides of the line~$(PQ)$.
Show that $[PX)$ does not intersect~$[QY)$. 
\end{thm}

}

\begin{thm}{Advanced exercise}\label{ex:angle-measures}
The following quantity 
$$\tilde\measuredangle ABC=
\begin{cases}
\pi&\text{if}\ \measuredangle ABC=\pi
\\
-\measuredangle ABC&\text{if}\ \measuredangle ABC<\pi
\end{cases}
$$
can serve as the angle measure; 
that is, the axioms hold if one exchanges $\measuredangle$ to $\tilde\measuredangle$ everywhere.

Show that $\measuredangle$ and $\tilde\measuredangle$ are the only possible angle measures on the plane. 

Show that without Axiom~\ref{def:birkhoff-axioms:2c}, this is no longer true.
\end{thm}

\section{Triangle with given sides}

Given $\triangle ABC$, set 
\begin{align*}
a&=BC,
&
b&=CA,
&
c&=AB.
\end{align*}
Without loss of generality, we may assume that 
\[a\le b \le c.\]
Then all three triangle inequalities for $\triangle ABC$
hold if and only if 
\[c\le a+b.\]
The following theorem states that this is the only restriction on $a$, $b$, and~$c$.

\begin{thm}[\abs]{Theorem}\label{thm:abc}
Assume that $0<a\le b\le c\le a+b$.
Then there is a triangle with sides $a$, $b$, and $c$;
that is, there is $\triangle ABC$ 
such that $a=BC$, $b=CA$, and $c=AB$.
\end{thm}

A proof of the following proposition is given at the end of the section.

\begin{thm}[\abs]{Proposition}\label{prop:C-cont}
Fix a real number $r>0$ 
and two distinct points $A$ and~$B$.
Then for 
any real number $\beta\in [0,\pi]$,
there is a unique point $C_\beta$ such that $BC_\beta=r$
and $\measuredangle ABC_\beta=\beta$.
Moreover, $\beta\mapsto C_\beta$ 
is a continuous map from $[0,\pi]$ to the plane.
\end{thm}

\parit{Proof of Theorem~\ref{thm:abc} modulo Proposition~\ref{prop:C-cont}.}\label{page:proof:thm:abc}
Fix the points $A$ and $B$ such that $AB=c$.
Given $\beta\in [0,\pi]$,
suppose that $C_\beta$ denotes the point in the plane such that $BC_\beta\z=a$ and $\measuredangle ABC=\beta$.

According to \ref{prop:C-cont},
the map
$\beta\mapsto C_\beta$ is continuous.
Therefore, the function $b\:\beta\mapsto AC_\beta$ is continuous
(formally, it follows from Exercise~\ref{ex:dist-cont} and Exercise~\ref{ex:comp+cont}).

Note that $b(0)=c-a$ and $b(\pi)=c+a$.
Since $c-a\le b\le c+a$,
by the intermediate value theorem (\ref{thm:intermidiate})
there is $\beta_0\in[0,\pi]$ such that
$b(\beta_0)=b$,
hence the result. 
\qeds

The proof of Proposition~\ref{prop:C-cont} relies on the following lemma.

Assume $r>0$ and $\pi>\beta>0$.
Consider the triangle $ABC$ such that 
$AB=BC=r$ and $\measuredangle ABC=\beta$.
The existence of such a triangle follows from Axiom~\ref{def:birkhoff-axioms:2a} and Exercise~\ref{ex:trig==}.
By Axiom~\ref{def:birkhoff-axioms:3},
the values
$\beta$ and~$r$ define the triangle $ABC$ up to the congruence.
In particular, the distance $AC$ depends only on $\beta$ and~$r$.

{

\begin{wrapfigure}{o}{22mm}
\vskip-5mm
\centering
\includegraphics{mppics/pic-32}
\end{wrapfigure}

Let
$$s(\beta,r)\df AC.$$

\begin{thm}[\abs]{Lemma}\label{lem:f(r,a)}
Given $r>0$ and $\epsilon>0$, there is $\delta>0$ such that
\[0<\beta<\delta\quad\Longrightarrow\quad s(r,\beta)<\epsilon.\]

\end{thm}

}

{

\begin{wrapfigure}{l}{36mm}
\vskip-8mm
\centering
\includegraphics{mppics/pic-34}
\end{wrapfigure}

\parit{Proof.}
Fix two points $A$ and $B$ such that $AB\z=r$.

Choose a point $X$ such that $\measuredangle ABX$ is positive.
Let $Y\in [AX)$ be the point such that $AY=\tfrac\epsilon2$;
it exists by Exercise~\ref{ex:trig==}.

Points $X$ and $Y$ lie on the same side of $(AB)$;
therefore, $\angle ABY$ is positive. 
Set $\delta\z=\measuredangle ABY$.

Suppose $0<\beta<\delta$;
by Axiom~\ref{def:birkhoff-axioms:2a}, we can choose $C$ so that $\measuredangle ABC\z=\beta$ and $BC\z=r$.
Furthermore, we can choose a half-line $[BZ)$ such that $\measuredangle ABZ\z=\tfrac12\cdot \beta$.

}

Points $A$ and $Y$ lie on opposite sides of~$(BZ)$ and  $\measuredangle ABZ\z\equiv -\measuredangle CBZ$.
In particular, $(BZ)$ intersects $[AY]$;
denote by $D$ the point of intersection.

Since $D$ lies between $A$ and $Y$, we have $AD<AY$.

Since $D$ lies on $(BZ)$ we have that $\measuredangle ABD\z\equiv -\measuredangle CBD$.
By Axiom~\ref{def:birkhoff-axioms:3},
$\triangle ABD\cong \triangle CBD$.
It follows that
\begin{align*}
s(r,\beta)&=AC\le
\\
&\le AD+DC=
\\
&=2\cdot AD< 
\\
&< 2\cdot AY=
\\
&=\epsilon.
\end{align*}
\qedsf

\parit{Proof of Proposition~\ref{prop:C-cont}.}
The existence and uniqueness of $C_\beta$ follow from Axiom~\ref{def:birkhoff-axioms:2a} and Exercise~\ref{ex:trig==}.

If $\beta_1\ne\beta_2$, then
$$C_{\beta_1}C_{\beta_2}=s(r,|\beta_1-\beta_2|).$$

By \ref{lem:f(r,a)}, the map $\beta\mapsto C_\beta$ is continuous.
\qeds

Given a positive real number $r$ and a point $O$,
the set $\Gamma$ of all points on distance $r$ from $O$ is called a \index{circle}\emph{circle} 
with \index{radius}\emph{radius} $r$ and \index{center}\emph{center}~$O$.

\begin{wrapfigure}{r}{40mm}
\centering
\includegraphics{mppics/pic-35}
\end{wrapfigure}

\begin{thm}{Exercise}\label{ex:intersecting-circles-3}
Show that two circles intersect if and only if 
\[|r_1-r_2|\le d\le r_1+r_2,\]
where $r_1$ and $r_2$ denote their radiuses, and $d$ is the distance between their centers.
\end{thm}

\chapter{Congruent triangles}\label{chap:cong}

\section{Side-angle-side}

Our next goal is to give conditions that guarantee the congruence of two triangles.

One such condition is given in Axiom~\ref{def:birkhoff-axioms:3}; it states that if two pairs of sides of two triangles are equal, and the included angles are equal up to sign, then the triangles are congruent.
This axiom is also called the {}\emph{side-angle-side congruence condition}, or briefly, \index{SAS congruence condition}\emph{SAS}.

\section{Angle-side-angle}

\begin{thm}[\abs]{ASA condition}\label{thm:ASA}\index{ASA congruence condition}
Assume that 
\begin{align*}
AB&=A'B',
&
\measuredangle A B C &= \pm\measuredangle A' B' C', 
&
\measuredangle C A B&=\pm\measuredangle C' A' B'
\end{align*}
 and $\triangle A' B' C'$ is nondegenerate.
Then 
$$\triangle A B C\cong\triangle A' B' C'.$$

\end{thm}

For degenerate triangles this statement does not hold.
For example, consider one triangle with sides $1$, $4$, $5$ 
and the other with sides $2$, $3$,~$5$.

\parit{Proof.} 
According to Theorem~\ref{thm:signs-of-triug},
either
\[
\measuredangle A B C = \measuredangle A' B' C'
\quad\text{and}\quad
\measuredangle C A B=\measuredangle C' A' B'
\eqlbl{eq:+angles}\]
or
\[
\measuredangle A B C = -\measuredangle A' B' C'
\quad\text{and}\quad
\measuredangle C A B=-\measuredangle C' A' B'.
\eqlbl{eq:-angles}\]
Let us assume that \ref{eq:+angles} holds; 
the case \ref{eq:-angles} is analogous.

\begin{wrapfigure}{o}{26mm}
\vskip-2mm
\centering
\includegraphics{mppics/pic-36}
\end{wrapfigure}

Let $C''$ be the point on the half-line $[A' C')$ such that $A' C''\z=A C$. 

By Axiom~\ref{def:birkhoff-axioms:3}, 
$\triangle A' B' C''\cong \triangle A B C$. 
Applying Axiom~\ref{def:birkhoff-axioms:3} again,
we get that
$$\measuredangle A' B' C'' = \measuredangle A B C=\measuredangle A' B' C'.$$
By Axiom~\ref{def:birkhoff-axioms:2a}, $[B'C')=[B C'')$. 
Hence
$C''$ lies on $(B' C')$ as well as on~$(A' C')$.

Since $\triangle A' B' C'$ is not degenerate, $(A' C')$ is distinct from~$(B' C')$.
Applying Axiom~\ref{def:birkhoff-axioms:1}, we get that $C''=C'$. 

Therefore, 
$\triangle A' B' C'=\triangle A' B' C''\cong\triangle A B C$.
\qeds

\section{Isosceles triangles}

A triangle with two equal sides is called \index{isosceles triangle}\emph{isosceles};
the remaining side is called the \index{base}\emph{base}.

\begin{thm}[\abs]{Theorem}\label{thm:isos}
Assume $\triangle A B C$ is an isosceles triangle with the base $[A B]$. 
Then 
$$\measuredangle A B C\equiv -\measuredangle B A C.$$
Moreover, the converse holds if $\triangle A B C$ is nondegenerate.
\end{thm}

The following proof is due to Pappus of Alexandria.

\begin{wrapfigure}[9]{r}{25mm}
\vskip-4mm
\centering
\includegraphics{mppics/pic-38}
\end{wrapfigure}

\parit{Proof.}
Note that
$$C A = C B,
\quad 
C B=C A,
\quad
\measuredangle A C B \equiv -\measuredangle B C A.$$
By Axiom~\ref{def:birkhoff-axioms:3},
$$\triangle C A B\cong\triangle C B A.$$
Applying the theorem on the signs of angles of triangles (\ref{thm:signs-of-triug}) and Axiom~\ref{def:birkhoff-axioms:3} again,
we get that 
$$\measuredangle B A C
\equiv -\measuredangle A B C.$$

To prove the converse, we assume that
$\measuredangle C A B \z\equiv - \measuredangle C B A$.
By the ASA condition (\ref{thm:ASA}), $\triangle C A B\z\cong\triangle CBA$.
Therefore,~$C A\z=C B$.
\qeds

A triangle with three equal sides is called \index{equilateral triangle}\emph{equilateral}. 

\begin{thm}{Exercise}\label{ex:equilateral}
Let $\triangle ABC$ be an equilateral triangle.
Show that 
\[\measuredangle ABC=\measuredangle BCA=\measuredangle CAB.\]

\end{thm}

\section{Side-side-side}

\begin{thm}[\abs]{SSS condition}\label{thm:SSS}\index{SSS congruence condition}
$\triangle A B C\cong\triangle A' B' C'$ if 
$$A' B'=A B,
\quad 
B' C'=B C,
\quad 
\text{and}
\quad 
C' A'=C A.$$

\end{thm}

Note that this condition is valid for degenerate triangles as well.

\parit{Proof.} 
Choose $C''$ so that $A' C''= A' C'$ and $\measuredangle B' A' C''= \measuredangle B A C$.
According to Axiom~\ref{def:birkhoff-axioms:3},
$$\triangle A' B' C''\cong\triangle A B C.$$

It will suffice to
prove that 
$$\triangle A' B' C'\cong\triangle A' B' C''.\eqlbl{eq:A'B'C'simA'B'C''}$$
The condition \ref{eq:A'B'C'simA'B'C''} holds trivially if $C''\z=C'$.
Thus, it remains to consider the case $C''\z\ne C'$.

\begin{wrapfigure}{o}{30mm}
\centering
\includegraphics{mppics/pic-40}
\end{wrapfigure}

Clearly, the corresponding sides of $\triangle A' B' C'$ and $\triangle A' B' C''$ are equal.
Hence the triangles
$\triangle C' A' C''$ and $\triangle C' B' C''$ are isosceles.
By Theorem~\ref{thm:isos}, we have 
\begin{align*}
 \measuredangle A' C'' C'&\equiv -\measuredangle A' C' C'',
\\
\measuredangle C' C'' B'&\equiv -\measuredangle C'' C' B'.
\end{align*}
Adding them, we get that
$$\measuredangle A' C'' B'
\equiv -\measuredangle A' C' B'.$$
Applying Axiom~\ref{def:birkhoff-axioms:3} again,
we get \ref{eq:A'B'C'simA'B'C''}.
\qeds

\begin{thm}[\abs]{Corollary}\label{cor:degenerate-trig}
If $AB+BC=AC$, then $B\in [AC]$.
\end{thm}

\parit{Proof.}
Since $AB+BC=AC$, we can choose $B'\in [AC]$ such that $AB=AB'$;
observe that $BC\z=B'C$.

We may assume that $AB>0$ and $BC>0$;
otherwise, $A=B$ or $B=C$, and the statement follows.
In this case, $B'$ lies between $A$ and $C$.
By Theorem~\ref{thm:straight-angle}, $\measuredangle AB'C=\pi$.

By SSS, 
\[\triangle ABC\cong \triangle AB'C.\]
Therefore, $\measuredangle ABC=\pi$.
By Theorem~\ref{thm:straight-angle}, $B$ lies between $A$ and $C$.
\qeds

\begin{thm}{Advanced exercise}\label{ex:SMS}
Consider two triangles $A B C$ and $A' B' C'$.
Let $M$ be the midpoint of $[A B]$, and
$M'$ be the midpoint of $[A' B']$.
Assume $C' A'=C A$, $C' B'= C B$, and $C' M'\z= C M$.
Prove that
\[\triangle A' B' C'\cong\triangle A B C.\]

\end{thm}

{

\begin{wrapfigure}[6]{r}{26mm}
\vskip-0mm
\centering
\includegraphics{mppics/pic-42}
\end{wrapfigure}

\begin{thm}{Exercise}\label{ex:isos-sides}
Let $\triangle A B C$ be an isosceles triangle with the base $[A B]$.
Suppose that point $A'$ lies between $B$ and $C$,
point $B'$ lies between $A$ and $C$,
and $C A'\z=C B'$.
Show that
\begin{enumerate}[(a)]
\item $\triangle A A' C\cong \triangle B B' C$;
\item $\triangle A B B'\cong \triangle B A A'$.
\end{enumerate}

\end{thm}

\begin{thm}{Exercise}\label{ex:ABC-motion}
Let $\triangle ABC$ be a nondegenerate triangle, and
let $f$ be a motion of the plane 
such that 
$$f(A)=A,
\quad 
f(B)=B,
\quad 
\text{and}
\quad
f(C)=C.$$

Show that $f$ is the identity map;
that is, $f(X)=X$ for every point $X$ on the plane.
\end{thm}

}






\section{On side-side-angle and side-angle-angle}

In each of the conditions SAS, ASA, and SSS we specify three corresponding parts of the triangles.
Now, let us discuss other triples of corresponding parts.

The first triple is called {}\emph{side-side-angle}, or briefly SSA;
it specifies two sides and a non-included angle.
This condition is not sufficient for congruence.
In other words, there exist two nondegenerate triangles $ABC$ and $A'B'C'$ such that
\[AB=A'B',\quad BC=B'C',\quad \measuredangle BAC\equiv\pm \measuredangle B'A'C',\]
but $\triangle ABC\not\cong\triangle A'B'C'$ and $AC\ne A'C'$.

\begin{wrapfigure}{r}{35mm}
\vskip-6mm
\centering
\includegraphics{mppics/pic-44}
\end{wrapfigure}

We will not use this negative statement, so there is no need to prove it formally.
You can guess an example from the picture.

The second triple is {}\emph{side-angle-angle}, or briefly SAA;
it specifies one side and two angles one of which is opposite to the side.
This provides a congruence condition; 
that is, if one of the triangles $ABC$ or $A'B'C'$ is nondegenerate and
$AB=A'B'$, $\measuredangle ABC\equiv\pm \measuredangle A'B'C'$, $\measuredangle BCA\z\equiv\pm \measuredangle B'C'A'$,
then $\triangle ABC\cong\triangle A'B'C'$.

The SAA condition will not be used directly in the sequel.
One proof of this condition can be obtained from ASA and the theorem on the sum of angles of a triangle, which is proved below (see~\ref{thm:3sum}). 
For a more direct proof, see Exercise~\ref{ex:SAA}.

Another triple is called {}\emph{angle-angle-angle}, or briefly AAA.
By Axiom~\ref{def:birkhoff-axioms:4}, it is not a congruence condition in the Euclidean plane.
However, in the hyperbolic plane it is; see \ref{thm:AAA}.

\section{Constructions}

The construction problems provide a valuable source of exercises in geometry,
which we will use further in the book.
Chapter~\ref{chap:car} will be devoted to the subject.

{

\begin{wrapfigure}{o}{33mm}
\vskip-8mm
\centering
\includegraphics{mppics/pic-445}
\bigskip
\includegraphics{mppics/pic-446}
\bigskip
\includegraphics{mppics/pic-447}
\bigskip
\includegraphics{mppics/pic-448}
\bigskip
\includegraphics{mppics/pic-449}
\end{wrapfigure}

A \index{ruler-and-compass construction}\emph{ruler-and-compass construction} is a construction with an idealized ruler and compass.
The idealized ruler can be used only to draw a line thru the given two points.
The idealized compass can be used only to draw a circle with a given center and radius.
That is, given three points $A$, $B$, and $O$
we can draw the circle of radius $AB$ centered at~$O$.
We may also mark new points in the plane,
as well as on the constructed lines, circles,
and their intersections (assuming that such points exist).

We may also look at the different sets of construction tools.
For example,
we may only use the ruler,
or we may invent a new tool,
say a tool that produces a midpoint for given two points,
or the reflection of one point across the other.

\begin{thm}{Problem}
Given two lines $\ell$ and $m$ and a point $M$ not on these lines,
construct a line segment with ends on $\ell$ and $m$ that is bisected by~$M$.
\end{thm}

The stages of the following construction are shown in the picture.

\parit{Construction.}
\begin{enumerate}[1.]
\item Draw a circle centered at $M$ that crosses line $\ell$ at two points, say $P$ and $Q$.
\item Draw the line $(PM)$ and mark by $P'$ its other point of intersection with the circle.
Similarly, draw $(QM)$ and mark the other intersection point by $Q'$.
\item Draw the line $\ell'=(P'Q')$ and mark by $A$ its point of intersection with $m$ (if it exists).
\item Draw the line $(MA)$;
mark by $B$ its point of intersection with $\ell$.
Then $[AB]$ is the required segment.
\end{enumerate}

}

It is often required to do \emph{problem analysis};
proving that the result of the construction satisfies the required conditions and also determine when the construction is possible and how many solutions may exist.
Let us do it for the example above.

\parit{Problem analysis.}
Notice that $P'$ and $Q'$ are reflections of $P$ and $Q$ across~$M$.
Therefore, the line $\ell'=(P'Q')$ is the reflection of $\ell=(PQ)$ across~$M$.

Since $A$ lies on $\ell'$,
the point $B$ is the reflection of $A$ across $M$.
Indeed, $B$ is the (necessarily unique) point of intersection of $(MA)$ and $\ell$, and $\ell'$ is the reflection of $\ell$.
By the definition of reflection, $M$ is the midpoint of $[AB]$.

The line $\ell'$ may have no points of intersection with $m$.
In this case, the problem has no solution.
Also, it might happen that $\ell'=m$;
in this case, we have an infinite set of solutions --- any point on $m$ can be marked as~$A$.
In all the remaining cases, we have exactly one solution.
\qeds

\begin{thm}{Exercise}\label{ex:motion}
Suppose $\triangle ABC$ and $\triangle A'B'C'$ are congruent nondegenerate triangles;
in particular, there is a motion $f$ such that
\[f\colon A\mapsto A',\quad f\colon B\mapsto B',\quad\text{and}\quad f\colon C\mapsto C'.\]

Given a point $X$, construct $X'=f(X)$.
\end{thm}

\chapter{Perpendicular lines}\label{chap:perp}

\section{Right, acute and obtuse angles}

\begin{itemize}
\item If $|\measuredangle A O B|=\tfrac\pi2$, we say that $\angle A O B$ is \index{right!angle}\emph{right};
\item If $|\measuredangle A O B|<\tfrac\pi2$, we say that $\angle A O B$ is 
\index{acute}\emph{acute};
\item If $|\measuredangle A O B|>\tfrac\pi2$, we say that $\angle A O B$ is \index{angle!right, acute, and obtuse angles}\index{obtuse}\emph{obtuse}.
\end{itemize}

\begin{wrapfigure}[2]{r}{25mm}
\vskip-25mm
\centering
\includegraphics{mppics/pic-46}
\end{wrapfigure}

Right angles will be marked with a small square, as shown.

If $\angle A O B$ is right, we also may say that $[O A)$ is \index{perpendicular}\emph{perpendicular} to $[O B)$.
This can be written as \index{38@$\perp$}$[O A)\z\perp [O B)$.

Theorem~\ref{thm:straight-angle} allows us to call two lines $(O A)$ and $(O B)$ {}\emph{perpendicular} if $[O A)\z\perp [O B)$; so, we can write $(O A)\z\perp (O B)$.

\begin{thm}{Exercise}\label{ex:acute-obtuce}
Assume that point $O$ lies between $A$ and $B$ and $X\ne O$.
Show that 
$\angle XOA$ is acute if and only if 
$\angle XOB$ is obtuse.
\end{thm}

\section{Perpendicular bisector}

Assume $M$ is the midpoint of the segment $[AB]$;
equivalently, $M\in(A B)$ and $A M \z= M B$.

The line $\ell$ thru $M$ and perpendicular to $(AB)$,
is called the \index{bisector!perpendicular bisector}\index{perpendicular!bisector}\emph{perpendicular bisector} to the segment~$[AB]$. 

\begin{thm}[\abs]{Theorem}\label{thm:perp-bisect}
Given distinct points $A$ and $B$,
all points that are equidistant from $A$ and $B$ and no
others lie on the perpendicular bisector of~$[A B]$.
\end{thm}

\parit{Proof.} Let $M$ be the midpoint of~$[A B]$.
Assume $P A= P B$ and $P\ne M$.
According to SSS (\ref{thm:SSS}),
$\triangle A M P \z\cong\triangle B M P$.
Hence
$\measuredangle A M P=\pm \measuredangle B M P$.
Since $A\not=B$, we have ``$-$'' in this formula.

\begin{wrapfigure}[8]{o}{34mm}
\centering
\includegraphics{mppics/pic-48}
\end{wrapfigure}

Furthermore,
\begin{align*}
\pi
&=
\measuredangle A M B
\equiv
\\
&\equiv\measuredangle A M P+\measuredangle P M B
\equiv
\\
&\equiv
2\cdot \measuredangle A M P.
\end{align*}
That is, $\measuredangle A M P
=
\pm
\tfrac\pi2$. 
Therefore, $P$ lies on the perpendicular bisector.

To prove the converse, 
suppose $P$ 
is any point on the perpendicular bisector of $[A B]$ and $P\z\ne M$.
Then $\measuredangle A M P=\pm \tfrac\pi2$, 
$\measuredangle B M P=\pm \tfrac\pi2$ and
$A M\z=B M$.
By SAS, $\triangle A M P\cong \triangle B M P$;
in particular, $A P\z= B P$.
\qeds

\begin{thm}{Exercise}\label{ex:pbisec-side}
Let $\ell$ be a perpendicular bisector of $[A B]$, and $X$ be an arbitrary point on the plane.

Show that 
$AX<BX$ if and only if $X$ and $A$ lie on the same side from~$\ell$.
\end{thm}

\begin{thm}{Exercise}\label{ex:side-angle}
Let $ABC$ be a nondegenerate triangle.
Show that 
\[AC>BC\iff|\measuredangle ABC|>|\measuredangle CAB|.\] 
\end{thm}

\begin{thm}{Exercise}\label{ex:pbisec-motion}
A motion has a fixed point $F$ and maps a point $X$ to another point $Y$.
Show that $F$ lies on the perpendicular bisector of $[XY]$.
\end{thm}

\section{Uniqueness of a perpendicular}

\begin{thm}[\abs]{Theorem}\label{perp:ex+un}
There is one and only one line that passes thru a given point $P$ and is perpendicular to a given line~$\ell$.
\end{thm}

According to the theorem above, 
there is a unique point $Q\in\ell$ such that $(QP)\perp\ell$.
This point $Q$ is called the \index{footpoint}\emph{footpoint} of $P$ on~$\ell$. 

\parit{Proof.} 
If $P\in\ell$, then both existence and uniqueness follow from Axiom~\ref{def:birkhoff-axioms:2}.

{

\begin{wrapfigure}{o}{30mm}
\vskip-4mm
\centering
\includegraphics{mppics/pic-50}
\end{wrapfigure}

\parit{Existence for $P\not\in\ell$.} 
Let $A$ and $B$ be two distinct points on~$\ell$. 
Choose $P'$ so that $AP'\z=AP$ and $\measuredangle  BAP' \equiv -\measuredangle   BAP$.
According to Axiom~\ref{def:birkhoff-axioms:3}, $\triangle A P' B\z\cong\triangle A P B$.
In particular, $A P= A P'$ and $B P= B P'$.

According to Theorem~\ref{thm:perp-bisect}, $A$ and $B$ lie on the perpendicular bisector of~$[P P']$.
In particular, $(P P')\perp (A B)=\ell$.

}

\parit{Uniqueness for $P\not\in\ell$.} 
From above we can choose a point $P'$ in such a way that $\ell$ forms the perpendicular bisector of~$[PP']$.

Assume $m\perp \ell$ and $m\ni P$.
Then $m$ is a perpendicular bisector of some segment $[Q Q']$ of $\ell$;
in particular, $P Q= P Q'$.

{

\begin{wrapfigure}{r}{37mm}
\centering
\includegraphics{mppics/pic-52}
\end{wrapfigure}

Since $\ell$ is the perpendicular bisector of $[P P']$,
we get that $PQ= P'Q$ and $PQ' \z= P'Q'$.
Therefore, 
$$P' Q=P Q=P Q'= P' Q'.$$
By Theorem~\ref{thm:perp-bisect}, 
$P'$ lies on the perpendicular bisector of $[QQ']$, which is~$m$.
By Axiom~\ref{def:birkhoff-axioms:1}, $m=(P P')$.
\qeds

}

\begin{thm}{Exercise}\label{ex:construction-perpendicular}
Describe a ruler-and-compass construction of a line thru a given point that is perpendicular to a given line.
\end{thm}

\section{Reflection across a line}

Assume the point $P$ and the line $(AB)$ are given.
To find the \index{reflection!across a line}\emph{reflection} $P'$ of $P$ across $(AB)$,
one drops a perpendicular from $P$ onto $(AB)$, 
and continues it to the same distance on the other side.

{

\begin{wrapfigure}{o}{37mm}
\vskip-2mm
\centering
\includegraphics{mppics/pic-54}
\end{wrapfigure}

According to Theorem~\ref{perp:ex+un}, $P'$ is uniquely determined by~$P$.
Note that $P=P'$ if and only if $P\in(AB)$.

\begin{thm}[\abs]{Proposition}\label{prop:reflection}
Assume $P'$ is the reflection of $P$ across $(AB)$.
Then $AP'=AP$, 
and if $A\ne P$, 
then
$\measuredangle BAP'\equiv -\measuredangle BAP$.
\end{thm}

\parit{Proof.} 
If $P\in (AB)$, 
then $P\z=P'$. 
By Corollary~\ref{cor:degenerate=pi}, $\measuredangle BAP\z=0$ or~$\pi$.
Hence the statement follows.

}

If $P\notin (AB)$, then~$P'\ne P$.
By the construction of $P'$, 
the line $(AB)$ is a perpendicular bisector of~$[PP']$.
Therefore, according to Theorem~\ref{thm:perp-bisect}, $AP'\z=AP$ and $BP'\z=BP$.
In particular, 
$\triangle ABP'\z\cong \triangle ABP$.
Therefore, $\measuredangle BAP'=\pm \measuredangle BAP$.

Since $P'\ne P$ and $AP'=AP$,
we get that $\measuredangle BAP'\ne \measuredangle BAP$.
That is, we are left with the case $\measuredangle BAP'\equiv-\measuredangle BAP$.
\qeds

\begin{thm}{Exercise}\label{ex:2-reflections}
Let $X$ and $Y$ be the reflections of $P$ across the lines $(AB)$ and $(BC)$ respectively.
Show that 
$$\measuredangle XBY\equiv 2\cdot \measuredangle ABC.$$

\end{thm}

\begin{thm}[\abs]{Corollary}\label{cor:reflection+angle}
A reflection across a line is a motion of the plane. 
Moreover,
$$\measuredangle Q'P'R'\equiv -\measuredangle QPR$$
if $\triangle P'Q'R'$ is the reflection of $\triangle PQR$.
\end{thm}

\parit{Proof.}
The composition of two reflections across the same line
is the identity map.
In particular, any reflection is a bijection.

Fix a line $(AB)$ and two points $P$ and $Q$;
denote their reflections across $(AB)$ by $P'$ and $Q'$.
Let us show that
$$P'Q'=PQ;\eqlbl{eq:P'Q'=PQ}$$
that is, the reflection is distance-preserving,

\begin{wrapfigure}{r}{33mm}
\centering
\vskip-15mm
\includegraphics{mppics/pic-56}
\end{wrapfigure}

Without loss of generality, we may assume that the points $P$ and $Q$ are distinct from $A$ and~$B$.
By Proposition~\ref{prop:reflection}, we get that
\begin{align*}
\measuredangle BAP'&\equiv -\measuredangle BAP,
&
\measuredangle BAQ'&\equiv -\measuredangle BAQ,
\\
AP'&=AP,
&
AQ'&=AQ.
\end{align*}
It follows that
\[\measuredangle P'AQ'\equiv -\measuredangle PAQ.\eqlbl{eq:P'AQ'=PAQ}\]
By SAS, 
$\triangle P'AQ'\cong\triangle PAQ$
and \ref{eq:P'Q'=PQ} follows.
Moreover, we also get that 
\[\measuredangle AP'Q'\equiv\pm\measuredangle APQ.\]
From \ref{eq:P'AQ'=PAQ} and the theorem on the signs of angles of triangles (\ref{thm:signs-of-triug}) we get
\[\measuredangle AP'Q'\equiv-\measuredangle APQ.\eqlbl{eq:AP'Q'=APQ}\]

Repeating the same argument for a pair of points $P$ and $R$,
we get that
$$\measuredangle AP'R'\equiv-\measuredangle APR.\eqlbl{eq:AP'R'=APR}$$
Subtracting \ref{eq:AP'R'=APR} from \ref{eq:AP'Q'=APQ},
we get that
$$\measuredangle Q'P'R'\equiv-\measuredangle QPR.$$
\qedsf

\section{Direct and indirect motions}
\label{direct motion}

\begin{thm}{Exercise}\label{ex:3-reflections}
Show that every motion of the plane can be presented as a
composition of at most three reflections across lines.
\end{thm}

A motion $X\mapsto X'$ is called \index{direct motion}\emph{direct} if 
$$\measuredangle Q'P'R'= \measuredangle QPR$$ 
for any triangle $PQR$;
if instead we always have 
$$\measuredangle Q'P'R'\equiv -\measuredangle QPR,$$ 
then the motion $f$ is called \index{indirect motion}\emph{indirect}.

By Corollary~\ref{cor:reflection+angle}, any reflection across a line is an indirect motion.
Note that the composition of two reflections is a direct motion.
More generally, the composition of two indirect motions is direct,
the composition of two direct motions is direct,
and the composition of direct and indirect motions is indirect.
In particular, the exercise above implies the following.

\begin{thm}{Proposition}\label{prop:direct-indirect}
Any motion of the plane is either direct or indirect.
\end{thm}

\section{The perpendicular is shortest}
\label{sec:perp<oblique}

\begin{thm}[\abs]{Lemma}\label{lem:perp<oblique}
Assume $Q$ is the footpoint of $P$ on the line~$\ell$.
Then 
$$PX>PQ$$
for every point $X$ on $\ell$ distinct from~$Q$.
\end{thm}

If $P$, $Q$, and $\ell$ are as above, 
then $PQ$ is called the \label{distance!from a point to a line}\index{distance!from a point to a line}\emph{distance from $P$ to~$\ell$}. 

\parit{Proof.}
If $P\in \ell$, 
then the result follows since $PQ=0$.
So, we can assume that $P\notin \ell$.

\begin{wrapfigure}{o}{24mm}
\centering
\vskip-4mm
\includegraphics{mppics/pic-58}
\vskip0mm
\end{wrapfigure}

Let $P'$ be the reflection of $P$ across the line~$\ell$.
Note that $Q$ is the midpoint of $[PP']$
and $\ell$ is the perpendicular bisector of $[PP']$.
Therefore
$$PX=P'X
\quad
\text{and}
\quad
PQ=P'Q=\tfrac12\cdot PP'.$$

The line $\ell$ meets $[PP']$ only at the point $Q$.
Therefore, $X\notin [PP']$.
By the triangle inequality and Corollary~\ref{cor:degenerate-trig}, we have 
$$PX+P'X>PP',$$
and hence the result: $PX>PQ$.
\qeds

\begin{thm}{Classroom exercise}\label{ex:right-acute}
Show that
\begin{enumerate}[(a)]
\item If $\triangle ABC$ has a right angle at $C$, then its other angles are acute.
\item If $\triangle ABC$ has an obtuse angle at $C$, then its other angles are acute.
\end{enumerate}
\end{thm}

This exercise is closely related to Proposition~\ref{prop:2sum}.

If a triangle has an obtuse or right angle, then it is called \index{triangle!right, acute, and obtuse triangle}\index{obtuse}\emph{obtuse} or \index{triangle!right triangle}\index{right!triangle}\emph{right}, respectively.
Otherwise, if all its angles are acute, it is called an \index{acute}\emph{acute triangle}.

\begin{thm}{Exercise}\label{ex:obtuce}
Suppose that point $X$ lies between $V$ and $W$.
Show that for every point $P$, we have $PX < PV$ or $PX < PW$.
\end{thm}

\begin{thm}{Exercise}\label{ex:PMQ}
Assume that points $P$ and $Q$ lie on one side of line~$\ell$.
Let $M\in \ell$ be the point that minimizes the sum $PM+MQ$.
Show that $M\in (P'Q)$, where $P'$ is the reflection of $P$ across $\ell$.
\end{thm}

\section{Circles}

Recall that a circle with radius $r$ and center $O$ is the set of all points on distance $r$ from $O$.
We say that a point $P$ lies \index{inside!a circle}\emph{inside} of the circle if $OP<r$; 
if $OP>r$, we say that $P$ lies \index{outside a circle}\emph{outside} of the circle.
\label{def:circle}

\begin{thm}{Exercise}\label{ex:inside-outside}
Let $\Gamma$ be a circle and $P\notin \Gamma$.
Assume a line $\ell$ is passing thru the point $P$
and intersects $\Gamma$ at two distinct points, $X$ and~$Y$.
Show that $P$ is inside $\Gamma$ if and only if $P$ lies between $X$ and~$Y$.
\end{thm}

A segment between two points on a circle is called a \index{chord}\emph{chord} of the circle.
A chord passing thru the center of the circle is called its \index{diameter}\emph{diameter}.

\begin{thm}{Exercise}\label{ex:chord-perp}
Assume two distinct circles $\Gamma$ and $\Gamma'$ have a common chord~$[A B]$.
Show that the line between centers of $\Gamma$ and $\Gamma'$ forms a perpendicular bisector of~$[A B]$.
\end{thm}

\begin{thm}{Exercise}\label{ex:center}
Describe a ruler-and-compass construction of a center
of a given circle.
\end{thm}

{

\begin{wrapfigure}[4]{r}{54mm}
\vskip-8mm
\centering
\includegraphics{mppics/pic-60}
\end{wrapfigure}

\begin{thm}[\abs]{Lemma}\label{lem:line-circle}
A line and a circle can have at most two points of intersection.
\end{thm}

\parit{Proof.} Assume $A$, $B$, and $C$ are distinct points that lie on a line $\ell$ and a circle $\Gamma$ with the center~$O$.
Then $OA=OB=OC$; in particular, $O$ lies on the perpendicular bisectors 
$m$ and $n$ to $[A B]$ and $[B C]$ respectively.
Note that the midpoints of $[AB]$ and $[BC]$ are distinct.
Therefore, $m$ and $n$ are distinct.
The latter contradicts the uniqueness of the perpendicular (Theorem~\ref{perp:ex+un}).
\qeds

}

\begin{thm}{Exercise}\label{ex:two-circ}
Show that two distinct circles can have at most two points of intersection.
\end{thm}

As a consequence of the above lemma, 
a line $\ell$ and a circle $\Gamma$ might have 2, 1, or 0 points of intersections.
In the first two cases, the line is called \index{secant line}\emph{secant} or \index{tangent!line}\emph{tangent} respectively;
if $P$ is the only point of intersection of $\ell$ and $\Gamma$,
we say that \textit{$\ell$ is tangent to $\Gamma$ at $P$}. 

Similarly, according to Exercise~\ref{ex:two-circ},
two distinct circles might have 2, 1, or 0 points of intersections.
If $P$ is the only point of intersection of circles $\Gamma$ and $\Gamma'$,
we say that \index{tangent!circles}\emph{$\Gamma$ is tangent to $\Gamma$ at $P$}; we also assume that a circle is tangent to itself at each of its points.

\begin{thm}[\abs]{Lemma}\label{lem:tangent}
Let $\ell$ be a line and $\Gamma$ be a circle centered at~$O$.
Assume $P$ is a common point of $\ell$ and~$\Gamma$. 
Then $\ell$ is tangent to $\Gamma$ at $P$ if and only if $(PO)\perp \ell$.
\end{thm}

\parit{Proof.}
Let $Q$ be the footpoint of $O$ on~$\ell$.

Assume~$P\ne Q$.
Let $P'$ be the reflection of $P$ across~$(OQ)$.
Note that $P'\in\ell$ and $(OQ)$ is the perpendicular bisector of~$[PP']$.
Therefore, $OP=OP'$.
Hence $P,P'\in \Gamma\cap \ell$;
that is, $\ell$ is secant to~$\Gamma$.

If $P=Q$, 
then according to Lemma~\ref{lem:perp<oblique},
$OP<OX$ for each point $X\in \ell$ distinct from~$P$.
Hence $P$ is the only point at the intersection $\Gamma\cap\ell$;
that is, $\ell$ is tangent to $\Gamma$ at~$P$. 
\qeds

\begin{thm}{Exercise}\label{ex:tangent-circles}
Let $\Gamma$ and $\Gamma'$ be two circles with centers at $O$ and $O'$, and radiuses $r$ and~$r'$, respectively.
Show that $\Gamma$ is tangent to $\Gamma'$ if and only if one of the following conditions holds:
\begin{enumerate}[(a)]
\item\label{ex:tangent-circles:a} The circles meet at a point~$P$ that lies on one line with $O$ and~$O'$.
\item\label{ex:tangent-circles:b} The circles meet at a point~$P$ and share a common tangent line at this point.
\item\label{ex:tangent-circles-2} $OO' = r + r'$ or $OO' = |r - r'|$.
\end{enumerate}

\end{thm}

{

\begin{wrapfigure}{r}{20mm}
\vskip-10mm
\centering
\includegraphics{mppics/pic-62}
\end{wrapfigure}

\begin{thm}{Exercise}\label{ex:tangent-circles-3}
Assume three circles have two points in common.
Prove that their centers lie on one line.
\end{thm}

}

\begin{thm}{Exercise}\label{ex:tangent}
Describe a ruler-and-compass construction of lines tangent to a given circle that pass thru a given point.
\end{thm}

\begin{thm}{Exercise}\label{ex:tangent-circle}
Given two circles $\Gamma_1$ and $\Gamma_2$ and a segment $[AB]$
describe a ruler-and-compass construction of a circle with the radius $AB$
that is tangent to each circle $\Gamma_1$ and~$\Gamma_2$.
\end{thm}

\chapter[Similar triangles]{Similar triangles}\label{chap:parallel}

Two triangles $A' B' C'$ and $A B C$ are called
\index{triangle!similar triangles}\index{similar triangles}\emph{similar} (briefly \index{30@$\sim$}$\triangle A' B' C'\z\sim\triangle A B C$) if (1) their sides are proportional; 
that is, 
$$A' B'
=
k\cdot A B,
\quad
B' C'=k\cdot B C
\quad
\text{and}
\quad
C' A'
=
k\cdot C A
\eqlbl{dist}
$$
for some coefficient $k>0$ (which is called the \index{similarity coefficient}\emph{similarity coefficient}), and (2) the corresponding angles are equal up to sign:
$$
\begin{aligned}
\measuredangle A' B' C'&=\pm\measuredangle A B C,
\\
\measuredangle B' C' A'&=\pm\measuredangle B C A,
\\
\measuredangle C' A' B'&=\pm\measuredangle CAB.
\end{aligned}
\eqlbl{angles}
$$

\parbf{Remarks.}
\begin{itemize}
\item According to \ref{thm:signs-of-triug}, in the above three equalities, the signs can be assumed to be the same.

\item The equations in \ref{dist} can be rewritten as
\[\frac{A'B'}{AB}=\frac{B'C'}{BC}=\frac{C'A'}{CA}.\]

\item If $\triangle A' B' C'\sim\triangle A B C$ with the coefficient $k=1$, 
 then $\triangle A' B' C'\z\cong\triangle A B C$.

\item Notice that ``$\sim$'' is an
\index{equivalence relation}\emph{equivalence relation}.
That is, 
\begin{enumerate}[(i)]
\item $\triangle A B C\sim\triangle A B C$
for every $\triangle A B C$.
\item If $\triangle A' B' C'\sim\triangle A B C$, then
$$\triangle A B C\sim\triangle A' B' C'.$$
\item If $\triangle A'' B'' C''\sim\triangle A' B' C'$ and $\triangle A' B' C'\z\sim\triangle A B C$, then 
$$\triangle A'' B'' C''\sim\triangle A B C.$$
\end{enumerate}
\end{itemize}

\section{Similarity conditions}

Using the new notation ``$\sim$'', we can reformulate Axiom~\ref{def:birkhoff-axioms:4}:

\begin{thm}{Reformulation of Axiom~\ref{def:birkhoff-axioms:4}}
If for  
$\triangle ABC$, 
$\triangle AB'C'$,
and $k>0$ we have
$B'\in [AB)$,
$C'\in [AC)$,
$AB'=k\cdot AB$ and
$AC'=k\cdot AC$,
then $\triangle ABC\sim\triangle AB'C'$.
\end{thm}

In other words, Axiom~\ref{def:birkhoff-axioms:4} provides 
a condition which guarantees that two triangles are similar.
Let us formulate three more such {}\emph{similarity conditions}.

\begin{thm}{Similarity conditions}\label{prop:sim}
Two triangles 
$\triangle ABC$ and $\triangle A'B'C'$
are similar if one of the following conditions holds:

(SAS)\index{SAS similarity condition}
For some constant $k>0$ we have
\[A B=k\cdot A' B',
\quad 
A C=k\cdot A' C',
\quad
\text{and}
\quad 
\measuredangle B A C=\pm\measuredangle B' A' C'.\]

(AA)\index{AA similarity condition} The triangle $A' B' C'$ is nondegenerate
and 
$$\measuredangle A B C
=
\pm\measuredangle A' B' C',
\quad 
\measuredangle B A C
=
\pm\measuredangle B' A' C'.$$

(SSS)\index{SSS similarity condition} For some constant $k>0$ we have
$$A B=k\cdot A' B',
\quad
A C=k\cdot A' C',
\quad
CB=k\cdot C'B'.$$

\end{thm}

Each of these conditions is proved by applying Axiom~\ref{def:birkhoff-axioms:4} with the SAS, ASA, and SSS congruence conditions respectively
(see Axiom~\ref{def:birkhoff-axioms:3} and the conditions \ref{thm:ASA}, \ref{thm:SSS}).

\parit{Proof.}
Let $k=\tfrac{AB}{A'B'}$.
Choose points $B''\in [A'B')$ and $C''\in [A'C')$,
so that $A'B''=k\cdot A'B'$ and $A'C''=k\cdot A'C'$.
By Axiom~\ref{def:birkhoff-axioms:4},
$\triangle A'B'C'\z\sim \triangle A'B''C''$.

Applying the SAS, ASA, or SSS congruence condition, depending on the case, 
we get that $\triangle A'B''C''\cong \triangle ABC$.
Hence the result.
\qeds

\begin{thm}{Exercise}\label{ex:mid-triangle}
Let $A'$, $B'$, and $C'$ be the midpoints of sides $[BC]$, $[CA]$, and $[AB]$ of $\triangle ABC$.
Show that $\triangle A'B'C'\sim\triangle ABC$ and find the similarity coefficient.
\end{thm}

\begin{thm}{Exercise}\label{ex:k*triangle}
Let $O$, $A$, $B$, $C$, $A'$, $B'$, and $C'$ be distict points.
Suppose that $A'\in [OA)$, $B'\in[OB)$, $C'\in [OC)$, $OA'=k\cdot OA$, $OB'=k\cdot OB$, and $OC'=k\cdot OC$ for some $k>0$.
Show that $\triangle A'B'C'\sim\triangle ABC$.
\end{thm}

A bijection $X\leftrightarrow X'$ from a plane to itself is called an \index{angle-preserving transformation}\emph{angle-preserving transformation} if 
\[\measuredangle ABC= \measuredangle A'B'C'\]
for every $\triangle ABC$ and its image $\triangle A'B'C'$.

(The term \index{transformation}\emph{transformation} is used for a bijection of space to itself that preserves a specified geometric structure.
For example, {}\emph{motions} are \textit{distance-preserving transformations}.)

\begin{thm}{Exercise}\label{ex:angle-preserving-euclid}
Show that every angle-preserving transformation of the plane multiplies all distances by a fixed constant.
\end{thm}

\section{The Pythagorean theorem}\index{Pythagorean theorem}

Recall that triangle is called \index{right!triangle}\index{triangle!right triangle}\emph{right} if one of its angles is right.
The side opposite the right angle is called the \index{hypotenuse}\emph{hypotenuse}. 
The sides adjacent to the right angle are called \index{leg}\emph{legs}.

\begin{thm}{Theorem}\label{thm:pyth}
Assume $\triangle ABC$ has a right angle at~$C$.
Then
$$AC^2+BC^2=AB^2.$$ 

\end{thm}

\parit{Proof.}
Let $D$ be the footpoint of $C$ on~$(AB)$.

\begin{wrapfigure}[4]{r}{40mm}
\vskip-2mm
\centering
\includegraphics{mppics/pic-66}
\end{wrapfigure}

According to Lemma~\ref{lem:perp<oblique},
\begin{align*}
AD&<AC<AB,
\shortintertext{and}
BD&<BC<AB.
\end{align*}
Therefore, $D$ lies between $A$ and $B$;
in particular, 
$$AD+BD=AB.\eqlbl{AD+BD=AB}$$

By the AA similarity condition, we have
$$\triangle ADC\sim\triangle ACB\sim \triangle CDB.$$
In particular,
$$
\frac{A D}{A C}=\frac{A C}{A B}
\quad
\text{and}
\quad
\frac{B D}{B C}=\frac{B C}{B A}.
\eqlbl{BCD}$$

Let us rewrite the two identities in \ref{BCD}:
\begin{align*}
AC^2=AB\cdot AD
\quad
\text{and}
\quad
BC^2=AB\cdot B D.
\end{align*}
Summing up these two identities and applying \ref{AD+BD=AB}, we get that
$$AC^2 +BC^2=AB\cdot (AD+ B D)=AB^2.$$
\qedsf

The idea in this proof appears in the Elements \cite[X.33]{euclid},
but the proof given there \cite[I.47]{euclid} is different; 
it uses the area method discussed in Chapter~\ref{chap:area}.

\begin{thm}{Exercise}\label{ex:pyth}
Assume $A$, $B$, $C$, and $D$ are as in the proof above.
Show that 
$CD^2=AD\cdot BD$.

\end{thm}

The following exercise is the converse to the Pythagorean theorem.

\begin{thm}{Exercise}\label{ex:pyth-conv}
Let $ABC$ be a triangle such that
$AC^2+BC^2=AB^2$.
Show that the angle at $C$ is right.
\end{thm}

\section{The method of similar triangles}

The proof of the Pythagorean theorem given above uses the {}\emph{method of similar triangles}.
To apply this method, one has to search for pairs of similar triangles and then use the proportionality of corresponding sides and/or equalities of corresponding angles.
Finding such pairs might be tricky at first. 

{

\begin{wrapfigure}{r}{25mm}
\vskip-6mm
\centering
\includegraphics{mppics/pic-68}
\end{wrapfigure}

\begin{thm}{Classroom exercise}\label{ex:two-pairs-sim}
Let $ABC$ be a nondegenerate triangle and the points $X$, $Y$, and $Z$ as shown in the picture.
Assume $\measuredangle CAY\z\equiv\measuredangle XBC$.
Find four pairs of similar triangles with these six points as the vertices
and prove their similarity.
\end{thm}

}

\begin{thm}{Exercise}\label{ex:ABC+D}
Suppose $\triangle ABC$ is nondegenerate.
Assume $D\in [AB]$ and $\measuredangle CAB=\measuredangle BCD$.
Show that $CD=\tfrac{b\cdot a}c$, where $a=BC$, $b=CA$, and $c=AB$. 
\end{thm}

{

\begin{wrapfigure}{r}{25mm}
\vskip-4mm
\centering
\includegraphics{mppics/pic-69}
\end{wrapfigure}

\begin{thm}{Exercise}\label{ex:right-perp-bi}
Assume that $\triangle ABC$ has a right angle at~$C$.
Let $M$ be the midpoint of $[AB]$.
Denote by $P$ and $Q$ the intersections of the perpendicular bisector of $[AB]$ with the half-lines
$[BC)$ and $[AC)$, respectively.
Show that
\[
4\cdot MP\cdot MQ = AB^2.
\]

\end{thm}

}


\section{Ptolemy's inequality}

A \index{quadrangle}\emph{quadrangle} is defined as an ordered quadruple of distinct points in the plane.
These four points are called \index{vertex!of quadrangle}\emph{vertices}.
The quadrangle $ABCD$ will also be denoted by \index{25@$\square$}$\square ABCD$.

Given a quadrangle $ABCD$,
the four segments $[AB]$, $[BC]$, $[CD]$, and $[DA]$ are called \index{side!of quadrangle}\emph{sides of $\square ABCD$};
the remaining two segments $[AC]$ and $[BD]$ are called \index{diagonal!of quadrangle}\emph{diagonals of $\square ABCD$}.

\begin{thm}{Ptolemy's inequality}\label{ptolemy-inq}
In every quadrangle, the product of diagonals cannot exceed the sum of the products of its opposite sides;
that~is, 
\[AC\cdot BD\le AB\cdot CD+ BC\cdot DA\]
for every $\square ABCD$.
\end{thm}

We will present a classical proof of this inequality using the method of similar triangles with additional construction.
This proof is given as an illustration --- it will not be used further in the sequel.

{

\begin{wrapfigure}{o}{33mm}
\vskip-3mm
\centering
\includegraphics{mppics/pic-70}
\end{wrapfigure}

\parit{Proof.}
Consider the half-line $[AX)$ such that $\measuredangle BAX=\measuredangle CAD$.
In this case, $\measuredangle XAD\z=\measuredangle BAC$ since adding $\measuredangle BAX$ or $\measuredangle CAD$ to the corresponding sides produces $\measuredangle BAD$.
We can assume that
\[AX=\frac{AB}{AC}\cdot AD.\]

}

In this case, we have
\begin{align*}\frac{AX}{AD}&=\frac{AB}{AC},
&
\frac{AX}{AB}&=\frac{AD}{AC}.\\
\shortintertext{Hence}
\triangle BAX&\sim \triangle CAD,
&
\triangle XAD&\sim\triangle BAC.\\
\shortintertext{Therefore}
\frac{BX}{CD}&=\frac{AB}{AC},
&
\frac{XD}{BC}&=\frac{AD}{AC},
\shortintertext{or}
AC\cdot BX&=AB\cdot CD,
&
AC\cdot XD&=BC\cdot AD.
\end{align*}
Adding these two equalities we get 
\[AC\cdot(BX+XD)=AB\cdot CD+BC\cdot AD.\]
It remains to apply the triangle inequality, $BD\le BX+XD$.
\qeds

Using this proof together with \ref{prop:inscribed-quadrangle}, one can show that the equality holds only if the vertices $A$, $B$, $C$, and $D$ appear on a line or a circle in the same cyclic order;
see also \ref{ptolemy-id} for another proof of the equality case.
Exercise~\ref{ex:ptolemy} suggests another proof of Ptolemy's inequality using complex coordinates.

\chapter{Parallel lines}\label{chap:angle-sum}

{

\begin{wrapfigure}{o}{19mm}
\vskip-4mm
\centering
\includegraphics{mppics/pic-72}
\end{wrapfigure}

Recall that in consequence of Axiom~\ref{def:birkhoff-axioms:1},
any two distinct lines $\ell$ and $m$ have either one point in common or none; see \ref{lem:line-line}.
In the first case they are \index{intersecting lines}\emph{intersecting} (briefly $\ell\nparallel m$); 
in the second case, $\ell$ and $m$ are said to be \index{parallel!lines}\emph{parallel} (briefly, \index{36@\hskip.4mm$\parallel$\hskip.4mm, \hskip.4mm$\nparallel$\hskip.4mm}$\ell\parallel m$);
in addition, a line is always regarded as parallel to itself.

}

To emphasize that two lines in a picture are parallel we will mark them with arrows of the same type.

\begin{thm}[\abs]{Proposition}\label{prop:perp-perp} Let $\ell$, $m$, and $n$ be three lines.
Assume that $n\perp m$ and $m\perp \ell$.
Then $\ell\parallel n$. 
\end{thm}

\parit{Proof.}
Assume the contrary; 
that is, $\ell\nparallel n$.
Then there is a point, say $Z$, of intersection of $\ell$ and~$n$.
Then by Theorem~\ref{perp:ex+un},
$\ell=n$.
Since every line is parallel to itself, we have that $\ell\parallel n$ --- a contradiction.
\qeds

\begin{thm}{Theorem}\label{thm:parallel}
For every point $P$ and every line $\ell$,
there is a unique line $m$
that passes thru $P$ and is parallel to~$\ell$.
\end{thm}

The above theorem has two parts, existence and uniqueness.
In the proof of uniqueness, we will use the method of similar triangles.

\parit{Proof; existence.} 
Apply Theorem~\ref{perp:ex+un} two times,
first to construct the line $n$ thru $P$ that is perpendicular to $\ell$,
and second to construct the line $n$ thru $P$ that is perpendicular to~$m$.
Then apply Proposition~\ref{prop:perp-perp}.

\parit{Uniqueness.}
If $P\in\ell$, then $m=\ell$ by the definition of parallel lines.
Now, assume $P\notin\ell$.

Let us construct the lines $n\ni P$ and $m\ni P$ as in the proof of existence; so $n\perp \ell$, $m\perp n$, and $m\parallel \ell$.

Assume there is another line $s\ni P$ parallel to~$\ell$.
Choose a point $Q\in s$ that lies with $\ell$ on the same side as~$m$.
Let $R$ be the footpoint of $Q$ on~$n$.

\begin{figure}[!ht]
\centering
\includegraphics{mppics/pic-74}
\end{figure}

Let $D$ be the point of intersection of $n$ and~$\ell$.
According to Proposition~\ref{prop:perp-perp} $(QR)\parallel m$. 
Therefore, $Q$, $R$, and $\ell$ lie on the same side of~$m$. 
In particular, $R\in [P D)$.

Choose $Z\in [P Q)$ such that 
$$\frac{PZ}{PQ}=\frac{PD}{PR};$$
it exists by Exercise~\ref{ex:trig==}.
By SAS similarity condition (or equivalently by Axiom~\ref{def:birkhoff-axioms:4})
we have that $\triangle RPQ\sim \triangle DPZ$;
therefore $(Z D)\perp(P D)$.
It follows that $Z$ lies on $\ell$ and $s$ --- a contradiction.\qeds

\begin{thm}{Corollary}\label{cor:parallel-1}
Assume $\ell$, $m$, and $n$ are lines
such that $\ell\parallel m$ and $m\parallel n$.
Then $\ell\parallel n$.
\end{thm}

\parit{Proof.}
Assume the contrary; that is, $\ell\nparallel n$.
Then there is a point $P\in \ell\cap n$.
By Theorem~\ref{thm:parallel},
$n=\ell$ --- a contradiction.
\qeds

By the definition, we have that $\ell\parallel m$ if and only if $m\z\parallel \ell$.
Therefore, according to the above corollary, ``$\parallel$'' is an 
\index{equivalence relation}\emph{equivalence relation}.
That is, for any lines $\ell$, $m$, and $n$ the following conditions hold:
\begin{enumerate}[(i)]
\item $\ell\parallel \ell$;
\item if $\ell\parallel m$, then $m\parallel \ell$;
\item if $\ell\parallel m$ and $m\parallel n$, then 
$\ell\parallel n$.
\end{enumerate}

\begin{thm}{Classroom exercise}\label{ex:perp-perp}
Let $k$, $\ell$, $m$, and $n$ be lines such that $k\perp \ell$, $\ell\perp m$, and $m\perp n$.
Show that $k\nparallel n$.
\end{thm}

\begin{thm}{Exercise}\label{ex:construction-parallel}
Perform a ruler-and-compass construction of a line thru a given point that is parallel to a given line.
\end{thm}

\section{Reflection across a point}

Recall that if $O$ is the midpoint of the line segment $[XX']$,
then we say that $X'$ is a \index{reflection across a point}\emph{reflection} of $X$ across a point $O$.
In addition, we assume that $O'=O$; that is, $O$ is a reflection of itself across itself.

The following statement is a refinement of Proposition~\ref{prop:point-reflection}.

\begin{thm}[\abs]{Proposition}\label{prop:point-reflection+}
A reflection across a point is a direct motion.
\end{thm}

\parit{Proof.}
Choose a point $O$.
By Proposition~\ref{prop:point-reflection}, the reflection across $O$ is a motion.
Note that any angle $\angle XOY$ is vertical to its reflection $\angle X'OY'$.
By \ref{prop:vert}, $\measuredangle XOY\z=\measuredangle X'OY'$.
Therefore, the reflection cannot be indirect.

By \ref{prop:direct-indirect}, any motion is either direct or indirect.
Therefore, the reflection across $O$ must be direct.
\qeds

\begin{thm}{Exercise}
Suppose $\angle AOB$ is right.
Show that the composition of reflections across the lines $(OA)$ and $(OB)$ is a reflection across $O$.

Use this statement and Corollary~\ref{cor:reflection+angle} to build another proof of~\ref{prop:point-reflection+}.
\end{thm}

{

\begin{wrapfigure}{r}{23mm}
\vskip-6mm
\centering
\includegraphics{mppics/pic-78}
\end{wrapfigure}

\begin{thm}{Theorem}\label{thm:parallel-point-reflection}
Let $\ell$ and $m$ be lines passing thru points $P$ and $Q$, respectively.
Suppose that $O$ is the midpoint of $[PQ]$.
Then $m \parallel \ell$ if and only if $m$ is the reflection of $\ell$ across $O$.
\end{thm}

}

\parit{Proof; ``if'' part.}
Assume $m$ is a reflection of $\ell$ across $O$.
Suppose $\ell\nparallel m$; that is $\ell$ and $m$ intersect at a single point $Z$.
Denote by $Z'$ be the reflection of $Z$ across $O$.

\begin{figure}[!ht]
\centering
\includegraphics{mppics/pic-80}
\end{figure}

The point $Z'$ lies on both lines $\ell$ and $m$.
It follows that $Z'=Z$ or equivalently $Z=O$.
In this case, $O\in \ell$ and therefore the reflection of $\ell$ across $O$ is $\ell$ itself;
that is, $\ell=m$ and in particular $\ell\parallel m$ --- a contradiction. 

\parit{``Only-if'' part.}
Let $\ell'$ be the reflection of $\ell$ across $O$.
According to the ``if'' part of the theorem, $\ell'\parallel \ell$.
Both lines $\ell'$ and $m$ pass thru $P$.
By uniqueness of parallel lines (\ref{thm:parallel}), if $m\parallel \ell$, then $\ell'=m$; whence the statement follows.
\qeds

\pagebreak

\section{The transversal property}

{

\begin{wrapfigure}{r}{25mm}
\centering
\vskip-10mm
\includegraphics{mppics/pic-82}
\end{wrapfigure}

If the line $t$ intersects each line $\ell$ and $m$ at one point, then we say that $t$ is a \index{transversal}\emph{transversal} to $\ell$ and~$m$.
For example, in the picture, line $(CB)$ is a transversal
to $(AB)$ and~$(CD)$.

}

\begin{thm}{Transversal property}\label{thm:parallel-2} 
Given points $A\ne B\ne C\ne D$,
we have $(AB)\parallel(C D)$ if and only if
$$2\cdot(\measuredangle A B C+\measuredangle B C D)\equiv 0.
\eqlbl{A B C + B C D}$$ 
Equivalently, 
$$\measuredangle A B C+\measuredangle B C D
\equiv 
0
\quad
\text{or}
\quad
\measuredangle A B C+\measuredangle B C D
\equiv
\pi.$$

Moreover, if $(AB)\ne(C D)$, then in the first case, 
$A$ and $D$ lie on opposite sides of $(BC)$;
in the second case,
$A$ and $D$ lie on the same sides of~$(BC)$.
\end{thm}

\parit{Proof; ``only-if'' part.}
Denote by $O$ the midpoint of $[BC]$.

Assume $(AB)\parallel(C D)$.
According to Theorem~\ref{thm:parallel-point-reflection},
$(CD)$ is a reflection of $(AB)$ across $O$.

\begin{wrapfigure}{r}{31mm}
\vskip-4mm
\centering
\includegraphics{mppics/pic-84}
\end{wrapfigure}

Let $A'$ be the reflection of $A$ across $O$.
Then $A'\in (CD)$ and by Proposition~\ref{prop:point-reflection+} we have that
\[\measuredangle ABO=\measuredangle A'CO.\eqlbl{A B O = A' C O}\]
Note that 
\[\measuredangle ABO\equiv\measuredangle ABC,
\qquad
\measuredangle A'CO\equiv-\measuredangle BCA'.\eqlbl{eq:A B O= A B C}
\]
Since $A'$, $C$, and $D$ lie on one line, Exercise~\ref{ex:ABCO-line} implies that 
\[2\cdot \measuredangle BCD\equiv 2\cdot \measuredangle BCA'.\eqlbl{eq:2BCD= 2BCA'}\]
Finally, \ref{A B O = A' C O}, \ref{eq:A B O= A B C}, and \ref{eq:2BCD= 2BCA'} imply \ref{A B C + B C D}.
\qeds

\parit{``If''-part.}
By Theorem~\ref{thm:parallel} there is a unique line $(CD)$ thru $C$ that is parallel to $(AB)$.
From the ``only-if'' part we know that  \ref{A B C + B C D} holds.

On the other hand, there is a \textit{unique} line $(CD)$ such that \ref{A B C + B C D} holds.
Indeed, suppose there are two such lines $(CD)$ and $(CD')$, then
$$2\cdot(\measuredangle A B C+\measuredangle B C D)\equiv 2\cdot(\measuredangle A B C+\measuredangle B C D')\equiv0.
$$ 
Therefore 
$2\cdot\measuredangle B C D\equiv 2\cdot\measuredangle B C D'$
and by Exercise~\ref{ex:ABCO-line},  $D'\in (CD)$, or equivalently the line $(CD)$ coincides with $(CD')$.

Therefore if \ref{A B C + B C D} holds, then $(CD)\parallel (AB)$.

\parit{Last statement.}
If $(AB)\ne(C D)$ and $A$ and $D$ lie on the opposite sides of $(BC)$, then $\angle ABC$ and $\angle BCD$ have opposite signs.
Therefore
\[-\pi\z<\measuredangle A B C+\measuredangle B C D<\pi.\]
Applying \ref{A B C + B C D}, we get $\measuredangle A B C+\measuredangle B C D=0$.

Similarly, if $A$ and $D$ lie on the same side of $(BC)$,
then $\angle ABC$ and $\angle BCD$ have the same sign.
Therefore
\[0<|\measuredangle A B C+\measuredangle B C D|<2\cdot\pi,\]
and \ref{A B C + B C D} implies that $\measuredangle A B C+\measuredangle B C D\z\equiv\pi$.
\qeds

\begin{thm}{Exercise}\label{ex:parallel-angles}
Suppose that $(AB) \parallel (A'B')$ and $(BC) \parallel (B'C')$.
Show that
$2\cdot \measuredangle ABC \equiv 2\cdot \measuredangle A'B'C'$.
Conclude that either
\[
\measuredangle ABC \equiv \measuredangle A'B'C' \quad \text{or} \quad \measuredangle ABC \equiv \pi + \measuredangle A'B'C',
\]
and draw a pair of angles to illustrate each case.
\end{thm}

\begin{thm}{Exercise}\label{ex:smililar+parallel}
Let $\triangle ABC$ be a nondegenerate triangle, and suppose $P$ lies between $A$ and $B$.
Suppose that a line $\ell$ passes thru $P$ and is parallel to $(AC)$.
Show that $\ell$ crosses the side $[BC]$ at another point, say $Q$, and 
\[\triangle ABC\sim\triangle PBQ.\]
In particular, $\tfrac{PB}{AB}=\tfrac{QB}{CB}$.
\end{thm} 

\begin{thm}{Exercise}\label{ex:trisection}
Trisect a given segment with a ruler and a compass.
\end{thm}

\section{Angles of triangles}

\begin{thm}{Theorem}\label{thm:3sum}
In each $\triangle A B C$, we have
$$\measuredangle A B C+ \measuredangle B C A + \measuredangle C A B \equiv \pi.$$

\end{thm}

\parit{Proof.} 
If $\triangle A B C$ is degenerate, then the equality follows from Corollary~\ref{cor:degenerate=pi}.
Now, assume that $\triangle A B C$ is nondegenerate.

Let $X$ be the reflection of $C$ across the midpoint $M$ of $[AB]$.
By Proposition~\ref{prop:point-reflection+}
$\measuredangle BAX\z=\measuredangle ABC$.
Note that $(AX)$ is a reflection of $(CB)$ across $M$;
therefore by Theorem~\ref{thm:parallel-point-reflection}, $(AX)\z\parallel (CB)$.

\begin{wrapfigure}{o}{23mm}
\centering
\vskip-0mm
\includegraphics{mppics/pic-86}
\end{wrapfigure}

Since $[BM]$ and $[MX]$ do not intersect $(CA)$,
the points $B$, $M$, and $X$ lie on the same side of $(CA)$.
Applying the transversal property for the transversal $(CA)$ to $(AX)$ and $(CB)$, we get that 
\[\measuredangle BCA+\measuredangle CAX\equiv \pi.\eqlbl{eq:ABC+CAB}\]

Since $\measuredangle BAX=\measuredangle ABC$,
we have 
\[\measuredangle CAX\equiv\measuredangle CAB+\measuredangle ABC.\]
The latter identity and \ref{eq:ABC+CAB} imply the theorem.\qeds

\begin{thm}{Classroom exercise}\label{ex:|3sum|}
Show that 
$$|\measuredangle A B C|+ |\measuredangle B C A| + |\measuredangle C A B| = \pi$$
for every $\triangle ABC$.
\end{thm} 

{

\begin{wrapfigure}{r}{30mm}
\vskip-0mm
\centering
\includegraphics{mppics/pic-88}
\vskip4mm
\includegraphics{mppics/pic-90}
\end{wrapfigure}

\begin{thm}{Exercise}\label{ex:pent}
Let $\triangle ABC$ be a nondegenerate triangle.
Assume there is a point $D\in [BC]$ 
such that 
\[\measuredangle BAD=\measuredangle DAC,
\quad
BA=AD=DC.\]
Show that $\measuredangle CBA= \measuredangle BAC$ and $AC=BC$.
\end{thm}

\begin{thm}{Exercise}\label{ex:right-isos}
Let $\triangle ABC$ be an isosceles nondegenerate triangle with the base~$[AC]$.
Suppose $D$ is a reflection of $A$ across $B$.
Show that $\angle ACD$ is right.
\end{thm}


}

\begin{thm}{Exercise}\label{ex:quadrangle}
Show that for every quadrangle $ABCD$, we have
$$\measuredangle ABC+\measuredangle BCD+\measuredangle CDA+\measuredangle DAB\equiv 0.$$

\end{thm}

\section{Parallelograms}

{

\begin{wrapfigure}{r}{26mm}
\vskip-10mm
\centering
\includegraphics{mppics/pic-94}
\end{wrapfigure}

A quadrangle $ABCD$ in the Euclidean plane is called \index{quadrangle!degenerate quadrangle}\index{degenerate!quadrangle}\emph{nondegenerate} if no three points from $A,B,C,D$ lie on one line.

}

A nondegenerate quadrangle  is called a \index{parallelogram}\emph{parallelogram}
if its opposite sides are parallel.

\begin{thm}{Lemma}\label{lem:parallelogram}
Any parallelogram is \index{central symmetry}\emph{centrally symmetric} with respect to a midpoint of one of its diagonals;
that is, the reflection across the midpoint maps the parallelogram to itself.

In particular, if $\square A B C D$ is a parallelogram, then
\begin{enumerate}[(a)]
\item\label{lem:parallelogram:midpoint} its diagonals $[AC]$ and $[BD]$ intersect each other at their midpoints;
\item $\measuredangle A B C= \measuredangle C D A$;
\item $AB=CD$.
\end{enumerate}
\end{thm}

{

\begin{wrapfigure}{r}{33mm}
\centering
\includegraphics{mppics/pic-96}
\end{wrapfigure}

\parit{Proof.} Let $\square A B C D$ be a parallelogram.
Denote by $M$ the midpoint of $[AC]$.

Since $(AB)\parallel (CD)$, Theorem~\ref{thm:parallel-point-reflection} implies that $(CD)$ is a reflection of $(AB)$ across $M$.
In the same way, $(BC)$ is a reflection of $(DA)$ across $M$.
Since $\square A B C D$ is nondegenerate, it follows that $D$ is a reflection of $B$ across $M$; in other words, $M$ is the midpoint of $[BD]$.

The remaining statements follow since reflection across $M$ is a direct motion of the plane (see \ref{prop:point-reflection+}).
\qeds

}

\begin{thm}{Exercise}\label{ex:4parallels}
Assume that point $P$ is distinct from vertices of a parallelogram $ABCD$.
Let $a$, $b$, $c$ and $d$ be lines thru points $A$, $B$, $C$ and $D$ such that 
$a\parallel (CP)$,
$b\parallel (DP)$,
$c\parallel (AP)$, and
$d\parallel (BP)$.
Show that the lines $a$, $b$, $c$ and $d$ intersect at one point.
\end{thm}

\begin{thm}{Exercise}\label{ex:romb}
Assume $ABCD$ is a quadrangle such that
\[AB=CD=BC=DA.\]
Show that $ABCD$ is a parallelogram.
\end{thm}

A quadrangle as in the exercise above is called a \index{rhombus}\emph{rhombus}.

A quadrangle $ABCD$ is called a \index{rectangle}\emph{rectangle} if the angles $ABC$, $BCD$, $CDA$, and $DAB$ are right.
According to the transversal property (\ref{thm:parallel-2}),
any rectangle is a parallelogram.

A rectangle with equal sides is called a \index{square}\emph{square}.

\begin{thm}{Exercise}\label{ex:rectangle}
Show that a parallelogram $ABCD$ is a rectangle
if and only if $AC=BD$.
\end{thm}

\begin{thm}{Exercise}\label{ex:romb2}
Show that a parallelogram $ABCD$ is a rhombus
if and only if $(AC)\perp (BD)$.
\end{thm}

Assume $\ell\parallel m$, and $X,Y\in m$.
Let $X'$ and $Y'$ denote the footpoints of $X$ and $Y$ on~$\ell$.
Note that $\square XYY'X'$ is a rectangle.
By Lemma~\ref{lem:parallelogram}, $XX'=YY'$.
That is, every point on $m$ lies at the same distance from $\ell$.
This distance is called the \index{distance!between parallel lines}\emph{distance between} $\ell$ and~$m$.

\begin{thm}{Exercise}\label{ex:inscribed-rhombus}
A rhombus $AXYZ$ with side length $s$ is inscribed in a nondegenerate triangle $ABC$ so that $X\in [AB]$, $Y\in [BC]$, and $Z\in[CA]$.
Show that $s=\sqrt{XB\cdot ZC}$.
\end{thm}

\section{The method of coordinates}

The following exercise is important;
it shows that our axiomatic definition agrees with the model described in Section~\ref{page:model}.

\begin{thm}{Exercise}\label{ex:coordinates} 
Let $\ell$ and $m$ be perpendicular lines in the Euclidean plane.
Given a point $P$, let $P_\ell$ and $P_m$ denote the footpoints of $P$ on $\ell$ and $m$ respectively.

\begin{enumerate}[(a)]
\item Show that for every $X\in \ell$ and $Y\in m$ there is a unique point $P$ such that $P_\ell=X$ and $P_m=Y$.
\end{enumerate}

\begin{enumerate}[(a)]\addtocounter{enumi}{1}
\item
Show that 
$PQ^2=P_\ell Q_\ell^2+P_mQ_m^2$
for every pair of points $P$ and~$Q$.
\end{enumerate}

\begin{enumerate}[(a)]\addtocounter{enumi}{2}
\item Conclude that the plane is isometric to $(\mathbb{R}^2,d_2)$; see \ref{ex:dist-square}.
\end{enumerate}

\end{thm}

\begin{wrapfigure}{r}{37mm}
\centering
\includegraphics{mppics/pic-98}
\end{wrapfigure}

Once this exercise is solved, we can apply 
the method of coordinates
to solve any problem in Euclidean plane geometry.
This method is powerful and universal;
it will be developed further in Chapter~\ref{chap:complex}.

\begin{thm}{Exercise}\label{ex:abc}
Use Exercise~\ref{ex:coordinates}
to give an alternative proof of Theorem~\ref{thm:abc} in the Euclidean plane.

That is, prove that given the real numbers $a$, $b$, and $c$ such that 
 $$0<a\le b\le c\le a+b,$$
there is a triangle $ABC$
such that $a=BC$, $b=CA$, and $c=AB$.
\end{thm} 

\begin{thm}{Exercise}\label{ex:line-coord}
Consider two distinct points $A=(x_A,y_A)$ and $B\z=(x_B,y_B)$ on the coordinate plane.
Show that the perpendicular bisector of $[AB]$ is described by the equation
\[2\cdot (x_B-x_A)\cdot x+2\cdot (y_B-y_A)\cdot y=x_B^2+y_B^2-x_A^2-y_A^2.\]

Conclude that a line can be defined as a subset of the coordinate plane of the following type:
\begin{enumerate}[(a)]
\item  Solutions of an equation $a\cdot x+b\cdot y=c$
for constants $a$, $b$, and $c$ such that $a\ne 0$ or $b\ne0$.
\item\label{ex:line-coord:parameter} The set of points $(a\cdot t+c,b\cdot t+d)$ for constants $a$, $b$, $c$, and $d$ such that $a\ne 0$ or $b\ne0$ and all $t\in \mathbb{R}$.
\end{enumerate}

\end{thm}

\section{Apollonian circles}\label{sec:Apollonian circle}
The exercises in this section illustrate the method of coordinates; they will not be used further in the sequel.

\begin{thm}{Exercise}\label{ex:circle-coord}
Show that for fixed real values $a$, $b$, and $c$ the equation 
\[x^2+y^2+a\cdot x+b\cdot y+c=0\]
describes a circle, a point, or an empty set.

Show that if it is a circle, then $(-\tfrac a2,-\tfrac b2)$ is its center,
and $r\z=\tfrac12\cdot \sqrt{a^2+b^2-4\cdot c}$ is its radius.
\end{thm}

\begin{thm}{Exercise}\label{ex:apolonnius}
Use the previous exercise to show that given a positive real number $k\ne1$,
the locus of points $M$ such that $AM=k\cdot BM$ 
for distinct points $A$ and $B$
is a circle. 
\end{thm}

\begin{figure}[!ht]
\centering
\includegraphics{mppics/pic-100}
\end{figure}

The circle in the exercise above is an example of the so-called \index{Apollonian circle}\emph{Apollonian circle with focuses $A$ and $B$}.
A few of these circles for different values of $k$ are shown in the picture;
for $k=1$, it becomes the perpendicular bisector of $[AB]$.

\begin{thm}{Exercise}\label{ex:apolonnius-construction}
Construct an Apollonian circle with given focuses $A$ and $B$ thru a given point $M$ using ruler and compass.
\end{thm}

\chapter{Triangle geometry}\label{chap:triangle}

Triangle geometry is the study of the properties of triangles, including associated centers and circles.

We discuss the most basic results in triangle geometry, 
mostly to show that we have developed sufficient machinery to prove things.

\section{Circumcircle and circumcenter}

\begin{thm}{Theorem}\label{thm:circumcenter}
Perpendicular bisectors to the sides of every nondegenerate triangle intersect at one point.
\end{thm}

The point of intersection of the perpendicular bisectors is called the \index{circumcenter}\emph{circumcenter}.
It is the center of the \index{circumcircle}\emph{circumcircle} of the triangle;
that is, a circle that passes thru all three vertices of the triangle.
The circumcenter of the triangle is usually denoted by~$O$.

\begin{wrapfigure}{o}{29mm}
\centering
\includegraphics{mppics/pic-102}
\end{wrapfigure}

\parit{Proof.}
Let $\triangle ABC$ be nondegenerate.
Let $\ell$ and $m$ be perpendicular bisectors to sides $[AB]$ and $[AC]$ respectively.

Assume $\ell$ intersects $m$, say, at $O$.

Let us apply Theorem~\ref{thm:perp-bisect}.
Since $O\in\ell$, we have $OA\z=OB$ and since $O\in m$, we have $OA\z=OC$.
It follows that $OB\z=OC$;
that is, $O$ lies on the perpendicular bisector of~$[B C]$.

It remains to show that $\ell\nparallel m$;
assume the contrary.
Since
$\ell\perp(AB)$ and $m\perp (AC)$, we get that $(AC)\z\parallel (AB)$ 
(see Exercise~\ref{ex:perp-perp}).
Therefore, by Theorem~\ref{perp:ex+un}, $(AC)=(AB)$;
that is, $\triangle ABC$ is degenerate --- a contradiction.
\qeds

\begin{thm}{Exercise}\label{ex:unique-cline}
Show that there is a unique circle that passes thru the vertices of a given nondegenerate triangle in the Euclidean plane.
\end{thm}

\section{Altitudes and orthocenter}

An \index{altitude}\emph{altitude} of a triangle is a line thru a vertex and perpendicular to the line containing the opposite side.
The term \index{altitude}\emph{altitude} may also be used for the distance from the vertex to its footpoint on the line containing the opposite side.

\begin{thm}{Theorem}\label{thm:orthocenter}
Three altitudes of every nondegenerate triangle intersect at a single point.
\end{thm}

The point of intersection of altitudes is called the \index{orthocenter}\emph{orthocenter}; 
it is usually denoted by~$H$.

{

\begin{wrapfigure}{o}{34mm}
\vskip-4mm
\centering
\includegraphics{mppics/pic-104}
\end{wrapfigure}

\parit{Proof.}
Fix a nondegenerate triangle $A B C$.
Consider three lines $\ell$, $m$, and $n$
such that 
\begin{align*}
\ell&\parallel(BC),
&
m&\parallel(CA),
&
n&\parallel(AB),
\\
\ell&\ni A,
&
m&\ni B,
&
n&\ni C.
\end{align*}
Since $\triangle A B C$ is nondegenerate,
no pair of the lines $\ell$, $m$, and $n$ is parallel.
Let $A'$, $B'$, and $C'$ be the points of intersection of
$m$ and $n$, $n$ and $\ell$, and $\ell$ and $m$, respectively.

}

Note that $\square A B C B'$, $\square B C A C'$, and $\square C A B A'$, are parallelograms.
Applying Lemma~\ref{lem:parallelogram} we get that $\triangle ABC$ is the \index{medial triangle}\emph{medial triangle} of $\triangle A' B' C'$;
that is, $A$, $B$, and $C$ are the midpoints of $[B' C']$, $[C' A']$, and $[A' B']$ respectively.

By Exercise~\ref{ex:perp-perp},
$(B' C')\parallel (BC)$,
the altitude from $A$ is perpendicular to $[B' C']$, 
and from above it bisects~$[B' C']$.

Hence the altitudes of $\triangle A B C$ 
are also perpendicular bisectors of $\triangle A' B' C'$.
Applying Theorem~\ref{thm:circumcenter}, we get that altitudes of $\triangle ABC$ intersect at one point.
\qeds

\begin{thm}{Exercise}\label{ex:orthic-4}
Assume that the orthocenter $H$ of $\triangle ABC$ is distinct from its vertices.
Show that $A$ is the orthocenter of $\triangle H B C$.
\end{thm}

\begin{thm}{Exercise}\label{ex:orthic-sim}
Let $A'$, $B'$, and $C'$ be the footpoints of the corresponding altitudes of an acute triangle $ABC$.
Show that
\[
\triangle ABC \sim \triangle A'B'C \sim \triangle AB'C' \sim \triangle A'BC'.
\]
\end{thm}

\section{Medians and centroid}

A median of a triangle is the segment joining a vertex to the midpoint of the opposing side. 

\begin{thm}{Theorem}\label{thm:centroid}
Three medians of every nondegenerate triangle intersect at a single point.
Moreover, the point of intersection divides each median in the ratio 2:1.
\end{thm}

The point of intersection of medians is called the \index{centroid}\emph{centroid} of the triangle; 
it is usually denoted by~$M$.
In the proof, we will apply exercises \ref{ex:chevinas} and \ref{ex:smililar+parallel}; their complete solutions are given in the hints.

\parit{Proof.}
Consider a nondegenerate triangle $A B C$.
Let $[A A']$ and $[B B']$ be its medians.
According to Exercise~\ref{ex:chevinas}, 
$[A A']$ and $[B B']$ have a point of intersection;
denote it by $M$.

\begin{wrapfigure}{o}{36mm}
\vskip-4mm
\centering
\includegraphics{mppics/pic-106}
\end{wrapfigure}

Draw a line $\ell$ thru $A'$ parallel to $(BB')$.
Applying Exercise~\ref{ex:smililar+parallel} for $\triangle BB'C$ and $\ell$, we get that $\ell$ crosses $[B'C]$ at a point, say $X$, and
\[\frac{CX}{CB'}=\frac{CA'}{CB}=\frac12;\]
that is, $X$ is the midpoint of $[CB']$.

Since $B'$ is the midpoint of $[AC]$ and $X$ is the midpoint of $[B'C]$, we get that 
\[\frac{AB'}{AX}=\frac23.\]

Applying Exercise~\ref{ex:smililar+parallel} for $\triangle XA'A$ and the line $(BB')$, we get that 
\[\frac{AM}{AA'}=\frac{AB'}{AX}=\frac23;\eqlbl{eq:2/3}\]
that is, $M$ divides $[AA']$ in the ratio 2:1.

Condition \ref{eq:2/3} uniquely defines $M$ on $[AA']$.
Repeating the same argument for medians $[AA']$ and $[CC']$, we get that they intersect at~$M$ as well,
hence the result.
\qeds

\begin{thm}{Exercise}\label{ex:midle}
Let $\square ABCD$ be a nondegenerate quadrangle
and $X$, $Y$, $V$, and~$W$ be midpoints of 
$[AB]$, $[BC]$, $[CD]$, and~$[DA]$ respectively.
Show that $\square XYVW$ is a parallelogram.
\end{thm}

\begin{thm}{Advaneced exercise}\label{ex:euler-line}
Show that the orthocenter $H$, centroid $M$, and circumcenter $O$ of every nondegenerate triangle $ABC$ lie on a single line.
Moreover, the centroid divides the segment $[HO]$ in the ratio $2:1$.
\end{thm}

The line in this exercise is called \emph{Euler's line}.

\section{Angle bisectors}

If $\measuredangle A B X\equiv-\measuredangle C B X$, 
then we say that the line $(BX)$ {}\emph{bisects} $\angle ABC$,
or the line $(BX)$ is а \index{bisector!angle bisector}\emph{bisector} of $\angle ABC$.
If $\measuredangle A B X\equiv\pi-\measuredangle C B X$, then the line $(BX)$ is called the \index{bisector!external bisector}\emph{external bisector} of $\angle ABC$.

\begin{wrapfigure}{o}{42mm}
\centering
\includegraphics{mppics/pic-108}
\end{wrapfigure}

If $\measuredangle ABA'=\pi$;
that is, if $B$ lies between $A$ and $A'$,
then the bisector of $\angle ABC$ is the external bisector of $\angle A' B C$ and the other way around.

Note that the bisector and the external bisector are uniquely defined by the angle.

\begin{thm}{Classroom exercise}\label{ex:perp-bisectors}
Show that for every angle, its bisector and external bisector are perpendicular.
\end{thm}

The bisectors of  $\angle ABC$, $\angle BCA$, and $\angle CAB$ of a nondegenerate triangle $A B C$
are called \index{bisector!of the triangle}\emph{bisectors of the triangle} $A B C$ at vertices $A$, $B$, and $C$ respectively.

\begin{thm}{Exercise}\label{ex:bisect=altitude}
Assume that, at one vertex of a nondegenerate triangle, its bisector coincides with its altitude.
Show that  the triangle is isosceles.
\end{thm}

\begin{thm}{Lemma}\label{lem:bisect-ratio}
Let $\triangle A B C$ be  a nondegenerate triangle.
Assume that the bisector at $A$ 
intersects $[BC]$ at~$D$.
Then 
$$\frac{AB}{AC}=\frac{DB}{DC}.
\eqlbl{bisect-ratio}$$

\end{thm}

\begin{wrapfigure}{r}{28mm}
\vskip-6mm
\centering
\includegraphics{mppics/pic-110}
\end{wrapfigure}

\parit{Proof.}
Let $\ell$ be a line passing thru $C$ that is parallel to~$(AB)$.
Note that the lines $\ell$ and $(AD)$ are not parallel;
let $E$ be their point of intersection.

Note also that $B$ and $C$ lie on opposite sides of~$(AD)$.
By the transversal property (\ref{thm:parallel-2}),
$$\measuredangle BAD=\measuredangle CED.\eqlbl{eq:<BAD=<CED}$$

Furthermore, the angles $ADB$ and $EDC$ are vertical; by \ref{prop:vert} we have
$$\measuredangle ADB=\measuredangle EDC.$$

By the AA similarity condition, 
$\triangle ABD\sim \triangle ECD$.
In particular, 
$$\frac{AB}{EC}=\frac{DB}{DC}.\eqlbl{eq:AB/EC=DB/DC}$$

Since $(AD)$ bisects $\angle BAC$, we get that
$\measuredangle BAD=\measuredangle DAC$.
Together with \ref{eq:<BAD=<CED},
it implies that 
$\measuredangle CEA=\measuredangle EAC$.
By Theorem~\ref{thm:isos}, $\triangle ACE$ is isosceles; 
that is, $$EC=AC.$$
Together with \ref{eq:AB/EC=DB/DC}, it implies \ref{bisect-ratio}.
\qeds 

\begin{thm}{Exercise}\label{ex:ext-disect}
Formulate and prove an analog of Lemma~\ref{lem:bisect-ratio} for the external bisector.
\end{thm}

\begin{thm}{Exercise}\label{ex:bisect=median} 
Assume that an angle bisector of a nondegenerate triangle bisects the opposite side. 
Show that the triangle is isosceles.
\end{thm}

{

\begin{wrapfigure}{r}{38mm}
\vskip-5mm
\centering
\includegraphics{mppics/pic-118}
\end{wrapfigure}

\begin{thm}{Exercise}\label{ex:bisector-parallel} 
Assume that the bisector at $A$ of a nondegenerate triangle $ABC$ intersects $[BC]$ at $D$;
a line thru $D$ and parallel to $(CA)$ intersects $(AB)$ at $E$;
a line thru $E$ and parallel to $(BC)$ intersects $(AC)$ at $F$.

Show that 
$AE=FC$.

\end{thm}

}

\section{The equidistant property}

Recall that the distance from a line $\ell$ to a point $P$ is defined as the distance from $P$ to its footpoint on $\ell$; see Section~\ref{sec:perp<oblique}. 

\begin{thm}[\abs]{Proposition}\label{prop:angle-bisect-dist}
Assume $\triangle ABC$ is not degenerate.
Then a point $X$ lies on a bisector or external bisector of $\angle ABC$
if and only if $X$ is equidistant from the lines $(AB)$ and $(BC)$.
\end{thm}

\parit{Proof.}
We can assume that $X$ does not lie on the union of $(AB)$ and $(BC)$.
Otherwise, the distance to one of the lines vanishes;
in this case, $X=B$ is the only point equidistant from the two lines.

Let $Y$ and $Z$ be the reflections of $X$ across $(AB)$ and $(BC)$ respectively.
Note that 
\[Y\ne Z.\]
Otherwise, both lines $(AB)$ and $(BC)$ are perpendicular bisectors of $[XY]$.
Hence $(AB)=(BC)$, but this is impossible since $\triangle ABC$ is not degenerate.

By Proposition~\ref{prop:reflection}, $XB=YB=ZB$.
Note that $X$ is equidistant from $(AB)$ and $(BC)$ if and only if $XY\z=XZ$.

{

\begin{wrapfigure}{r}{24mm}
\centering
\includegraphics{mppics/pic-112}
\end{wrapfigure}

Applying SSS and then SAS, we get that
$$\begin{aligned}
XY&=XZ
\\
&\hskip0.4mm\Updownarrow
\\
\triangle XBY&\cong\triangle XBZ
\\
&\hskip0.4mm\Updownarrow
\\
\measuredangle XBY&\equiv \pm \measuredangle XBZ.
\end{aligned}
$$

Since $Y\ne Z$, we get that $\measuredangle XBY\ne \measuredangle XBZ$.
Therefore $X$ is equidistant from $(AB)$ and $(BC)$ if and only if
\[\measuredangle XBY\equiv -\measuredangle XBZ.
\eqlbl{eq:iff-chain}\]

}

By Proposition~\ref{prop:reflection}, $A$ lies on the bisector of $\angle XBY$,
and $C$ lies on the bisector of $\angle XBZ$; that is,
\begin{align*}
2\cdot \measuredangle XBA&\equiv \measuredangle XBY,
&
2\cdot \measuredangle XBC&\equiv \measuredangle XBZ.
\end{align*}
By \ref{eq:iff-chain},
\[2\cdot \measuredangle XBA\equiv -2\cdot \measuredangle XBC.\]
The last identity means either
\[
\measuredangle XBA+\measuredangle XBC\equiv 0
\quad
\text{or}
\quad
\measuredangle XBA+\measuredangle XBC\equiv \pi
\]
--- hence the result.
\qeds

\section{The incenter}

\begin{thm}[\abs]{Theorem}\label{thm:incenter}
Angle bisectors of every nondegenerate triangle intersect at one point.
\end{thm}

The point of intersection of bisectors is called the \index{incenter}\emph{incenter} of the triangle, 
usually denoted by~$I$.
The point $I$ lies at the same distance from each side.
In particular, it is the center of a circle tangent to each side of the triangle.
This circle is called 
the \index{incircle}\emph{incircle} and its radius is called 
the \index{inradius}\emph{inradius} of the triangle.

\parit{Proof.} 
Let $\triangle ABC$ be a nondegenerate triangle.

Points $B$ and $C$ lie on opposite sides of the bisector of $\angle BAC$.
Hence this bisector intersects $[BC]$ at a point, say~$A'$.

Analogously, there is $B'\in[AC]$ 
such that $(BB')$ bisects $\angle ABC$.

Applying Exercise~\ref{ex:chevinas},
we get that $[AA']$ and $[BB']$ intersect.
Suppose that $I$ denotes the point of intersection.

{

\begin{wrapfigure}{o}{28mm}
\vskip0mm
\centering
\includegraphics{mppics/pic-114}
\end{wrapfigure}

Let $X$, $Y$, and $Z$ be the footpoints of $I$ on  $(B C)$, $(C A)$, and $(A B)$ respectively.
Applying Proposition~\ref{prop:angle-bisect-dist}, we get that
$$I Y=I Z=I X.$$
From the same lemma, we get that $I$ lies on the bisector or on the exterior bisector of $\angle B C A$.

The line $(C I)$ intersects $[B B']$;
points $B$ and $B'$ lie on opposite sides of~$(C I)$.
Therefore, the angles $I C B'$ and $I C B$ have opposite signs.
Note that $\angle I C A=\angle I C B'$.
Therefore, $(C I)$ cannot be the exterior bisector of $\angle B C A$.
Hence the result.
\qeds

}

{

\begin{wrapfigure}{r}{30mm}
\centering
\vskip-2mm
\includegraphics{mppics/pic-116}
\end{wrapfigure}

\begin{thm}{Exercise}\label{ex:2x=b+c-a}
Assume sides $[B C]$, $[C A]$, and $[A B]$ of $\triangle A B C$ are tangent to the incircle at $X$, $Y$, and $Z$ respectively. 
Show that 
$$AY=AZ= \tfrac12\cdot(A B+ A C- B C).$$

\end{thm}

Assume that footpoints $A'$, $B'$, and $C'$ of the altitudes of a given triangle $ABC$ are all mutually distinct.
Then $\triangle A'B'C'$ is called an \index{triangle!orthic triangle}\index{orthic triangle}\emph{orthic triangle} of $\triangle ABC$.

}

\begin{thm}{Exercise}\label{ex:orthic-triangle}
Prove that an orthocenter of an acute triangle coincides with an incenter of its orthic triangle.

What should be an analog of this statement for an obtuse triangle?
\end{thm}

\begin{thm}{Exercise}\label{ex:bisector-incenter}
Let $I$ be the intersection of angle bisectors at $A$ and $B$ of a nondegenerate triangle $ABC$.
Denote by $D$ the intersection of the angle bisector at $A$ with the side $[BC]$.
Show that $\frac{AI}{DI}=\frac{b+c}{a}$,
where $a=BC$, $b=CA$, and $c=AB$.

Use it to build another proof of \ref{thm:incenter}.
\end{thm}

\addtocontents{toc}{\protect\contentsline{part}{\protect\numberline{}Inversive geometry}{}{}}

\chapter{Inscribed angles}\label{chap:inscribed-angle}

\section{Angle between a tangent line and a chord}

\begin{thm}{Theorem}\label{thm:tangent-angle}
Let $\Gamma$ be a circle with a center $O$.
Assume that a line $(XQ)$ is tangent to $\Gamma$ at $X$
and $[XY]$ is a chord of~$\Gamma$.
Then 
$$2\cdot\measuredangle QXY
\equiv\measuredangle X O Y.
\eqlbl{eq:tangent-angle}$$
Equivalently, 
$$\measuredangle QXY
\equiv
\tfrac12\cdot\measuredangle X O Y
\quad 
\text{or}
\quad
\measuredangle QXY
\equiv
\tfrac12\cdot\measuredangle X O Y+\pi.$$

\end{thm}

\begin{wrapfigure}{o}{33mm}
\centering
\includegraphics{mppics/pic-120}
\end{wrapfigure}

\parit{Proof.}
By Lemma~\ref{lem:tangent}, $(OX)\z\perp(XQ)$;
therefore,
$$\measuredangle QXY+\measuredangle YXO \equiv\pm\tfrac\pi2.$$

Note that $\triangle XOY$ is isosceles.
Therefore,
$$\measuredangle YXO=\measuredangle OYX.$$

Applying Theorem~\ref{thm:3sum}
to $\triangle XOY$,
we get
\begin{align*}
\pi&\equiv\measuredangle YXO+\measuredangle OYX+\measuredangle XOY\equiv
\\
&\equiv 2\cdot \measuredangle YXO+\measuredangle XOY.
\end{align*}

It follows that
$$2\cdot\measuredangle QXY
\equiv \pi -2\cdot \measuredangle YXO
\equiv\measuredangle X O Y.
$$
\qedsf

\section{Inscribed angles}\label{sec:inscribed}

We say that a triangle is \index{inscribed triangle}\emph{inscribed} in the circle $\Gamma$ if all its vertices lie on~$\Gamma$.

\begin{thm}{Theorem}\label{thm:inscribed-angle}
Let $X$ and $Y$ be distinct points on a circle $\Gamma$ centered at~$O$.
Then
$\triangle X P Y$ is inscribed in $\Gamma$ if and only if
$$2\cdot\measuredangle X P Y\equiv\measuredangle X O Y.
\eqlbl{eq:inscribed-angle}$$
Equivalently, if and only if
$$\measuredangle XPY\equiv\tfrac12\cdot\measuredangle X O Y
\quad
\text{or}
\quad
\measuredangle XPY\equiv\pi+\tfrac12\cdot\measuredangle X O Y.$$

\end{thm}

\begin{wrapfigure}{o}{33mm}
\vskip-6mm
\centering
\includegraphics{mppics/pic-122}
\vskip4mm
\includegraphics{mppics/pic-124}
\vskip4mm
\includegraphics{mppics/pic-126}
\end{wrapfigure}

\parit{Proof; the ``only if'' part.}
Let $(PQ)$ be the tangent line to $\Gamma$ at~$P$.
By Theorem~\ref{thm:tangent-angle},
\begin{align*}
2\cdot\measuredangle QPX&\equiv\measuredangle POX,
&
2\cdot\measuredangle QPY&\equiv\measuredangle POY.
\end{align*}
Subtracting one identity from the other, we get~\ref{eq:inscribed-angle}.

\parit{``If'' part.}
Assume that \ref{eq:inscribed-angle} holds for some $P\notin \Gamma$.
Note that $\measuredangle X O Y\ne 0$. 
Therefore, $\measuredangle X P Y\ne 0$ nor $\pi$;
that is, $\triangle PXY$ is nondegenerate.

The line $(PX)$ is tangent to $\Gamma$ at the point $X$, or it intersects $\Gamma$ at another point.
In the latter case, suppose that $P'$ denotes this point of intersection. 

In the first case, by Theorem~\ref{thm:tangent-angle}, we have
\begin{align*}
2\cdot \measuredangle PXY&\equiv \measuredangle XOY\equiv 
 2\cdot\measuredangle  XPY.
\end{align*}
Applying the transversal property (\ref{thm:parallel-2}), we get that
$(XY)\parallel (PY)$, which is impossible since $\triangle PXY$ is nondegenerate.

In the second case, 
applying the ``if'' part and that  $P$, $X$, and $P'$ lie on one line (see Exercise~\ref{ex:ABCO-line}) we get that 
\begin{align*}
2\cdot \measuredangle P'PY&\equiv
2\cdot \measuredangle XPY\equiv 
 \measuredangle  XOY\equiv
 \\
&\equiv 2\cdot\measuredangle  XP'Y\equiv
 2\cdot\measuredangle  PP'Y.
\end{align*}
Again, by the transversal property,
$(PY)\z\parallel (P'Y)$, which is impossible since $\triangle PXY$ is nondegenerate.
\qeds

{

\begin{wrapfigure}{r}{34mm}
\centering
\includegraphics{mppics/pic-128}
\vskip3mm
\includegraphics{mppics/pic-130}
\vskip3mm
\includegraphics{mppics/pic-132}
\end{wrapfigure}

\begin{thm}{Exercise}\label{ex:inscribed-angle}
Let $X$, $X'$, $Y$, and $Y'$ be distinct points on the circle $\Gamma$.
Assume $(XX')$ meets $(YY')$ at a point~$P$.
Show that 
\begin{enumerate}[(a)]
\item $2\cdot \measuredangle XPY\equiv\measuredangle XOY+\measuredangle X'OY'$;
\item\label{ex:inscribed-angle:b} $\triangle PXY\sim \triangle PY'X'$;
\item\label{ex:inscribed-angle:power} $PX\cdot PX'=|OP^2-r^2|$, where $O$ is the center and $r$ is the radius of $\Gamma$.
\end{enumerate}

\end{thm}

(The value $OP^2-r^2$ is called the \index{power of a point}\emph{power} of the point $P$ with respect to the circle $\Gamma$. 
Part \textit{(\ref{ex:inscribed-angle:power})} of the exercise makes it a useful tool to study circles, but we are not going to consider it further in the book.) 

\begin{thm}{Exercise}\label{ex:inscribed-hex}
Three chords $[XX']$, $[YY']$, and $[ZZ']$
of a circle $\Gamma$ intersect at a point $P$.
Show that 
$$XY'\cdot ZX'\cdot YZ'=X'Y\cdot Z'X\cdot Y'Z.$$

\end{thm}

\begin{thm}{Exercise}\label{ex:altitudes-circumcircle}
Let $\Gamma$ be a circumcircle of an acute triangle $A B C$.
Denote by $A'$ and $B'$ the second points of intersection of the altitudes from $A$ and $B$ with~$\Gamma$.
Show that $\triangle A' B' C$ is isosceles.
\end{thm}

}

\begin{thm}{Exercise}\label{ex:two-chords}
Let $[XY]$ and $[X'Y']$ be two parallel chords of a circle.
Show that $XX'=YY'$.
\end{thm}

\begin{thm}{Exercise}
Watch ``Why is pi here? And why is it squared? A geometric answer to the Basel problem'' by Grant Sanderson. (It is available on \href{https://youtu.be/d-o3eB9sfls}{YouTube}.) 
Prepare one question.
\end{thm}

\section{Points on a circle}

Recall that the diameter of a circle is a chord that passes thru the center.
If $[XY]$ is the diameter of a circle with center $O$, then $\measuredangle X O Y\z=\pi$. 
Hence Theorem~\ref{thm:inscribed-angle} implies the following:

\begin{thm}{Corollary}\label{cor:right-angle-diameter}
Suppose $\Gamma$ is a circle with the diameter~$[AB]$.
A triangle $ABC$ has a right angle at $C$ if and only if $C\in\Gamma$.
\end{thm}

The following two exercises can be used in construction problems.

\begin{thm}{Exercise}\label{ex:tangent-construction-inscribed}
Let $\Gamma$ be a circle with center $O$, and let $\Delta$ be a circle with a diameter $[OP]$.
Assume that $\Delta$ and $\Gamma$ intersect at a point $X$.
Prove that the line $(PX)$ is tangent to $\Gamma$.
\end{thm}

\begin{thm}{Exercise}\label{ex:perp-construction-inscribed}
Let $\Gamma$ be a circle with diameter $[PQ]$.
Assume that a line $\ell$ passes thru $Q$ and intersects $\Gamma$ at another point $X$.
Prove that $\ell \perp (PX)$.
\end{thm}

\begin{thm}{Exercise}\label{ex:altitude+circles}
Assume two circles intersect at two distinct points $A$ and $B$, and have diameters $[AX]$ and $[AY]$.
Prove that the points $X$, $Y$, and $B$ lie on one line.
\end{thm}

\begin{thm}{Exercise}\label{ex:perpendicular-ruler}
Given a line $\ell$, a circle $\Gamma$ with its center on $\ell$,
and a point $P$ that lies on neither $\ell$ nor $\Gamma$,
construct the perpendicular from $P$ to $\ell$ using only a ruler.
\end{thm}

\begin{wrapfigure}{r}{32mm}
\vskip-4mm
\centering
\includegraphics{mppics/pic-133}
\end{wrapfigure}

\begin{thm}{Exercise}\label{ex:tnagents+midpoint}
Suppose that lines $\ell$, $m$, and $n$ pass thru a point $P$;
the lines $\ell$ and $m$ are tangent to a circle $\Gamma$ at $L$ and $M$;
the line $n$ intersects $\Gamma$ at two points $X$ and $Y$.
Let $N$ be the midpoint of $[XY]$.
Show that the points $P$, $L$, $M$, and $N$ lie on one circle.
\end{thm}

We say that a quadrangle $ABCD$ is 
\index{quadrangle!inscribed quadrangle}\emph{inscribed in circle $\Gamma$}
if all the points $A$, $B$, $C$, and $D$ lie on $\Gamma$.

\begin{thm}{Corollary}\label{cor:inscribed-quadrangle}
A nondegenerate quadrangle $ABCD$ is inscribed in a circle if and only if 
\[2\cdot\measuredangle ABC\equiv 2\cdot\measuredangle ADC.\]

\end{thm}

\parit{Proof.}
Since $\square ABCD$ is nondegenerate, so is $\triangle ABC$.
Let $O$ and $\Gamma$ denote the circumcenter and circumcircle of $\triangle ABC$ (they exist by Exercise~\ref{ex:unique-cline}).

{

\begin{wrapfigure}[10]{o}{25mm}
\vskip-1mm
\centering
\includegraphics{mppics/pic-134}
\end{wrapfigure}

According to Theorem~\ref{thm:inscribed-angle},
$$
2\cdot\measuredangle ABC
\equiv
\measuredangle AOC.
$$
From the same theorem, $D\in\Gamma$ if and only if 

$$
2\cdot\measuredangle ADC
\equiv\measuredangle AOC,
$$
hence the result.
\qeds

}

{

\begin{thm}{Exercise}\label{ex:VVAA}
Let $\triangle A B C$ be a nondegenerate triangle,
$A'$ and $B'$ be footpoints of altitudes from $A$ and~$B$ respectively.
Show that the four points $A$, $B$, $A'$, and $B'$ lie on one circle.
What is the center of this circle?
\end{thm}

\begin{wrapfigure}{r}{33mm}
\vskip-0mm
\centering
\includegraphics{mppics/pic-136}
\end{wrapfigure}

\begin{thm}{Exercise}\label{ex:secant-circles}
Consider two circles $\Gamma$ and $\Gamma'$ that intersect at two distinct points $A$ and~$B$.
Assume $[XY]$ and $[X'Y']$ are the chords of $\Gamma$ and $\Gamma'$ respectively,
such that $A$ lies between $X$ and $X'$ and $B$ lies between $Y$ and~$Y'$.
Show that $(XY)\parallel (X'Y')$.
\end{thm}

}

\begin{thm}{Advanced exercise}\label{ex:perim+angle+side}
Perform a ruler-and-compass construction of $\triangle ABC$, given its perimeter $p$, $\beta=\measuredangle ABC$, and $b=AC$.
\end{thm}

\section{The method of additional circle}

\begin{thm*}{Problem}
 Assume that two chords $[AA']$ and $[BB']$ intersect at the point $P$ inside their circle.
Let $X$ be a point such that both angles $XAA'$ and $XBB'$ are right.
Show that $(XP)\perp(A'B')$.
\end{thm*}

\parit{Solution.}
Suppose that the lines $(A'B')$ and $(XP)$ meet at $Y$.

Both angles $XAA'$ and $XBB'$ are right;
therefore
\[2\cdot\measuredangle XAA'
\equiv
2\cdot\measuredangle XBB'.\]
By Corollary~\ref{cor:inscribed-quadrangle},  $\square XAPB$ is inscribed.
Applying this theorem again we get that
\[2\cdot\measuredangle AXP
\equiv
2\cdot\measuredangle ABP.\]

\begin{wrapfigure}[9]{o}{32mm}
\vskip-7mm
\centering
\includegraphics{mppics/pic-138}
\end{wrapfigure}

Since $\square ABA'B'$ is inscribed, 
\[2\cdot\measuredangle ABB'
\equiv
2\cdot\measuredangle AA'B'.\]

It follows that 
\[2\cdot\measuredangle AXY
\equiv
2\cdot\measuredangle AA'Y.\]
By the same theorem, $\square XAYA'$ is inscribed,
and
therefore, 
\[2\cdot\measuredangle XAA'
\equiv
2\cdot\measuredangle XYA'.\]
Since $\angle XAA'$ is right, 
so is $\angle XYA'$. 
That is, $(XP)\perp(A'B')$.
\qeds

\begin{thm}{Exercise}\label{ex:inaccuracy}
Find an inaccuracy in the solution of the problem and try to fix it.
\end{thm}

The method used in the solution 
is called the \textit{method of additional circle}
--- the circumcircles of the quadrangles $XAPB$ and $XAYA'$ 
 above can be considered as \textit{additional constructions}. 

{

\begin{wrapfigure}{r}{28mm}
\vskip-0mm
\centering
\includegraphics{mppics/pic-140}
\end{wrapfigure}

\begin{thm}{Exercise}\label{ex:equilateral-2}
Assume three lines $\ell$, $m$, and $n$ intersect at point $O$ and form six equal angles at~$O$. 
Let $X$ be a point distinct from $O$.
Let $L$, $M$, and $N$ denote the footpoints of perpendiculars from $X$ on $\ell$, $m$, and $n$, respectively.
Show that $\triangle LMN$ is equilateral.
\end{thm}

}

\begin{thm}{Exercise}\label{ex:median-angle}
Let $[AA']$ and $[BB']$ be the medians of a nondegenerate triangle $ABC$.
Assume that $\measuredangle A'AC = \measuredangle CBB'$.
Show that $AC = BC$.
\end{thm}

\begin{thm}{Advanced exercise}\label{ex:simson}
Assume that a point $P$ lies on a circumcircle of $\triangle ABC$.
Show that three footpoints of $P$ on the lines $(AB)$, $(BC)$, and $(CA)$ lie on one line
(this is the so-called \index{Simson line}\emph{Simson line} of $P$).
\end{thm}

\section{Arcs of circlines}

A subset of a circle bounded by two points is called a circular arc.

More precisely,
suppose $A$, $B$, and $C$ are distinct points on a circle $\Gamma$.
The \index{circular arc}\emph{circular arc}~$ABC$ is the subset that includes the points $A$, $C$,
as well as all the points on $\Gamma$ that lie with $B$ on the same side of $(AC)$.

Points $A$ and $C$ are called 
\index{endpoint of arc}\emph{endpoints} of the circular arc $ABC$. 
There are precisely two circular arcs of $\Gamma$ with the given endpoints; they are \index{opposite arc}\emph{opposite} to each other.

\begin{wrapfigure}[7]{r}{33mm}
\vskip-6mm
\centering
\includegraphics{mppics/pic-142}
\end{wrapfigure}

Suppose $X$ is another point on $\Gamma$.
By \ref{cor:inscribed-quadrangle} we have
that $2\cdot\measuredangle AXC\equiv 2\cdot\measuredangle ABC$;
that is,
\[\measuredangle AXC\equiv\measuredangle ABC
\quad\text{or}\quad
\measuredangle AXC\equiv\measuredangle ABC+\pi.\]

Recall that $X$ and $B$ lie on the same side of $(AC)$ if and only if $\angle AXC$ and $\angle ABC$ have the same sign (see Exercise~\ref{ex:signs-PXQ-PYQ}).
It follows that 
\begin{itemize}
\item $X$ lies on the arc $ABC$ if and only if $\measuredangle AXC\equiv\measuredangle ABC$;
\item $X$ lies on the arc opposite to $ABC$ if $\measuredangle AXC\equiv\measuredangle ABC+\pi$.
\end{itemize}

A circular arc $ABC$ is defined if $\triangle ABC$ is not degenerate.
If $\triangle ABC$ is degenerate, then arc $ABC$ is defined as a subset of the line bounded by $A$ and $C$ that contains $B$;
in this case we say that arc $ABC$ is \index{degenerate arc}\emph{degenerate}.

\begin{wrapfigure}{r}{43mm}
\vskip-2mm
\centering
\includegraphics{mppics/pic-144}
\vskip4mm
\includegraphics{mppics/pic-146}
\end{wrapfigure}

More precisely, if $B$ lies between $A$ and $C$, then the {}\emph{arc} $ABC$ is defined as 
the line segment $[AC]$.
If $B'$ lies on the extension of $[AC]$,
then the arc $AB'C$ is defined as a union of disjoint half-lines $[AX)$ and $[CY)$ in $(AC)$.
In this case, the arcs $ABC$ and $AB'C$ are called opposite to each other.

\begin{wrapfigure}[2]{r}{43mm}
\vskip-0mm
\centering
\includegraphics{mppics/pic-148}
\end{wrapfigure}

In addition, every half-line $[AB)$ will be regarded as an arc.
If $A$ lies between $B$ and $X$, then $[AX)$ will be called opposite to $[AB)$.
This arc has only one endpoint $A$.

It will be convenient to use the notion of 
\index{circline}\emph{circline},
which means \textit{circle or line}.
For example, every arc is a subset of a circline;
we also may use the term \index{circline!arc}\emph{circline arc} if we want to emphasize that the arc might be degenerate.
Note that for any three distinct points $A$, $B$, and $C$ there is a unique circline arc $ABC$.

The following statement summarizes the discussion above.

\begin{thm}{Proposition}\label{prop:arcs}
Let $ABC$ be a circline arc and $X$ be a point distinct from $A$ and $C$.
Then 
\begin{enumerate}[(a)]
\item $\measuredangle AXC=\measuredangle ABC$ if and only if $X$ lies on the arc $ABC$;
\item $\measuredangle AXC\equiv\measuredangle ABC+\pi$ if and only if $X$ lies on the arc opposite to $ABC$.
\end{enumerate}
\end{thm}

\begin{thm}{Exercise}\label{ex:3x120}
Given an acute triangle $ABC$,
perform a ruler-and-compass of a point $Z$ such that
$\measuredangle AZB
= \measuredangle BZC
= \measuredangle CZA
=\pm\tfrac23\cdot\pi$.
\end{thm}

\begin{thm}{Exercise}\label{ex:a+b=c}
A point $P$ lies on the circumcircle of an equilateral triangle $ABC$,
and $PA\le PB\le PC$.
Show that $PA+PB=PC$.
\end{thm}

A quadrangle $ABCD$ is 
\index{quadrangle!inscribed quadrangle}\emph{inscribed}
if all the points $A$, $B$, $C$, and $D$ lie on a circline $\Gamma$.
If the arcs $ABC$ and $ADC$ are opposite, then we say that the points $A$, $B$, $C$, and $D$ appear on $\Gamma$ in the same \index{cyclic order}\emph{cyclic order}.

This definition allows us to formulate the following refinement of Corollary~\ref{cor:inscribed-quadrangle} which includes degenerate quadrangles.
It follows directly from \ref{prop:arcs}.

\begin{thm}{Proposition}\label{prop:inscribed-quadrangle}
A quadrangle $ABCD$ is inscribed in a circline if and only if 
\[\measuredangle ABC+\measuredangle CDA\equiv 0
\quad\text{or}\quad
\measuredangle ABC+\measuredangle CDA\equiv\pi.\]
Moreover, the second identity holds if and only if the points $A,B,C,D$ appear on the circline in the same cyclic order.
\end{thm}

\section{Tangent half-lines}

{

\begin{wrapfigure}[6]{r}{30mm}
\vskip-14mm
\centering
\includegraphics{mppics/pic-150}
\end{wrapfigure}

Suppose $ABC$ is an arc of a circle $\Gamma$.
A half-line $[AX)$ is called 
\index{tangent!half-line}\emph{tangent} 
to the arc $ABC$ at $A$
if the line $(AX)$ is tangent to $\Gamma$, and the points $X$ and $B$ lie on the same side of the line $(AC)$.

}

If the arc is formed by the line segment $[AC]$, then the half-line $[AC)$ is considered to be tangent at $A$.
If the arc is formed by a union of two half-lines $[AX)$ and $[BY)$ in $(AC)$,
then the half-line $[AX)$ is considered to be tangent to the arc at $A$.

\begin{thm}{Proposition}\label{prop:arc(angle=tan)}
A half-line $[AX)$ is tangent to an arc $ABC$ if and only if 
$\measuredangle ABC+\measuredangle CAX\equiv \pi$.
\end{thm}

\parit{Proof.}
For a degenerate arc $ABC$, 
the statement is evident.
Assume it is nondegenerate.

Note that the tangent half-line to the arc $ABC$ at $A$ is uniquely defined.
Furthermore, there is a unique half-line $[AX)$ such that the equation in the proposition holds.
Therefore, it is sufficient to prove the ``only-if'' part.

\begin{wrapfigure}{o}{24mm}
\vskip-0mm
\centering
\includegraphics{mppics/pic-152}
\vskip4mm
\end{wrapfigure}

If $[AX)$ is tangent to the arc $ABC$, then by \ref{thm:tangent-angle}
and \ref{thm:inscribed-angle},
we get that
$$2\cdot \measuredangle ABC+2\cdot\measuredangle CAX\equiv 0.$$
Therefore, either 
$$\measuredangle ABC+\measuredangle CAX
\equiv 
\pi,
\quad
\text{or}
\quad
\measuredangle ABC+\measuredangle CAX
\equiv 
0.$$

By the definition of a tangent half-line,
$X$ and $B$ lie on the same side of~$(AC)$.
By \ref{cor:half-plane} and \ref{thm:signs-of-triug}, the angles $CAX$, $CAB$, and $ABC$ 
have the same sign.
In particular,
$\measuredangle ABC\z+\measuredangle CAX\not\equiv 0$;
that is, we are left with the case 
$\measuredangle ABC+\measuredangle CAX\equiv \pi$.
\qeds

\begin{thm}{Exercise}\label{ex:tangent-arc}
Show that there is a unique arc 
with given endpoints $A$ and $C$, 
that is tangent to a given half-line~$[AX)$.
\end{thm}

\begin{thm}{Exercise}\label{ex:tangent-lim}
Let $[AX)$ be a tangent half-line to an arc $ABC$.
Assume $Y$ is a point on the arc $ABC$ that is distinct from $A$.
Show that $\measuredangle XAY$ approaches $0$ as $AY$ approaches $0$.

\end{thm}

{

\begin{wrapfigure}{r}{33mm}
\vskip-10mm
\centering
\includegraphics{mppics/pic-154}
\end{wrapfigure}

\begin{thm}{Exercise}\label{ex:two-arcs}
Given two circular arcs $AB_1C$ and $AB_2C$, let $[AX_1)$ and $[AX_2)$ be half-lines tangent to the arcs $AB_1C$ and $AB_2C$ at $A$, and $[CY_1)$ and $[CY_2)$ be half-lines tangent to the arcs $AB_1C$ and $AB_2C$ at~$C$.
Show that
$$\measuredangle X_1AX_2\equiv -\measuredangle Y_1CY_2.$$

\end{thm}

\chapter{Inversion}\label{chap:inversion}

Let $\Omega$ be the circle with center $O$ and radius~$r$.
The \index{inversion}\emph{inversion} of a point $P$ across $\Omega$ is the point $P'\in[OP)$ such that
$$OP\cdot OP'=r^2.$$
In this case, the circle $\Omega$ will be called the 
\index{inversion!circle of inversion}\emph{circle of inversion}, 
and its center $O$ is called the \index{inversion!center of inversion}\emph{center of inversion}.

The inverse of $O$ is undefined.

If $P$ is inside $\Omega$, then $P'$ is outside and the other way around. 
Furthermore, $P=P'$ if and only if $P\in \Omega$.

{

\begin{wrapfigure}{r}{37mm}
\vskip-6mm
\centering
\includegraphics{mppics/pic-158}
\end{wrapfigure}

Note that the inversion maps $P'$ back to~$P$.

\begin{thm}{Exercise}\label{ex:constr-inversion}
Let $\Omega$ be a circle centered at~$O$.
Suppose that a line $(PT)$ is tangent to $\Omega$ at~$T$.
Let $P'$ be a footpoint of $T$ on $(OP)$.

Show that $P'$ is the inverse of $P$ across $\Omega$.
\end{thm}

}

\begin{thm}{Lemma}\label{lem:inversion-sim}
Let $\Gamma$ be a circle centered at~$O$.
Assume $A'$ and $B'$ are the inverses of $A$ and $B$ across~$\Gamma$.
Then 
$$\triangle O A B\sim\triangle O B' A'.$$
Moreover,
$$\begin{aligned}
\measuredangle AOB&\equiv -\measuredangle B'OA',
\\
\measuredangle OBA&\equiv -\measuredangle OA'B',
\\
\measuredangle BAO&\equiv -\measuredangle A'B'O.
\end{aligned}\eqlbl{eq:-angle}$$

\end{thm}

\parit{Proof.}
Let $r$ be the radius of the circle of the inversion.

\begin{wrapfigure}[12]{r}{32mm}
\centering
\includegraphics{mppics/pic-160}
\end{wrapfigure}

By the definition of an inversion, 
\begin{align*}
OA\cdot OA'=OB\cdot OB'=r^2.
\end{align*}
Therefore, 
$$\frac{OA}{OB'}=\frac{OB}{OA'}.$$

Clearly,
$$\measuredangle AOB= \measuredangle A'OB'\equiv -\measuredangle B'OA'.\eqlbl{eq:-AOB}$$
From SAS, we get that
$$\triangle O A B\z\sim\triangle O B' A'.$$

Applying Theorem~\ref{thm:signs-of-triug} and \ref{eq:-AOB},
we get \ref{eq:-angle}.
\qeds

\begin{thm}{Exercise}%
\label{ex:appolo-circ}
Let $P'$ be an inverse of $P$ across a circle $\Gamma$.
Assume that $P\ne P'$.
Show that the value $\frac{PX}{P'X}$ is the same for all $X\in\Gamma$.
\end{thm}

The converse to the exercise above also holds.
Namely, given a positive real number $k\ne 1$ 
and two distinct points $P$ and $P'$
the locus of points $X$ such that $\frac{PX}{P'X}=k$ forms a circle which is called the \index{Apollonian circle}\emph{Apollonian circle} (see also Section \ref{sec:Apollonian circle}).
In this case, $P'$ is the inverse of $P$ across the Apollonian circle.

\begin{thm}{Exercise}%
\label{ex:incenter+inversion=orthocenter}
Let $A'$, $B'$, and $C'$ be images of $A$, $B$, and $C$ 
under the inversion across the incircle of $\triangle A B C$.
Show that the incenter of $\triangle A B C$ 
is the orthocenter of $\triangle A' B' C'$.
\end{thm}

\begin{thm}{Exercise}\label{ex:consturuction-of-inversion}
Construct an inverse of a given point across a given circle using ruler and compass.
\end{thm}

\section{The cross-ratio}

The following theorem lists quantities that do not change after inversion.

\begin{thm}{Theorem}\label{lem:inverse-4-angle}
Let $ABCD$ and $A'B'C'D'$  be two quadrangles
such that the points $A'$, $B'$, $C'$, and $D'$ are the inverses of $A$, $B$, $C$, and $D$ respectively.
Then 
\begin{enumerate}[(a)]
\item\label{lem:inverse-4-angle:cross-ratio} $$\frac{AB\cdot CD}{BC\cdot DA}= \frac{A'B'\cdot C'D'}{B'C'\cdot D'A'}.$$
\item\label{lem:inverse-4-angle:angle} 
$$\measuredangle ABC+\measuredangle CDA\equiv -(\measuredangle A'B'C'+\measuredangle C'D'A').$$
\item\label{lem:inverse-4-angle:inscribed}
If the quadrangle $ABCD$ is inscribed, 
then so is $\square A'B'C'D'$.
\end{enumerate}
\end{thm}

The number $\frac{AB \cdot CD}{BC \cdot DA}$ from \textit{(\ref{lem:inverse-4-angle:cross-ratio})} is called the \index{cross-ratio}\emph{cross-ratio} of the four points $A$, $B$, $C$, and $D$.

\parit{Proof; (\ref{lem:inverse-4-angle:cross-ratio}).}
Let $O$ be the center of the inversion.
According to Lemma~\ref{lem:inversion-sim},
$\triangle AOB\z\sim \triangle B'OA'$.
Therefore, 
\begin{align*}
&&\frac{AB}{A'B'} &=\frac{OA}{OB'}.
\intertext{Analogously,}
\frac{BC}{B'C'}&=\frac{OC}{OB'},&
\frac{CD}{C'D'}&=\frac{OC}{OD'},&
\frac{DA}{D'A'}&=\frac{OA}{OD'}.
\end{align*}

Therefore, 
\begin{align*}
\frac{AB}{A'B'}\cdot \frac{CD}{C'D'}&=\frac{OA}{OB'}\cdot\frac{OC}{OD'}=
\frac{OC}{OB'}\cdot\frac{OA}{OD'}
=\frac{BC}{B'C'}\cdot\frac{DA}{D'A'}.
\end{align*}
Hence \textit{(\ref{lem:inverse-4-angle:cross-ratio})} follows.

\parit{(\ref{lem:inverse-4-angle:angle}).}
According to Lemma~\ref{lem:inversion-sim},
\[\begin{aligned}
\measuredangle ABO&\equiv -\measuredangle B'A'O,
&
\measuredangle OBC&\equiv -\measuredangle OC'B',\\
\measuredangle CDO&\equiv -\measuredangle D'C'O,
&
\measuredangle ODA&\equiv -\measuredangle OA'D'.
\end{aligned}\eqlbl{eq:angle=-angle}\]
By Axiom~\ref{def:birkhoff-axioms:2b},
\begin{align*}
\measuredangle ABC&\equiv\measuredangle ABO+\measuredangle OBC,
&
\measuredangle D'C'B'&\equiv\measuredangle D'C'O+\measuredangle OC'B',
\\
\measuredangle CDA&\equiv\measuredangle CDO+\measuredangle ODA,
&
\measuredangle B'A'D'&\equiv\measuredangle B'A'O+\measuredangle OA'D'.
\end{align*}
Therefore, 
summing the four identities in \ref{eq:angle=-angle}, we get that
\begin{align*}
\measuredangle ABC+\measuredangle CDA
&\equiv -(\measuredangle D'C'B'+\measuredangle B'A'D').
\intertext{Applying Axiom~\ref{def:birkhoff-axioms:2b} and Exercise~\ref{ex:quadrangle}, we get that}
\measuredangle A'B'C'+\measuredangle C'D'A'
&\equiv -(\measuredangle B'C'D'+\measuredangle D'A'B')\equiv
\\
&\equiv \measuredangle D'C'B'+\measuredangle B'A'D'.
\end{align*}
Hence \textit{(\ref{lem:inverse-4-angle:angle})} follows.

\parit{(\ref{lem:inverse-4-angle:inscribed}).}
Follows by \textit{(\ref{lem:inverse-4-angle:angle})} and Corollary~\ref{cor:inscribed-quadrangle}.
\qeds

\section{Inversive plane and circlines}

Let $\Omega$ be a circle with the center $O$ and the radius~$r$.
Consider the inversion across~$\Omega$.

Recall that the inverse of $O$ is undefined.
To deal with this problem it is useful to add an extra point to the plane, which will be called the \index{point!at infinity}\emph{point at infinity}, denoted as~\index{50@$\infty$}$\infty$.
We can assume that $\infty$ is the inverse of $O$ and the other way around.

The Euclidean plane with the added point at infinity is called the \index{plane!inversive plane}\index{inversive!plane}\emph{inversive plane}.

We will always assume that every line and half-line contains~$\infty$.

Recall that
\index{circline}\emph{circline}
means \textit{circle or line}.
Therefore we may say 
\textit{``if a circline contains $\infty$, then it is a line''} or \textit{``a circline that does not contain $\infty$  is a circle''}.

According to Theorem~\ref{thm:circumcenter}, 
for every $\triangle ABC$ there is a unique circline that passes thru $A$, $B$, and~$C$
(if $\triangle ABC$ is degenerate, then this is a line, and if not it is a circle).

\begin{thm}{Theorem}\label{thm:inverse-cline}
In the inversive plane, the inverse of a circline is a circline.
\end{thm}

Exercise~\ref{ex:inversions-inversive} gives a partial converse to the theorem.

\parit{Proof.}
Let $O$ be the center of the inversion, and let $r$ be its radius.

Choose a circline $\Gamma$ and three distinct points $A$, $B$, and $C$ on it.
We may assume that $A \ne O$, $B \ne O$, and $C \ne O$.
(If~$\triangle ABC$ is nondegenerate, then $\Gamma$ is the circumcircle of $\triangle ABC$; if $\triangle ABC$ is degenerate, then $\Gamma$ is the line passing thru $A$, $B$, and~$C$.)

Let $A'$, $B'$, and $C'$ be the inverses of $A$, $B$, and $C$, respectively,
and let $\Gamma'$ be the circline passing thru $A'$, $B'$, and~$C'$.

Assume $D$ is a point in the inversive plane, distinct from $A$, $B$, $C$, $O$, and~$\infty$.
Let $D'$ be the inverse of~$D$.
By Theorem~\ref{lem:inverse-4-angle}\textit{\ref{lem:inverse-4-angle:inscribed}},
$D' \in \Gamma'$ if and only if $D \in \Gamma$.

It remains to prove that
$\infty \in \Gamma$ $\Leftrightarrow$ $O \in \Gamma'$, which is equivalent to
$O \in \Gamma$ $\Leftrightarrow$ $\infty \in \Gamma'$.

The condition $\infty \in \Gamma$ means that $\Gamma$ is a line;
equivalently, for any $\epsilon > 0$, the circline $\Gamma$ contains a point $P \ne \infty$ such that $OP > r/\epsilon$.
For the inversion $P' \in \Gamma'$ of $P$, we have $OP' = r^2 / OP < r \cdot \epsilon$.
That is, the circline $\Gamma'$ contains points arbitrarily close to~$O$.
It follows that $O \in \Gamma'$.
In other words,
$\infty \in \Gamma$
$\Rightarrow$
$O \in \Gamma'$.

The implication
$\infty \in \Gamma$
$\Leftarrow$
$O \in \Gamma'$ is proved in the same way.
We can choose $P' \ne O$ on $\Gamma'$ such that $OP' < r \cdot \epsilon$.
Its inverse $P$ lies on $\Gamma$, and $OP = r^2 / OP' > r/\epsilon$.
Since $\epsilon>0$ is arbitrary, we have $\infty \in \Gamma$.
\qeds

\begin{thm}{Exercise}\label{ex:inv-center not=center-inv}
Assume that a circle $\Gamma'$ is an inverse of a circle $\Gamma$ centered at $Q$.
Suppose $Q'$ denotes the inverse of~$Q$.
Show that $Q'$ is not the center of~$\Gamma'$.
\end{thm}

{

\begin{wrapfigure}{r}{37mm}
\vskip-8mm
\centering
\includegraphics{mppics/pic-162}
\end{wrapfigure}

Assume that a \index{circumtool}\emph{circumtool} is a geometric construction tool 
that produces a circline passing thru three given points.

\begin{thm}{Exercise}\label{ex:circumtool}
Show that with only a circumtool,
it is impossible to construct the center of a given circle.
\end{thm}

}

\begin{thm}{Exercise}\label{ex:tangent-circ->parallels}
Show that for every pair of tangent circles in the inversive plane, there is an inversion that sends them to a pair of parallel lines.
\end{thm}

\begin{thm}{Theorem}\label{thm:inverse}
Consider an inversion of the inversive plane across a circle $\Omega$ centered at~$O$. 
Then 
\begin{enumerate}[(a)]
\item\label{thm:inverse:line-line}
A line passing thru $O$ is inverted into itself.
\item\label{thm:inverse:line} 
A line not passing thru $O$ is inverted into a circle that passes thru $O$, and the other way around.
\item\label{thm:inverse:circle} 
A circle not passing thru $O$ 
is inverted into a circle not passing thru~$O$. 
\end{enumerate}
\end{thm}

\parit{Proof.}
In the proof, we use Theorem~\ref{thm:inverse-cline} without mentioning it.

\parit{(\ref{thm:inverse:line-line}).}
If a line passes thru $O$, then it contains both $\infty$ and~$O$.
Therefore, its inverse also contains $\infty$ and~$O$.
In particular, the image is a line passing thru~$O$.

\parit{(\ref{thm:inverse:line}).}
Since each line $\ell$ passes thru $\infty$, its image $\ell'$ has to contain~$O$.
If the line does not contain $O$, 
then $\ell'\not\ni \infty$;
that is, $\ell'$ is not a line.
Therefore, $\ell'$ is a circle that passes thru~$O$. 

\parit{(\ref{thm:inverse:circle}).}
If the circle $\Gamma$ does not contain $O$, 
then its image $\Gamma'$ does not contain~$\infty$.
Therefore, $\Gamma'$ is a circle.
Since  $\Gamma\not\ni\infty$ we get that $\Gamma' \not\ni O$.
Hence the result.
\qeds

\section{The method of inversion}

Here is an application of inversion,
which we include as an illustration;
we will not use it further in the book.

\begin{thm}{Ptolemy's identity}\label{ptolemy-id}
Let $ABCD$ be an inscribed quadrangle.
Assume that points $A$, $B$, $C$, and $D$ appear on the circline in the same order.
Then 
$$ AB\cdot CD+ BC\cdot DA=AC\cdot BD.$$

\end{thm}

\parit{Proof.}
Assume the points $A$, $B$, $C$, and $D$ lie on one line in this order.

\begin{wrapfigure}{o}{39mm}
\centering
\includegraphics{mppics/pic-164}
\end{wrapfigure}

Set $x\z=AB$, $y=BC$, $z\z=CD$.
Note that
$$x\cdot z+y\cdot (x+y+z)=(x+y)\cdot(y+z).$$
Since $AC\z=x+y$, $BD=y+z$, and $DA\z=x+y+z$,
it proves the identity.

\begin{wrapfigure}{o}{39mm}
\centering
\includegraphics{mppics/pic-166}
\end{wrapfigure}

It remains to consider the case when the quadrangle $ABCD$ is inscribed in a circle, say~$\Gamma$. 

The identity can be rewritten as 
$$\frac{AB\cdot DC}{ BD\cdot CA}+ \frac{BC\cdot AD}{CA\cdot DB}=1.$$
On the left-hand side we have two cross-ratios.
According to Theorem~\ref{lem:inverse-4-angle}\textit{\ref{lem:inverse-4-angle:cross-ratio}}, the left-hand side does not change if we apply an inversion to each point.

Consider an inversion across a circle centered at point $O$ that lies on $\Gamma$ between $A$ and~$D$.
By 
Theorem~\ref{thm:inverse},
this inversion maps $\Gamma$ to a line.
This reduces the problem to the case when $A$, $B$, $C$, and $D$ lie on one line, which was already considered.
\qeds

{

\begin{wrapfigure}{r}{27mm}
\vskip-4mm
\centering
\includegraphics{mppics/pic-168}
\vskip4mm
\includegraphics{mppics/pic-170}
\end{wrapfigure}

In the proof above, we rewrite Ptolemy's identity in a form that is invariant with respect to inversion 
and then apply an inversion which makes the statement evident.
The solution of the following exercise is based on the same idea;
one has to make a right choice of inversion.

\begin{thm}{Exercise}\label{ex:4-circles}
Assume that four circles are mutually tangent to each other.
Show that four (among six) of their points of tangency lie on one circline.
\end{thm}

\begin{thm}{Advanced exercise}\label{ex:inverse}
Assume that three circles are tangent to each other and to two parallel lines as shown in the picture.

Show that $(AB)$ is tangent to two circles at~$A$.
\end{thm}

}

\section{Perpendicular circles}

Assume two circles $\Gamma$ and $\Omega$ intersect at two points $X$ and~$Y$.
Let $\ell$ and $m$ be the tangent lines at $X$ to $\Gamma$ and $\Omega$ respectively.
Analogously, let $\ell'$ and $m'$ be the tangent lines at $Y$ to $\Gamma$ and~$\Omega$.

From Exercise~\ref{ex:two-arcs}, we can deduce that  
 $\ell\perp m$ if and only if $\ell'\perp m'$.

We say that the circle $\Gamma$ is {}\emph{perpendicular} to the circle $\Omega$ 
(briefly \index{38@$\perp$}$\Gamma\perp \Omega$)
if they intersect and the lines tangent to the circles at one point 
(and therefore, both points) 
of intersection are perpendicular.

Similarly, we say that the circle $\Gamma$ is perpendicular to the line $\ell$ (briefly $\Gamma\perp \ell$)
if $\Gamma\cap\ell\ne \emptyset$ and $\ell$ is perpendicular to the tangent lines to $\Gamma$ at one point (and therefore, both points) of intersection.
According to Lemma~\ref{lem:tangent}, 
this happens only if the line $\ell$ passes thru the center of~$\Gamma$.

Now we can talk about \index{perpendicular!circlines}\emph{perpendicular circlines}.

\begin{thm}{Theorem}\label{thm:perp-inverse}
Assume $\Gamma$ and $\Omega$ are distinct circles. 
Then $\Omega\perp\Gamma$ if and only if the circle $\Gamma$ coincides with its inversion across~$\Omega$.
\end{thm}

\begin{wrapfigure}[7]{o}{29mm}
\vskip-4mm
\centering
\includegraphics{mppics/pic-172}
\end{wrapfigure}

\parit{Proof.} 
Suppose that $\Gamma'$ denotes the inverse of~$\Gamma$.

\parit{``Only if'' part.}
Let $O$ be the center of $\Omega$
and $Q$ be the center of~$\Gamma$.
Let $X$ and $Y$ denote the points of intersections of  $\Gamma$ and~$\Omega$.
By Lemma~\ref{lem:tangent}, $\Gamma\perp\Omega$ if and only if $(OX)$ and $(OY)$ are tangent to~$\Gamma$.

Since $O\ne X$, Lemma \ref{lem:perp<oblique} implies that $O$ lies outside of~$\Gamma$.
By Theorem \ref{thm:inverse}\textit{\ref{thm:inverse:circle}}, $\Gamma'$ is a circle.

Note that $\Gamma'$ is also tangent to $(OX)$ and $(OY)$ at $X$ and $Y$ respectively. 
It follows that $X$ and $Y$ are the footpoints of the center of $\Gamma'$ on $(OX)$ and $(OY)$.
Therefore, both $\Gamma'$ and $\Gamma$ have the center~$Q$.
Finally, $\Gamma'=\Gamma$, since both circles pass thru~$X$.

\parit{``If'' part.}
Assume $\Gamma=\Gamma'$.

Since $\Gamma\ne \Omega$, there is a point $P$ that lies on $\Gamma$, but not on~$\Omega$.
Let $P'$ be the inverse of $P$ across~$\Omega$.
Since $\Gamma=\Gamma'$, we have that $P'\in \Gamma$.
In particular, the half-line $[OP)$ intersects $\Gamma$ at two points.
By Exercise~\ref{ex:inside-outside}, 
 $O$ lies outside of~$\Gamma$.

As $\Gamma$ has points inside and outside of $\Omega$,
the circles $\Gamma$ and $\Omega$ intersect.
The latter follows from Exercise~\ref{ex:intersecting-circles-3}.

Let $X$ be a point of their intersection.
We need to show that $(OX)$ is tangent to $\Gamma$, which means $X$ is the only intersection point of $(OX)$ and~$\Gamma$.

Assume $Z$ is another point of intersection of $(OX)$ and~$\Gamma$.
Since $O$ is outside of $\Gamma$, 
the point $Z$ lies on the half-line $[OX)$.

Suppose that $Z'$ denotes the inverse of $Z$ across~$\Omega$.
Clearly, the three points $Z, Z', X$ lie on $\Gamma$ and $(OX)$, which contradicts Lemma~\ref{lem:line-circle}.
\qeds 

It is convenient to define the 
\index{inversion!inversion across the line}\emph{inversion across the line} $\ell$
as the reflection across $\ell$.
This way we can talk about \index{inversion!inversion across the circline}\emph{inversion across an arbitrary circline}.

\begin{thm}{Corollary}\label{cor:perp-inverse-clines}
Let $\Omega$  and $\Gamma$ be distinct circlines in the inversive plane.
Then
the inversion across $\Omega$ sends $\Gamma$ to itself if and only if $\Omega\perp\Gamma$.
\end{thm}

\parit{Proof.}
By Theorem~\ref{thm:perp-inverse}, it is sufficient to consider the case when $\Omega$ or $\Gamma$ is a line.

Assume $\Omega$ is a line, so the inversion across $\Omega$ is a reflection.
In this case, the statement follows from Corollary~\ref{cor:reflection+angle}.

If $\Gamma$ is a line, 
then the statement follows from Theorem~\ref{thm:inverse}.
\qeds

\begin{thm}{Corollary}\label{cor:perp-inverse}
Let $P$ and $P'$ be two distinct points
such that $P'$ is the inverse of $P$ across the circle~$\Omega$.
If a circline $\Gamma$ passes thru $P$ and~$P'$, then $\Gamma\perp\Omega$.
\end{thm}

\parit{Proof.} 
Without loss of generality, we may assume that $P$ is inside and $P'$ is outside~$\Omega$.
By Theorem~\ref{thm:abc}, $\Gamma$ intersects $\Omega$.
Suppose that $A$ denotes a point of intersection.

Suppose that $\Gamma'$ denotes the inverse of~$\Gamma$.
Since $A$ is a self-inverse, the points $A$, $P$, and $P'$ lie on~$\Gamma'$.
By Exercise~\ref{ex:unique-cline},
$\Gamma'=\Gamma$
and by Theorem~\ref{thm:perp-inverse}, $\Gamma\perp\Omega$.
\qeds

\begin{thm}{Corollary}\label{cor:h-line} 
Let $P$ and $Q$ be two distinct points inside a circle~$\Omega$.
Then there is a unique circline $\Gamma$ perpendicular to $\Omega$ that passes thru $P$ and~$Q$.  
\end{thm}

\parit{Proof.}
Let $P'$ be the inverse of the point $P$ across the circle~$\Omega$.
According to Corollary~\ref{cor:perp-inverse},
if a circline that passes thru $P$ and $Q$ is perpendicular to $\Omega$, then it passes thru~$P'$, and the converse holds as well.

Note that $P'$ lies outside of~$\Omega$.
Therefore, the points $P$, $P'$, and $Q$ are distinct.

According to Exercise~\ref{ex:unique-cline},
there is a unique circline passing thru $P$, $Q$, and~$P'$.
Hence the result.
\qeds

\begin{thm}{Exercise}\label{ex:inscribed+inv}
Let $P$, $Q$, $P'$, and $Q'$ be points in the Euclidean plane.
Assume $P'$ and $Q'$ are inverses of $P$ and $Q$ respectively.
Show that the quadrangle $PQP'Q'$ is inscribed.
\end{thm}

\begin{thm}{Exercise}\label{ex:centers-of-perp-circles}
Let $\Omega_1$ and $\Omega_2$ be two perpendicular circles with centers at $O_1$ and $O_2$ respectively.
Show that the inverse of $O_1$ across $\Omega_2$ 
coincides with 
the inverse of $O_2$ across~$\Omega_1$.
\end{thm}

\begin{thm}{Exercise}\label{ex:4-th-perp-circ}
Three distinct circles --- $\Omega_1$, $\Omega_2$, and $\Omega_3$, intersect at two points --- $A$ and~$B$.
Assume that a circle $\Gamma$ is perpendicular to $\Omega_1$ and $\Omega_2$.
Show that $\Gamma\perp\Omega_3$.
\end{thm}

Let us consider two new construction tools:
the \index{circumtool}\emph{circumtool} that constructs a circline thru three given points, 
and the \index{inversor}\emph{inversor} --- a tool that constructs an inverse of a given point across a given circline.

\begin{thm}{Exercise}\label{ex:construction-perp-clines}
Given  two circles $\Omega_1$, $\Omega_2$ and a point $P$ that does not lie on the circles,
use only circumtool and inversor to construct a circline $\Gamma$ thru $P$, 
and perpendicular to both $\Omega_1$ and $\Omega_2$.
\end{thm}

\begin{thm}{Advanced exercise}\label{ex:3-construction-perp-clines}
Given  three disjoint circles $\Omega_1$, $\Omega_2$, and $\Omega_3$,
use only circumtool and inversor to construct a circline $\Gamma$ that is perpendicular to each circle $\Omega_1$, $\Omega_2$, and $\Omega_3$.

Think what to do if two of the circles intersect.
\end{thm}

\section{Angles after inversion}

\begin{thm}{Proposition}
In the inversive plane, an inverse of an arc is an arc.

\end{thm}

\parit{Proof.} 
Consider four distinct points $A$, $B$, $C$, and $D$; 
let $A'$, $B'$, $C'$, and $D'$  be their inverses.
We need to show that $D$ lies on the arc $ABC$ if and only if $D'$ lies on the arc $A'B'C'$.
According to Proposition~\ref{prop:arcs},
the latter is equivalent to the following:
$$\measuredangle ADC= \measuredangle ABC
\quad
\iff
\quad  
\measuredangle A'D'C'= \measuredangle A'B'C'.$$
The latter follows from Theorem~\ref{lem:inverse-4-angle}\textit{\ref{lem:inverse-4-angle:angle}}.
\qeds

The following theorem states that the angle between arcs changes only its sign after the inversion.

{

\begin{wrapfigure}{r}{55mm}
\vskip-6mm
\centering
\includegraphics{mppics/pic-174}
\end{wrapfigure}

\begin{thm}{Theorem}\label{thm:angle-inversion}
Let $AB_1C_1$, $AB_2C_2$ be two arcs in the inversive plane,
and arcs $A'B_1'C_1'$, $A'B_2'C_2'$ be their inverses.
Let $[AX_1)$ and $[AX_2)$ be the half-lines tangent to $AB_1C_1$ and  $AB_2C_2$ at $A$,
and
$[A'Y_1)$ and $[A'Y_2)$ be the half-lines tangent to $A'B_1'C_1'$ and  $A'B_2'C_2'$ at~$A'$.
Then
$$\measuredangle X_1AX_2\equiv-\measuredangle Y_1A'Y_2.$$

\end{thm}

}

{

\begin{wrapfigure}[6]{o}{21mm}
\vskip-6mm
\centering
\includegraphics{mppics/pic-176}
\end{wrapfigure}

The \index{angle!between arcs}\emph{angle between arcs} can be defined as the angle between their tangent half-lines at the common endpoint.
Therefore under inversion, the angles between arcs are preserved up to sign.

From Exercise~\ref{ex:tangent-lim}, it follows that the angle between arcs with the common endpoint $A$ is the limit of $\measuredangle P_1AP_2$ where $P_1$ and $P_2$ are points approaching $A$ along the corresponding arcs. 
This observation can be used to define the angle between a pair of curves emerging from one point.
It turns out that under inversion, angles between curves are also preserved up to sign.

}

\parit{Proof.}
By Proposition~\ref{prop:arc(angle=tan)},
\begin{align*}
\measuredangle X_1AX_2&\equiv\measuredangle X_1AC_1+\measuredangle C_1AC_2+\measuredangle C_2AX_2\equiv
\\
&\equiv(\pi-\measuredangle C_1B_1A)+\measuredangle C_1AC_2+(\pi-\measuredangle AB_2 C_2)\equiv
\\
&
\equiv -(\measuredangle C_1B_1A+\measuredangle AB_2 C_2 +\measuredangle C_2 A C_1)\equiv
\\
&\equiv 
-(\measuredangle C_1B_1A+\measuredangle AB_2 C_1)
-(\measuredangle C_1B_2C_2 +\measuredangle C_2 A C_1).
\intertext{The same way, we get}
\measuredangle Y_1A'Y_2
&\equiv-(\measuredangle C_1'B_1'A'+\measuredangle A'B_2' C_1')
-(\measuredangle C_1'B_2'C_2' +\measuredangle C_2' A' C_1').
\end{align*}

By Theorem~\ref{lem:inverse-4-angle}\textit{\ref{lem:inverse-4-angle:angle}},
\begin{align*}
\measuredangle C_1B_1A+\measuredangle AB_2 C_1&\equiv-(\measuredangle C_1'B_1'A'+\measuredangle A'B_2' C_1'),
\\
\measuredangle C_1B_2C_2 +\measuredangle C_2 A C_1&\equiv-(\measuredangle C_1'B_2'C_2' +\measuredangle C_2' A' C_1'),
\end{align*}
and hence the result.\qeds

\begin{thm}{Corollary}\label{cor:invese-comp}
Let $P$ be an inverse of a point $Q$ across a circle $\Gamma$.
Assume that $P'$, $Q'$, and $\Gamma'$ 
are the inverses of  $P$, $Q$, and $\Gamma$ across another circle $\Omega$.
Then $P'$ is the inverse  of $Q'$ across~$\Gamma'$.
\end{thm}

{

\begin{wrapfigure}{r}{45mm}
\vskip-6mm
\centering
\includegraphics{mppics/pic-178}
\end{wrapfigure}

\parit{Proof.}
If $P=Q$, then $P'=Q'\z\in\Gamma'$. 
Therefore, $P'$ is the inverse of $Q'$ across~$\Gamma'$.

It remains to consider the case $P\z\ne Q$. 
Let $\Delta_1$ and $\Delta_2$ be two distinct circles that intersect at $P$ and~$Q$.
According to Corollary~\ref{cor:perp-inverse}, 
$\Delta_1\perp\Gamma$ and $\Delta_2\perp\Gamma$.

Let $\Delta_1'$ and $\Delta_2'$ denote the inverses of $\Delta_1$ and $\Delta_2$ across~$\Omega$.
Clearly, $\Delta_1'$ meets $\Delta_2'$ at $P'$ and~$Q'$.

By Theorem~\ref{thm:angle-inversion},  $\Delta_1'\perp\Gamma'$ and $\Delta_2'\z\perp\Gamma'$.
By Corollary~\ref{cor:perp-inverse-clines}, $P'$ is the inverse of $Q'$ across~$\Gamma'$.
\qeds

}

\addtocontents{toc}{\protect\contentsline{part}{\protect\numberline{}Non-Euclidean geometry}{}{}}

\chapter{Neutral plane}\label{chap:non-euclid}

Let us remove Axiom~\ref{def:birkhoff-axioms:4} from our axiomatic system (Section~\ref{sec:axioms}).
This way we define a new object called the 
\index{plane!neutral plane}\index{neutral plane}\emph{neutral plane} or \index{plane!absolute plane}\index{absolute plane}\emph{absolute plane}.
(In a neutral plane, Axiom~\ref{def:birkhoff-axioms:4} may or may not hold.)

Every theorem in neutral geometry holds in Euclidean geometry.
In other words, the Euclidean plane is an example of a neutral plane. 
In the next chapter, we will construct an example of a neutral plane that is not Euclidean.

In this book, 
Axiom~\ref{def:birkhoff-axioms:4} was used starting from Chapter~\ref{chap:parallel}.
Therefore all the statements before hold in neutral geometry.
This makes all the discussed results
about
half-planes,
signs of angles,
congruence conditions,
perpendicular lines,
and reflections 
true in neutral geometry.
Recall that a statement is marked with ``$\a$''\label{a-mark} (for example, ``\textbf{Theorem.\abs}'') if it holds in every neutral plane, and the same proof works.

Let us give an example of a theorem in neutral geometry that admits a simpler proof in Euclidean geometry. 

\begin{thm}{Hypotenuse-leg congruence condition}\label{thm:hypotenuse-leg}
Assume that triangles $ABC$ and $A'B'C'$
have right angles at $C$ and $C'$ respectively, 
$AB\z=A'B'$ and $AC\z=A'C'$.
Then $\triangle ABC\cong\triangle A'B'C'$.
\end{thm}

\parit{Euclidean proof.} 
By the Pythagorean theorem $BC=B'C'$.
Then the statement follows from the SSS congruence condition.
\qeds

The proof of the Pythagorean theorem used properties of similar triangles, which in turn used Axiom~\ref{def:birkhoff-axioms:4}. 
Therefore this proof does not work in a neutral plane.

\parit{Neutral proof.}
Suppose that $D$ denotes the reflection of $A$ across $(BC)$
and $D'$ denotes the reflection of $A'$ across $(B'C')$.
Note that 
$$
AD=2\cdot AC=2\cdot A'C'=A'D',\qquad
BD=BA=B'A'=B'D'.
$$

{

\begin{wrapfigure}{r}{28mm}
\vskip-0mm
\centering
\includegraphics{mppics/pic-180}
\end{wrapfigure}

By the SSS congruence condition (\ref{thm:SSS}), 
we get that $\triangle ABD\cong \triangle A'B'D'$.

The statement follows since $C$ is the midpoint of $[AD]$
and $C'$ is the midpoint of $[A'D']$.  
\qeds

\begin{thm}{Exercise}\label{ex:tangent-angle-neutral}
Read the proof of Theorem~\ref{thm:tangent-angle} and identify the first statement that does not work in the netral plane.
\end{thm}

}

\begin{thm}{Exercise}\label{ex:abs-bisect=median}
Give a proof of Exercise~\ref{ex:bisect=median}
that works in the neutral plane. 
\end{thm}

\begin{thm}{Exercise}\label{ex:abs-inscibed}
Let $ABCD$ be an inscribed quadrangle in the neutral plane.
Show that
$$\measuredangle ABC+\measuredangle CDA\equiv \measuredangle BCD+\measuredangle DAB.$$

\end{thm}

One cannot use Corollary~\ref{cor:inscribed-quadrangle} to solve the exercise above since it uses Theorems~\ref{thm:tangent-angle} and \ref{thm:inscribed-angle},
which in turn uses Theorem~\ref{thm:3sum}.

\section{Two angles of a triangle}

In this section, we will prove a weaker form of Theorem~\ref{thm:3sum}
which holds in every neutral plane.

\begin{thm}{Proposition}\label{prop:2sum}
Let $\triangle ABC$ be a nondegenerate triangle in the neutral plane.
Then 
$$|\measuredangle CAB|+|\measuredangle ABC|< \pi.$$

\end{thm}

\begin{wrapfigure}{r}{23mm}
\vskip-1mm
\centering
\includegraphics{mppics/pic-182}
\end{wrapfigure}

According to \ref{thm:signs-of-triug}, the angles $ABC$, $BCA$, and $CAB$
have the same sign.
Therefore, in the Euclidean plane, the theorem follows immediately from Theorem~\ref{thm:3sum}.

\parit{Proof.}
Let $X$ be the reflection of $C$ across the midpoint $M$ of $[AB]$.
Applying \ref{ex:between}, we get
\[|\measuredangle CAX|=|\measuredangle CAB|+|\measuredangle ABC|.\eqlbl{eq:CAX=CAB+ABC}\]

By Proposition~\ref{prop:point-reflection+}
$\measuredangle BAX\z=\measuredangle ABC$.

Note that $\measuredangle CAX\z{\not\equiv} \pi$; otherwise, $X$ would lie on $(AC)$.
Therefore the identity \ref{eq:CAX=CAB+ABC} implies that
\[|\measuredangle CAB|+|\measuredangle ABC|=|\measuredangle CAX|<\pi.\]
\qedsf

\begin{thm}{Exercise}\label{ex:parallel-abs}
Assume $A$, $B$, $C$, and $D$ are points in a neutral plane
such that 
$$2\cdot \measuredangle ABC+2\cdot\measuredangle BCD\equiv 0.$$
Show that $(AB)\parallel (CD)$.
\end{thm}

Note that one cannot apply the transversal property (\ref{thm:parallel-2}) here.

\begin{thm}{Exercise}\label{ex:SAA}
Prove the \index{SAA congruence condition}\emph{side-angle-angle} congruence condition in neutral geometry.

In other words, let $ABC$ and $A'B'C'$ be two triangles in a neutral plane;
suppose that $\triangle A'B'C'$ is nondegenerate.
Show that $\triangle ABC\z\cong \triangle A'B'C'$
if 
$$AB=A'B',
\quad  
\measuredangle ABC=\pm\measuredangle A'B'C',
\quad 
\text{and}
\quad
\measuredangle BCA=\pm\measuredangle B'C'A'.$$

\end{thm}

In the Euclidean plane, the above exercise follows from ASA and the theorem on the sum of angles of a triangle (\ref{thm:3sum}).
However, Theorem~\ref{thm:3sum} cannot be used here, since its proof uses Axiom~\ref{def:birkhoff-axioms:4}.
Later (Theorem~\ref{thm:3sum-h}) 
we will show that Theorem~\ref{thm:3sum} does not hold in a neutral plane.

\begin{thm}{Exercise}\label{ex:chev<side}
Assume that point $D$ lies between vertices $A$ and $B$ of $\triangle ABC$ in the neutral plane.
Show that 
$$CD<CA
\quad
\text{or}
\quad
CD<CB.$$

\end{thm}

\section{Three angles of a triangle}

\begin{thm}{Proposition}\label{prop:angle-side}
Let $\triangle ABC$ and $\triangle A'B'C'$ be two triangles in the neutral plane
such that $AC=A'C'$ and $BC=B'C'$.
Then 
$$AB\z<A'B'
\quad
\text{if and only if}
\quad 
|\measuredangle ACB|<|\measuredangle A'C'B'|.$$

\end{thm}

\begin{wrapfigure}{o}{38mm}
\vskip-4mm
\centering
\includegraphics{mppics/pic-184}
\end{wrapfigure}

\parit{Proof.}
Without loss of generality, we may assume that $A=A'$, $C=C'$, and $\measuredangle ACB,\measuredangle ACB'\ge 0$.
In this case, we need to show that 
$$AB<AB'
\ 
\iff
\  
\measuredangle ACB<\measuredangle ACB'.$$

Choose a point $X$ so that 
$$\measuredangle ACX=\tfrac12\cdot(\measuredangle ACB+\measuredangle ACB').$$
Note that 
\begin{itemize}
\item $(CX)$ bisects $\angle BCB'$.
\item $(CX)$ is the perpendicular bisector of $[BB']$.
\item $A$ and $B$ lie on the same side of $(CX)$ if and only if $$\measuredangle ACB<\measuredangle ACB'.$$
\end{itemize}
From Exercise~\ref{ex:pbisec-side}, $A$ and $B$ lie on the same side of $(CX)$ if and only if $AB<AB'$.
Hence the result.
\qeds

\begin{thm}{Theorem}\label{thm:3sum-a}
Let $\triangle ABC$ be a triangle in the neutral plane.
Then 
$$|\measuredangle ABC|+|\measuredangle BCA|+|\measuredangle CAB|\le \pi.$$

\end{thm}

The following proof is due to Adrien-Marie Legendre \cite{legendre}, 
earlier proofs were given by Giovanni Saccheri \cite{saccheri}
and Johann Lambert \cite{lambert}.

\parit{Proof.} 
Set 
\begin{align*}
a&=BC,
&
b&=CA,
&
c&=AB,
\\
\alpha&=\measuredangle CAB,
&
\beta&=\measuredangle ABC,
&
\gamma&=\measuredangle BCA.
\end{align*}
Without loss of generality, we may assume that $\alpha,\beta,\gamma\ge 0$.

\begin{figure}[!ht]
\centering
\includegraphics{mppics/pic-186}
\end{figure}

Fix a positive integer~$n$.
Consider points $A_0$, $A_1,\dots,A_n$ on the half-line
$[BA)$, such that $BA_i=i\cdot c$ for each~$i$.
(In particular, $A_0=B$ and $A_1=A$.)
Construct points $C_1$, $C_2,\dots,C_n$,
so that
$\measuredangle A_iA_{i-1}C_i\z=\beta$ and $A_{i-1}C_i=a$ for each~$i$.

By SAS, we have constructed $n$ congruent triangles 
\begin{align*}
\triangle ABC&=\triangle A_{1}A_0C_1\cong\triangle A_2A_{1}C_2\cong
\dots
\cong\triangle A_nA_{n-1}C_n.
\end{align*}

Set $d=C_1C_2$ and $\delta=\measuredangle C_2A_1C_1$.
Note that 
$$\alpha+\beta+\delta=\pi.
\eqlbl{eq:gamma'}$$
By Proposition~\ref{prop:2sum}, we get that $\delta\ge 0$.

By construction
\begin{align*}
\triangle A_1C_1C_2&\cong\triangle A_{2}C_2C_3\cong\dots
\cong\triangle A_{n-1}C_{n-1}C_n.
\end{align*}
In particular, $C_iC_{i+1}=d$ 
for each~$i$.

Applying the triangle inequality several times, we get that
\begin{align*}
n\cdot c&=A_0A_n\le 
\\
&\le A_0C_1+C_1C_2+\dots+C_{n-1}C_n+C_nA_n=
\\
&=a+(n-1)\cdot d+b.
\end{align*}

In particular, 
\[c\le  d+\tfrac1n\cdot (a+b-d).\]
Since  $n$ is an arbitrary positive integer,
the latter implies
$c\le d$.
By Proposition~\ref{prop:angle-side}, 
it is equivalent to 
\[\gamma\le \delta.\] 
From \ref{eq:gamma'}, 
the theorem follows.
\qeds

\begin{thm}{Exercise}\label{ex:neutral-quadrangle}
Let $ABCD$ be a quadrangle in the neutral plane.
Suppose that $\angle DAB$ and $\angle ABC$ are right.
Show that $AB\le CD$.
\end{thm}

\section{The defect}\label{The defect}

The \index{defect of triangle}\emph{defect of triangle} $\triangle ABC$ is defined as 
$$\defect(\triangle ABC)
\df 
\pi-|\measuredangle ABC|-|\measuredangle BCA|-|\measuredangle CAB|.$$

Theorem~\ref{thm:3sum-a} states that \textit{the defect of each triangle in a neutral plane has to be nonnegative}.
According to Theorem~\ref{thm:3sum}, every triangle in
the Euclidean plane has zero defect.

{

\begin{wrapfigure}{o}{24mm}
\vskip-6mm
\centering
\includegraphics{mppics/pic-188}
\vskip4mm
\includegraphics{mppics/pic-189}
\end{wrapfigure}

\begin{thm}{Classroom exercise}\label{ex:defect}
Let $\triangle ABC$ be a nondegenerate triangle in the neutral plane.
Assume $D$ lies between $A$ and~$B$.
Show that 
$$\defect(\triangle ABC)=\defect(\triangle ADC)+\defect(\triangle DBC).$$

\end{thm}

\begin{thm}{Exercise}\label{ex:defect=} Let $ABC$ be a nondegenerate triangle in the neutral plane.
Suppose $X$ is a reflection of $C$ across a midpoint $M$ of $[AB]$.
Show that 
$$\defect(\triangle ABC)=\defect(\triangle AXC).$$
\end{thm}

}

\vskip-2mm

\begin{thm}{Exercise}\label{ex:neutral-rectangle}
Let $ABCD$ be a \index{rectangle}\emph{rectangle} in the neutral plane;
that is, $ABCD$ is a quadrangle with all right angles.
Show that $AB\z=CD$. 
\end{thm}

\begin{thm}{Advanced exercise}\label{ex:neutral-rectangle+}
Show that if the neutral plane has a rectangle, then all its triangles have zero defect.
\end{thm}

\section{Proving that something cannot be proved}
\label{sec:unprovable}

Many attempts were made to prove that every theorem in Euclidean geometry holds in neutral geometry.
The latter is equivalent to the statement that Axiom~\ref{def:birkhoff-axioms:4} is a \textit{theorem} in neutral geometry.

Some of these attempts were accepted as proof for long periods until a mistake was found.

Many statements in neutral geometry are equivalent to Axiom~\ref{def:birkhoff-axioms:4}.
It means that if we exchange Axiom~\ref{def:birkhoff-axioms:4} for each of these statements, then we will obtain an equivalent axiomatic system.

The following theorem provides a short list of such statements.
We are not going to prove it in the book.

\begin{thm}{Theorem}\label{thm:=IV}
The neutral plane is Euclidean if and only if one of the following equivalent conditions holds:
\begin{enumerate}[(a)]
\item\label{thm:=IV:main} 
There is a line $\ell$ 
and a point $P\notin\ell$ 
such that there is only one line passing thru $P$ 
and parallel to~$\ell$.
\item 
Every nondegenerate triangle can be circumscribed.
\item
There exists a pair of distinct lines that lie at a bounded distance from each other.
\item
There is a triangle with an arbitrarily large inradius.
\item\label{thm:=IV:defect}
There is a nondegenerate triangle with zero defect.
\item\label{thm:=IV:rectangle}
There exists a \index{rectangle}\emph{rectangle}; that is, a quadrangle with all right angles.
\end{enumerate}
\end{thm}

It is hard to imagine a neutral plane that does not satisfy some of the properties above.
That is partly the reason for a large number of false proofs;
each used one of such statements by accident.

Let us formulate the negation of \textit{(\ref{thm:=IV:main})} above as a new axiom;
we label it h-$\!$\ref{def:birkhoff-axioms:4} as a \textit{hyperbolic version} of Axiom~\ref{def:birkhoff-axioms:4}.

\begin{framed}
\begin{description}
\item[{\rm h-$\!$\ref{def:birkhoff-axioms:4}.}]\label{def:hyperbolic-4a}  
For every line $\ell$ and every point $P\notin\ell$
there are at least two lines that pass thru $P$ and parallel to~$\ell$.
\end{description}
\end{framed}

By Theorem~\ref{thm:parallel}, a neutral plane that satisfies Axiom~h-$\!$\ref{def:birkhoff-axioms:4} is not Euclidean. 
Moreover, according to Theorem~\ref{thm:=IV} (which we do not prove) 
in every non-Euclidean neutral plane, Axiom~h-$\!$\ref{def:birkhoff-axioms:4} holds.

It opens a way to look for a proof by contradiction.
Simply exchange  Axiom~\ref{def:birkhoff-axioms:4} to Axiom~h-$\!$\ref{def:birkhoff-axioms:4}
 and start to prove theorems in the obtained axiomatic system.
In the case if we arrive at a contradiction, 
we prove Axiom~\ref{def:birkhoff-axioms:4} in a neutral plane.
This idea was growing since the $5^\text{th}$ century;
the most notable results were obtained by Giovanni Saccheri \cite{saccheri}.

The system of axioms \ref{def:birkhoff-axioms:0}--\ref{def:birkhoff-axioms:3} and h-$\!$\ref{def:birkhoff-axioms:4} defines a new geometry which is now called \index{hyperbolic!geometry}\emph{hyperbolic} or \index{Lobachevsky  geometry}\emph{Lobachevsky  geometry}.
The more this geometry was developed,
it became more and more believable that there is no contradiction;
that is, the system of axioms \ref{def:birkhoff-axioms:0}--\ref{def:birkhoff-axioms:3}, and h-$\!$\ref{def:birkhoff-axioms:4} is \index{consistent}\emph{consistent}.
In fact, the following theorem holds true:

\begin{thm}{Theorem}\label{thm:consistent}
Hyperbolic geometry is consistent if and only if so is Euclidean geometry.
\end{thm}

The claims
that hyperbolic geometry has no contradiction can be found in the private letters of
Carl Friedrich Gauss, 
Ferdinand  Schweikart, 
and Franz Taurinus.%
\footnote{The oldest surviving letters were the Gauss letter to Christian Gerling in 1816 
and the yet more convincing letter dated 1818 
of Schweikart sent to Gauss via Gerling.}
They all seem to be afraid to state it in public.
For instance, in 1818 Gauss writes to Gerling:

\smallskip

\begin{quotation}{\it
\dots I am happy that you have the courage to express yourself as if you recognized the possibility that our parallels theory along with our entire geometry could be false.
But the wasps whose nest you disturb will fly around your head.}
\end{quotation}

\smallskip

Nikolai Lobachevsky came to the same conclusion independently.
Unlike the others, he dared to state it in public and in print \cite{lobachevsky}.
That cost him serious trouble.
A couple of years later, also independently, János Bolyai published his work \cite{bolyai}.

It seems that Lobachevsky was the first who had a proof of Theorem~\ref{thm:consistent} altho its formulation required rigorous axiomatics which was not developed at his time.
Later, Beltrami gave a cleaner proof of the ``if'' part of the theorem.
It was done by modeling points, lines, distances, and angle measures of one geometry using some other objects in another geometry.
The same idea was used earlier by Lobachevsky \cite[\S34]{lobachevsky-1840}; 
he modeled the Euclidean plane in the hyperbolic space.

The proof of Beltrami is the subject of the next chapter. 



\section{Curvature}

In a letter from 1824 Gauss writes: 

\begin{quotation}{\it
The assumption that the sum of the three angles is less than $\pi$ leads to a curious geometry, 
quite different from ours but completely consistent, 
which I have developed to my entire satisfaction, 
so that I can solve every problem in it with the exception of a determination of a constant, which cannot be designated a priori. 
The greater one takes this constant, the nearer one comes to Euclidean geometry, 
and when it is chosen indefinitely large the two coincide.
The theorems of this geometry appear to be paradoxical and, 
to the uninitiated, absurd; but calm, steady reflection reveals that they contain nothing at all impossible. 
For example, the three angles of a triangle become as small as one wishes, if only the sides are taken sufficiently large; 
yet the area of the triangle can never exceed a definite limit, regardless how great the sides are taken, 
nor indeed can it ever reach it.}
\end{quotation} 

In modern terminology, the constant that Gauss mentions 
can be expressed as $1/\sqrt{-k}$, 
where $k\le 0$, is the so-called \index{curvature}\emph{curvature} of the neutral plane, which we are about to introduce.

The identity in Exercise~\ref{ex:defect} suggests that the defect of a triangle should be proportional to its area.%
\footnote{The area in the neutral plane is discussed briefly at the end of Chapter~\ref{chap:area},
but the reader could also refer to an intuitive understanding of area measurement.}

In fact, for every neutral plane, there is a nonpositive real number $k$
such that 
$$k\cdot\area(\triangle ABC)+\defect(\triangle ABC)=0$$
for every $\triangle ABC$.
This number $k$ is called the \index{curvature}\emph{curvature} of the plane.

For example, by Theorem~\ref{thm:3sum}, the Euclidean plane has zero curvature.
By Theorem~\ref{thm:3sum-a}, the curvature of every neutral plane is nonpositive.

It turns out that up to isometry, the neutral plane is characterized by its curvature;
that is, two neutral planes are isometric if and only if they have the same curvature. 

In the next chapter, we will construct a {}\emph{hyperbolic plane};
this is, an example of a neutral plane with curvature $k=-1$.

Any neutral plane, distinct from Euclidean,
can be obtained by rescaling the metric on the hyperbolic plane.
Indeed,
if we rescale the metric by a positive factor $c$,
the area changes by factor $c^2$, while the defect stays the same.
Therefore, taking $c=\sqrt{-k}$,
we can get the neutral plane of the given curvature $k<0$.
In other words, all the non-Euclidean neutral planes become identical
if we use $r=1/\sqrt{-k}$ as the unit of length.

\medskip

In Chapter~\ref{chap:sphere}, we discuss spherical geometry.
Altho spheres are not neutral planes,
the spherical geometry is a close relative of Euclidean and hyperbolic geometries.

Nondegenerate spherical triangles have negative defects.
Moreover, 
if $r$ is the radius of the sphere, then
$$\tfrac1{r^2}\cdot\area(\triangle ABC)+\defect(\triangle ABC)=0$$
for every spherical triangle $ABC$; compare to Lemma~\ref{lem:area-spher-triangle}.
In other words, 
the sphere of radius $r$ has the curvature $k=\tfrac1{r^2}$.

\chapter{Hyperbolic plane}\label{chap:poincare}

\begin{figure}[!ht]
\vspace*{-9mm}
\centering
\includegraphics[scale=0.25]{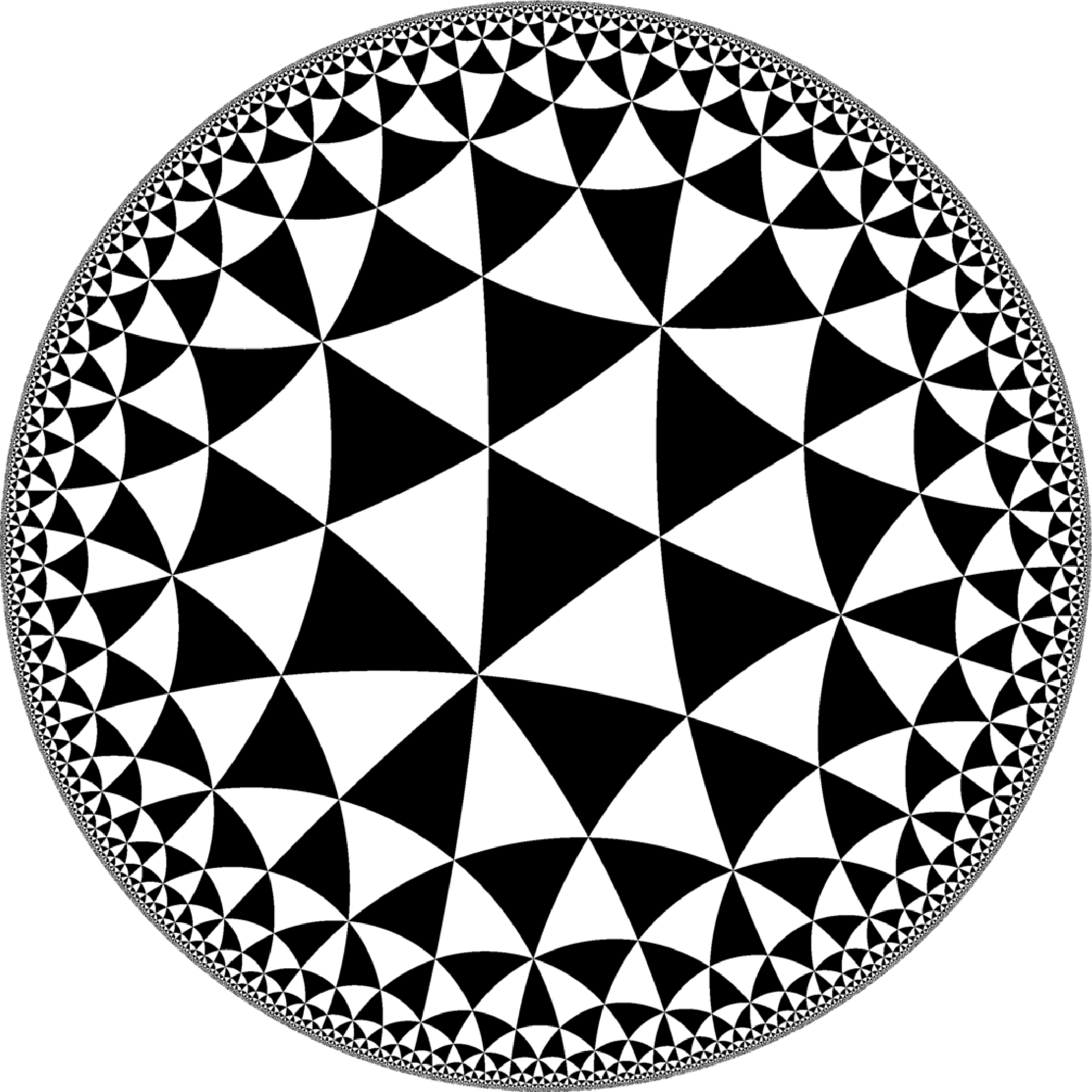}
\end{figure}

In this chapter, we use inversive geometry 
to construct a model of a hyperbolic plane --- a neutral plane that is not Euclidean.

Namely, we construct the so-called \index{conformal disc model}\emph{conformal disc model} of the hyperbolic plane.
This model was discovered by Eugenio Beltrami \cite{beltrami},
and it is commonly referred to as the {}\emph{Poincar\'e disc model}. 

The figure above shows the conformal disc model of the hyperbolic plane which is divided into congruent triangles with angles $\tfrac\pi3$, $\tfrac\pi3$, and~$\tfrac\pi4$.

\section{The conformal disc model}
\label{sec:conformal-model}

In this section, we give new names for certain objects in the Euclidean plane
which will represent lines, angle measures, and distances in the hyperbolic plane.

\parbf{Hyperbolic plane.}
Let us fix a circle on the Euclidean plane 
and call it \index{absolute}\emph{absolute}.
The set of points inside the absolute will be referred to as the \index{hyperbolic!plane}\index{plane!hyperbolic plane}\emph{hyperbolic plane} (or \index{plane!h-plane}\index{h-plane}\emph{h-plane}).

The points on the absolute do \textit{not} belong to the h-plane.
The points in the h-plane will be also called \index{h-point}\emph{h-points}.

We will often assume that the absolute is a unit circle.

\parbf{Hyperbolic lines.}
The intersections of the h-plane with circlines that are perpendicular to the absolute are called {}\emph{hyperbolic lines} or \index{h-line}\emph{h-lines}.

\begin{wrapfigure}{o}{48mm}
\centering
\includegraphics{mppics/pic-190}
\end{wrapfigure}

By Corollary~\ref{cor:h-line}, there exists a unique h-line that passes thru any two distinct h-points $P$ and~$Q$.
This h-line will be denoted by~\index{62@$(PQ)_h$, $[PQ)_h$,$[PQ]_h$}$(PQ)_h$.

The arcs of hyperbolic lines will be called {}\emph{hyperbolic segments} or \index{h-segment}\emph{h-segments}.
An h-segment with endpoints $P$ and $Q$ will be denoted as~$[PQ]_h$.

The portion of an h-line on one side of a point will be called a {}\emph{hyperbolic half-line} (or \index{h-half-line}\emph{h-half-line}).
More precisely, an h-half-line is the intersection of the h-plane with an arc that is perpendicular to the absolute and has exactly one of its endpoints in the h-plane.
An h-half-line starting at $P$ and passing thru $Q$ will be denoted as~$[PQ)_h$.

If $\Gamma$ is the circline containing the h-line $(PQ)_h$, then the points ($A$ and $B$ in the picture) where $\Gamma$ intersects the absolute are called
\index{point!ideal point}\index{ideal!point}\emph{ideal points} of~$(PQ)_h$.
(The ideal points of an h-line do not belong to the h-line.)

An ordered triple of h-points, say $(P,Q,R)$, will be called an {}\emph{h-triangle $PQR$} and denoted by \index{21@$\triangle_h$}$\triangle_h P Q R$.

It is important to clarify that at this point, that \textit{so far an h-line $(PQ)_h$ is just a subset of the h-plane},
but soon we will introduce h-distance 
and show that $(PQ)_h$ is a line for the h-distance in the sense of the Definition~\ref{def:line}. 

\begin{thm}{Exercise}\label{ex:ideal-line-unique}
Show that an h-line is uniquely determined by its ideal points.
\end{thm}

\begin{thm}{Exercise}\label{ex:1ideal-line-unique}
Show that an h-line is uniquely determined by one of its ideal points and one h-point on it.
\end{thm}

\begin{thm}{Classroom exercise}\label{ex:line/h-line}
Show that an h-segment $[PQ]_h$ coincides with the Euclidean segment $[PQ]$
if and only if the line $(PQ)$ passes thru the center of the absolute.
\end{thm}

\parbf{Hyperbolic distance.}\label{h-dist}
Consider two distinct h-points $P$ and $Q$,
and let $A$ and $B$ denote the ideal points of $(PQ)_h$.
Without loss of generality, we may assume that on the Euclidean circline containing the h-line $(PQ)_h$, the points $A,P,Q,B$ appear in the same order.

Consider the function 
$$\delta(P,Q)\df\frac{AQ\cdot PB}{AP\cdot QB}.$$
The right-hand side is a cross-ratio;
by Theorem~\ref{lem:inverse-4-angle} it is invariant under inversion.
We set $\delta(P,P)=1$ for any h-point~$P$.
Let us define h-distance as the logarithm of $\delta$; that is,
$$PQ_h\df\ln[\delta(P,Q)].$$

The proof that $PQ_h$ is a metric on the h-plane will be given later.
For now, it is just a function that returns a real value $PQ_h$ for any pair of h-points $P$ and~$Q$.

\begin{thm}{Exercise}\label{ex:h-dist-eq}
Let $O$ be the center of the absolute.
Assume that h-points $O$, $X$, and $Y$ lie on an h-line in the same order, and $OX\z=XY$.
Prove that $OX_h<XY_h$.
\end{thm}

\parbf{Hyperbolic angles.}\label{h-angle measure}
Consider three h-points $P$, $Q$, and $R$
such that $P\ne Q$ and $R\ne Q$.
The \index{hyperbolic!angle}\emph{hyperbolic angle $PQR$} (briefly $\angle_h PQR$)\index{11@$\angle_h$, $\measuredangle_h$} consists of the ordered pair of h-half-lines $[QP)_h$ and $[QR)_h$.

Let $[QX)$ and $[QY)$ be (Euclidean) half-lines 
that are tangent to $[QP]_h$ and $[QR]_h$ 
at~$Q$.
Then the \index{angle!measure!hyperbolic angle measure}\index{hyperbolic!angle measure}\emph{hyperbolic angle measure} (or \index{h-angle measure}\emph{h-angle measure}) of $\angle_h PQR$ is denoted by
$\measuredangle_h PQR$ and is defined as
$\measuredangle XQY$.

\begin{thm}{Exercise}\label{ex:h-perp-unique}
Let $P$ be an h-point that does not lie on an h-line~$\ell$.
Show that there is a unique h-line thru $P$ that is perpendicular to~$\ell$.
\end{thm}

\section{Plan of the proof}

We have introduced all the {}\emph{h-notions} needed in the formulation of the axioms \ref{def:birkhoff-axioms:0}--\ref{def:birkhoff-axioms:3} and h-\ref{def:birkhoff-axioms:4}.
It remains to show that all these axioms hold; 
this will be done by the end of this chapter.

Once our proofs are complete, we will have a model that serves as an example of a neutral plane.
After that we can use axiomatic approach in the h-plane;
for example, Exercise~\ref{ex:h-perp-unique} can be proved the same way as Theorem~\ref{perp:ex+un}.

Most importantly we will prove the ``if''-part of Theorem~\ref{thm:consistent}.

Indeed, any statement in hyperbolic geometry can be restated in the Euclidean plane using the introduced h-notions.
Consequently, if the system of axioms \ref{def:birkhoff-axioms:0}--\ref{def:birkhoff-axioms:3}, and h-\ref{def:birkhoff-axioms:4} leads to a contradiction, then so does the system axioms \ref{def:birkhoff-axioms:0}--\ref{def:birkhoff-axioms:4}.

\section{Auxiliary statements}

One may compare the conformal model with a telescope --- it makes it possible to see the h-plane from the Euclidean plane.
Continuing this analogy further, we may say that the following lemma will be used to \textit{aim} the telescope at any particular point in the h-plane.

\begin{thm}{Lemma}\label{lem:P-->O} 
Consider an h-plane with a unit circle $\Omega$ as the absolute.
Let $O$ be the center of $\Omega$ and $P$ be another h-point.
Then there is a circle $\Gamma$ perpendicular to the absolute such that $O$ is the inverse of $P$ across~$\Gamma$.

Moreover, $\Gamma$ has radius $\tfrac{\sqrt{1-OP^2}}{OP}$,
and its center $P'$ is the inverse of $P$ across the absolute. 
\end{thm}

\begin{wrapfigure}[8]{o}{45mm}
\vskip-6mm
\centering
\includegraphics{mppics/pic-192}
\end{wrapfigure}

\parit{Proof.}
Let $P'$ be the inversion of $P$ across the absolute.
Choose a point $T$ on the absolute such that $\angle OPT$ is right.
Let $\Gamma$ be the circle thru $T$ with center $P'$.

By \ref{lem:inversion-sim}, $\triangle OPT\sim \triangle OTP'$.
It follows that $\angle OTP'$ is right and $\Gamma\perp\Omega$.

By AA, $\triangle TPP'\sim \triangle OTP'$
(the angle at $P'$ is shared, and $\angle TPP'$ and $\angle OTP'$ are right). 
Therefore $\tfrac{OP'}{P'T}=\tfrac{P'T}{PP'}$; so, $O$ is the inversion of $P$ across~$\Gamma$.

Note that $OP'\cdot OP=OT^2$.
By the Pythagorean theorem $P'T^2\z=P'P^2+PT^2$ and $PT^2=OT^2-PO^2$.
Since $P'P=OP'-OP$ and $OT=1$, we get $P'T=\tfrac{\sqrt{1-OP^2}}{OP}$,
and the second statement follows.
\qeds

Assume $\Gamma$ is a circline that is perpendicular to the absolute.
Consider the inversion $X\mapsto X'$ across $\Gamma$;
if $\Gamma$ is a line, set $X\mapsto X'$ to be the reflection across~$\Gamma$.

The following observation says that the map $X\mapsto X'$ respects all the notions introduced in the previous section.
Together with the lemma above, it implies that in any problem that is formulated entirely \textit{in h-terms} we can assume that a given h-point lies in the center of the absolute.

\begin{thm}{Main observation}\label{thm:main-observ}
The map $X\mapsto X'$ described above is a bijection from the h-plane to itself. 
Moreover, for any h-points $P$, $Q$, $R$ such that $P\ne Q$ and $Q\ne R$, the following conditions hold:
\begin{enumerate}[(a)]
\item\label{h-line-to-hline} The h-line $(PQ)_h$, h-half-line $[PQ)_h$, and h-segment $[PQ]_h$ are transformed into $(P'Q')_h$, $[P'Q')_h$, and $[P'Q']_h$ respectively.
\item\label{h-reflect} $\delta(P',Q')=\delta(P,Q)$ and $P'Q'_h=PQ_h$.
\item\label{h-angle-mes} 
$\measuredangle_h P'Q'R'\equiv-\measuredangle_h PQR$.
\end{enumerate}

\end{thm}

It is instructive to compare this observation with Proposition~\ref{prop:reflection}.

\parit{Proof.}
According to Theorem~\ref{thm:perp-inverse}, the map sends the absolute to itself. 
Therefore the points on $\Gamma$ do not move.
It follows that points inside of the absolute remain inside after the mapping.
Whence the $X\mapsto X'$ is a bijection from the h-plane to itself.

Part~\textit{(\ref{h-line-to-hline})} follows from \ref{thm:inverse-cline} and \ref{thm:angle-inversion}.

Part~\textit{(\ref{h-reflect})} follows from Theorem~\ref{lem:inverse-4-angle}.

Part~\textit{(\ref{h-angle-mes})} follows from Theorem~\ref{thm:angle-inversion}.
\qeds

\begin{thm}{Lemma}\label{lem:O-h-dist}
Assume that the absolute is a unit circle centered at~$O$.
Given an h-point $P$, set $x=OP$ and $y=OP_h$.
Then
\begin{align*}
y&=\ln\frac{1+x}{1-x}
&
&\text{and}
&
x&=\frac{e^y-1}{e^y+1}.
\end{align*}
 
\end{thm}

Observe that according to the lemma, $OP_h\to \infty$ as $OP\to 1$.
That is, if $P$ the approaches absolute in the Euclidean sense, it escapes to infinity in the h-sense.

\begin{wrapfigure}[6]{o}{38mm}
\vskip-4mm
\centering
\includegraphics{mppics/pic-194}
\end{wrapfigure}

\parit{Proof.}
The h-line $(OP)_h$ lies on a diameter of the absolute.
If $A$ and $B$ are the ideal points as in the definition of the h-distance, then
\begin{align*}
OA&=OB=1,
\\ 
PA&=1+x,
\\
PB&=1-x.\end{align*}
In particular,
\begin{align*}
y&=\ln \frac{AP\cdot BO}{PB\cdot OA}=\ln\frac{1+x}{1-x}.
\end{align*}

Taking the exponential function of the left and the right-hand side and applying obvious algebra manipulations, we get that
$$x=\frac{e^y-1}{e^y+1}.$$
\qedsf

\begin{thm}{Lemma}\label{lem:h-tiangle=}
Assume that points $P$, $Q$, and $R$ appear on an h-line in the same order.
Then 
$$PQ_h+QR_h=PR_h.$$ 

\end{thm}

\parit{Proof.}
Note that
$$PQ_h+QR_h=PR_h$$
is equivalent to 
\[\delta(P,Q)\cdot\delta(Q,R)=\delta(P,R).\eqlbl{eq:deltaPQR}\]

Let $A$ and $B$ be the ideal points of~$(PQ)_h$. 
Without loss of generality, we can assume that the points $A$, $P$, $Q$, $R$, and $B$ appear in the same order on the circline containing $(PQ)_h$.
Then
\begin{align*}
\delta(P,Q)\cdot\delta(Q,R)
&=
\frac{AQ\cdot BP}{QB\cdot PA}\cdot\frac{AR\cdot BQ}{RB\cdot QA}=
\\
&=\frac{AR\cdot BP}{RB\cdot PA}=
\\
&=\delta(P,R).
\end{align*}
Hence \ref{eq:deltaPQR} follows.
\qeds

Let $P$ be an h-point and $\rho>0$.
The set of all h-points $Q$ such that $PQ_h=\rho$ is called an \index{h-circle}\emph{h-circle} with the center $P$ and the \index{h-radius}\emph{h-radius} $\rho$.

\begin{thm}{Lemma}\label{lem:h-circle=circle}
Any h-circle is a Euclidean circle that lies completely in the h-plane.

More precisely for any h-point $P$ and $\rho\ge 0$
there is a $\hat\rho\ge 0$ and a point $\hat P$ such that 
$$PQ_h= \rho
\quad 
\iff
\quad
\hat PQ= \hat\rho$$
for every h-point~$Q$.

Moreover, if $O$ is the center of the absolute, then 
\begin{enumerate}
\item $\hat O=O$ for any $\rho$ and
\item $\hat P\in (OP)$ for every $P\ne O$.
\end{enumerate}

\end{thm}

\begin{wrapfigure}{o}{33mm}
\vskip-4mm
\centering
\includegraphics{mppics/pic-196}
\end{wrapfigure}

\parit{Proof.}
According to Lemma~\ref{lem:O-h-dist}, 
$OQ_h\z= \rho$ if and only if $$OQ= \hat\rho=\frac{e^\rho-1}{e^\rho+1}.$$
Therefore, the locus of h-points $Q$ such that $OQ_h= \rho$ is a Euclidean circle, 
denote it by $\Delta_\rho$.

If $P\ne O$, then by Lemma~\ref{lem:P-->O} and the main observation (\ref{thm:main-observ})
there is an inversion that respects all h-notions and sends $O\mapsto P$.

Let $\Delta_\rho'$ be the inverse of $\Delta_\rho$.
Since the inversion preserves the h-distance,
$PQ_h=\rho$ if and only if $Q\z\in\Delta_\rho'$.

According to Theorem~\ref{thm:inverse-cline}, $\Delta_\rho'$ is a Euclidean circle.
Let $\hat P$ and $\hat\rho$ denote the Euclidean center and radius of $\Delta_\rho'$.

Finally, notice that $\Delta_\rho'$ reflects to itself across $(OP)$;
that is, the center $\hat P$ lies on~$(OP)$.
\qeds

\begin{thm}{Exercise}\label{ex:h-circle=circle}
Describe a nondegenerate h-triangle $\triangle_hPQR$ that does not have an h-circumcircle;
that is, its vertices $P$, $Q$, and $R$ do not lie on an h-circle or h-line.
\end{thm}

\section[Axioms]{Axioms}
\subsection*{Axiom~\ref{def:birkhoff-axioms:0}}

Evidently, the h-plane contains at least two points.
Therefore, to show that Axiom~\ref{def:birkhoff-axioms:0} holds in the h-plane, we need to show that the h-distance defined in Section~\ref{sec:conformal-model} is a metric;
that is, the conditions \textit{(\ref{def:metric-space:a})}--\textit{(\ref{def:metric-space:d})} 
in Definition~\ref{def:metric-space} hold for the h-distance.

The following claim says that the h-distance satisfies conditions \textit{(\ref{def:metric-space:a})} 
and \textit{(\ref{def:metric-space:b})}.

\begin{thm}{Claim}
For any pair of h-points $P$ and $Q$, we have
$PQ_h\ge 0$
and $PQ_h=0$ if and only if $P=Q$.
\end{thm}

\parit{Proof.}
According to Lemma~\ref{lem:P-->O}
and the main observation (\ref{thm:main-observ}), 
we can assume that $Q$ is the center of the absolute.
In this case
$$
\delta(Q,P)=\frac{1+QP}{1-QP}\ge 1$$
and therefore
$$QP_h=\ln[\delta(Q,P)]\ge 0.$$
Moreover, the equalities hold if and only if $P=Q$.
\qeds

The following claim says that the h-distance meets~\ref{def:metric-space}\textit{\ref{def:metric-space:c}}.

\begin{thm}{Claim}
For any h-points $P$ and $Q$, we have
$PQ_h=QP_h$.
\end{thm}

\parit{Proof.}
Let $A$ and $B$ be ideal points of $(PQ)_h$ and
$A,P,Q,B$ appear on the circline containing $(PQ)_h$ in the same order.

{

\begin{wrapfigure}{o}{33mm}
\vskip-5mm
\centering
\includegraphics{mppics/pic-198}
\end{wrapfigure}

Then
\begin{align*}
PQ_h
&=\ln\frac{AQ\cdot BP}{QB\cdot PA}
=
\\
&=\ln\frac{BP\cdot AQ}{PA\cdot QB}=
\\
&=QP_h.
\end{align*}
\qedsf

}

The following claim shows, in particular, that
the triangle inequality 
(which is condition \ref{def:metric-space}\textit{\ref{def:metric-space:d}})
holds for $h$-distance.

\begin{thm}{Claim}\label{clm:h-dist+trig-inq}
Given a triple of h-points $P$, $Q$, and $R$,
we have
\[PQ_h+QR_h\ge PR_h.\]
Moreover, the equality holds if and only if $P$, $Q$, and $R$ lie on one h-line in the same order.
\end{thm}

\parit{Proof.}
Without loss of generality, we may assume that $P$ is the center of the absolute
and 
$0<QR_h\le PQ_h$.

Let $\Delta$ be the h-circle with the center $Q$ and h-radius $QR_h$.
Choose points $S$ and $T$ on the intersection of $(PQ)$ with~$\Delta$ so that $P$, $S$, $Q$, and $T$ appear on the h-line in the same order.
The latter is possible by Lemma~\ref{lem:h-tiangle=}, since $QS_h\z=QT_h\z=QR_h\z\le PQ_h$.

{

\begin{wrapfigure}{o}{35mm}
\vskip-0mm
\centering
\includegraphics{mppics/pic-200}
\end{wrapfigure}

According to Lemma~\ref{lem:h-circle=circle}, $\Delta$ is a Euclidean circle;
let $\hat Q$ be its Euclidean center.

Note that $\hat QS\z=\hat QT\z=\hat QR$.
By the Euclidean triangle inequality,
$$PT
=
P\hat Q+\hat Q R
\ge 
PR,
\eqlbl{RT>RQ}$$
and the equality holds if and only if $T=R$. 

By Lemma~\ref{lem:O-h-dist},
\begin{align*}
PT_h&=\ln\frac{1+PT}{1-PT},\\
PR_h&=\ln\frac{1+PR}{1-PR}.
\end{align*}
Note that the function $x\mapsto\ln\frac{1+x}{1-x}$ is increasing for $0\le x<1$.
Therefore, \ref{RT>RQ} implies
$$PT_h\ge PR_h;$$
moreover, the equality holds if and only if $T\z=R$.

}

Finally, applying Lemma~\ref{lem:h-tiangle=} again, 
we get that
$$PT_h=PQ_h+QR_h.$$
Hence the claim follows.
\qeds

\subsection*{Axiom~\ref{def:birkhoff-axioms:1}}

Once the following claim is proved,
Axiom~\ref{def:birkhoff-axioms:1} 
follows from Corollary~\ref{cor:h-line}.

\begin{thm}{Claim}
A subset of the h-plane is an h-line if and only if it forms a line for the h-distance in the sense of Definition~\ref{def:line}.
\end{thm}

\parit{Proof.}
Let $\ell$ be an h-line.
Applying the main observation (\ref{thm:main-observ}) we can assume that $\ell$ contains the center of the absolute.
In this case, $\ell$ is an intersection of a diameter of the absolute and the h-plane.
Let $A$ and $B$ be the endpoints of the diameter.

\begin{wrapfigure}{o}{29mm}
\vskip-3mm
\centering
\includegraphics{mppics/pic-201}
\vskip-2mm
\end{wrapfigure}

Consider the map $\iota\:\ell\to \mathbb{R}$ defined as
$$\iota(X)=\ln \frac{AX}{XB}.$$
Note that $\iota\:\ell\to \mathbb{R}$ is a bijection.

Furthermore, if $X,Y\in \ell$ and the points $A$, $X$, $Y$, and $B$ appear on $[AB]$ in the same order, then
\[\iota(Y)-\iota(X)=\ln \frac{AY}{YB}-\ln \frac{AX}{XB}=\ln \frac{AY\cdot BX}{YB\cdot XB}=XY_h.\]

We have shown that every h-line is a line for h-distance.
The converse follows from Claim~\ref{clm:h-dist+trig-inq}.
\qeds

\subsection*{Axiom~\ref{def:birkhoff-axioms:2}}

The first part of Axiom~\ref{def:birkhoff-axioms:2} follows directly from the definition of the h-angle measure (defined in Section~\ref{sec:conformal-model}).
It remains to show that $\measuredangle_h$ satisfies the conditions \ref{def:birkhoff-axioms:2a}, \ref{def:birkhoff-axioms:2b}, and \ref{def:birkhoff-axioms:2c} (see Section~\ref{sec:axioms}).

The following two claims say that
$\measuredangle_h$ satisfies
 \ref{def:birkhoff-axioms:2a} and \ref{def:birkhoff-axioms:2b}.

\begin{thm}{Claim}\label{clm:h2a}
Given an h-half-line $[O P)_h$ and $\alpha\in(-\pi,\pi]$, there is a unique h-half-line $[O Q)_h$ such that $\measuredangle_h P O Q= \alpha$.
\end{thm}

\begin{thm}{Claim}\label{clm:h2b}
For any h-points $P$, $Q$, and $R$ distinct from an h-point $O$, we have
$$\measuredangle_h P O Q+\measuredangle_h Q O R
\equiv\measuredangle_h P O R.$$

\end{thm}

\parit{Proof of \ref{clm:h2a} and \ref{clm:h2b}.}
Applying the main observation, 
we may assume that $O$ is the center of the absolute.
In this case, for any h-point $P\z\ne O$, the h-half-line
$[OP)_h$ is the intersection of the Euclidean half-line $[OP)$ with the h-plane.
Hence \ref{clm:h2a} and \ref{clm:h2b} 
follow from the axioms \ref{def:birkhoff-axioms:2a} and \ref{def:birkhoff-axioms:2b} of the Euclidean plane.
\qeds

The following claim says that
$\measuredangle_h$ satisfies
 \ref{def:birkhoff-axioms:2c}.

\begin{thm}{Claim}\label{clm:h2c}
The function 
$$\measuredangle_h\:(P,Q,R)\mapsto\measuredangle_h P Q R$$
is continuous at every triple of points $(P,Q,R)$
such that $Q\ne P$, $Q\ne R$, and $\measuredangle_h P Q R\ne\pi$.
\end{thm}

\parit{Proof.}
Suppose that $O$ denotes the center of the absolute.
We can assume that $Q$ is distinct from~$O$;
the latter follows from the main observation.

\begin{wrapfigure}{o}{45mm}
\vskip-4mm
\centering
\includegraphics{mppics/pic-199}
\end{wrapfigure}

Let $Z$ be the inverse of $Q$ across the absolute;
denote by $\Gamma$ the circle with the center at~$Z$ that is perpendicular to the absolute.
According to Lemma~\ref{lem:P-->O},
point $O$ is the inverse of $Q$ across~$\Gamma$.

Let $P'$ and $R'$ be the inversions across $\Gamma$ of the points $P$ and $R$ respectively.
The point $P'$ is completely determined by $Q$ and $P$.
Moreover, the map $(Q,P)\mapsto P'$ is continuous at any pair of h-points $(Q,P)$ such that $Q\ne O$.
The same is true for the map $(Q,R)\mapsto R'$.

According to the main observation 
$$\measuredangle_h P Q R\equiv -\measuredangle_h P' O R'.$$
Since $\measuredangle_h P' O R'=\measuredangle P' O R'$ and 
the maps $(Q,P)\mapsto P'$, $(Q,R)\mapsto R'$ are continuous,
the claim follows from the corresponding axiom of the Euclidean plane.
\qeds

\subsection*{Axiom~\ref{def:birkhoff-axioms:3}}

The following claim says that Axiom~\ref{def:birkhoff-axioms:3} holds in the h-plane.

\begin{thm}{Claim}
In the h-plane, we have
$\triangle_h P Q R 
\cong
\triangle_h P' Q' R'$
if and only if 
\begin{align*}
Q' P'_h&=Q P_h, & Q' R'_h&= Q R_h &&\text{and}
&\measuredangle_h P' Q' R'&=\pm\measuredangle P Q R.
\end{align*}
 
\end{thm}

\parit{Proof.}
Applying the main observation, 
we can assume that $Q$ and $Q'$ coincide with the center of the absolute; in particular, $Q=Q'$.
In this case, 
$$\measuredangle P' Q R'=\measuredangle_h P' Q R'=\pm\measuredangle_h P Q R=\pm\measuredangle P Q R.$$
Since 
$$Q P_h=Q P'_h\quad \text{and}\quad Q R_h=Q R'_h,$$
Lemma~\ref{lem:O-h-dist} implies that the same holds for the Euclidean distances;
that is,
$$Q P=Q P'
\quad
\text{and}
\quad
Q R=Q R'.$$
By SAS,
there is a motion of the Euclidean plane that sends
$Q$ to itself,
$P$ to $P'$, 
and $R$ to $R'$.

Note that the center of the absolute is fixed by the corresponding motion.
It follows that this motion gives also a motion of the h-plane;
in particular, the h-triangles 
$\triangle_h P Q R$ and $\triangle_h P' Q R'$ are h-congruent.
\qeds

\subsection*{Axiom h-$\!$\ref{def:birkhoff-axioms:4}}

Finally, we need to check that Axiom~h-$\!$\ref{def:birkhoff-axioms:4} in Section~\ref{sec:unprovable} holds;
that is, we need to prove the following claim.

{

\begin{wrapfigure}{r}{31mm}
\vskip-4mm
\centering
\includegraphics{mppics/pic-202}
\end{wrapfigure}

\begin{thm}{Claim}
For every h-line $\ell$ and every h-point $P\notin\ell$ there are at least two h-lines that pass thru $P$
and have no points of intersection with~$\ell$.
\end{thm}

\parit{Instead of proof.}
Applying the main observation we can assume that $P$ is the center of the absolute.

The remaining part of the proof can be guessed from the picture.
\qeds

}

\begin{thm}{Exercise}\label{ex:3-h-lines}
\begin{enumerate}[(a)]

\item Draw three three mutually parallel h-lines
such that every pair of these three lines lies on one side of the remaining h-line.

\item Draw h-lines $\ell$, $m$, and $n$ such that $\ell\parallel m$, $m\parallel n$, but $\ell\nparallel n$.
Conclude that the parallelness is not an equivalence relation for h-lines.

\item Draw $h$-lines $k$, $\ell$, $m$, and $n$, such that $k\perp \ell$, $\ell\perp m$, $m\perp n$, and $k\parallel n$.
\end{enumerate}
\end{thm}
 
\section{Hyperbolic trigonometry}
\label{sec:hyp-trig}

This section provides formulas for h-distance using \index{hyperbolic!functions}\emph{hyperbolic functions}.
One of these formulas will be used in the proof of the hyperbolic Pythagorean theorem (\ref{thm:pyth-h-poincare}).

Recall that $\cosh$, $\sinh$, and $\tanh$ denote \index{ch@$\cosh$}\index{hyperbolic!cosine}\emph{hyperbolic cosine}, \index{sh@$\sinh$}\index{hyperbolic!sine}\emph{hyperbolic sine}, and \index{th@$\tanh$}\index{hyperbolic!tangent}\emph{hyperbolic tangent}\label{hyperbolic tangent};
that is, the functions defined by
\[\cosh x\df \tfrac{e^x+e^{-x}}2,
 \quad
 \sinh x\df \tfrac{e^x-e^{-x}}2,
\]
\[\tanh x\df \tfrac{\sinh x}{\cosh x}.
\]

These hyperbolic functions are analogous to sine, cosine, and tangent. 

\begin{thm}{Exercise}\label{ex:hyp-fun}
Prove the following identities:
\[\cosh' x=\sinh x;\quad \sinh'x=\cosh x;\quad (\cosh x)^2-(\sinh x)^2=1.\]
\end{thm}

\begin{thm}{Exercise}\label{ex:O-h-dist}
Assume that the absolute is a unit circle centered at~$O$.
Show that 
\[OP=\tanh(\tfrac12 \cdot OP_h).\]
for every h-point $P$.
\end{thm}

\begin{thm}{Double-argument identities}\label{double-argument}
The identities
\begin{align*}
\cosh (2\cdot x)&=(\cosh x)^2+(\sinh x)^2 
&&\text{and}&
\sinh (2\cdot x)&=2\cdot\sinh x\cdot \cosh x
\end{align*}
hold for any real value $x$.
\end{thm}

\parit{Proof.}
\begin{align*}
(\sinh x)^2+(\cosh x)^2
&=(\tfrac{e^x-e^{-x}}2)^2+(\tfrac{e^x+e^{-x}}2)^2=
\\
&=\tfrac{e^{2\cdot x}+e^{-2\cdot x}}2=
\\
&=\cosh (2\cdot x);
\\
2\cdot\sinh x\cdot \cosh x
&=2\cdot(\tfrac{e^x-e^{-x}}2)\cdot(\tfrac{e^x+e^{-x}}2)=
\\
&=\tfrac{e^{2\cdot x}-e^{-2\cdot x}}2=
\\
&=\sinh (2\cdot x).
\end{align*}
\qedsf

\begin{thm}{Advanced exercise}\label{ex:cosh}
Let $P$ and $Q$ be two h-points distinct from the center of the absolute.
Denote by $P'$ and $Q'$ the inverses of $P$ and $Q$ across the absolute.

\begin{wrapfigure}[20]{r}{40mm}
\centering
\includegraphics{mppics/pic-204}
\end{wrapfigure}

Show that 
\medskip
\begin{enumerate}[(a)]
\item\label{ex:cosh/2} 
$\displaystyle{\cosh[\tfrac12\cdot PQ_h]=\sqrt{\frac{PQ'\cdot P'Q}{PP'\cdot QQ'}};}$
\medskip
\item\label{ex:coshsinh} 
$\displaystyle{\sinh[\tfrac12\cdot PQ_h]=\sqrt{\frac{PQ\cdot P'Q'}{PP'\cdot QQ'}};}$
\medskip
\item\label{ex:coshtanh} 
$\displaystyle{\tanh[\tfrac12\cdot PQ_h]=\sqrt{\frac{PQ\cdot P'Q'}{PQ'\cdot P'Q}};}$
\medskip
\item\label{ex:coshcosh} 
$\displaystyle{\cosh PQ_h=\frac{PQ\cdot P'Q'+PQ'\cdot P'Q}{PP'\cdot QQ'}.}$
\end{enumerate}

\end{thm}

\chapter{Geometry of the h-plane}\label{chap:h-plane}

In this chapter, we study the geometry of the plane described by the conformal disc model, briefly the h-plane.

We can work with this model inside the Euclidean plane. 
Additionally, we have the option to apply the axioms of neutral geometry since they all hold in the h-plane; the latter is proved in the previous chapter.

\section{The angle of parallelism}

Let $P$ be a point off an h-line~$\ell$. 
Drop a perpendicular $(PQ)_h$ from $P$ to $\ell$;
let $Q$ be its footpoint.
Let $\phi$ be the smallest value such that the h-line $(PZ)_h$ with $|\measuredangle_h Q P Z|=\phi$ does not intersect~$\ell$.

The value $\phi$ is called the \index{angle!angle of parallelism}\emph{angle of parallelism} at $P$ to~$\ell$.
Clearly, $\phi$ depends only on the h-distance $s=PQ_h$.
Furthermore, $\phi(s)\to \pi/2$ as $s\to 0$, 
and $\phi(s)\to0$ as $s\to\infty$.
(In Euclidean geometry, the angle of parallelism is identically equal to~$\pi/2$.)

\begin{thm}{Exercise}\label{ex:lambert-parallelism}
Suppose that $\square_hABCD$ has right h-angles at $A$, $B$, and~$C$.
Show that $|\measuredangle_h CDA|>\phi$, where $\phi$ is the angle of parallelism at $D$ with respect to~$(AB)_h$.
\end{thm}

\begin{wrapfigure}{o}{34mm}
\vskip-3mm
\centering
\includegraphics{mppics/pic-206}
\end{wrapfigure}

If $\ell$, $P$, and $Z$ are as above, then the h-line $m=(PZ)_h$ is called \index{asymptotically parallel lines}\emph{asymptotically parallel} to~$\ell$.
In other words, two h-lines are asymptotically parallel if they share one ideal point.
(In hyperbolic geometry, the term parallel lines is often used for asymptotically parallel lines; we do not follow this convention.)

Given $P\not\in\ell$, there are exactly two asymptotically parallel lines thru $P$ to $\ell$; 
the remaining parallel lines are called \index{parallel lines!ultra parallel lines}\emph{ultra parallel}.

In the picture, the two solid h-lines passing thru $P$ are asymptotically parallel to~$\ell$;
the dashed h-line is ultra parallel to~$\ell$.

\begin{thm}{Exercise}\label{ex:ultra-parallel}
Show that two distinct h-lines $\ell$ and $m$ are ultraparallel if and only if they have a {}\emph{common perpendicular};
that is, there is an $h$-line $n$ such that $n\perp \ell$ and $n\perp m$.
\end{thm}

\begin{thm}{Proposition}\label{prop:angle-parallelism}
Let $Q$ be a footpoint of $P$ on an h-line~$\ell$.
Then
\[h=\tfrac12\cdot\ln \frac{1+\cos\phi}{1-\cos\phi}
\qquad\text{and}\qquad
\cos\phi=\frac{e^{2\cdot h}-1}{e^{2\cdot h}+1},\]
where $h=PQ_h$, and $\phi$ is the angle of parallelism at $P$ with respect to~$\ell$.

In particular, if $P\notin\ell$ and $\beta\z=|\measuredangle_h XPY|$ for some points $X,Y\in\ell$, then 
\[
h<\tfrac12\cdot\ln \frac{1+\cos\tfrac\beta2}{1-\cos\tfrac\beta2}
\qquad\text{and}\qquad
\cos\tfrac\beta2>\frac{e^{2\cdot h}-1}{e^{2\cdot h}+1}.\]

\end{thm}

\begin{wrapfigure}{o}{50mm}
\vskip-6mm
\centering
\includegraphics{mppics/pic-208}
\end{wrapfigure}

\parit{Proof.} Applying a motion of the h-plane if necessary,
we may assume $P$ is the center of the absolute.
Then the h-lines thru $P$ are the intersections of Euclidean lines with the h-plane.

Let $A$ and $B$ denote the ideal points of~$\ell$.
Without loss of generality, we may assume that $\angle APB$ 
is positive.
In this case, 
$$\phi=\measuredangle QPB=\measuredangle APQ=\tfrac12 \cdot\measuredangle APB.$$

Let $Z$ be the center of the circle $\Gamma$ containing the h-line~$\ell$.
Set $X$ to be the point of intersection of the Euclidean segment $[AB]$ and the line~$(PQ)$.

Note that, $PX=\cos\phi$.
Therefore, by Lemma~\ref{lem:O-h-dist},
$$PX_h=\ln \tfrac{1+\cos\phi}{1-\cos\phi}.$$

Note that both angles $PBZ$ and $BXZ$ are right.
Since the angle $PZB$ is shared, $\triangle ZBX\sim \triangle ZPB$.
In particular, 
$$ZX\cdot ZP=ZB^2;$$
that is, $X$ is the inverse of $P$ across~$\Gamma$.

{

\begin{wrapfigure}{o}{36mm}
\vskip-0mm
\centering
\includegraphics{mppics/pic-209}
\end{wrapfigure}

The inversion across $\Gamma$ is the reflection of the h-plane across~$\ell$. 
Therefore
\begin{align*}
h&=PQ_h=QX_h=
\\
&=\tfrac12\cdot PX_h=
\\
&=\tfrac12\cdot\ln \tfrac{1+\cos\phi}{1-\cos\phi}.
\end{align*}
Taking exponent of the left and right hand side and simplifying, we get the second identity.

}

The last statement follows since $\phi\z>\tfrac\beta2$ and the function 
\[\phi\mapsto  \tfrac12\cdot\ln \tfrac{1+\cos\phi}{1-\cos\phi}\] 
is decreasing in the interval $(0,\tfrac\pi2]$.
\qeds

\begin{thm}{Exercise}\label{ex:right-angle-parallelism} 
Assume $\triangle_hABC$ has a right h-angle at $C$.
Show that the h-distance from $C$ to $[AB]_h$ cannot exceed 1.
\end{thm}

\begin{thm}{Exercise}\label{ex:small-angle}
Let $ABC$ be an equilateral h-triangle with side $100$.
Show that 
\[|\measuredangle_h ABC|<\frac1{10\,000\,000\,000}.\]
Sketch $\triangle_h ABC$ in the conformal model.
\end{thm}

\section{Inradius of h-triangle}

\begin{thm}{Theorem}\label{thm:h-inradius}
The inradius of every h-triangle
is less than $\tfrac12\cdot\ln3$.
\end{thm}

\parit{Proof.}
Let $I$ and $r$ be the h-incenter and h-inradius of $\triangle_hXYZ$.

Note that the h-angles 
$XIY$, 
$YIZ$, and 
$ZIX$
have the same sign.
Without loss of generality, we can assume that all of them are positive.
Therefore,
\[\measuredangle_hXIY+ 
\measuredangle_hYIZ+ 
\measuredangle_hZIX=2\cdot\pi,
\]
and we can assume that
$\measuredangle_hXIY\ge\tfrac23\cdot\pi$;
if not relabel $X$, $Y$, and~$Z$.

{

\begin{wrapfigure}{o}{44mm}
\centering
\vskip-4mm
\includegraphics{mppics/pic-210}
\end{wrapfigure} 

Since $r$ is the h-distance from $I$ to $(XY)_h$,
Proposition~\ref{prop:angle-parallelism} implies that
\begin{align*}r&<\tfrac12\cdot\ln \tfrac{1+\cos\frac\pi3}{1-\cos\frac\pi3}=
\\
&=\tfrac12\cdot\ln\frac{1+\tfrac12}{1-\tfrac12}=
\\
&=\tfrac12\cdot\ln 3.
\end{align*}
\qedsf

}

\begin{thm}{Exercise}\label{ex:side-sup}
Let $\square_h ABCD$ be a quadrangle in the h-plane 
such that the h-angles at $A$, $B$, and $C$ are right and $AB_h=BC_h$.
Find the optimal upper bound for~$AB_h$.
\end{thm}

\section{Circles, horocycles, and equidistants}

According to Lemma~\ref{lem:h-circle=circle},
every h-circle is a Euclidean circle that lies entirely within the h-plane.
Furthermore, every h-line is an intersection of the h-plane with a circle
perpendicular to the absolute.

In this section, we will describe the 
h-geometric meaning of the intersections 
of other circles with the h-plane.

You will see that all these intersections have a {}\emph{perfectly round shape} in the h-plane.

One may think of these curves as the trajectories of a car with a fixed position of the steering wheel.
In the Euclidean plane, 
this way you either travel along a circle or a line.

In the hyperbolic plane, the picture is different.
If you turn the steering wheel to the far right, you will travel along a circle.
If you turn it less, at a certain position of the wheel, you will never return to the same point, but the path will be different from a line.
If you turn the wheel a bit further, then you start to travel along a path that stays at a fixed distance from an h-line.

\parbf{Equidistants of h-lines.}
Consider the h-plane with the absolute~$\Omega$.
Assume a circle $\Gamma$ intersects $\Omega$ at two distinct points, $A$ and~$B$. 
Suppose that $g$ denotes the intersection of $\Gamma$ with the h-plane.

\begin{wrapfigure}{o}{46mm}
\vskip-0mm
\centering
\includegraphics{mppics/pic-212}
\end{wrapfigure}

Let us draw an h-line $m$ with the ideal points $A$ and~$B$.
According to Exercise~\ref{ex:ideal-line-unique}, $m$ is uniquely defined.

Consider an h-line $\ell$ perpendicular to~$m$;
let $\Delta$ be the circle containing~$\ell$.

Note that $\Delta\perp \Gamma$.
Indeed,
according to Corollary~\ref{cor:perp-inverse-clines}, $m$ and $\Omega$ invert to themselves in~$\Delta$.
It follows that $A$ is the inverse of $B$ across~$\Delta$.
Finally, by Corollary~\ref{cor:perp-inverse}, we get that $\Delta\perp \Gamma$.

Therefore, inversion across $\Delta$ sends both $m$ and $g$ to themselves.
For any two points $P',P\in g$ there is a choice of $\ell$ and $\Delta$ as above such that
$P'$ is the inverse of $P$ across $\Delta$.
By the main observation (\ref{thm:main-observ}) the inversion across $\Delta$ is a motion of the h-plane. Therefore, all points of $g$ lie at the same distance from~$m$.

In other words, $g$ is the set of points that lie at a fixed h-distance and on the same side of~$m$.

Such a curve $g$ is called 
\index{equidistant}\emph{equidistant} to h-line~$m$.%
\footnote{It can also be called \index{hypercycle}\emph{hypercycle} with central line $m$.}
In Euclidean geometry, the equidistant from a line is a line;
apparently, in hyperbolic geometry, the picture is different.

\begin{thm}{Exercise}\label{ex:equidistant-reflection}
Let $g$ be the set of all h-reflections of a given h-point $P$ across the points on a given h-line $m$.
Show that $g$ is equidistant to $m$.
\end{thm}

\parbf{Horocycles.}
If the circle $\Gamma$ touches the absolute from inside at one point $A$, then the complement $h=\Gamma\backslash\{A\}$ lies in the h-plane.
This set is called a \index{horocycle}\emph{horocycle}.
It also has a perfectly round shape in the sense described above.

\begin{wrapfigure}{o}{33mm}
\vskip-0mm
\centering
\includegraphics{mppics/pic-214}
\end{wrapfigure}

The shape of a horocycle is between the shapes of circles and equidistants to h-lines.
A horocycle might be considered as a limit of circles 
thru a fixed point, say $P$,
with the centers $O_n$ running to infinity along an h-line $\ell$.
The same horocycle is a limit of equidistants thru $P$ to the sequence of h-lines $m_n$ passing thru $O_n$ and perpendicular to $\ell$.

Since any three points lie on a circline, we have that any nondegenerate h-triangle is inscribed in an h-circle, horocycle, or equidistant.

\begin{thm}{Exercise}\label{ex:right-trig-horocycle}
Find the leg of an isosceles right h-triangle inscribed in a horocycle.
\end{thm}

\section{Hyperbolic triangles}

\begin{thm}{Theorem}\label{thm:3sum-h}
Any nondegenerate hyperbolic triangle has a positive defect.
\end{thm}

\begin{wrapfigure}{o}{35mm}
\centering
\includegraphics{mppics/pic-216}
\end{wrapfigure}

\parit{Proof.}
Choose an h-triangle $ABC$.
According to Theorem~\ref{thm:3sum-a},
$$\defect(\triangle_hABC)\ge 0.\eqlbl{eq:defect<0}$$
It remains to show that in the case of equality, $\triangle_hABC$ degenerates.

Without loss of generality, we may assume that $A$ is the center of the absolute;
in this case, 
$\measuredangle_h CAB\z=\measuredangle CAB$.
Yet we may assume that 
$$\measuredangle_h CAB,
\quad 
\measuredangle_h ABC,
\quad
\measuredangle_h BCA,
\quad
\measuredangle ABC,
\quad
\measuredangle BCA\ge 0.$$

Let $D$ be an arbitrary point in $[CB]_h$ distinct from $B$ and~$C$.
From Proposition~\ref{prop:arc(angle=tan)}, we have
$$\measuredangle ABC-\measuredangle_h ABC \equiv 
\pi-\measuredangle CDB
\equiv \measuredangle BCA-\measuredangle_h BCA.$$

From Exercise~\ref{ex:|3sum|}, we get that
$$\defect(\triangle_hABC)=2\cdot(\pi-\measuredangle CDB).$$
Therefore, if we have equality in \ref{eq:defect<0}, then $\measuredangle CDB=\pi$.
In particular, the h-segment $[BC]_h$ coincides with the Euclidean segment~$[BC]$.
By Exercise~\ref{ex:line/h-line},
the latter can happen only if the h-line $(BC)_h$ passes thru the center of the absolute ($A$);
that is, if $\triangle_hABC$ degenerates.
\qeds

The following theorem states, in particular, that nondegenerate hyperbolic triangles are congruent if their corresponding angles are equal.
In particular, in hyperbolic geometry, similar triangles have to be congruent.

{\sloppy 
\begin{thm}{AAA congruence condition}\label{thm:AAA}\index{AAA congruence condition}
Two nondegenerate h-triangles
 $ABC$ and $A'B'C'$
 are congruent if
$\measuredangle_hABC\z=\pm\measuredangle_hA'B'C'$,
$\measuredangle_hBCA\z=\pm\measuredangle_hB'C'A'$
and 
$\measuredangle_hCAB=\pm\measuredangle_hC'A'B'$.
\end{thm}

}

\parit{Proof.}
If $AB_h=A'B'_h$, then the theorem follows from ASA.
Assume $AB_h\ne A'B'_h$. 
Without loss of generality, we may assume that $AB_h\z<A'B'_h$.

\begin{wrapfigure}{r}{32mm}
\centering
\vskip-5mm
\includegraphics{mppics/pic-218}
\end{wrapfigure}

Let us choose $B''\in [A'B')_h$ and $C''\in [A'C')_h$ such that 
$A'B''_h=AB_h$ and $A'C''_h=AC_h$.
By SAS, $\triangle_hA'B''C''\cong\triangle_hABC$;
it follows that $\measuredangle_h A'B''C''=\pm\measuredangle_h A'B'C'$.
Since angles of a triangle have the same signs (\ref{thm:signs-of-triug}), we have
\[\measuredangle_h A'B''C''=\measuredangle_h A'B'C'.\eqlbl{A'B''C''congA'B'C'}\]

Let us show that points $B''$ and $C''$ lie on the sides $[A'B']_h$ and $[A'C']_h$ respectively.
Indeed, according to Exercise~\ref{ex:parallel-abs}, \ref{A'B''C''congA'B'C'} implies that $(B''C'')_h\z\parallel(B'C')_h$.
In particular, $B''$ and $C''$ lie on one side of $(B'C')_h$.
Since $AB_h<A'B'_h$ and $B''\in [A'B')_h$, we have $B''$ and $A'$ lie on one side of $(B'C')_h$.
It follows that $C''$ and $A'$ lie on one side of $(B'C')_h$;
therefore, $C''\in [A'C']_h$.

Since $\triangle_h ABC\cong\triangle_h A'B''C''$, we have
$$\defect(\triangle_h ABC)=\defect(\triangle_h A'B''C'').
\eqlbl{eq:defect=defect}$$

Applying Exercise~\ref{ex:defect} twice, we get that
$$\begin{aligned}
\defect(\triangle_h A'B'C')
&=
\defect(\triangle_h A'B''C'')
+
\\
&+\defect(\triangle_h B''C''C')+\defect(\triangle_h B''C'B').
\end{aligned}
\eqlbl{eq:defect+defect}$$
By Theorem~\ref{thm:3sum-h}, all the defects have to be positive.
Therefore
$$\defect(\triangle_h A'B'C')
>\defect(\triangle_h ABC).$$
On the other hand, by assumption we have
$$\begin{aligned}
\defect(\triangle_h A'B'C')
&= \pi -|\measuredangle_hA'B'C'|-|\measuredangle_hB'C'A'|-|\measuredangle_hC'A'B'|=
\\
&=\pi -|\measuredangle_hABC|-|\measuredangle_hBCA|-|\measuredangle_hCAB|=
\\
&=\defect(\triangle_h ABC)
 \end{aligned}$$
--- a contradiction.
\qeds

Recall that \index{angle-preserving transformation}\emph{angle-preserving transformation} is a bijection from an h-plane to itself such that 
\[\measuredangle_h ABC= \measuredangle_h A'B'C'\]
for every $\triangle_h ABC$ and its image $\triangle_h A'B'C'$.

\begin{thm}{Exercise}\label{ex:angle-preserving-hyp}
Show that every angle-preserving transformation of the h-plane is a motion.
\end{thm}

\section{The conformal interpretation}

Let us give another interpretation of the h-distance.

\begin{thm}{Lemma}\label{lem:conformal}
Consider the h-plane with the unit circle centered at~$O$ as the absolute.
Fix a point $P$ and let $Q$ be another point in the h-plane.
Then
$$\frac{PQ_h}{PQ}\to \frac{2}{1-OP^2}.$$
as $Q\to P$.
\end{thm}

The above formula tells us that the h-distance from $P$ to a nearby point $Q$ is almost proportional to the Euclidean distance
with the coefficient $\tfrac{2}{1-OP^2}$. 
The value $\lambda(P)=\tfrac{2}{1-OP^2}$ is called the \index{conformal factor}\emph{conformal factor} of the h-metric.

The value $\tfrac1{\lambda(P)}=\tfrac12\cdot(1-OP^2)$
can be interpreted as the {}\emph{speed limit} at the given point~$P$. 
In this case, the h-distance is the minimal time needed to travel from one point of the h-plane to another point.

\parit{Proof.}
Set $x=PQ$ and $y=PQ_h$.

If $P=O$, then by Lemma~\ref{lem:O-h-dist} we have
$$\frac{y}{x}=\frac{\ln \tfrac{1+x}{1-x}}{x}\to 2\eqlbl{eq:O=P}$$
as $x\to0$.%
\footnote{In other words, the function $x\mapsto \ln \tfrac{1+x}{1-x}$ has derivative 2 at $x=0$.} 

Suppose $P\ne O$; let $Z$ denotes the inverse of $P$ across the absolute.
Let $\Gamma$ be the circle with the center $Z$ 
perpendicular to the absolute.

According to the main observation (\ref{thm:main-observ}) and Lemma~\ref{lem:P-->O}, 
the inversion across $\Gamma$ is a motion of the h-plane, and it sends $P$ to~$O$.
In particular, $OQ'_h=PQ_h$, where $Q'$ denotes the inverse of $Q$ across~$\Gamma$.

\begin{wrapfigure}{o}{40mm}
\centering
\includegraphics{mppics/pic-220}
\end{wrapfigure}

Set $x'=OQ'$.
According to Lemma~\ref{lem:inversion-sim},
$$\frac{x'}{x}=\frac{OZ}{ZQ}.$$
Since $Z$ is the inverse of $P$ across the absolute, we have that $PO\cdot OZ=1$.
Therefore, 
$$\frac{x'}{x}\to \frac{OZ}{ZP}=\frac{1}{1-OP^2}$$
as $x\to 0$.

According to \ref{eq:O=P}, $\frac{y}{x'}\to 2$ as $x'\to 0$.
Therefore
$$\frac{y}{x}=\frac{y}{x'}\cdot \frac{x'}{x}\to \frac{2}{1-OP^2}$$
as $x\to 0$.\qeds

Here is an application of the lemma above.

\begin{thm}{Proposition}\label{prop:circum}
The circumference of an h-circle of the h-radius $r$ is 
$$2\cdot\pi\cdot\sinh r,$$
where \index{sh@$\sinh$}$\sinh r$ denotes the \index{hyperbolic!sine}\emph{hyperbolic sine} of $r$;
that is,
$$\sinh r\df \frac{e^r-e^{-r}}{2}.$$

\end{thm}

Before we proceed with the proof, let us discuss the same problem in the Euclidean plane.

The circumference of a circle in the Euclidean plane
can be defined as the limit of perimeters of regular $n$-gons inscribed in the circle as $n\to \infty$.

Namely, let us fix $r>0$.
Given a positive integer $n$, consider $\triangle AOB$
such that
$\measuredangle AOB=\tfrac{2\cdot\pi}{n}$ and $OA=OB=r$.
Set $x_n=AB$.
Note that $x_n$ is the side of a regular $n$-gon inscribed in the circle with radius $r$. 
Therefore, the perimeter of the $n$-gon is $n\cdot x_n$.

\begin{wrapfigure}[11]{o}{47mm}
\centering
\vskip-0mm
\includegraphics{mppics/pic-222}
\end{wrapfigure}

The circumference of the circle with radius $r$ 
can then be defined as the limit
$$\lim_{n\to\infty} n\cdot x_n=2\cdot\pi\cdot r.\eqlbl{eq:2pir}$$
(This limit can be taken as the definition of~$\pi$.)

In the following proof, we will replicate the same construction in the h-plane.

\parit{Proof.}
Without loss of generality, we can assume that the center $O$ of the circle is also the center of the absolute.

By Lemma~\ref{lem:O-h-dist}, 
the h-circle with the h-radius $r$ is the Euclidean circle with center $O$ and radius 
$$a=\frac{e^r-1}{e^r+1}.$$

Let $x_n$ and $y_n$ denote the side lengths of the regular $n$-gons inscribed in the circle in the Euclidean and hyperbolic plane respectively.

Note that $x_n\to0$ as $n\to\infty$.
By Lemma~\ref{lem:conformal},
\begin{align*}
\lim_{n\to\infty}\frac{y_n}{x_n}
&=\frac{2}{1-a^2}.
\end{align*}

Applying \ref{eq:2pir},
we get that the circumference of the h-circle can be found the following way:
\begin{align*}
\lim_{n\to\infty}n\cdot y_n
&=\frac{2}{1-a^2}\cdot\lim_{n\to\infty}n\cdot x_n=
\\
&=\frac{4\cdot\pi\cdot a}{1-a^2}=
\\
&=\frac{4\cdot\pi\cdot\left(\frac{e^r-1}{e^r+1}\right)}{1-\left(\frac{e^r-1}{e^r+1}\right)^2}=
\\
&=2\cdot\pi\cdot\frac{e^{r}-e^{-r}}{2}=
\\
&=2\cdot\pi\cdot\sinh r.
\end{align*}
\qedsf

\begin{thm}{Exercise}\label{ex:circum}
Let $\circum_h(r)$ denote the circumference of the h-circle with the h-radius~$r$.
Show that 
$$\circum_h(r+1)>2\cdot \circum_h(r)$$
for all $r>0$.
\end{thm}

\section{The Pythagorean theorem}\index{Pythagorean theorem}

Recall that $\cosh$ denotes \index{ch@$\cosh$}\index{hyperbolic!cosine}\emph{hyperbolic cosine};
that is, the function defined by
$$\cosh x\df \tfrac{e^x+e^{-x}}2.$$

\begin{thm}{Hyperbolic Pythagorean theorem}\label{thm:pyth-h-poincare}
Assume that h-triangle $ACB$ has right angle at~$C$.
Then
\[\cosh c=\cosh a\cdot\cosh b,
\eqlbl{eq:thm:pyth-h-poincare}\]
where $a=BC_h$, $b=CA_h$, and $c=AB_h$.
\end{thm}

The formula \ref{eq:thm:pyth-h-poincare} will be proved by means of direct calculations.
Before giving the proof, let us discuss the limit cases of this formula.

Note that $\cosh x$ can be written using the Taylor expansion
\[\cosh x=1+\tfrac1{2}\cdot x^2+\tfrac1{24}\cdot x^4+\dots\]

It follows that if $a$, $b$, and $c$ are small, then
\begin{align*}
1+\tfrac1{2}\cdot c^2&\approx \cosh c=\cosh a\cdot\cosh b\approx
\\
&\approx(1+\tfrac1{2}\cdot a^2)\cdot (1+\tfrac1{2}\cdot b^2)
\approx 
\\
&\approx
1+\tfrac1{2}\cdot (a^2+b^2).
\end{align*}
In other words, the original Pythagorean theorem (\ref{thm:pyth}) is a limit case of the hyperbolic Pythagorean theorem for small triangles.

For large $a$ and $b$ the terms $e^{-a}$, $e^{-b}$, and $e^{-a-b+\ln 2}$ are neglectable.
In this case, we have the following approximations:
\begin{align*}
\cosh a\cdot\cosh b&\approx \tfrac{e^a}2\cdot\tfrac{e^b}2=
\\
&=\frac{e^{a+b-\ln 2}}{2}\approx
\\
&\approx \cosh(a+b-\ln 2).
\end{align*}
Therefore $c\approx a+b-\ln 2$. 

\begin{thm}{Exercise}\label{ex:c+1>a+b}
Assume $\triangle_h ACB$ has a right angle at~$C$.
Set $a=BC_h$, $b=CA_h$, and $c=AB_h$.
Show that
\[c+\ln 2>a+b.\]

\end{thm}

In the proof of the hyperbolic Pythagorean theorem, we use the following formula from Exercise~\ref{ex:cosh}\textit{\ref{ex:coshcosh}}:
\[\cosh AB_h=\frac{AB\cdot A'B'+AB'\cdot A'B}{AA'\cdot BB'},\]
here $A$, $B$ are h-points distinct from the center of the absolute and $A'$, $B'$ are their inversions across the absolute.
This formula is derived in the hints.

\begin{wrapfigure}{o}{43mm}
\centering
\vskip-6mm
\includegraphics{mppics/pic-224}
\end{wrapfigure}

\parit{Proof of \ref{thm:pyth-h-poincare}.}
We assume that the absolute is a unit circle.
By the main observation (\ref{thm:main-observ}) we can assume that $C$ is the center of the absolute.
Let $A'$ and $B'$ denote the inverses of $A$ and $B$ across the absolute.

Set $x=BC$, $y=AC$.
By Lemma~\ref{lem:O-h-dist}
\begin{align*}
a&=\ln \tfrac{1+x}{1-x},
&
b&=\ln \tfrac{1+y}{1-y}.
\end{align*}
Therefore
\[\begin{aligned}
\cosh a&=\tfrac12\cdot (\tfrac{1+x}{1-x}+\tfrac{1-x}{1+x})=
&&&&&
\cosh b&=\tfrac12\cdot (\tfrac{1+y}{1-y}+\tfrac{1-y}{1+y})=
\\
&=\frac{1+x^2}{1-x^2},
&&&&&
&=\frac{1+y^2}{1-y^2}.
\end{aligned}
\eqlbl{cosha+coshb}
\]

Note that 
\begin{align*}
B'C&=\tfrac1x,
&
A'C&=\tfrac1y.
\shortintertext{Therefore}
BB'&=\tfrac1x-x,
&
AA'&=\tfrac1y-y.
\intertext{Since the triangles $ABC$, $A'BC$, $AB'C$, $A'B'C$ are right, the original Pythagorean theorem (\ref{thm:pyth}) implies}
AB&=\sqrt{x^2+y^2},
&
AB'&=\sqrt{\tfrac1{x^2}+y^2},
\\
A'B&=\sqrt{x^2+\tfrac1{y^2}},
&
A'B'&=\sqrt{\tfrac1{x^2}+\tfrac1{y^2}}.
\end{align*}

According to Exercise~\ref{ex:cosh}\textit{\ref{ex:coshcosh}},
\[
\begin{aligned}
\cosh c&= \frac{AB\cdot A'B'+AB'\cdot A'B}{AA'\cdot BB'}=
\\
&=
\frac{\sqrt{x^2+y^2}\cdot \sqrt{\tfrac1{x^2}
+\tfrac1{y^2}}+\sqrt{\tfrac1{x^2}+y^2}\cdot \sqrt{x^2+\tfrac1{y^2}}}
{(\tfrac1y-y)\cdot (\tfrac1x-x)}=
\\
&=\frac{x^2+y^2+1+x^2\cdot y^2}
{(1-y^2)\cdot (1-x^2)}=
\\
&=\frac{1+x^2}
{1-x^2}\cdot
\frac{1+y^2}
{1-y^2}.
\end{aligned}
\eqlbl{coshc}
\]
Finally, \ref{cosha+coshb} and \ref{coshc} imply \ref{eq:thm:pyth-h-poincare}.
\qeds

\addtocontents{toc}{\protect\contentsline{part}{\protect\numberline{}Additional topics}{}{}}

\chapter{Affine geometry}\label{chap:trans}

\section{Affine transformations}

\emph{Affine geometry} studies the so-called \index{incidence structure}\emph{incidence structure} of the Euclidean plane.
The incidence structure sees only which points lie on which lines and nothing else;
it does not directly see distances, angle measures, and many other things.

A bijection from the Euclidean plane to itself is called an \index{affine transformation}\emph{affine transformation} if it maps lines to lines;
that is, the image of any line is a line.
So, one can say that affine geometry studies the properties of the Euclidean plane that are preserved under affine transformations.

\begin{thm}{Classroom exercise}\label{ex:affine-par}
Show that an affine transformation of the Euclidean plane sends every pair of parallel lines to a pair of parallel lines.
\end{thm}

The observation below follows since the lines are defined using the metric only.

\begin{thm}{Observation}
Any motion of the Euclidean plane is an affine transformation.
\end{thm}

The following exercise provides more general examples of affine transformations.

\begin{thm}{Exercise}\label{ex:afine-linear}
Show that the following maps of a coordinate plane to itself define affine transformations:
\begin{enumerate}[(a)]
\item\label{ex:afine-linear:shear} \index{shear map}\emph{Shear map} defined by $(x,y)\mapsto (x+k\cdot y,y)$ for a constant~$k$.
\item\label{ex:afine-linear:scaling}
Rescaling defined by $(x,y)\mapsto (a\cdot x,a\cdot y)$ for a constant $a\ne 0$.
\item $x$-scaling and $y$-scaling defined respectively by 
\[(x,y)\mapsto (a\cdot x,y),\quad\text{and}\quad(x,y)\mapsto (x,a\cdot y)\]
for a constant $a\ne 0$.
\item\label{affine-general-formula} A transformation defined by
\[(x,y)\mapsto(a\cdot x+b\cdot y+r,c\cdot x+d\cdot y+s)\]
for constants $a,b,c,d,r,s$ such that the matrix $(\begin{smallmatrix}a&b\\c&d\end{smallmatrix})$ is invertible. 
\end{enumerate}
\end{thm}

From the fundamental theorem of affine geometry (\ref{thm:fundamental-theorem-of-affine-geometry}), it will follow that any affine transformation can be written in the form \textit{(\ref{affine-general-formula})}.

Recall that points are \index{collinear}\emph{collinear} if they lie on a single line.

\begin{thm}{Exercise}\label{ex:collinear=affine}
Suppose $P\mapsto P'$ is a bijection of the Euclidean plane that maps collinear triples of points to collinear triples.
Show that $P\mapsto P'$ maps noncollinear triples to noncollinear ones.

Conclude that $P\mapsto P'$ is an affine transformation.
\end{thm}

\begin{thm}{Exercise}\label{ex:circle=affine}
Let $\alpha$ be a bijection of the Euclidean plane that maps every circle to a circle.
Show that $\alpha$ is an affine transformation.
\end{thm}

\section{Constructions}

\begin{wrapfigure}{r}{34mm}
\vskip-15mm
\centering
\includegraphics{mppics/pic-287}
\end{wrapfigure}

Let us consider geometric constructions with a ruler and a \index{parallel!tool}\emph{parallel tool};
the latter allows us to draw a line thru a given point parallel to a given line.
(One may think of the tool on the figure.
It consists of two straight edges joined by two arms that can move while remaining parallel to each other.)

By Exercisers~\ref{ex:affine-par}, any construction with these two tools is invariant with respect to affine transformations.
For example, 
to solve the following exercise,
it is sufficient to prove that the midpoint of a given segment can be constructed with a ruler and a parallel tool.

\begin{thm}{Exercise}\label{ex:midpoint-affine}
Let $M$ be the  midpoint of a segment $[AB]$ in the Euclidean plane.
Assume that an affine transformation sends the points $A$, $B$, and $M$
to $A'$, $B'$, and $M'$ respectively.
Show that $M'$ is the midpoint of~$[A'B']$.
\end{thm}

The following exercise will be used in the proof of Proposition~\ref{prop:affine-linear}.

\begin{thm}{Exercise}\label{ex:R-hom}
Assume that points with coordinates $(0,0)$, $(1,0)$, $(a,0)$, and $(b,0)$ are given.
Using a ruler and a parallel tool, construct points with coordinates $(a+b,0)$, $(a-b,0)$, $(a\cdot b,0)$, and $(\tfrac a b,0)$.
\end{thm}

\begin{thm}{Exercise}\label{ex:center-circ-affine}
Construct the center of a given circle using the ruler and the parallel tool.
\end{thm}

The shear map (described in \ref{ex:afine-linear}\textit{\ref{ex:afine-linear:shear}}) can change angles between lines almost arbitrarily.
This observation can be used to prove the impossibility of some constructions;
here is one example:

\begin{thm}{Exercise}\label{ex:affine-perp}
Show that one cannot construct a line perpendicular to a given line with the ruler and the parallel tool.
\end{thm}

\section{Fundamental theorem of affine geometry}

In this section, we assume knowledge of vector algebra; namely, multiplication by a real number, addition, and the parallelogram rule.

\begin{thm}{Exercise}\label{ex:parallelogram-rule}
Show that affine transformations map parallelograms to parallelograms.
Conclude that if $P\mapsto P'$ is an affine transformation, then
\[\overrightarrow{XY}=\overrightarrow{AB},
\quad\text{if and only if}\quad
\overrightarrow{X'Y'}=\overrightarrow{A'B'}.\]

\end{thm}

\begin{thm}{Proposition}\label{prop:affine-linear}
Let $P\mapsto P'$ be an affine transformation of the Euclidean plane.
Then, for every triple of points $O$, $X$, $P$, we have
\[\overrightarrow{OP}=t\cdot \overrightarrow{OX}
\quad\text{if and only if}\quad
\overrightarrow{O'P'}=t\cdot \overrightarrow{O'X'}.
\eqlbl{eq:OP=tOX}\]

\end{thm}

\parit{Proof.}
Observe that the affine transformations described in Exercise~\ref{ex:afine-linear}, as well as all motions, satisfy the condition \ref{eq:OP=tOX}.
Therefore a given affine transformation $P\mapsto P'$ satisfies \ref{eq:OP=tOX} if and only if its composition with motions and rescalings satisfies \ref{eq:OP=tOX}.

Applying this observation, we can reduce the problem to its partial case.
Namely, we may assume that $O'=O$, $X'=X$, the point $O$ is the origin of a coordinate system, and $X$ has coordinates $(1,0)$.

In this case, $\overrightarrow{OP}=t\cdot \overrightarrow{OX}$ if and only if $P=(t,0)$.
Since $O$ and $X$ are fixed, the transformation maps the $x$-axis to itself.
That is, $P'=(f(t),0)$ for a function $t\mapsto f(t)$,
or, equivalently, $\overrightarrow{O'P'}=f(t)\cdot \overrightarrow{O'X'}$.
It remains to show that 
\[f(t)=t
\eqlbl{eq:f(x)=x}\]
for any $t$.

Since $O'=O$ and $X'=X$, we get that $f(0)=0$ and $f(1)=1$.
Furthermore, according to Exercise~\ref{ex:R-hom}, we have that 
$f(x\cdot y)=f(x)\cdot f(y)$ and $f(x+y)=f(x)+f(y)$ for any $x,y\in\mathbb{R}$.
By the algebraic lemma (proved below, see \ref{lem:R-auto}), these conditions imply \ref{eq:f(x)=x}.
\qeds

\begin{thm}{Fundamental theorem of affine geometry}\label{thm:fundamental-theorem-of-affine-geometry}\index{Fundamental theorem of affine geometry}
Suppose that an affine transformation $P\mapsto P'$ maps a nondegenerate triangle $OXY$ to a triangle $O'X'Y'$.
Then $\triangle O'X'Y'$ is nondegenerate, and
\[\overrightarrow{OP}=x\cdot\overrightarrow{OX}+y\cdot\overrightarrow{OY}\quad\text{if and only if}\quad\overrightarrow{O'P'}=x\cdot\overrightarrow{O'X'}+y\cdot\overrightarrow{O'Y'}.\]
\end{thm}

\parit{Proof.}
Since an affine transformation maps lines to lines, the triangle
$O'X'Y'$ is nondegenerate.

Consider points $V$ and $W$ defined by
\[\overrightarrow{OV}=x\cdot\overrightarrow{OX},
\qquad
\overrightarrow{OW}=y\cdot\overrightarrow{OY}.
\]

{

\begin{wrapfigure}{o}{33mm}
\centering
\vskip-0mm
\includegraphics{mppics/pic-225}
\end{wrapfigure}

Note that 
$\overrightarrow{WP}=\overrightarrow{OV}$.
Applying Exercise~\ref{ex:parallelogram-rule} and the proposition, we get
\begin{align*}\overrightarrow{O'P'}&=\overrightarrow{O'W'}+\overrightarrow{W'P'}=
\\
&=\overrightarrow{O'W'}+\overrightarrow{O'V'}=
\\
&=x\cdot\overrightarrow{O'X'}+y\cdot\overrightarrow{O'Y'}.
\end{align*}
\qedsf

}

\begin{thm}{Exercise}\label{ex:affine-continuous}
Show that every affine transformation is continuous.
\end{thm}

The following exercise provides the converse to Exercise \ref{ex:afine-linear}\textit{\ref{affine-general-formula}}.

\begin{thm}{Exercise}\label{ex:affine-coordinates}
Show that every affine transformation can be written in coordinates as
$(x,y)\mapsto(a\cdot x+b\cdot y+r,\ c\cdot x+d\cdot y+s)$
for constants $a,b,c,d,r,s$ such that the matrix $(\begin{smallmatrix}a&b\\c&d\end{smallmatrix})$ is invertible. 
\end{thm}

\begin{thm}{Exercise}\label{ex:preserved-circle}
Let $\Gamma$ be a circle with center $O$.
Suppose $\alpha\:P\mapsto P'$ is an affine transformation that maps $\Gamma$ to itself.
Show that $\alpha$ is a motion that fixes $O$;
that is, $O'=O$.
\end{thm}

A bijection from the inversive plane to itself is called an \index{inversive!transformation}\emph{inversive transformation} if it maps circlines to circlines.
By \ref{thm:inverse-cline},  inversions are  inversive transformations.

\begin{thm}{Advaneced exercise}\label{ex:inversions-inversive}
Show that every inversive transformation is a composition of inversions and reflections.
\end{thm}

\section{An algebraic lemma}

The following lemma was used in the proof of Proposition~\ref{prop:affine-linear}.

\begin{thm}{Lemma}\label{lem:R-auto}
Assume $f\:\mathbb{R}\to\mathbb{R}$ is a function such that for any $x,y\in\mathbb{R}$ we have
\begin{enumerate}[(a)]
\item\label{lem:R-auto:a} $f(1)=1$,
\item\label{lem:R-auto:b} $f(x+y)=f(x)+f(y)$,
\item\label{lem:R-auto:c} $f(x\cdot y)=f(x)\cdot f(y)$.
\end{enumerate}

Then $f$ is the identity function; that is,
$f(x)=x$ for any $x\in \mathbb{R}$.
\end{thm}

Note that we did not assume that $f$ is continuous.

A function $f$ satisfying these three conditions
is called a \index{field automorphism}\emph{field automorphism}.
Therefore, the lemma states that the identity function is the only automorphism of the field of real numbers.
For the field of complex numbers, the conjugation $z\mapsto\bar z$ (see Section~\ref{sec:complex-conjugation}) gives an example of a nontrivial automorphism.

\begin{thm}{Exercise}\label{ex:f(1)=1}
Suppose that $f\:\mathbb{R}\to\mathbb{R}$ satisfies the condition (\ref{lem:R-auto:c}) in the lemma.
Show that 
\begin{enumerate}[(a)]
 \item if $f(1)\ne 1$, then $f(x)=0$ for any $x$;
  \item if $f(0)\ne 0$, then $f(x)=1$ for any $x$.
\end{enumerate}
\end{thm}

\parit{Proof.}
By \textit{(\ref{lem:R-auto:b})} we have
$f(0)+f(0)=f(0+0)=f(0)$.
Therefore,
\[f(0)=0.
\eqlbl{eq:0=0}\]

Applying \textit{(\ref{lem:R-auto:b})} again, we get that $0=f(0)=f(x)+f(-x)$.
Therefore, 
\[f(-x)=-f(x)
\quad
\text{for any}
\quad
x\in \mathbb{R}.
\eqlbl{eq:f-x}\] 

Applying \textit{(\ref{lem:R-auto:a})} and \textit{(\ref{lem:R-auto:b})} recursively, we get that
\begin{align*}
f(2)&=f(1)+f(1)=1+1=2;\\
f(3)&=f(2)+f(1)=2+1=3;\\
&\dots
\end{align*}
Together with \ref{eq:f-x},
the latter implies that 
$f(n)=n$
for every integer
$n$. 
By~\textit{(\ref{lem:R-auto:c})},
$f(m)=f(\tfrac mn)\cdot f(n)$.
Therefore
$$f(\tfrac mn)=\tfrac mn \eqlbl{eq:m/n}$$
for every rational number~$\tfrac mn$.

Assume~$a\ge 0$.
Then the equation $x^2=a$ has a real solution $x\z=\sqrt{a}$.
Suppose $y=f(x)$.
By~\textit{(\ref{lem:R-auto:c})}, $f(a)=y^2$,
and hence $f(a)\ge 0$.

\raggedcolumns\setlength{\multicolsep}{.5mm}
\setlength{\columnseprule}{1pt}
\begin{multicols}{2}
That is,
\[a\ge 0\quad\Longrightarrow\quad f(a)\ge 0.\eqlbl{a>0=>b>0}\]

\columnbreak

By \ref{eq:f-x}, 
we also get 
\[a\le 0\quad \Longrightarrow\quad f(a)\le 0.\eqlbl{a<0=>b<0}\]
\end{multicols}
\setlength{\columnseprule}{0pt}

Now assume $f(a)\ne a$ for some $a\in\mathbb{R}$.
Then there is a rational number $\tfrac{m}{n}$ that lies between $a$ and $f(a)$;
that is, 
the numbers 
\[x\z=a-\tfrac{m}{n}\quad\text{and}\quad y=f(a)-\tfrac{m}{n}\]
have opposite signs.

By \ref{eq:m/n},
\begin{align*}
y+\tfrac{m}{n}&=f(a)=
\\
&=f(x+\tfrac{m}{n})=
\\
&=f(x)+f(\tfrac{m}{n})=
\\
&=f(x)+\tfrac{m}{n};
\end{align*}
that is, $f(x)=y$.
By \ref{a>0=>b>0} and \ref{a<0=>b<0} the values $x$ and $y$ cannot have opposite signs --- a contradiction.
\qeds

\section{Three theorems}

We discussed several statements in Euclidean geometry that have an affine nature;
that is, their assumptions and conclusions survive under affine transformations.
The examples include 
\ref{thm:parallel-point-reflection},
\ref{lem:parallelogram}\ref{lem:parallelogram:midpoint},
\ref{ex:4parallels},
\ref{thm:centroid},
and \ref{ex:midle}.
In this section, we present three more examples of that type.

Suppose that $A$, $B$, and $X$ are distinct points on one line.
Note that 
\[\overrightarrow{AX}=t\cdot \overrightarrow{BX},\]
where $t=\pm\frac{AX}{BX}$.

Choose an affine transformation $P\mapsto P'$.
By \ref{prop:affine-linear}, we have 
\[\overrightarrow{A'X'}=t\cdot \overrightarrow{B'X'}.\]
It follows that $\frac{A'X'}{B'X'}=\frac{AX}{BX}$;
that is, \textit{if $A$, $B$, and $X$ are distinct points on one line, then the ratio $\frac{AX}{BX}$ is preserved under affine transformations}.

\begin{thm}{Menelaus's theorem}\index{Menelaus's theorem}
Let $ABC$ be a nondegenerate triangle.
Suppose a line $\ell$ crosses the lines $(BC)$, $(CA)$, and $(AB)$ at three distinct points points $A'$, $B'$, and $C'$.
Then 
\[\frac{AC'}{BC'}\cdot\frac{BA'}{CA'}\cdot \frac{CB'}{AB'}=1.\]
\end{thm}

As we saw, the ratios $\frac{AC'}{BC'}$, $\frac{BA'}{CA'}$, $\frac{CB'}{AB'}$ are preserved under affine transformations.
Therefore Menelaus's theorem belongs to affine geometry.
However, in the proof we can (and will) use Euclidean-type arguments.

\parit{Proof.}
Let $X$, $Y$ and $Z$ be the footpoints of $A$, $B$ and $C$ on $\ell$.
By the AA similarity condition, $\triangle AXC'\sim \triangle BYC'$.
\begin{figure}[!ht]
\centering
\includegraphics{mppics/pic-227}
\vskip-10mm
\end{figure}
Therefore, 
\[\frac{AC'}{BC'}=\frac{AX}{BY}.\]
The same way we get
\[\frac{BA'}{CA'}=\frac{BY}{CZ}
\quad\text{and}\quad
\frac{CB'}{AB'}=\frac{CZ}{AX}.\]
Therefore,
\[\frac{AC'}{BC'}\cdot\frac{BA'}{CA'}\cdot \frac{CB'}{AB'}=\frac{AX}{BY}\cdot \frac{BY}{CZ}\cdot \frac{CZ}{AX}=1.\]
\qedsf

{

\begin{thm}{Exercise}\label{thm:ceva-affine}\label{ex:ceva-affine}
Let $ABC$ be a nondegenerate triangle.
Suppose three distinct points points $A'$, $B'$, and $C'$ lie on the lines $(BC)$, $(CA)$, and $(AB)$ respectively.
Assume that the lines $(AA')$, $(BB')$ and $(CC')$ meet at a point $X$.

\begin{wrapfigure}{r}{31mm}
\centering
\vskip-2mm
\includegraphics{mppics/pic-229}
\end{wrapfigure}

Use Menelaus's theorem to prove the following.

\begin{enumerate}[(a)]
\item \index{Ceva's theorem}\textbf{Ceva's theorem:}
\[\frac{AC'}{BC'}\cdot\frac{BA'}{CA'}\cdot \frac{CB'}{AB'}=1.\]
\item \index{Van Aubel's theorem}\textbf{Van Aubel's theorem:} If $A'$ lies between $B$ and $C$, then
\[\frac{AC'}{BC'}+\frac{AB'}{CB'}=\frac{AX}{A'X}.\]
\end{enumerate}
\end{thm}

}

\parbf{Signed versions.}
These three theorems can be formulated in a more precisely way by taking the ratios $\frac{AC'}{BC'}$, $\frac{BA'}{CA'}$, and $\frac{CB'}{AB'}$ with sign.
Namely, given distinct points $A$, $B$, and $X$ on a line, let us define the \index{signed ratio}\emph{signed ratio} $(\frac{AX}{BX})\df\pm\frac{AX}{BX}$, where the sign is chosen to meet the equation
$\overrightarrow{AX}=\left(\frac{AX}{BX}\right)\cdot \overrightarrow{BX}$.
Equivalently, $(\frac{AX}{BX})<0$ if and only if $X$ lies between $A$ and $B$.

The more exact identities of Menelaus's and Ceva's theorems will be
\[\left(\tfrac{AC'}{BC'}\right)\cdot\left(\tfrac{BA'}{CA'}\right)\cdot \left(\tfrac{CB'}{AB'}\right)=1
\quad\text{and}\quad
\left(\tfrac{AC'}{BC'}\right)\cdot\left(\tfrac{BA'}{CA'}\right)\cdot \left(\tfrac{CB'}{AB'}\right)=-1
\]
respectively.
In Van Aubel's theorem we have
\[\left(\tfrac{AC'}{BC'}\right)+\left(\tfrac{AB'}{CB'}\right)=\left(\tfrac{AX}{A'X}\right)\]
and this identity holds without the assumption that $A'$ lies between $B$ and~$C$.

Moreover, the converse of each of the three signed versions also holds.

\chapter{Projective geometry}\label{chap:proj}

\section{The projective completion}

In the Euclidean plane, two distinct lines might have one or zero points of intersection 
(in the latter case the lines are parallel).
We aim to extend the Euclidean plane by ideal points so that any two distinct lines will have exactly one point of intersection.

\begin{wrapfigure}{o}{31mm}
\centering
\includegraphics{mppics/pic-226}
\vskip4mm
\includegraphics{mppics/pic-228}
\end{wrapfigure}

A collection of lines in the Euclidean plane is called \index{concurrent}\emph{concurrent} if they all intersect at a single point 
or if they are all pairwise parallel.
A maximal set of concurrent lines in the plane is called a \index{pencil}\emph{pencil}.
There are two types of pencils: 
\emph{central pencils} contain all lines passing thru a fixed point called the \index{center!center of the pencil}\emph{center of the pencil}
and  
\emph{parallel pencils} contain pairwise parallel lines.

Any two lines completely determine the pencil containing both.

A point $P$ in the Euclidean plane uniquely defines a central pencil with the center at $P$.
Let us add one \index{point!ideal point}\index{ideal!point}\emph{ideal point} for each parallel pencil,
assuming that all these ideal points lie on one \index{line!ideal line}\index{ideal!line}\emph{ideal line}.
We also assume that the ideal line belongs to each parallel pencil.

We obtain the so-called \index{projective!plane}\emph{projective plane} (or  the \index{projective!completion}\emph{projective completion} of the original plane). 
It comes with an incidence structure --- we say that three points lie on one line if the corresponding pencils contain a common line.
Projective geometry studies this incidence structure.

A parallel pencil contains the ideal line and the lines represented by equations of the form $y=m\cdot x+b$ with a fixed slope $m$.
If $m=\infty$, we assume that the lines are given by equations of the form $x=a$.
Therefore the projective completion contains all the points in the $(x,y)$-plane along with the ideal line containing one ideal point $P_m$ for every slope $m\z\in\mathbb{R}\cup\{\infty\}$.

\section{Euclidean space}

Let us revisit the construction of the metric $d_2$ (Exercise~\ref{ex:dist-square}), but this time in three-dimensional space. 

Consider $\mathbb{R}^3$, which represents the set of all triples $(x,y,z)$ of real numbers.
Let  $A=(x_A,y_A,z_A)$ and $B=(x_B,y_B,z_B)$ be arbitrary points in $\mathbb{R}^3$.
Define the metric on $\mathbb{R}^3$ as follows:
$$AB
\df
\sqrt{(x_A-x_B)^2+(y_A-y_B)^2+(z_A-z_B)^2}.$$
The obtained metric space is called \index{Euclidean!space}\emph{Euclidean space}.

A subset of points in $\mathbb{R}^3$ is called a \index{plane!plane in the space}\emph{plane} if it can be
described by an equation of the form
$$a\cdot x+b\cdot y+c\cdot z+d=0,$$ 
where $a$, $b$, $c$, and $d$ are constants, and at least one of the values $a$, $b$ or $c$ is distinct from zero.

It is straightforward to show the following:
\begin{itemize}
 \item Any plane in the Euclidean space is isometric to the Euclidean plane.
 \item Any three points in the Euclidean space lie on a plane.
 \item The  intersection of two distinct planes (if it is nonempty) forms a line in each of these planes.
\end{itemize}

These statements enable us to generalize many results from Euclidean plane geometry to Euclidean space.

\section{Space model}

Let us identify the Euclidean plane with a plane $\Pi$ in the Euclidean space $\mathbb{R}^3$ that does not pass thru the origin $O$.
Denote by $\hat\Pi$ the projective completion of $\Pi$.

Denote by $\Phi$ the set of all lines in $\mathbb{R}^3$ thru $O$.
Let us define a bijection $P\leftrightarrow \dot P$ between $\hat \Pi$ and $\Phi$.
If $P\in \Pi$, we take the line $\dot P=(OP)$;
if $P$ is an ideal point of $\hat \Pi$ defined by a parallel pencil of lines, we take the line $\dot P$ thru $O$ and parallel to the lines in this pencil.

Furthermore, let $\Psi$ be the set of all planes in $\mathbb{R}^3$ thru $O$.
In a similar fashion, we can define a bijection $\ell\leftrightarrow \dot \ell$ between lines in $\hat \Pi$ and $\Psi$.
If a line $\ell$ is not ideal, we take the plane $\dot \ell$ that contains $\ell$ and $O$;
if the line $\ell$ is ideal, we take $\dot \ell$ to be the plane passing thru $O$ and parallel to $\Pi$ (that is, $\dot\ell\cap\Pi=\emptyset$).

\begin{thm}{Observation}\label{obs:bijections}
Let $P$ and $\ell$ be a point and a line in the projective plane.
Then 
\[P\in \ell \quad\iff\quad \dot P\subset \dot \ell,\]
where $\dot P$ and $\dot \ell$ are defined by the constructed bijections.
\end{thm}

\section{The perspective projection}
\label{sec:perspective-projection}

Let $O$ be a point that does not lie on planes $\Pi$ and $\Pi'$ in $\mathbb{R}^3$.

\begin{wrapfigure}{o}{43mm}
\centering
\vskip-4mm
\includegraphics{mppics/pic-230}
\end{wrapfigure}

A \index{perspective projection}\emph{perspective projection from $\Pi$ to $\Pi'$ with center $O$} maps a point $P\z\in \Pi$
to the intersection point $P'$ of the plane $\Pi'$ and the line $(OP)$.

In general, perspective projection is not a bijection between the planes.
Indeed, if the line $(OP)$ is parallel to $\Pi'$ 
(that is, if $(OP)\cap\Pi'=\emptyset$)
then the perspective projection of $P\in \Pi$ is undefined.
Also, if $(OP')\parallel \Pi$ 
for $P'\in \Pi'$,
then the point $P'$ is not an image of the perspective projection.

Denote by $\hat \Pi$ and $\hat \Pi'$ the projective completions of $\Pi$ and $\Pi'$ respectively. 
The perspective projection is a restriction of the composition of two bijections $\hat \Pi\leftrightarrow\Phi \leftrightarrow\hat \Pi'$ constructed in the previous section.
By Observation~\ref{obs:bijections}, the perspective projection can be extended to a bijection $\hat \Pi\leftrightarrow\hat \Pi'$ that sends lines to lines.%
\footnote{A similar story happened with inversion.
Inversion is not defined at its center;
moreover, the center is not an the inverse of any point.
To deal with this problem, we introduced the inversive plane, 
which is the Euclidean plane extended by one ideal point.
The same strategy works for the perspective projection $\Pi\to\Pi'$, but this time we need to add an ideal line.}

For example, suppose $O$ is the origin of the $(x,y,z)$-coordinate space,
and the planes $\Pi$ and $\Pi'$ are defined by the equations
$z=1$ and $x=1$ respectively.
Then the perspective projection from $\Pi$ to $\Pi'$
can be written in coordinates as
\[(x,y,1)\mapsto (1,\tfrac yx,\tfrac 1x).\]
Indeed the coordinates have to be proportional;
points on $\Pi$ have a unit $z$-coordinate, 
and points on $\Pi'$ have a unit $x$-coordinate.

The perspective projection maps one plane to another.
However, we can identify the two planes by fixing a coordinate system in each.
In this case, we get a partially defined map from the plane to itself.
We will keep the name {}\emph{perspective transformation} for such maps.

For the described perspective projection, we may get the map 
\[\beta\:(x,y)\mapsto (\tfrac 1x,\tfrac yx).
\eqlbl{eq:(x,y)-perspective}\]
This map is undefined on the line $x=0$.
Also, points on this line are not images of points under the perspective projection.

For example, to define an extension of the perspective projection $\beta$ in \ref{eq:(x,y)-perspective},
we have to observe that 
\begin{itemize}
\item The pencil of vertical lines $x=a$ is mapped to itself.
\item The ideal points defined by pencils of lines $y=m\cdot x+ b$ are mapped to the point $(0,m)$, and the other way around --- point $(0,m)$ is mapped to the ideal point defined by the  pencil of lines $y=m\cdot x+ b$.
\end{itemize}

\section{Projective transformations}

A bijection from the projective plane to itself 
that sends lines to lines 
is called \index{projective!transformation}\emph{projective transformation}.

Any affine transformation defines  a projective transformation on the corresponding projective plane.
We will call such projective transformations \index{affine transformation}\emph{affine}; 
these are projective transformations that send the ideal line to itself.

The extended perspective projection, as  discussed in the previous section, 
provides another source of examples of projective transformations.

\begin{thm}{Theorem}\label{thm:moving}
Given a line $\ell$ in the projective plane, there is a perspective projection that sends $\ell$ to the ideal line.

Moreover, a perspective transformation is either affine or, in a suitable coordinate system, it can be written as a composition of the extension of perspective projection 
\[\beta\:(x,y)\mapsto (\tfrac xy,\tfrac 1y)\]
and an affine transformation.
\end{thm}

\parit{Proof.}
Choose an $(x,y)$-coordinate system so that the line $\ell$ is defined by the equation $y=0$. 
Then the extension of $\beta$ gives the needed transformation.

Fix a projective transformation $\gamma$.
If $\gamma$ sends the ideal line to itself,
then it must to be affine. 
In this case, the theorem is proved.

Suppose $\gamma$ sends the ideal line to a line $\ell$.
Choose a perspective projection $\beta$ as above.
The composition $\beta\circ\gamma$ sends the ideal line to itself.
That is, $\alpha=\beta\circ\gamma$ is affine.
Note that $\beta$ is self-inverse; therefore 
\[\gamma=\beta\circ\beta\circ\gamma=\beta\circ\alpha\]
--- hence the result.
\qeds

\begin{thm}{Exercise}\label{ex:proj-cross-ratio}
Let $P\mapsto P'$ be (a) an affine transformation, (b) the perspective projection defined by $(x,y)\mapsto (\tfrac xy,\tfrac 1y)$, or (c) an arbitrary projective transformation.
Suppose $P_1,P_2,P_3,P_4$ are distinct points on one line.
Show that 
\[\frac{P_1P_2\cdot P_3P_4}{P_2P_3\cdot P_4P_1}=\frac{P'_1P'_2\cdot P'_3P'_4}{P'_2P'_3\cdot P'_4P'_1}\]
if none of these points lies on the ideal line.
In other words, each of these maps preserves the cross-ratio for quadruples of points on one line.

\end{thm}

{

\begin{wrapfigure}{o}{60mm}
\vskip-0mm
\centering
\includegraphics{mppics/pic-231}
\end{wrapfigure}

\begin{thm}{Advanced exercise}\label{ex:proj-cross-ratio=1}
Let $A$, $B$, $C$, $P$, $Q$, $V$, $W$, $X$, $Y$ be as in the picture.
\begin{enumerate}[(a)]
\item\label{ex:proj-cross-ratio=1:=} Show that 
\[\frac{AX\cdot BY}{AY\cdot BX}=\frac{AP\cdot CQ}{AQ\cdot CP}.\]
\item\label{ex:proj-cross-ratio=1:1} 
Use part (\ref{ex:proj-cross-ratio=1:=}) to show that
\[AX\cdot BY=AY\cdot BX.\]
\end{enumerate}

\end{thm}

}

\section{Moving points to infinity}

{

\begin{wrapfigure}{r}{37mm}
\vskip-8mm
\centering
\includegraphics{mppics/pic-234}
\end{wrapfigure}

Theorem~\ref{thm:moving} allows us declare that a given line is ideal.
In other words, we can choose a preferred affine plane by removing one line from the projective plane.
This construction provides a method for solving problems in projective geometry 
which will be illustrated by the following classical example:

\begin{thm}{Desargues' theorem}\label{thm:desargues}\index{Desargues' theorem}
Consider three concurrent lines $(AA')$, $(BB')$, and $(CC')$ in the projective plane.
Let $X$, $Y$, and $Z$ be the intersection points of $(BC)$ and $(B'C')$, $(CA)$ and $(C'A')$,
and $(AB)$ and $(A'B')$, respectively.
Then the points $X$, $Y$, and $Z$ are collinear.
\end{thm}

}

\parit{Proof.}
We may assume that the line $(XY)$ is ideal.
If not, apply a perspective projection that sends the line $(XY)$ to the ideal line.

\begin{wrapfigure}{o}{40mm}
\vskip-0mm
\centering
\includegraphics{mppics/pic-236}
\end{wrapfigure}

That is, we can assume that 
\[(BC)\z\parallel (B'C')\quad\text{and}\quad(CA)\z\parallel (C'A'),\]
and we need to show that 
\[(AB)\z\parallel(A'B').\]

Assume that the lines $(AA')$, $(BB')$, and $(CC')$ intersect at point~$O$.
Since $(BC)\z\parallel (B'C')$, 
the transversal property (\ref{thm:parallel-2}) implies that $\measuredangle OBC\z= \measuredangle OB'C'$ and $\measuredangle OCB\z= \measuredangle OC'B'$.
By the AA similarity condition, $\triangle OBC\z\sim\triangle OB'C'$.
In particular,
\[\frac{OB}{OB'}=\frac{OC}{OC'}.\]

In the same way, we get that $\triangle OAC\z\sim\triangle OA'C'$ and
\[\frac{OA}{OA'}=\frac{OC}{OC'}.\]
Therefore, 
\[\frac{OA}{OA'}=\frac{OB}{OB'}.\]
By the SAS similarity condition, 
we get that $\triangle OAB\sim\triangle OA'B'$;
in particular, $\measuredangle OAB=\pm\measuredangle OA'B'$.

Note that $\measuredangle AOB=\measuredangle A'OB'$.
Therefore, 
\[\measuredangle OAB=\measuredangle OA'B'.\]
By the transversal property (\ref{thm:parallel-2}), we have
$(AB)\parallel (A'B')$.

The case $(AA')\parallel(BB')\parallel(CC')$ is done similarly.
In this case, the quadrangles $B'BCC'$ and $A'ACC'$ are parallelograms.
Therefore, 
\[BB'=CC'=AA'.\]
Hence $\square B'BAA'$ is a parallelogram and $(AB)\parallel (A'B')$.
\qeds

Here is another classical theorem of projective geometry.

\begin{thm}{Pappus' theorem}\label{thm:pappus}\index{Pappus' theorem}
Assume that two triples of points $A$, $B$, $C$,
and $A'$, $B'$, $C'$ are collinear.
Let $X$, $Y$, and $Z$ be the intersection points of $(BC')$ and $(B'C)$, $(CA')$ and $(C'A)$,
and $(AB')$ and $(A'B)$, respectively.
Then the points $X$, $Y$, $Z$ are collinear.
\end{thm}

Pappus' theorem can be proved the same way as Desargues' theorem.

\parit{Idea of the proof.}
Applying a perspective projection, we can assume that $Y$ and $Z$ lie on the ideal line.
It remains to show that $X$ lies on the ideal line.

In other words, assuming that $(AB')\parallel (A'B)$ and $(AC')\parallel (A'C)$, we need to show that $(BC')\parallel(B'C)$.

\begin{figure}[!ht]
\centering
\includegraphics{mppics/pic-238}
\hskip15mm
\includegraphics{mppics/pic-240}
\end{figure}

\begin{thm}{Exercise}\label{ex:pappus}
Finish the proof of Pappus' theorem using the idea described above.
\end{thm}

The following exercise gives a partial converse to Pappus' theorem.

\begin{thm}{Exercise}\label{ex:pappus-converse}
Given two triples of points $A$, $B$, $C$,
and $A'$, $B'$, $C'$,
suppose that distinct points $X$, $Y$, and $Z$ are the unique intersection points of $(BC')$ and $(B'C)$, $(CA')$ and $(C'A)$,
and $(AB')$ and $(A'B)$, respectively.
Assume that the triples $A$, $B$, $C$,
and $X$, $Y$, $Z$ are collinear.
Show that the triple $A'$, $B'$, $C'$ is collinear.
\end{thm}

\begin{thm}{Exercise}\label{ex:desargues-construction}
Solve the following construction problem
\begin{enumerate}[(a)]
\item\label{ex:desargues-construction:desargues} using Desargues' theorem;
\item\label{ex:desargues-construction:pappus} using Pappus' theorem.
\end{enumerate}
\parbf{Problem.}
Suppose a parallelogram and a line $\ell$ are given.
Assume that $\ell$ crosses all sides (or their extensions) of the parallelogram at different points. 
Construct another line parallel to $\ell$ with a ruler only.
\end{thm}

\section{Duality}

Assume $P\leftrightarrow p$ is a bijection between the set of lines and the set of points of a plane is given.
\begin{figure}[!ht]
\centering
\includegraphics{mppics/pic-242}
\hskip15mm
\includegraphics{mppics/pic-244}
\caption*{Dual configurations.}
\end{figure}
That is,
given a point $P$, we denote by $p$ the corresponding line;
and the other way around, 
given a line $\ell$ we denote by $L$ the corresponding point. 

The bijection between points and lines is called \index{duality}\emph{duality}\label{page:duality}%
\footnote{The standard definition of duality is more general; we consider a special case which is also called \index{polarity}\emph{polarity}.}
if 
\[P\in \ell
\quad
\iff
\quad 
p\ni L.\]
for any point $P$ and line~$\ell$.

{

\begin{wrapfigure}{r}{40mm}
\vskip-4mm
\centering
\includegraphics{mppics/pic-245}
\end{wrapfigure}

\begin{thm}{Exercise}\label{ex:dual-configurations}
Consider the configuration of lines and points in the picture.

\begin{enumerate}[(a)]
\item\label{ex:dual-configurations:infty} Redraw the picture after moving the line $q$ to infinity; mark the points and lines respectively.

\item\label{ex:dual-configurations:dual} Start with a generic quadrangle $KLMN$ and extend it to a dual picture; label the lines and points using the convention described above.
\end{enumerate}

\end{thm}

\begin{thm}{Exercise}\label{ex:dual-euclid}
Show that the Euclidean plane does not admit a duality. 
\end{thm}

}

\begin{thm}{Theorem}\label{thm:dual}
The projective plane admits a duality.
\end{thm}

\parit{Proof.}
Consider a plane $\Pi$ and a point $O\notin\Pi$ in the space;
suppose that $\hat \Pi$ denotes the corresponding projective plane.

Recall that $\Phi$ and $\Psi$ denote the set of all lines and planes passing thru~$O$.
By \ref{obs:bijections}, there are bijections $P\leftrightarrow\dot P$  between points of $\hat\Pi$ and $\Phi$ and $\ell\leftrightarrow\dot\ell$ between lines in $\hat\Pi$ and $\Psi$ such that 
$P\in\ell$ if and only if $\dot P\subset \dot \ell$.

It remains to construct a bijection $\dot \ell \leftrightarrow \dot L$
between $\Phi$ and $\Psi$ 
such that 
\[\dot P\subset \dot \ell
\quad
\iff
\quad
\dot p\supset \dot L
\eqlbl{iff-dual}\]
for any two lines $\dot P$ and $\dot L$ passing thru~$O$.

Set $\dot \ell$ to be the plane thru $O$ 
that is perpendicular to~$\dot L$.
Both conditions \ref{iff-dual} are equivalent to $\dot P\perp \dot L$;
hence the result follows.
\qeds

\begin{thm}{Exercise}\label{ex:dula-coordinates}
Consider the Euclidean plane with $(x,y)$-coordinates; suppose that $O$ denotes the origin.
Given a point $P\ne O$ with coordinates $(a,b)$ consider the line $p$ 
given by the equation 
$a\cdot x+b\cdot y=1$.

Show that the correspondence $P$ to $p$ can be extended to a duality of the projective plane.

Which line corresponds to $O$?

Which point corresponds to the line  $a\cdot x\z+b\cdot y=0$?
\end{thm}

Duality says that lines and points have the same rights in terms of incidence.
In particular, we can formulate an equivalent dual statement for  any statement in projective geometry.
For example, the dual statement of  ``the points $X$, $Y$, and $Z$ lie on one line $\ell$''
would be the ``lines $x$, $y$, and $z$ intersect at one point $L$''.
Let us formulate the dual statement for Desargues' theorem~\ref{thm:desargues}.

\begin{thm}{Dual Desargues' theorem}\label{thm:dual-desargues}\index{Desargues' theorem}
Let $X$, $Y$, and $Z$ be the unique intersection points of $(BC)$ and $(B'C')$, $(CA)$ and $(C'A')$,
and $(AB)$ and $(A'B')$, respectively.
If the points $X$, $Y$, and~$Z$ are collinear, then the lines $(AA')$, $(BB')$, and $(CC')$ are concurrent.
\end{thm}

In this theorem, the points $X$, $Y$, and $Z$ 
are dual to the lines $(AA')$, $(BB')$, and $(CC')$ in the original formulation, and the other way around.

Once Desargues' theorem is proved, applying duality (\ref{thm:dual})
we get the dual Desargues' theorem.
Note that the dual Desargues' theorem is the converse to the original Desargues' theorem~\ref{thm:desargues}.
So, we get a stronger version of the theorem for free: \textit{the points $X$, $Y$, and~$Z$ are collinear {}\emph{if and only if} the lines  $(AA')$, $(BB')$, and $(CC')$ are concurrent.}

\begin{thm}{Exercise}\label{ex:dual-pappus}
Formulate the dual Pappus' theorem (see \ref{thm:pappus}).
\end{thm}

\begin{thm}{Exercise}\label{ex:dual-desargues-construction} 
Solve the following construction problem
\begin{enumerate}[(a)]
\item\label{ex:dual-desargues-construction:desargues} using dual Desargues' theorem;
\item\label{ex:dual-desargues-construction:pappus} using Pappus' theorem or its dual.
\end{enumerate}
\parbf{Problem.}
Given two parallel lines, construct a third parallel line thru a given point with a ruler only.
\end{thm}

\section{The polar and the pole}

In this section, we describe a powerful trick that can be used in the constructions with a ruler.

Assume $\Gamma$ is a circle in the plane and $P\notin \Gamma$.
\begin{figure}[!ht]
\centering
\includegraphics{mppics/pic-290}
\end{figure}
Draw two lines $x$ and $y$ thru $P$ that intersect $\Gamma$ at two pairs of points $X$, $X'$ and $Y$, $Y'$.
Let $Z=(XY)\cap(X'Y')$ and $Z'=(XY')\cap(X'Y)$.
Consider the line $p=(ZZ')$.

\begin{thm}{Claim}\label{clm:polar}
The constructed line $p=(ZZ')$ does not depend on the choice of lines $x$ and $y$.

Moreover, 

\begin{enumerate}[(a)]
\item $P\leftrightarrow p$ can be extended to a duality such that any point $P$ on the circle $\Gamma$ corresponds to a line $p$ that is tangent to $\Gamma$ at~$P$.
\item If $O$ is the center of $\Gamma$, then the inversion of $P$ across $\Gamma$ is the footpoint of $O$ on~$p$.
\end{enumerate}
 
\end{thm}

We will not prove this claim, but the proof is not hard.
If $P$ lies outside of $\Gamma$, it can be done by moving $P$ to infinity while keeping $\Gamma$ fixed as a set.
If $P$ lies inside $\Gamma$, it can be done by moving $P$ to the center of $\Gamma$.
The existence of corresponding projective transformations follows from the idea in Exercise~\ref{ex:cone}.

The line $p$ is called the \index{polar}\emph{polar} of the point $P$ with respect to~$\Gamma$.

The point $P$ is called the \index{pole}\emph{pole} of the line $p$ with respect to~$\Gamma$.

\begin{thm}{Exercise}\label{ex:revert}
Reverse the described construction.
That is, given a circle $\Gamma$ and a line $p$ that is not tangent to $\Gamma$, construct a point $P$ such that the described construction for $P$ and $\Gamma$ produces the line $p$.
\end{thm}

\begin{thm}{Exercise}\label{ex:tangent ruler}
Let $p$ be a polar line of a point $P$ with respect to the circle~$\Gamma$.
Suppose $p$ intersects $\Gamma$ at points $V$ and~$W$.
Show that the lines $(PV)$ and $(PW)$ are tangent to~$\Gamma$.

Provide a ruler-only construction of the tangent lines to the given circle $\Gamma$ thru the given point $P\notin\Gamma$.
\end{thm}

\begin{thm}{Exercise}\label{ex:concentric-circ}
Assume two concentric circles $\Gamma$ and $\Gamma'$ are given.
Construct the common center of $\Gamma$ and $\Gamma'$ using only a ruler.
\end{thm}

\begin{thm}{Exercise}\label{ex:proj-perp}
Assume a line $\ell$ and a circle $\Gamma$ with its center $O$ are given.
Suppose $O\notin \ell$.
Construct a perpendicular from $O$ on $\ell$ using only a ruler.
\end{thm}

\chapter{Spherical insights}
\label{chap:sphere}

Spherical geometry studies the surface of a unit sphere; it has applications in cartography, navigation, and astronomy.

It is closely related to both Euclidean and hyperbolic geometry.
Since a sphere is easier to visualize, studying spherical geometry can aid in understanding hyperbolic geometry.
We will discuss several such aspects, but we do not intend to provide a complete treatment of spherical geometry.

\section{The sphere}

Recall that Euclidean space is the set \( \mathbb{R}^3 \) of all triples \( (x,y,z) \) of real numbers,
where the distance between two points
\( A = (x_A, y_A, z_A) \) and \( B = (x_B, y_B, z_B) \)
is given by the formula:
\[
AB \df \sqrt{(x_A - x_B)^2 + (y_A - y_B)^2 + (z_A - z_B)^2}.
\]

Planes in this space are defined as the set of solutions to the equation
\[
a \cdot x + b \cdot y + c \cdot z + d = 0
\]
for some constants \( a \), \( b \), \( c \), and \( d \),
where at least one of \( a \), \( b \), or \( c \) is nonzero.
Every such plane is isometric to the Euclidean plane.

A \index{sphere}\emph{sphere} with center \( O \) and radius \( r \) is the set of all points in space that are at a distance \( r \) from \( O \).

Let \(\Sigma\) be a unit sphere with center \( O \).
A \index{great circle}\emph{great circle} of \(\Sigma\) is the intersection of a sphere with a plane passing through \( O \).
In spherical geometry, great circles serve as the analog of lines, even though they do not satisfy Definition~\ref{def:line}.
A great circle divides the sphere into two \index{half-sphere}\emph{half-spheres}, analogous to half-planes in Euclidean geometry.

Let \( A \) and \( B \) be two distinct, non-antipodal points on the unit sphere.
Then there is a unique great circle passing through \( A \) and \( B \), which we denote by \( (AB)_s \).
The shorter of the two arcs of \( (AB)_s \) with endpoints \( A \) and \( B \) will be denoted by \( [AB]_s \).
The \index{distance!spherical}\index{spherical distance}\emph{spherical distance} between \( A \) and \( B \) (denoted briefly as \( AB_s \)) is defined as the absolute value of the angle measure \( \measuredangle AOB \), measured in the plane containing \( O \), \( A \), and \( B \).

Note that each great circle is a circle of radius \( \tfrac{\pi}{2} \) on $\Sigma$.

\section{Area of spherical triangle}

Fix a unit sphere \( \Sigma \).
Similar to the Euclidean plane, we can define a spherical triangle, denoted \( \triangle_s ABC \), along with its sides \( [AB]_s \), \( [BC]_s \), and \( [CA]_s \).
We will always assume that the triangle has no antipodal vertices.
The angle measures
\[
\alpha = |\measuredangle_s CAB|,\quad
\beta = |\measuredangle_s ABC|,\quad
\text{and} \quad
\gamma = |\measuredangle_s BCA|
\]
are defined as the angles between the tangent half-lines to the corresponding arcs.
We will consider only the absolute values of angle measures, altho their signs could also be defined.

A solid spherical triangle \( \solidtriangle_s ABC \) is defined as the set that includes the sides of \( \triangle_s ABC \) along with the rigin inside;
that is, the intersection of three half-spheres:
the first bounded by \( (AB)_s \) and containing \( C \),
the second bounded by \( (BC)_s \) and containing \( A \),
and the third bounded by \( (CA)_s \) and containing \( B \).
This definition closely follows Section~\ref{Solid triangles}.

The concept of area on a sphere will be briefly discussed in Section \ref{Neutral planes and spheres},
but an intuitive understanding of area measurement should be sufficient for the proof of the following lemma.

{

\begin{wrapfigure}{r}{46 mm}
\vskip-4mm
\centering
\includegraphics{mppics/pic-251}
\vskip0mm
\end{wrapfigure}

\begin{thm}{Lemma}\label{lem:area-spher-triangle}
For every solid spherical triangle $\Delta$, we have
\[\area\Delta = \alpha + \beta + \gamma - \pi,\]
where $\alpha$, $\beta$, and $\gamma$ are the absolute angle measures of $\Delta$.
\end{thm}

\parit{Informal proof.}
We will use that
\[\area\Sigma = 4\cdot\pi \eqlbl{eq:area(S2)}\]
without proof.

}

Note that the area of a spherical slice \( S_\alpha \) between two meridians meeting at an angle \( \alpha \) is proportional to \( \alpha \).
Since \( S_\pi \) is a hemisphere, from \ref{eq:area(S2)}, we get \( \area S_\pi = 2\cdot\pi \).
Therefore, the proportionality coefficient is 2; that is,
\[\area S_\alpha = 2\cdot \alpha
\eqlbl{eq:area(Sa)}\]
for any \( \alpha \).

Extending the sides of \( \Delta \), we obtain six slices: two \( S_\alpha \), two \( S_\beta \), and two \( S_\gamma \).
They cover most of the sphere once,
but the triangle \( \Delta \) and its centrally symmetric copy are each covered three times.
It follows that
\begin{align*}
2\cdot \area S_\alpha + 2\cdot \area S_\beta + 2\cdot \area S_\gamma = \area\Sigma + 4\cdot\area\Delta.
\end{align*}
It remains to apply \ref{eq:area(S2)} and \ref{eq:area(Sa)}.
\qeds

Recall that the defect of a triangle with absolute angle measures $\alpha$, $\beta$, and $\gamma$ is defined as $\pi - \alpha - \beta - \gamma$.
According to Section~\ref{The defect}, \textit{the defect of every triangle in a neutral plane is nonnegative}.

\begin{thm}{Exercise}\label{ex:defect-sphere}
Show that the defect of every spherical triangle is nonpositive.
\end{thm}

\section{The central projection}

Any map from a sphere (or its part) to a plane distorts the surface in some way.
In this section we describe the so-called central projection; it maps great circles to straight lines on the plane.
This property makes it convenient in navigation;
while the spherical shapes are distorted, the shortest path on the plane corresponds to a shortest path on the sphere.

The central projection is analogous to the projective model of the hyperbolic plane discussed in Chapter~\ref{chap:klein}.

\medskip

Let $\Sigma$ be the unit sphere centered at the origin~$O$.
Suppose that $\Pi^+$ denotes the plane defined by the equation $z=1$.
This plane is parallel to the $(x,y)$-plane and passes thru 
the north pole $N\z=(0,0,1)$ of~$\Sigma$.

The northern half-sphere of $\Sigma$
is the subset of points $(x,y,z)\z\in \Sigma$ such that $z>0$; it will be denoted by~$\Sigma^+$.

{

\begin{wrapfigure}{o}{40mm}
\centering
\vskip-7mm
\includegraphics{mppics/pic-254}
\end{wrapfigure}

Given a point $P\in \Sigma^+$, consider the half-line $[OP)$. 
Suppose that $P'$ denotes the intersection of $[OP)$ and~$\Pi^+$.
If $P=(x,y,z)$, then $P'=(\tfrac xz,\tfrac yz,1)$.
It follows that $P\leftrightarrow P'$ is a bijection between $\Sigma^+$ and~$\Pi^+$.

}

The described bijection $\Sigma^+\leftrightarrow \Pi^+$ is called the \index{central projection}\emph{central projection} of 
the half-phere~$\Sigma^+$.

The central projection sends the intersections of the great circles with $\Sigma^+$ to the lines in~$\Pi^+$.
This follows since the great circles are intersections of $\Sigma$ with planes passing thru the origin,
and the lines in $\Pi^+$ are the intersection of $\Pi^+$ with these planes.

The following exercise 
is analogous to Exercise~\ref{ex:h-median}.

\begin{thm}{Exercise}\label{ex:s-medians}
Let $\triangle_sABC$ be a nondegenerate spherical triangle.
Assume that the plane $\Pi^+$ is parallel to the plane passing thru $A$, $B$, and~$C$.
Let $A'$, $B'$, and $C'$ denote the central projections of $A$, $B$, and~$C$.
\begin{enumerate}[(a)]
\item\label{ex:s-medians:a} Show that the midpoints of $[A'B']$, $[B'C']$, and $[C'A']$
are central projections of the midpoints of $[AB]_s$, $[BC]_s$, and $[CA]_s$ respectively.
\item\label{ex:s-medians:b} Use part (\ref{ex:s-medians:a}) to show that the medians of a spherical triangle intersect at one point.
\end{enumerate}

\end{thm}

\begin{thm}{Exercise}\label{ex:s-altitudes}
Let $P\leftrightarrow P'$ be the central projection described above
and $N$ be the north pole; so, $N'=N$.
Show that $|\measuredangle_s NPQ|=\tfrac\pi2$ if and only if $|\measuredangle N'P'Q'|=\tfrac\pi2$.
\end{thm}

\section{The stereographic projection}

In this section, we define stereographic projection;
it is closely related to the conformal model of the hyperbolic plane, which is discussed in Chapter~\ref{chap:poincare}.

{

\begin{wrapfigure}{r}{48mm}
\vskip-6mm
\centering
\includegraphics{mppics/pic-252}
\caption*{The plane thru $P$, $O$, and~$S$.}
\end{wrapfigure}

Consider the unit sphere $\Sigma$ 
centered at the origin $(0,0,0)$.
This sphere can be described by the equation $x^2+y^2\z+z^2\z=1$. 

Suppose that $\Pi$ denotes the $(x,y)$-plane;
it is defined by the equation $z \z= 0$.
Clearly, $\Pi$
runs thru the center of~$\Sigma$.

Let $N = (0, 0, 1)$ and $S\z=(0, 0, -1)$ denote the ``north'' and ``south'' poles of $\Sigma$;
these are the points on the sphere that have extremal distances to~$\Pi$.
Suppose that $\Omega$ denotes the ``equator'' of $\Sigma$;
it is the intersection $\Sigma\cap\Pi$.

}

For a given point $P\ne S$ on $\Sigma$,
consider the line $(SP)$ in the space. 
This line intersects $\Pi$ at exactly one point, denoted by~$P'$. 
Set $S'=\infty$.

The map $\xi_s\: P\mapsto P'$ is called the \index{stereographic projection}\emph{stereographic projection from $\Sigma$ to $\Pi$ with respect to the south pole}.
The inverse of this map $\xi^{-1}_s\: P'\z\mapsto P$ is called the {}\emph{stereographic projection from $\Pi$ to $\Sigma$ with respect to the south pole}.

In the same way, one can define the
{}\emph{stereographic projections $\xi_n$ and $\xi^{-1}_n$ with respect to the north pole}~$N$.

Note that $P=P'$ if and only if $P\in\Omega$.

Exercise~\ref{ex:stereographic-inversion} below states that the stereographic projection preserves 
the angles between arcs;
more precisely, \textit{the absolute value of the angle measure} between arcs on the sphere.
This is a useful property in cartography;
the curves on the sphere meet at the same angles as their stereographic projections.

\section[The inversion]{Inversion across a sphere}

The inversion across a sphere is defined the same way as the inversion across a circle.

Formally, let $\Sigma$ be the sphere with center $O$ and radius~$r$.
The \index{inversion!inversion across a sphere}\emph{inversion} across $\Sigma$ of a point $P$ is the point $P'\in[OP)$ such that
$$OP\cdot OP'=r^2.$$
In this case, the sphere $\Sigma$  will be called the 
\index{inversion!sphere of inversion}\emph{sphere of inversion},
and its center is called the \index{inversion!center of inversion}\emph{center of inversion}.

We also add $\infty$ to the space and assume that the center of inversion is mapped to $\infty$ and the other way around. 
The space $\mathbb{R}^3$ with the point $\infty$ will be called \index{inversive!space}\emph{inversive space}.

The inversion of space shares many properties with the inversion of the plane.
Most importantly, analogs of theorems \ref{lem:inverse-4-angle}, \ref{thm:inverse-cline}, \ref{thm:angle-inversion} can be summarized as follows:

\begin{thm}{Theorem}\label{thm:inversion-3d}
The inversion across the sphere has the following properties:
\begin{enumerate}[(a)]
\item\label{thm:inversion-3d:a} Inversion maps a sphere or a plane into a sphere or a plane.
\item\label{thm:inversion-3d:b} Inversion maps a circline into a circline. 
\item\label{thm:inversion-3d:cross-ratio} Inversion preserves the cross-ratio;
that is, if $A'$, $B'$, $C'$, and $D'$ are the inverses of the points $A$, $B$, $C$, and $D$ respectively,
then
$$\frac{AB\cdot CD}{BC\cdot DA}= \frac{A'B'\cdot C'D'}{B'C'\cdot D'A'}.$$
\item Inversion maps arcs into arcs.
\item\label{thm:inversion-3d:angle}
Inversion preserves the absolute value of the angle
measure between tangent half-lines to the arcs.
\end{enumerate}
\end{thm}

In the following exercises, we apply the theorem above to tie the stereographic projection with the inversion. 
We assume that $\Sigma$, $\Pi$, $\Omega$, $O$, $S$, $N$ and $\xi_s$ are as in the previous section.

{

\begin{wrapfigure}{r}{32mm}
\vskip-0mm
\centering
\includegraphics{mppics/pic-253}
\end{wrapfigure}

\begin{thm}{Exercise}\label{ex:stereographic-inversion}
Show that the stereographic projections 
$\xi_s\: \Sigma\to\Pi$ and $\xi^{-1}_s\: \Pi\to\Sigma$
are the restrictions to $\Sigma$ and $\Pi$ respectively, of the inversion across the sphere $\Upsilon$ with the center $S$ and radius $\sqrt{2}$.

Conclude that the stereographic projection preserves 
the angles between arcs;
more precisely, the absolute value of the angle measure between arcs on the sphere.
\end{thm}

\begin{thm}{Exercise}\label{ex:great-circ}
Show that the stereographic projection $\xi_s\:\Sigma\to\Pi$
sends the great circles to plane circlines that intersect $\Omega$ at opposite points.
\end{thm}

}

\begin{thm}{Exercise}\label{ex:conform-sphere}
Fix a point $P\in \Pi$  and let $Q$ be another point in~$\Pi$.
Let $P'$ and $Q'$ denote their stereographic projections to~$\Sigma$.
Set $x=PQ$ and $y=P'Q'_s$.
Show that
$$\lim_{x\to 0}\, \frac{y}{x}=\frac{2}{1+OP^2}.$$
\end{thm}

The last exercise is analogous to Lemma~\ref{lem:conformal}.

The proof of \ref{thm:inversion-3d} resembles the corresponding proofs in plane geometry.
Let us give couple of a hints to the reader who wants to reconstruct its proof.
To prove \ref{thm:inversion-3d}\textit{\ref{thm:inversion-3d:a}}, one needs the following lemma;
its proof is left to the reader.

\begin{thm}{Lemma}
Let $\Sigma$ be a subset of the Euclidean space
that contains at least two points.
Fix a point $O$ in the space.

Then $\Sigma$ is 
a sphere 
if and only if
for every plane $\Pi$ passing thru $O$,
the intersection $\Pi\cap \Sigma$ is either an empty set,
a one-point set, or a circle.
\end{thm}  

The following observation reduces part~\textit{(\ref{thm:inversion-3d:b})} to part~\textit{(\ref{thm:inversion-3d:a})}.

\begin{thm}{Observation}
Any circle in the space is an intersection of two spheres.
\end{thm}

{

\begin{wrapfigure}{o}{25mm}
\centering
\includegraphics{mppics/pic-250}
\end{wrapfigure}

Let us define a \index{circular cone}\emph{circular cone} as a set formed by line segments from a fixed point, called the \index{tip of cone}\emph{tip} of the cone, to all the points on a fixed circle, called the \index{base!of cone}\emph{base} of the cone;
we always assume that the base does not lie in the same plane as the tip.
We say that the cone is \index{right!circular cone}\emph{right}
if the center of the base circle is the footpoint of the tip on the base plane;
otherwise, we call it \index{oblique circular cone}\emph{oblique}.

}

\begin{thm}{Exercise}\label{ex:cone}
Let $K$ be an oblique circular cone.
Show that there is a plane $\Pi$ that is not parallel to the base plane of $K$ such that the intersection $\Pi\cap K$ is a circle.
\end{thm}

\section{The Pythagorean theorem}

Here is an analog of the Pythagorean theorems (\ref{thm:pyth} and \ref{thm:pyth-h-poincare}) in spherical geometry.

\begin{thm}{Spherical Pythagorean theorem}\label{thm:s-pyth}\index{Pythagorean theorem}
Let $\triangle_sABC$ be a spherical triangle with a right angle at~$C$.
Set $a=BC_s$, $b=CA_s$, and $c=AB_s$.
Then
$$\cos c=\cos a\cdot\cos b.$$

\end{thm}

In the proof, we will use the notion of the scalar product which we are about to discuss.

Let $v_A\z=(x_A,y_A,z_A)$ and $v_B=(x_B,y_B,z_B)$ denote the position vectors of points $A$ and $B$.
The \index{scalar product}\emph{scalar product} of the two vectors $v_A$ and $v_B$ in $\mathbb{R}^3$
is defined as 
$$\langle v_A,v_B\rangle
\df
x_A\cdot x_B+y_A\cdot y_B+z_A\cdot z_B.\eqlbl{eq:scal-def}$$

Assume both vectors $v_A$ and $v_B$ are nonzero;
suppose that $\phi$ denotes the angle measure between them.
Then the scalar product can be expressed the following way:
$$\langle v_A,v_B\rangle=|v_A|\cdot|v_B|\cdot\cos\phi,
\eqlbl{eq:scal-angle}$$
where 
\begin{align*}
|v_A|&=\sqrt{x_A^2+y_A^2+z_A^2},
&
|v_B|&=\sqrt{x_B^2+y_B^2+z_B^2}.
\end{align*}

Now, assume that the points $A$ and $B$ 
lie on the unit sphere $\Sigma$ in $\mathbb{R}^3$ centered at the origin.
In this case, $|v_A|=|v_B|=1$.
By \ref{eq:scal-angle} we get that
$$\cos AB_s=\langle v_A,v_B\rangle.
\eqlbl{eq:scalar-s-dist}$$

\parit{Proof of the spherical Pythagorean theorem.}
Since the angle at $C$ is right,
we can choose the coordinates in $\mathbb{R}^3$ so that 
$v_C\z=(0,0,1)$, $v_A$ lies in the $(x,y)$-plane, so $v_A\z=(x_A,0,z_A)$,
and $v_B$ lies in the $(y,z)$-plane, so $v_B=(0,y_B,z_B)$.

{

\begin{wrapfigure}{r}{40mm}
\vskip-6mm
\centering
\includegraphics{mppics/pic-248}
\end{wrapfigure}

Applying, \ref{eq:scalar-s-dist},
we get that
\begin{align*}
z_A&=\langle v_C,v_A\rangle
=\cos b,
\\
z_B&=\langle v_C,v_B\rangle
=\cos a.
\end{align*}

Applying, \ref{eq:scal-def} and \ref{eq:scalar-s-dist}, we get that
\begin{align*}
\cos c &=\langle v_A,v_B\rangle=
\\
&=x_A\cdot 0+0\cdot y_B+z_A\cdot z_B=
\\
&=\cos b\cdot\cos a.
\end{align*}
\qedsf

}

\begin{thm}{Exercise}\label{ex:2(pi/4)=pi/3}
Show that 
if $\triangle_sABC$ is a spherical triangle with a right angle at $C$,
and $AC_s=BC_s=\tfrac\pi4$, then $AB_s=\tfrac\pi3$.
\end{thm}

\section{Imaginary distance}

Recall that  
\[
\cosh x=\frac {e^{x}+e^{-x}}2
\qquad\text{and}\qquad
\cos x=\frac {e^{i\cdot x}+e^{-i\cdot x}}2.
\]
The first formula is the definition of the hyperbolic cosine (see \ref{sec:hyp-trig});
the second one is called Euler's formula (see \ref{sec:Euler's formula}).
It follows that 
\[\cosh x=\cos (i\cdot x)
\qquad\text{and}\qquad
\cos x=\cosh (i\cdot x).\]

Let us compare the formulas in hyperbolic and spherical Pythagorean theorems (see \ref{thm:s-pyth} and \ref{thm:pyth-h-poincare}):
\begin{align*}
\cosh c&=\cosh a\cdot \cosh b,
&
\cos c&=\cos a\cdot \cos b.
\end{align*}
Observe that if we change $a$, $b$, and $c$ to $i\cdot a$, $i\cdot b$, and $i\cdot c$,
then the first formula transforms into the second one, and the other way around.

This is not a coincidence;
the same holds for all analytic formulas ---
changing every distance $d$ to $i\cdot d$ transforms a valid spherical formula into a valid hyperbolic formula;
the angle measures need no change.
This magic substitution was found by Franz Taurinus \cite{taurinus}; we are not going to prove it.

\begin{thm}{Advanced exercise}\label{ex:taurinus}
Consider a spherical triangle $ABC$; set
\begin{align*}
a&=BC_s,& b&=CA_s,& c&=AB_s,
\\
\alpha&=|\measuredangle_sCAB|,& \beta&=|\measuredangle_sABC|,& \gamma&=|\measuredangle_sBCA|.
\end{align*}
Use the Taurinus substitution to write the hyperbolic analogs to the following formulas in spherical geometry:

\begin{enumerate}[(a)]
\item The spherical cosine rule:
\[\cos c=\cos a \cdot \cos b+\sin a\cdot \sin b\cdot \cos\gamma.\]
\item The dual spherical cosine rule:
\[\cos \gamma=-\cos \alpha \cdot \cos \beta+\sin \alpha\cdot \sin \beta \cdot \cos c.\]
\item
The spherical sine rule:
\[\frac{\sin \alpha}{\sin a}=\frac{\sin \beta}{\sin b}=\frac{\sin \gamma}{\sin c}.\]
\end{enumerate}
 
\end{thm}

\chapter{Projective model}\label{chap:klein}

The \textit{projective model} is another model of the hyperbolic plane discovered by Eugenio Beltrami; it is often called the {}\emph{Klein model}.
The projective and conformal models are saying exactly the same thing but in two different languages.
Some problems  are easier in projective model, and others in the conformal model.
Therefore, it is worth knowing both. 

\section{A special bijection on the h-plane}
\label{sec:special-bijection}

Consider the conformal disc model.
Assume its absolute is a unit circle $\Omega$ centered at~$O$.
Choose a coordinate system $(x,y)$ on the plane with the origin at $O$, 
so the circle $\Omega$ is described by the equation $x^2+y^2=1$.

\label{pic:stereographic_projection-klein}
\begin{wrapfigure}{o}{48mm}
\centering
\vskip-3mm
\includegraphics{mppics/pic-256}
\caption*{Plane thru $P$, $O$, and $S$.}
\vskip-3mm
\end{wrapfigure}

Let us think that our plane is the coordinate $(x,y)$-plane in the Euclidean space; denote it by $\Pi$.
Let $\Sigma$ be the unit sphere centered at $O$;
it is described by the equation 
$x^2+y^2+z^2=1$.
Set $S\z=(0,0,-1)$ and $N\z=(0,0,1)$; 
these are the south and north poles of~$\Sigma$.

Consider the stereographic projection $\Pi\to\Sigma$ from $S$;
given a point $P\in\Pi$ denote by $P'$ its image in $\Sigma$.
Note that the  h-plane is mapped to the {}\emph{northern half-sphere},
that is, to the set of points $(x,y,z)$ in $\Sigma$ described by the inequality~$z>0$.

For a point $P'\in \Sigma$ consider its footpoint $\hat P$
on $\Pi$;
this is the closest point to~$P'$.

The composition of these two maps $P\z\leftrightarrow P'\z\leftrightarrow\hat P$
gives a bijection from the h-plane to itself.
Furthermore, $P=\hat P$
 if and only if  $P\in \Omega$ or $P=O$.


\begin{thm}{Lemma}\label{lem:P-hat-chord}
Let $P\leftrightarrow\hat P$ be the bijection of the h-plane described above.
Suppose $(PQ)_h$ is an h-line with the ideal points $A$ and~$B$.
Then $\hat P,\hat Q\in[AB]$.

Moreover, 
$$\frac{A\hat Q\cdot B\hat P}{\hat QB\cdot \hat PA}
=
\left(\frac{AQ\cdot BP}{QB\cdot PA}\right)^2.
\eqlbl{eq:lem:P-hat-chord}$$
In particular, if $A,P,Q,B$ appear in the same order, then
$$PQ_h=\tfrac12\cdot\ln\frac{A\hat Q\cdot B\hat P}{\hat QB\cdot \hat PA}.$$
\end{thm}

\parit{Proof.}
Consider the stereographic projection $\Pi\to \Sigma$ from the south pole~$S$.
It fixes $A$ and $B$;
denote by $P'$ and $Q'$ the images of $P$ and~$Q$;

According to Theorem~\ref{thm:inversion-3d}\textit{\ref{thm:inversion-3d:cross-ratio}},
$$\frac{AQ\cdot BP}{QB\cdot PA}=\frac{AQ'\cdot BP'}{Q'B\cdot P'A}.\eqlbl{eq:(AB;PQ)=(AB;P'Q')}$$

By Theorem~\ref{thm:inversion-3d}\textit{\ref{thm:inversion-3d:angle}}, 
each circline in $\Pi$ that is perpendicular to $\Omega$ 
is mapped to a circle in $\Sigma$ that is still perpendicular to~$\Omega$.
It follows that the stereographic projection sends $(PQ)_h$ to the intersection of the northern half-sphere of $\Sigma$ with a plane perpendicular to~$\Pi$.

Denote this plane by $\Lambda$;
it contains the points $A$, $B$, $P'$, $\hat P$ and the circle $\Gamma=\Sigma\cap\Lambda$.
(It also contains $Q'$ and $\hat Q$ but we will not use these points for a while.)

{

\begin{wrapfigure}{r}{30mm}
\vskip-0mm
\centering
\includegraphics{mppics/pic-258}
\caption*{Plane $\Lambda$.}
\end{wrapfigure}

Note that 
\begin{itemize}
\item 
$A,B,P'\in\Gamma$,
\item $[AB]$ is a diameter of $\Gamma$,
\item $(AB)=\Pi\cap\Lambda$,
\item $\hat P\in [AB]$
\item $(P'\hat P)\perp (AB)$.
\end{itemize}

Since $[AB]$ is the diameter of $\Gamma$, 
by Corollary~\ref{cor:right-angle-diameter},
the angle $AP'B$ is right. 
Hence $\triangle A\hat PP'\z\sim \triangle AP'B\z\sim \triangle P'\hat PB$.
In particular
$$\frac{AP'}{BP'}=\frac{A\hat P}{P'\hat P}=\frac{P'\hat P}{B\hat P}.$$
Therefore
$$\frac{A\hat P}{B\hat P}=\left(\frac{AP'}{BP'}\right)^2
\quad\text{and similarly}\quad
\frac{A\hat Q}{B\hat Q}=\left(\frac{AQ'}{BQ'}\right)^2.
\eqlbl{eq:AP/BP}$$
Finally, 
\ref{eq:(AB;PQ)=(AB;P'Q')}+\ref{eq:AP/BP} imply \ref{eq:lem:P-hat-chord}.

}

The last statement follows from \ref{eq:lem:P-hat-chord} and the definition of the h-distance.
Indeed,
\begin{align*}
PQ_h&\df\ln\frac{A Q\cdot B P}{QB\cdot PA}=
\\
&=\ln\left(\frac{A \hat Q\cdot B \hat P}{\hat QB\cdot \hat PA}\right)^{\frac12}=
\\
&=\tfrac12\cdot\ln\frac{A \hat Q\cdot B \hat P}{\hat QB\cdot \hat PA}.
\end{align*}
\qedsf

{

\begin{wrapfigure}[8]{r}{36mm}
\centering
\vskip-8mm
\includegraphics{mppics/pic-260}
\end{wrapfigure}

\begin{thm}{Exercise}\label{ex:hex}
Let $\Gamma_1$, $\Gamma_2$, and $\Gamma_3$ 
be three circles perpendicular to a circle~$\Omega$.
Let $[A_1B_1]$, $[A_2B_2]$, and $[A_3B_3]$ denote
the common chords of $\Omega$ and $\Gamma_1$, $\Gamma_2$, $\Gamma_3$ respectively.
Show that the chords $[A_1B_1]$, $[A_2B_2]$, and $[A_3B_3]$ intersect at one point inside $\Omega$ if and only if $\Gamma_1$, $\Gamma_2$, and $\Gamma_3$ intersect at two points.
\end{thm}

\begin{thm}{Exercise}\label{ex:P<->hatP}
Show that $2\cdot OP_h=O\hat P_h$, where $O$ is the center of absolute, and $P\leftrightarrow \hat P$ is the bijection of the h-plane described above.
\end{thm}

}

\section{The projective model}
\label{sec:proj-model}

The following picture illustrates the bijection $P\mapsto \hat P$ of the h-plane described in the previous section --- 
\begin{figure}[!ht]
\centering
\includegraphics{mppics/pic-262}
\end{figure}
if you take the picture on the left and apply the map $P\z\mapsto \hat P$,
you get the picture on the right.
The pictures are {}\emph{conformal} and \index{projective!model}\emph{projective models} of the hyperbolic plane respectively.
The map $P\mapsto \hat P$ serves as a \textit{translation} between the models.

In the projective model, things look different;
some become simpler, and others become more complex.

\parbf{Lines.}
In the projective model, the h-lines are represented as chords of the absolute, more precisely, chords without their endpoints.

This observation can be used to transfer statements about lines and points from the Euclidean plane to the h-plane.
As an example let us state a version of Pappus' theorem for h-plane.

\begin{thm}{Hyperbolic Pappus' theorem}\label{thm:pappus-h}\index{Pappus' theorem}
Assume that two triples of h-points $A$, $B$, $C$,
and $A'$, $B'$, $C'$ in the h-plane are h-collinear.
Suppose that the h-points $X$, $Y$, and $Z$ are uniquelly defined by
\begin{align*}
X&\in(BC')_h\cap(B'C)_h,
&
Y&\in(CA')_h\cap(C'A)_h,
&
Z&\in(AB')_h\cap(A'B)_h.
\end{align*}
Then the points $X$, $Y$, $Z$ are h-collinear.
\end{thm}

In the projective model, this statement follows immediately from the original Pappus' theorem \ref{thm:pappus}.
The same can be done for Desargues' theorem \ref{thm:desargues}.
The same argument shows that the construction of a tangent line with a ruler only described in Exercise~\ref{ex:tangent ruler} works in the h-plane as well.
These statements are not at all evident in the conformal model.

\parbf{Circles and equidistants.}
In the projective model, the h-circles and equidistants are ellipses and their open arcs.
This follows since the stereographic projection sends circles on the plane to circles on the unit sphere and the footpoint projection of the circle back to the plane is an ellipse.
(One may define an \index{ellipse}\emph{ellipse} as a footpoint projection of a circle.)

\parbf{Distance.}
Consider a pair of h-points $P$ and $Q$.
Let $A$ and $B$ be the ideal points of the h-line in the projective model;
that is, $A$ and $B$ are the intersections of the Euclidean line $(PQ)$ with the absolute.

\begin{wrapfigure}{o}{36mm}
\vskip-2mm
\centering
\includegraphics{mppics/pic-264}
\vskip-2mm
\end{wrapfigure}

Then by Lemma~\ref{lem:P-hat-chord},
$$PQ_h=\tfrac12\cdot\ln\frac{AQ\cdot BP}{QB\cdot PA},\eqlbl{eq:proj-h-dist}$$
assuming the points $A, P, Q, B$ appear on the line in the same order.

\parbf{Angles.}
The hyperbolic angle measure in the projective model is very different from Euclidean, making it hard to determine from the picture. 
\label{klein-angles}
For example, all the intersecting h-lines on the picture 
are perpendicular.

To find the angle measure,
one may apply a motion of the h-plane that moves 
the vertex of the angle to the center of the absolute;
once done, the hyperbolic and Euclidean angles have the same measure.
In particular, if $O$ is the center of the absolute, then 
$$\measuredangle_hAOB=\measuredangle AOB.$$

\begin{thm}{Observation}\label{obs:h-p-perp}
If $O$ is the center of the absolute, then 
\[(OA)\z\perp (AB)\qquad\Longleftrightarrow\qquad(OA_h)\perp (AB)_h.\]

\end{thm}

\parit{Proof.}
The Euclidean reflection across $(OA)$ induces the h-reflection across $(OA)_h$.
Let $B'$ be the reflection of $B$.
Observe that $(OA)\perp (AB)$ $\Longleftrightarrow$ $A\in(BB')$ $\Longleftrightarrow$  $A\in(BB')_h$ $\Longleftrightarrow$  $(OA_h)\perp (AB)_h$.
\qeds

\parbf{Motions.}
The motions of the h-plane in the conformal and projective models are relevant to inversive transformations and projective transformations in the same way.
Namely: 
\begin{itemize}
\item Inversive transformations that preserve the h-plane describe motions of the h-plane in the conformal model.
\item Projective transformations that preserve the h-plane describe motions of the h-plane in the projective model.\footnote{The idea described in the solution of Exercise~\ref{ex:cone}.
The sketch of proof of \ref{thm:circle-center-proj} can be used to construct many projective transformations of this type.}
\end{itemize}

The following exercise is a hyperbolic analog of Exercise~\ref{ex:s-medians}. 

\begin{thm}{Exercise}\label{ex:h-median}
Let $P$ and $Q$ be points in the h-plane that lie at the same distance from the center of the absolute.
Observe that in the projective model, the h-midpoint of $[PQ]_h$ coincides with the Euclidean midpoint of $[PQ]_h$.

Conclude that if an h-triangle is inscribed in an h-circle, then its medians meet at one point.

Recall that an h-triangle might also be inscribed in a horocycle or an equidistant.
Think about how to prove the statement in this case.
\end{thm}

\begin{thm}{Exercise}\label{ex:h-altitudes}
Show that the altitudes of a hyperbolic triangle either intersect at one point or are pairwise disjoint. 
\end{thm}

{

\begin{wrapfigure}[4]{r}{34mm}
\vskip-11mm
\centering
\includegraphics{mppics/pic-266}
\end{wrapfigure}

\begin{thm}{Exercise}\label{ex:klein-perp}
Let $\ell$ and $m$ be  h-lines in the projective model.
Let $s$ and $t$ denote the Euclidean lines tangent to the absolute
at the ideal points of $\ell$. 
Show that
if the lines $s$, $t$, and the extension of $m$ meet at one point, then $\ell\perp m$. 
\end{thm}

}

\begin{thm}{Exercise}\label{ex:klein-for-angle-parallelism}
Use the projective model to derive the formula for the angle of parallelism (Proposition~\ref{prop:angle-parallelism}). 
\end{thm}

\begin{thm}{Exercise}\label{ex:klein-inradius}
Use the projective model to find an inradius of an ideal triangle.
\end{thm}

\begin{thm}{Advanced exercise}\label{ex:pyth-h-proj}
Prove the hyperbolic Pythagorean theorem (\ref{thm:pyth-h-poincare}) using the following idea.
\end{thm}

Recall that the hyperbolic Pythagorean theorem\index{Pythagorean theorem} (\ref{thm:pyth-h-poincare}) states that
\[\cosh c=\cosh a\cdot\cosh b,
\eqlbl{eq:hyp-pyth-proj}\]
where $a\z=BC_h$, $b=CA_h$, and $c=AB_h$ and
$\triangle_hACB$ is an h-triangle with a right angle at~$C$.

\begin{wrapfigure}{r}{34mm}
\vskip-6mm
\centering
\includegraphics{mppics/pic-268}
\end{wrapfigure}

\parit{Idea.}
We can assume that $A$ is the center of the absolute.
By \ref{obs:h-p-perp} the Euclidean triangle $ABC$ is right.
Set 
$s=BC$, $t =CA$, $u\z= AB$.
According to the Euclidean Pythagorean theorem (\ref{thm:pyth}), we have
$$u^2=s^2+t^2.\eqlbl{eq:hyp-proj}$$
It remains to express $a$, $b$, and $c$ using $s$, $u$, and $t$ and show that \ref{eq:hyp-proj} implies~\ref{eq:hyp-pyth-proj}.

\section{Bolyai's construction}

Assume we need to construct a line thru $P$ that is asymptotically parallel to the given line $\ell$ in the h-plane.

If $A$ and $B$ are ideal points of $\ell$ in the projective model, 
then we could simply draw the Euclidean line $(PA)$.
However, the ideal points do not lie in the h-plane; therefore there is no way to use them in the construction.

In the following construction we assume that you know a ruler-and-compass construction of the perpendicular line; see Exercise~\ref{ex:construction-perpendicular}.

\begin{thm}{Bolyai's construction}
\begin{enumerate}
\item Drop a perpendicular from $P$ to~$\ell$; denote it by~$m$.
Let $Q$ be the footpoint of $P$ on~$\ell$.
\item Erect a perpendicular from $P$ to~$m$; denote it by~$n$.
\item Mark a point $R$ on $\ell$ distinct from $Q$.
\item Drop a perpendicular from $R$ to~$n$; denote it by~$k$. 
\item Draw the circle $\Gamma$ with center $P$ and the radius $QR$. 
Mark a point of intersection of $\Gamma$ with~$k$; denote it by $T$.
\item The line $(PT)_h$ is asymptotically parallel to~$\ell$.
\end{enumerate}
\end{thm}

\begin{thm}{Exercise}\label{ex:Boyai-in-Euclid}
Explain what happens if one performs the Bolyai construction in the Euclidean plane.
\end{thm}

The following proposition implies that Bolyai's construction works.

\begin{thm}{Proposition}\label{prop:boyai}
Assume $P$, $Q$, $R$, $S$, $T$ are points in the h-plane
such that 
$S\in (RT)_h$,
$(PQ)_h\perp (QR)_h$,
$(PS)_h\perp(PQ)_h$,
$(RT)_h\perp (PS)_h$ and 
$(PT)_h$ and $(QR)_h$ are asymptotically parallel.
Then $QR_h=PT_h$.
\end{thm}

\parit{Proof.}
We will use the projective model.
Without loss of generality, we may assume that $P$ is the center of the absolute.
By~\ref{obs:h-p-perp},
the corresponding Euclidean lines are also perpendicular;
that is, $(PQ)\perp (QR)$, $(PS)\perp(PQ)$, and $(RT)\z\perp (PS)$.

Let $A$ be the common ideal point of $(QR)_h$ and $(PT)_h$.
Let $B$ and $C$ denote the remaining ideal points of $(QR)_h$ and $(PT)_h$ respectively.

The Euclidean lines $(PQ)$, $(TR)$, and $(CB)$ are parallel.
Therefore,  
\[\triangle AQP\sim \triangle ART \sim\triangle ABC.\]

\begin{wrapfigure}{o}{55mm}
\vskip-0mm
\centering
\includegraphics{mppics/pic-270}
\vskip2mm
\end{wrapfigure}

In particular,
\[\frac{AC}{AB}=\frac{AT}{AR}=\frac{AP}{AQ}.\]
It follows that
\[\frac{AT}{AR}=\frac{AP}{AQ}=\frac{BR}{CT}=\frac{BQ}{CP}.\]
In particular, 
\[\frac{AT\cdot CP}{TC\cdot PA}=\frac{AR\cdot BQ}{RB\cdot QA}.\]
Applying the formula for h-distance \ref{eq:proj-h-dist}, we get that $QR_h=PT_h$.
\qeds

\begin{wrapfigure}{r}{35mm}
\vskip-6mm
\centering
\includegraphics{mppics/pic-271}
\vskip2mm
\end{wrapfigure}

\begin{thm}{Advanced exercise}\label{ex:common-perp}
Look at the picture and recover the construction of a common perpendicular $n$ to the h-lines $\ell$ and $m$.

Prove that it works.
\end{thm}

\chapter{Complex coordinates}\label{chap:complex}

In this chapter, we give an interpretation of inversive geometry using complex coordinates.
The results of this chapter lead to a deeper understanding of both concepts.

\section{Complex numbers}

Informally,
a complex number has form 
$$z=x+i\cdot y,
\eqlbl{eq:z=x+iy}$$ 
where $x$ and $y$ 
are real numbers and $i^2=-1$. 

The set of complex numbers 
will be denoted by~$\mathbb{C}$.
If $x$, $y$, and $z$ are as in \ref{eq:z=x+iy}, 
then $x$ is called the \index{real!part}\emph{real part} and $y$ is the \index{imaginary!part}\emph{imaginary part} of the complex number~$z$.
Briefly, it is written as 
\[x=\Re z
\quad
\text{and}
\quad 
y=\Im z.\]

On a more formal level, the expression $x + i\cdot y$ 
is just a convenient way 
to write the pair $(x,y)$.
Furthermore, addition and multiplication of the pairs are defined as follows:
\[
\begin{aligned}
(x_1+i\cdot y_1) + (x_2+i\cdot y_2) 
&\df (x_1+x_2) + i\cdot(y_1+y_2);
\\
(x_1+i\cdot y_1)\cdot(x_2+i\cdot y_2) 
&\df 
(x_1\cdot x_2-y_1\cdot y_2) + i\cdot(x_1\cdot y_2+y_1\cdot x_2).
\end{aligned}
\eqlbl{eq:comlex+x}
\] 

\section{Complex coordinates}

Remember that we can view the Euclidean plane as the set of all pairs of real numbers $(x,y)$ equipped with the metric 
$$AB=\sqrt{(x_A-x_B)^2+(y_A-y_B)^2},$$
where $A=(x_A,y_A)$ and $B=(x_B,y_B)$.

We can pack the coordinates $(x,y)$ of a point in a single complex number $z=x+i\cdot y$.
This creates a direct link between points in the Euclidean plane and~$\mathbb{C}$ (the set of complex numbers).
Given a point $Z=(x,y)$, 
the complex number $z=x+ i\cdot y$ is referred to as the \index{complex coordinate}\emph{complex coordinate} of~$Z$.

If $O$, $E$, and $I$ are points in the plane 
with complex coordinates $0$, $1$, and $i$ respectively, then $\measuredangle EOI=\pm\tfrac\pi2$.
We will assume that $\measuredangle EOI=\tfrac\pi2$;
if not, one needs to adjust the direction of the $y$-coordinate. 

\section{Conjugation and absolute value}
\label{sec:complex-conjugation}

Let $z=x+i\cdot y$; 
that is, $z$ is a complex number with real part $x$ and imaginary part $y$.
If $y=0$, we say that the complex number $z$ is \index{real!complex number}\emph{real}; if $x=0$, we say that $z$ is \index{imaginary!number}\emph{imaginary}.
The set of points with real (imaginary) complex coordinates forms a line in the plane, known as the real (or imaginary) line. 
The real line will be denoted as $\mathbb{R}$.

The complex number
\[\bar z\df x-i\cdot y\] is called the \index{complex conjugate}\emph{complex conjugate} of $z=x+i\cdot y$.
Let $Z$ and $\bar Z$ be the points in the plane with the complex coordinates $z$ and $\bar z$ respectively.
The point $\bar Z$ is the reflection of $Z$ across the real line.

It is straightforward to check that
$$\begin{aligned}
x&=\Re z=\frac{z+\bar z}2,
&
y&=\Im z=\frac{z-\bar z}{i\cdot2},
\quad\text{and}
&
x^2+y^2&=z\cdot\bar z.
\end{aligned}\eqlbl{eq:conj-1}$$

The last formula in \ref{eq:conj-1} allows us to express the quotient $\tfrac{w}{z}$ of two complex numbers $w$ and $z=x+i\cdot y$:
$$\frac{w}{z}=\tfrac{1}{z\cdot\bar z}\cdot w\cdot\bar z=\tfrac{1}{x^2+y^2}\cdot w\cdot\bar z.$$

\label{page:cojugation=authomorphism}
Note that
\begin{align*}
\overline {z+ w}&=\bar z+\bar w,
&
\overline {z- w}&=\bar z-\bar w,
&
\overline {z\cdot w}&=\bar z\cdot\bar w,
&
\overline {z/w}&=\bar z/\bar w.
\end{align*}
That is, complex conjugation
\textit{respects}
all the arithmetic operations.

The value 
\begin{align*}
|z|&\df\sqrt{x^2+y^2}=\sqrt{(x+i\cdot y)\cdot(x-i\cdot y)}
=
\sqrt{z\cdot\bar z}
\end{align*}
is called the
\index{absolute value}\emph{absolute value} of $z$.
If $|z|=1$, then $z$ is called a \index{unit complex number}\emph{unit complex number}.

\begin{thm}{Exercise}\label{ex:|zw|}
Show that $|v\cdot w|=|v|\cdot |w|$ for all $v,w\in\mathbb{C}$.
\end{thm}

Suppose that $Z$ and $W$ are points with complex coordinates $z$ and $w$.
Note that
$$ZW=|z-w|.\eqlbl{eq:C-dist}$$
The triangle inequality for points with complex coordinates $0$, $v$, and $v+w$ implies that
\[|v+w|\le |v|+|w|\]
for all $v,w\in\mathbb{C}$;
this inequality is also called the \index{triangle!inequality}\emph{triangle inequality}.

\begin{thm}{Exercise}\label{ex:ptolemy}
Use the identity 
\[u\cdot (v-w)+v\cdot (w-u)+w\cdot(u-v)=0\]
for $u,v,w\in\mathbb{C}$ and the triangle inequality
to prove Ptolemy's inequality (\ref{ptolemy-inq}).
\end{thm}

\section{Euler's formula}\label{sec:Euler's formula}

{

\begin{wrapfigure}{r}{36mm}
\vskip-15mm
\centering
\includegraphics{mppics/pic-272}
\end{wrapfigure}

Let $\alpha$ be a real number.
The following identity 
$$e^{i\cdot\alpha}=\cos\alpha+i\cdot\sin\alpha
\eqlbl{eq:euler}$$
is called \index{Euler's formula}\emph{Euler's formula}.
In particular, $e^{i\cdot\pi}\z=-1$ and $e^{i\cdot\frac\pi2}\z=i$.

Geometrically, Euler’s formula can be understood as follows: Imagine
$O$ and $E$ as points with complex coordinates $0$ and $1$ respectively.
Assume 
\[OZ=1\quad \text{and}\quad \measuredangle EOZ\z\equiv \alpha,\]
then $e^{i\cdot\alpha}$ is the complex coordinate of $Z$.
In particular, the complex coordinate of every point on the unit circle centered at~$O$
can be uniquely expressed as $e^{i\cdot\alpha}$ for some $\alpha\in(-\pi,\pi]$.

}

\parbf{Why is this true?}
The proof of Euler's identity depends on how you define the exponential function.

If you have never applied the exponential function to an imaginary number, 
you may take the right-hand side in \ref{eq:euler} 
as the definition of $e^{i\cdot\alpha}$.
In this case, nothing has to be proved,
but it is better to check that $e^{i\cdot\alpha}$ satisfies familiar identities.
Primarily,
$$e^{i\cdot \alpha}\cdot e^{i\cdot \beta}= e^{i\cdot(\alpha+\beta)}.$$
The latter can be proved using \ref{eq:comlex+x} and the following well-known trigonometric formulas:
\begin{align*}
\cos(\alpha+\beta)&=\cos\alpha\cdot\cos\beta-\sin\alpha\cdot\sin\beta,
\\
\sin(\alpha+\beta)&=\sin\alpha\cdot\cos\beta+\cos\alpha\cdot\sin\beta.
\end{align*}

If you are familiar with the power series for sine, cosine, and the exponential function, then the following might convince you that identity \ref{eq:euler} holds:
\begin{align*}
 e^{i\cdot \alpha } &{}= 1 + i\cdot \alpha  + \frac{(i\cdot \alpha )^2}{2!} + \frac{(i\cdot \alpha  )^3}{3!} + \frac{(i\cdot \alpha )^4}{4!} + \frac{(i\cdot  \alpha )^5}{5!} +  \cdots =
 \\
&= 1 + i\cdot \alpha  - \frac{\alpha ^2}{2!} - i\cdot\frac{ \alpha ^3}{3!} + \frac{\alpha ^4}{4!} + i\cdot\frac{ \alpha ^5}{5!} -  \cdots =
\\
&= \left( 1 - \frac{\alpha ^2}{2!} + \frac{\alpha ^4}{4!}  - \cdots \right) +  i\cdot\left( \alpha  - \frac{\alpha ^3}{3!} + \frac{\alpha ^5}{5!} -  \cdots \right) =
\\
&= \cos \alpha  +  i\cdot\sin \alpha.
\end{align*}

\section{Argument and polar coordinates}

As before, we assume that $O$ and $E$ are the points with complex coordinates $0$ and $1$ respectively.

Let $Z$ be a point distinct from $O$.
Set $\rho=OZ$ and $\theta=\measuredangle EOZ$.
The pair $(\rho,\theta)$ is called the \index{polar!coordinates}\emph{polar coordinates} of~$Z$.

If $z$ is the complex coordinate of $Z$, then $\rho=|z|$. 
The value $\theta$ is called the \index{argument}\emph{argument} of $z$
(briefly, $\theta=\arg z$).
In this case, 
$$z=\rho\cdot e^{i\cdot\theta}=\rho\cdot(\cos\theta+i\cdot\sin\theta).\qquad\qquad\qquad\qquad $$

\begin{wrapfigure}[4]{o}{30mm}
\vskip-14mm
\centering
\includegraphics{mppics/pic-274}
\end{wrapfigure}

Note that 
\begin{align*}
\arg (z\cdot w)&\equiv \arg z+\arg w,
\shortintertext{and}
\arg \tfrac z w&\equiv \arg z-\arg w
\end{align*}
if $z\ne0 $ and $w\ne0$.
In particular, if $Z$, $V$, and $W$ are points with complex coordinates $z$, $v$, and $w$ respectively, then
$$
\begin{aligned}
\measuredangle VZW
&=\arg\left(\frac{w-z}{v-z}\right)\equiv
\\
&\equiv \arg(w-z)-\arg(v-z)
\end{aligned}
\eqlbl{eq:angle-arg}$$
if $\measuredangle VZW$ is defined.

\begin{thm}{Exercise}\label{ex:3-sum-C}
Use the formula \ref{eq:angle-arg} to show that  
$$\measuredangle ZVW+\measuredangle VWZ+\measuredangle WZV\equiv \pi$$
for every triangle $ZVW$ in the Euclidean plane.
\end{thm}

\begin{thm}{Exercise}\label{ex:C-sim}
Suppose that points $O$, $E$, $V$, $W$, and $Z$ have complex coordinates $0$, $1$, $v$, $w$, and $z=v\cdot w$ respectively.
Show that 
\[\triangle OEV\sim \triangle OWZ.\]

\end{thm}

The following theorem is a translation of Corollary~\ref{cor:inscribed-quadrangle} into complex-number language.

\begin{thm}{Theorem}\label{thm:inscribed-quadrangle-C}
Let $\square UVWZ$ be a quadrangle and $u$, $v$, $w$, and $z$ be the complex coordinates of its vertices. 
Then $\square UVWZ$ is inscribed 
if and only if the number
$$\frac{(v-u)\cdot(z-w)}{(v-w)\cdot(z-u)}$$ 
is real.
\end{thm}

The value $\frac{(v-u)\cdot(w-z)}{(v-w)\cdot(z-u)}$ is called the 
\index{cross-ratio!complex cross-ratio}\emph{complex cross-ratio} of $u$, $w$, $v$, and $z$; 
it will be denoted by \index{64@$(u,v;w,z)$}$(u,w;v,z)$.

\begin{thm}{Exercise}\label{ex:real-cross-ratio}
Observe that the complex number $z\ne 0$ is real if and only if $\arg z=0$ or $\pi$;
in other words, $2\cdot\arg z\equiv 0$.

Use this observation to show that Theorem~\ref{thm:inscribed-quadrangle-C}
is indeed a reformulation of Corollary~\ref{cor:inscribed-quadrangle}.
\end{thm}

\section{The method of complex coordinates}

The following problem illustrates the method of complex coordinates.

\begin{thm}{Problem}\label{prob:2right-tringles}
Let $\triangle OPV$ and $\triangle OQW$ be isosceles right triangles such that 
\[\measuredangle VPO=\measuredangle OQW=\tfrac\pi2\] 
and $M$ be the midpoint of $[VW]$.
Assume $P$, $Q$, and $M$ are distinct points.
Show that  $\triangle PMQ$ is an isosceles right triangle.
\end{thm}

\parit{Solution.}
Choose complex coordinates so that $O$ is the origin.
Denote by $v, w, p, q, m$ the complex coordinates of the remaining points respectively.

Since $\triangle OPV$ and $\triangle OQW$ are isosceles and $\measuredangle VPO=\measuredangle OQW=\tfrac\pi2$,
\ref{eq:C-dist} and \ref{eq:angle-arg} imply that
\begin{align*}
v-p&=i\cdot p,
&
q-w&=i\cdot q.
\end{align*}

\begin{wrapfigure}{o}{37mm}
\centering
\includegraphics{mppics/pic-276}
\end{wrapfigure}

Therefore
\begin{align*}
m
&=\tfrac12\cdot(v+w)=
\\
&=\tfrac{1+i}2\cdot p+\tfrac{1-i}2\cdot q.
\end{align*}

By straightforward computations, we get that
\[p-m=i\cdot (q-m).\]
In particular, $|p-m|=|q-m|$ and  $\arg\frac{p-m}{q-m}=\tfrac\pi2$;
that is, $PM=QM$ and $\measuredangle QMP =\tfrac\pi2$.  
\qeds

{

\begin{wrapfigure}{r}{36mm}
\vskip-4mm
\centering
\includegraphics{mppics/pic-278}
\end{wrapfigure}

\begin{thm}{Exercise}\label{ex:3-squares}
Consider three squares with common sides as shown in the picture.
Use the method of complex coordinates to show that 
\[\measuredangle EOA+\measuredangle EOB+\measuredangle EOC=\pm\tfrac\pi2.\]

\end{thm}

}

\begin{thm}{Exercise}\label{ex:6-circles}
Check the following identity with six complex cross-ratios:
\[(u,w;v,z)\cdot(u',w';v',z')=\frac{(v,w';v',w)\cdot(z,u';z',u)}{(u,v';u',v)\cdot(w,z';w',z)}.\]
Use it together with Theorem~\ref{thm:inscribed-quadrangle-C} to prove that if
$\square UVWZ$, $\square UVV'U'$, $\square VWW'V'$, $\square WZZ'W'$, and $\square ZUU'Z'$
are inscribed, then  $\square U'V'W'Z'$ is inscribed as well.

\end{thm}

\begin{minipage}{.47\textwidth}
\centering
\includegraphics{mppics/pic-280}
\end{minipage}
\hfill
\begin{minipage}{.47\textwidth}
\centering
\includegraphics{mppics/pic-282}
\end{minipage}

\medskip

\begin{thm}{Exercise}\label{ex:4-sim}
Suppose that points $U$, $V$, and $W$ lie on one side of line $(AB)$ and 
$\triangle UAB\sim \triangle BVA \sim \triangle ABW$.
Denote by $a$, $b$, $u$, $v$, and $w$ the complex coordinates of $A$, $B$, $U$, $V$, and $W$ respectively.
\begin{enumerate}[(a)]
 \item Show that $\tfrac{u-a}{b-a}=\tfrac{b-v}{a-v}=\tfrac{a-b}{w-b}=\tfrac{u-v}{w-v}$.
 \item Conclude that $\triangle UAB\sim \triangle BVA \sim \triangle ABW\sim \triangle UVW$.
\end{enumerate}
 
\end{thm}

\section{Fractional linear transformations}

\begin{thm}{Exercise}\label{ex:movie}
Watch the video ``Möbius transformations revealed'' by Douglas Arnold and Jonathan Rogness.
(It is available on \href{http://youtu.be/JX3VmDgiFnY}{YouTube}.)
\end{thm}

The complex plane $\mathbb{C}$ extended by one ideal number $\infty$ 
is called the \index{extended complex plane}\emph{extended complex plane}.
It is denoted by $\hat{\mathbb{C}}$, so $\hat{\mathbb{C}}=\mathbb{C}\cup\{\infty\}$

A \index{fractional linear transformation}\emph{fractional linear transformation} or \index{M\"obius transformation}\emph{M\"obius transformation} of  $\hat{\mathbb{C}}$ is a function of one complex variable $z$
that can be written as
$$f(z) = \frac{a\cdot z + b}{c\cdot z + d},$$
where the coefficients $a$, $b$, $c$, $d$ are complex numbers satisfying $a\cdot d \z- b\cdot c \not= 0$.
(If $a\cdot d - b\cdot c = 0$ the function defined above is a constant; it is not a fractional linear transformation.) 

In case $c\not=0$, we assume that
$$f(-d/c) = \infty
\quad
\text{and}
\quad
f(\infty) = a/c;$$
and if $c=0$ we assume
$f(\infty) = \infty$.

\section{Elementary transformations}

The following three types of fractional linear transformations are called \index{elementary transformation}\emph{elementary}:
\begin{enumerate}
\item $z\mapsto z+w,$
\item $z\mapsto w\cdot z$ for $w\ne0,$
\item $z\mapsto \frac1z.$
\end{enumerate}
 
\parbf{The geometric interpretations.}
Suppose that $O$ denotes the point with the complex coordinate~$0$.

The first map $z\mapsto z+w,$ corresponds to the so-called 
\index{parallel!translation}\emph{parallel translation} 
of the Euclidean plane; its geometric meaning should be evident.

The second map is called the \index{rotational homothety}\emph{rotational homothety} with the center at~$O$.
That is, the point $O$ maps to itself,
and any other point $Z$ maps to a point $Z'$ such that $OZ'=|w|\cdot OZ$ and $\measuredangle ZOZ'=\arg w$.

The third map can be described as a composition of the inversion across the unit circle centered at $O$ and the reflection across $\mathbb{R}$ 
(the composition can be taken in any order).

Indeed, $\arg z\equiv -\arg \tfrac1z$.
Therefore, 
$$\arg z=\arg (1/\bar z);$$
that is, if the points $Z$ and $Z'$ have complex coordinates $z$ and $1/\bar z$,
then $Z'\in[OZ)$.
Clearly, $OZ=|z|$ and $OZ'=|1/\bar z|=\tfrac{1}{|z|}$.
Therefore, $Z'$ is the inverse of $Z$ across the unit circle centered at~$O$.
Finally, $\tfrac1z\z=\overline{(1/\bar z)}$ is the complex coordinate of
the reflection of $Z'$ across $\mathbb{R}$.

\begin{thm}{Proposition}\label{prop:mob-comp}
A map $f\:\hat{\mathbb{C}}\to\hat{\mathbb{C}}$ is a fractional linear transformation if and only if $f$ can be expressed as a composition of elementary transformations. 
\end{thm}

\parit{Proof; the ``only if'' part.}
Fix a fractional linear transformation
\begin{align*}
f(z) &= \frac{a\cdot z + b}{c\cdot z + d}.
\shortintertext{Assume $c\ne 0$. Then}
f(z) &= \frac ac-\frac{a\cdot d-b\cdot c}{c\cdot(c\cdot z + d)} =
\\
&= \frac ac-\frac{a\cdot d-b\cdot c}{c^2}\cdot \frac1{z + \frac dc}.
\end{align*}
That is, 
$$f(z)=f_4\circ f_3\circ f_2\circ f_1 (z),
\eqlbl{eq:moebius-compose}$$
where $f_1$, $f_2$, $f_3$, and $f_4$ are the following elementary transformations:
\begin{align*}
f_1(z)&= z+\tfrac dc,
&
f_2(z)&= \tfrac1z,
\\
f_3(z)&= - \tfrac{a\cdot d-b\cdot c}{c^2} \cdot z,
&
f_4(z)&= z+\tfrac ac.
\end{align*}

If $c=0$, then
\[f(z) = \frac{a\cdot z + b}{ d}.\]
In this case, $f(z)=f_2\circ f_1 (z)$,
where 
\begin{align*}
f_1(z)&= \tfrac ad\cdot z,
&
f_2(z)= z+\tfrac bd.
\end{align*}

\parit{``If'' part.}
We need to show that by composing elementary transformations,
we can only get fractional linear transformations.
It is sufficient to check that the composition of a fractional linear transformation
$$f(z) = \frac{a\cdot z + b}{c\cdot z + d}.$$
with every elementary transformation $z\mapsto z+w$, $z\mapsto w\cdot z$, or $z\mapsto \tfrac1z$ is a fractional linear transformation.

The latter is done by means of direct calculations.
\begin{align*}
\frac{a\cdot (z+w) + b}{c\cdot (z+w) + d}
&=
\frac{a\cdot z + (b+a\cdot w)}{c\cdot z + (d+c\cdot w)},
\\
\frac{a\cdot (w\cdot z) + b}{c\cdot (w\cdot z) + d}
&=
\frac{(a\cdot w)\cdot z + b}{(c\cdot w)\cdot z + d},
\\
\frac{a\cdot \frac1z + b}{c\cdot \frac1z + d}
&=
\frac{b\cdot z + a}{d\cdot z + c}.
\end{align*}
\qedsf

\begin{thm}{Corollary}\label{cor:cline-Moeb}
The fractional linear transformations map circlines to circlines.
\end{thm}

\parit{Proof.}
By Proposition~\ref{prop:mob-comp},
it is sufficient to check that each elementary transformation sends a circline to a circline.

For the first and second elementary transformations, the latter is evident.

As noted above,
the map $z\mapsto\tfrac1z$ is a composition of inversion and reflection.
By Theorem~\ref{thm:inverse}, the inversion sends a circline to a circline.
Hence the result.
\qeds

\begin{thm}{Exercise}\label{ex:inverse-Mob}
Show that if $f$ is a fractional linear transformation, then so is its inverse $f^{-1}$.
\end{thm}

\begin{thm}{Exercise}\label{ex:3-point-Mob}
Given distinct values $z_0,z_1,z_\infty\in \hat{\mathbb{C}}$,
construct a fractional linear transformation $f$ such that 
$f(z_0)=0$,
$f(z_1)=1$,
and 
$f(z_\infty)\z=\infty$.
Show that such a transformation is unique.
\end{thm}

\begin{thm}{Exercise}\label{ex:inversion-Mob}
Show that every inversion is a composition of complex conjugation and a fractional linear transformation.

Use \ref{ex:inversions-inversive} to conclude that every inversive transformation is either a fractional linear transformation or a complex conjugate of a fractional linear transformation.
\end{thm}

\section{The complex cross-ratio}

Let $u$, $v$, $w$, and $z$ be four distinct complex numbers.
Recall that 
the complex number
$$
\frac{(u-w)\cdot(v-z)}{(v-w)\cdot(u-z)}$$
is called the \index{cross-ratio!complex cross-ratio}\emph{complex cross-ratio} of $u$, $v$, $w$, and $z$; 
it is denoted by $(u,v;w,z)$.

If one of the numbers $u$, $v$, $w$, $z$ is $\infty$, 
then the complex cross-ratio has to be defined by taking the appropriate limit; in other words, we assume that $\frac\infty\infty=1$.
For example,
$$(u, v; w, \infty)=\frac{(u-w)}{(v-w)}.$$

Assume that $U$, $V$, $W$, and  $Z$ are the points with complex coordinates  
$u$, $v$, $w$, and $z$ respectively.
Note that 
\begin{align*}
\frac{UW\cdot VZ}{VW\cdot UZ}&=|(u,v;w,z)|,
\\
\measuredangle WUZ +\measuredangle ZVW&=\arg\frac{u-w}{u-z}+\arg\frac{v-z}{v-w}\equiv \arg(u,v;w,z).
\end{align*}
These equations give the following reformulation of Theorem~\ref{lem:inverse-4-angle}.

\begin{thm}{Theorem}\label{lem:inverse-4-angle-C}
Let $UWVZ$ and $U'W'V'Z'$  be two quadrangles 
such that the points $U'$, $W'$, $V'$, and $Z'$ are inverses of $U$, $W$, $V$, and $Z$ respectively.
Assume $u$, $w$, $v$, $z$, $u'$, $w'$, $v'$, and $z'$ are the complex coordinates of $U$, $W$, $V$, $Z$, $U'$, $W'$, $V'$, and $Z'$ respectively.

Then 
$$(u',v';w',z')=\overline{(u,v;w,z)}.$$

\end{thm}

The following exercise is a generalization of the theorem above.
It has a short solution using Proposition~\ref{prop:mob-comp}.

\begin{thm}{Exercise}\label{ex:C-cross-ratio}
Show that the complex cross-ratio is {}\emph{invariant} under fractional linear transformations. 

That is, if a fractional linear transformation maps four distinct complex numbers $u, v, w, z$ to complex numbers $u', v', w', z'$ respectively, then
$$
(u',v';w',z')
=
(u,v;w,z).
$$

\end{thm}

\section{The Schwarz--Pick theorem}

The following theorem shows 
that the metric in the conformal disc model naturally appears in other branches of mathematics.
We do not give its proof, but it can be found in any textbook on geometric complex analysis.

Let $\mathbb{D}$ be the \index{disc}\emph{open unit disc} in the complex plane centered at~$0$;
that is, a complex number $z$
belongs to $\mathbb{D}$ if and only if $|z|<1$.

Let us use the disc $\mathbb{D}$ as an h-plane in the conformal disc model;
the h-distance between $z, w\in\mathbb{D}$ will be denoted by $d_h(z,w)$;
that is,
\[d_h(z,w)\df ZW_h,\]
where $Z$ and $W$ are h-points with complex coordinates $z$ and $w$ respectively.

A function $f\:\mathbb{D}\to \mathbb{C}$ is called \index{holomorphic function}\emph{holomorphic} if for every $z\in \mathbb{D}$
there is a complex number $s$ such that
\[f(z+w)=f(z)+s\cdot w+o(|w|).\]
In other words, $f$ is {}\emph{complex-differentiable}
at any $z\in\mathbb{D}$.
The complex number $s$ is called the {}\emph{derivative} of $f$ at $z$ (briefly $s=f'(z)$).

\begin{thm}{Schwarz--Pick theorem}\index{Schwarz--Pick theorem}
Assume $f\: \mathbb{D}\to \mathbb{D}$ is a holomorphic function.
Then 
\[d_h(f(z),f(w))\le d_h(z,w)\]
for all $z,w\in \mathbb{D}$.

If the equality holds for one pair of distinct numbers $z,w\in \mathbb{D}$, then it holds for each pair.
In this case, $f$ is a fractional linear transformation as well as a motion of the h-plane.
\end{thm}

\begin{thm}{Exercise}\label{ex:schwarz-moebius}
Show that if a fractional linear transformation $f$ appears in the equality case of the Schwarz--Pick theorem, then it can be written as 
\[f(z)=\frac{v\cdot z+\bar w}{w\cdot z+\bar v}.\]
where $v$ and $w$ are complex constants such that $|v|>|w|$.
\end{thm}

Recall that hyperbolic tangent $\tanh$ is defined in Section~\ref{sec:hyp-trig}.

\begin{thm}{Exercise}\label{ex:schwarz-tanh}
Show that 
\[\tanh [\tfrac12\cdot d_h(z,w)]=\left|\frac{z-w}{1-z\cdot\bar w}\right|.\]
Conclude that the inequality in Schwarz--Pick theorem can be rewritten as
\[\left|\frac{z'-w'}{1-z'\cdot\bar w'}\right|\le\left|\frac{z-w}{1-z\cdot\bar w}\right|,\]
where
$z'=f(z)$ and $w'=f(w)$.
\end{thm}

\begin{thm}{Exercise}\label{ex:schwarz}
Show that the Schwarz lemma stated below 
follows from the Schwarz--Pick theorem.
\end{thm}

\begin{thm}{Schwarz lemma}
Let $f\: \mathbb{D}\to \mathbb{D}$ be a holomorphic function
and $f(0)=0$.
Then 
$|f(z)|\le |z|$
for every $z\in \mathbb{D}$.

Moreover, if equality holds for some $z\ne 0$, then there is a unit complex number $u$ 
such that 
$f(z)=u\cdot z$
for every $z\in\mathbb{D}$.
\end{thm}

\chapter{Geometric constructions}
\label{chap:car}

Geometric constructions were introduced at the end of Chapter~\ref{chap:perp} and have been used thruout the text.
They possess significant pedagogical value as an introduction to mathematical proofs.

In this chapter, we briefly explore the theory behind geometric constructions.
Let us also mention the computer game ``Euclidea'' \cite{euclidea};
it is an addictive way to become an expert in the field.

\section{Classical problems}

The solutions to the following two problems are quite nontrivial.

\begin{thm}{Problem of Brahmagupta} 
Construct an inscribed quadrangle with given sides.
\end{thm}

{

\begin{wrapfigure}[7]{r}{22mm}
\vskip-7mm
\centering
\includegraphics{mppics/pic-284}
\end{wrapfigure}

Several solutions to the following problem are presented in \cite{hadamard}.
 
\begin{thm}{Problem of Apollonius} Construct a circle that is tangent to three given circles.
\end{thm}

}

\section{Impossible constructions}

Impossible construction problems cannot be solved in principle; 
that is, the required ruler-and-compass construction does not exist.
The following problems have been known for about two thousand years;
their impossibility was proved only in the $19^\text{th}$ century.
The method used in the proofs is indicated in the next section.

\begin{thm*}{Doubling the cube}
Given a cube, construct a side of a cube with twice the volume.
\end{thm*}

Given a segment of length $a$,
construct a segment of length~$\sqrt[3]{2}\cdot a$.

\begin{thm*}{Squaring the circle}
Construct a square with the same area as a given disc.
\end{thm*}

If $r$ is the radius of the given circle, we need to construct a segment of length~$\sqrt{\pi}\cdot r$. 

\begin{thm*}{Angle trisection} 
Divide a given angle into three equal angles.
\end{thm*}

Moreover, ruler-and-compass constructions cannot trisect an angle of size~$\tfrac\pi3$. 
The existence of such a construction would imply the constructability of a regular 9-gon which is prohibited by the following famous result:

A \index{regular $n$-gon}\emph{regular $n$-gon} inscribed in a circle with center $O$ is a sequence of points $A_1\dots A_n$ on the circle such that 
\[\measuredangle A_nOA_1=\measuredangle A_1OA_2=\dots=\measuredangle A_{n-1}OA_n=\pm\tfrac2n\cdot \pi.\]
The points $A_1,\dots, A_n$ are \index{vertex!of a regular $n$-gon}\emph{vertices},
the segments $[A_1A_2], \z\dots, [A_nA_1]$ are \index{side!of a regular $n$-gon}\emph{sides} 
and the remaining segments $[A_iA_j]$ are \index{diagonal!of a regular $n$-gon}\emph{diagonals} of the $n$-gon.

The construction of a regular $n$-gon, therefore, is reduced to the construction of an angle with size $\tfrac2n\cdot \pi$.

\begin{thm}{Gauss--Wantzel theorem}\index{Gauss--Wantzel theorem}
A regular $n$-gon can be constructed with a ruler and a compass 
if and only if 
$n$ is the product of a power of $2$ and distinct Fermat primes.
\end{thm}

A \index{Fermat prime}\emph{Fermat prime} is a prime number of the form $2^k+1$ for an integer~$k$.
Only five Fermat primes are known  today: $3$, $5$, $17$, $257$, and $65537$.
For example, 
\begin{itemize}
\item one can construct a regular 340-gon since $340=2^2\cdot 5\cdot 17$ and $5$, as well as $17$, are Fermat primes;
\item one cannot construct a regular 7-gon since $7$ is not a Fermat prime;
\item one cannot construct a regular 9-gon; 
altho $9=3\cdot 3$ is a product of two Fermat primes, 
these primes are not distinct.
\end{itemize}

\section{Constructible numbers}

Let us give an intuitive definition of ruler-and-compass constructions; a formal definition is subtler than one might think \cite{engeler}.

In classical ruler-and-compass constructions, the initial configuration can be completely described by a finite number of points;
each line is defined by two points on it and each circle is described by its center and a point on it.

Similarly, the result of a construction can be described by a finite collection of points.

We may always assume that the initial configuration has at least two points;
if not, we could add one or two points to the configuration.
Moreover, by scaling the whole plane, we can assume that the first two points in the initial configuration lie at a unit distance from each other.

In this case, we can choose a  coordinate system, 
so that the points $(0,0)$ and $(1,0)$ are among the initial points;
so the initial configuration of $n$ points is described by 
$2\cdot n-4$ numbers --- their coordinates $x_3,y_3,\z\dots,x_n,y_n$. 

\medskip

It turns out that the coordinates of every point constructed with a ruler and compass
can be written using the numbers $x_3,y_3,\z\dots,x_n,y_n$ and the four arithmetic operations ``$+$'', ``$-$'', ``$\cdot$'', ``$/\,$''
along with the square root ``$\sqrt{\phantom{a}}\,$''.

For example, assume we want to find the points $X_1=(x_1,y_1)$ and $X_2\z=(x_2,y_2)$ of the intersections of 
a line passing thru $A=(x_A,y_A)$ and $B\z=(x_B,y_B)$ and
the circle with center $O=(x_O,y_O)$ that passes thru the point $W=(x_W,y_W)$.
Let us write the equations of the circle and the line in the coordinates $(x,y)$:
$$
\left\{
\begin{aligned}
(x-x_O)^2+(y-y_O)^2&=(x_W-x_O)^2+(y_W-y_O)^2,
\\
(x-x_A)\cdot(y_B-y_A)&=(y-y_A)\cdot(x_B-x_A).
\end{aligned}
\right.
$$
The coordinates $(x_1,y_1)$ and $(x_2,y_2)$ of the points $X_1$ and $X_2$ solve this system.
Expressing $y$ from the second equation and substituting the result in the first one, gives us a quadratic equation in $x$, 
which can be solved using ``$+$'', ``$-$'', ``$\cdot$'', ``$/\,$''
and  ``$\sqrt{\phantom{a}}\,$'' only.

The same can be performed for the intersection of two circles. 
The intersection of two lines is even simpler; 
it is described as a solution of two linear equations and can be expressed using only four arithmetic operations;
the square root ``$\sqrt{\phantom{a}}\,$'' is not needed.

\medskip

On the other hand, it is easy to construct segments of lengths $a+b$ and $a-b$ from two given segments of lengths $a>b$.

\begin{wrapfigure}{r}{37mm}
\vskip-4mm
\centering
\includegraphics{mppics/pic-286}
\end{wrapfigure}

To perform ``$\cdot$'', ``$/\,$''
and ``$\sqrt{\phantom{a}}\,$'' consider the following picture:
let $[AB]$ be a diameter of a circle; 
fix a point $C$ on the circle and let $D$ be the footpoint of $C$ on~$[AB]$.
By Corollary~\ref{cor:right-angle-diameter}, the angle $ACB$ is right.
Therefore 
$$\triangle ABC\sim\triangle ACD\sim \triangle CBD.$$
It follows that $AD\cdot BD=CD^2$. 

Using this picture, one should guess a solution to the following exercise.

\begin{thm}{Exercise}\label{ex:a2/b}
Given two line segments with lengths $a$ and $b$, provide a ruler-and-compass give construction of segments with lengths $\tfrac {a^2}b$ and $\sqrt{a\cdot b}$.
\end{thm}

Taking $1$ for $a$ or $b$ above, we can produce 
$\sqrt a$, $a^2$, $\tfrac1b$.
Combining these constructions we can produce
$a\cdot b=(\sqrt{a\cdot b})^2$,
$\tfrac ab=a\cdot\tfrac 1b$.
In other words, we produced a \index{ruler-and-compass calculator}\emph{ruler-and-compass calculator} that can do ``$+$'', ``$-$'', ``$\cdot$'', ``$/\,$'', and the square root ``$\sqrt{\phantom{a}}\,$''.

The discussion above sketches a proof of the following theorem:
 
\begin{thm}{Theorem}\label{thm:constructible-numbers}
Suppose $A_1=(0,0)$, $A_2=(1,0)$, $A_3\z=(x_3,y_3),\z\dots,A_n\z=(x_n,y_n)$ is 
an initial configuration of points.
A point $X\z=(x,y)$ can be obtained by a ruler-and-compass construction
if and only if both coordinates $x$ and $y$ can be expressed from the integer numbers and $x_3$, $y_3$, $x_4$, $y_4,\z\dots,x_n,y_n$ using the arithmetic operations ``$+$'', ``$-$'', ``$\cdot$'', ``$/\,$'', and ``$\sqrt{\phantom{a}}\,$''.
\end{thm}

The numbers that can be expressed from the given numbers using the arithmetic operations and the square root ``$\sqrt{\phantom{a}}\,$'' are called \index{constructible numbers}\emph{constructible};
if the list of given numbers is not given, then we can only use the integers.

{\sloppy
The theorem above translates every ruler-and-compass construction problem into a purely algebraic language.
Let us give some examples.
\begin{itemize}
\item The impossibility of the doubling-cube problem states that $\sqrt[3]{2}$ is not a constructible number.
That is, $\sqrt[3]{2}$ cannot be expressed thru integers using
``$+$'', ``$-$'', ``$\cdot$'', ``$/\,$'', and ``$\sqrt{\phantom{a}}\,$''.

\item The impossibility of squaring the circle states that 
$\sqrt{\pi}$, or equivalently $\pi$, is not a constructible number.

\item The impossibility of angle trisection states that $\cos\tfrac\alpha3$ is not a constructible number from $\cos\alpha$.

\item The \index{Gauss--Wantzel theorem}\emph{Gauss--Wantzel theorem} says for which integers $n$ the number 
$\cos\tfrac{2\cdot\pi}n$ is constructible.
\end{itemize} 
Some of these statements might look evident, 
but rigorous proofs require some knowledge of abstract algebra (namely, the \textit{field theory})
which is out of the scope of this book. 
In the next section, we discuss similar but simpler examples of impossible constructions with an unusual tool.

}

\begin{thm}{Exercise}\label{ex:5-gon}
\begin{enumerate}[(a)]
 \item\label{ex:5-gon:a} Show that a diagonal of a regular pentagon is $\tfrac{1+\sqrt5}2$ times larger than its side.
 \item\label{ex:5-gon:b} Based on (\ref{ex:5-gon:a}), describe a ruler-and-compass construction of a regular pentagon.
\end{enumerate}
\end{thm}

\section{Set-square constructions}

A \index{set-square}\emph{set-square} (or 45°-{}\emph{set-square}) is a construction tool shown in the picture ---
it can produce a line thru a given point
that makes the angles
$\tfrac\pi2$ or $\pm\tfrac\pi4$ 
to a given line, and it can also can be used as a ruler.

\begin{thm}{Exercise}\label{ex:trisect-set-square}
Trisect a given segment with a set-square.
\end{thm}

The following theorem is an analog of Theorem~\ref{thm:constructible-numbers} for set-square constructions.

\begin{thm}{Theorem}\label{thm:set-square-constructible-numbers}
Suppose $A_1 = (0, 0)$, $A_2 = (1, 0)$, $A_3 = (x_3, y_3), \z\dots, A_n = (x_n, y_n)$ is
an initial configuration of points.
Then a point $X \z= (x, y)$ can be constructed using a set-square if and only if both coordinates $x$ and $y$
can be expressed using only integers and the numbers $x_3$, $y_3$, $x_4$, $y_4$, $\dots$, $x_n$, $y_n$,
along with the arithmetic operations ``$+$'', ``$-$'', ``$\cdot$'', and ``$/$''.
\end{thm}

\begin{wrapfigure}{о}{26mm}
\vskip-2mm
\centering
\includegraphics{mppics/pic-288}
\end{wrapfigure}

The proof of this theorem is close to Theorem~\ref{thm:constructible-numbers}.
(The ``if'' part nearly follows from Exercise~\ref{ex:R-hom}.
The ``only-if'' part is proved by induction on the number of elementary constructions; one needs to write an equation for each line in a set-square construction and verify that an intersection point of such lines satisfies the theorem.)

Unlike Theorem~\ref{thm:constructible-numbers} it can be applied directly to show the impossibility of some constructions with a set-square --- no need to delve into the field theory.

If all the coordinates $x_3,y_3,\dots,x_n,y_n$ are rational numbers, then the theorem above implies that with a set-square, one can only construct points with rational coordinates.
A point with both rational coordinates is called \index{rational point}\emph{rational},
and if at least one of the coordinates is irrational, then the point is called \index{irrational point}\emph{irrational}.

\begin{thm}{Exercise}\label{ex:equilateral triangle}
Show that it is impossible to construct an equilateral triangle with a given base using
the set-square.
\end{thm}

{

\begin{wrapfigure}{r}{26mm}
\vskip-8mm
\centering
\includegraphics{mppics/pic-289}
\end{wrapfigure}

\begin{thm}{Exercise}\label{ex:set-square-bisect}
Show that it is impossible to bisect a given angle with the set-square only.
\end{thm}

\begin{thm}{Advanced exercise}\label{ex:90-60-30}
Consider another tool --- a 30°-set-square that can produce a line thru a given point
that makes the angles
$\tfrac\pi2$, $\pm\tfrac\pi3$, $\pm\tfrac\pi6$
to a given line and can be used as a ruler.

{Show that it is impossible to construct a square with the 30°-set-square.}
\end{thm}

}

\section{Verifications}
\label{sec:verification}

Suppose we need to verify that a given configuration is defined by a certain property. 
Is it possible to do this task by geometric construction with the given tools?
We assume that we can \textit{verify} that two constructed points coincide.

Evidently, if a configuration is constructible, then it is \index{verifiable construction}\emph{verifiable}%
\footnote{Adopting the terminology of computability theory, we may also say that such a construction is \index{decidable construction}\emph{decidable}.} --- simply repeat the construction and check if the result is the same.
Some nonconstructible configurations are verifiable.
For example, it does not pose a problem to verify that the given angle is trisected while it is impossible to trisect a given angle with a ruler and compass.
A regular 7-gon provides another example of that type --- it is easy to verify, while Gauss--Wantzel theorem states that it is impossible to construct with a ruler and compass.

Since we did not prove the impossibility of angle trisection and the Gauss--Wantzel theorem, the following example might be more satisfactory.
It is based on Exercise~\ref{ex:equilateral triangle} which states that it is impossible to construct an equilateral with a set-square only.

\begin{thm}{Exercise}\label{ex:equilateral triangle-verify}
Describe a set-square construction verifying that 
\begin{enumerate}[(a)]
\item\label{ex:verify:triangle} a triangle is equilateral.
\item\label{ex:verify:bisector} a line bisects an angle.
\end{enumerate}
\end{thm}

This observation leads to a source of impossible constructions in a stronger sense --- those that are even not verifiable.

The following example is closely related to Exercise~\ref{ex:circumtool}.
Recall that a \index{circumtool}\emph{circumtool} produces a circle passing thru given three points
or a line if these points lie on one line;
the \index{inversor}\emph{inversor} --- a tool that constructs an inverse of a given point in a given circline.

\begin{thm}{Problem}\label{prob:center-inversor+circumtool}
Show that with a circumtool and inversor,
it is impossible to verify that a given point is the center of a given circle~$\Gamma$.
In particular, it is impossible to construct the center with a circumtool only.
\end{thm}

\parbf{Remark.}
In geometric constructions, we allow choosing \index{free points}\emph{free points}, such as any point on the plane, or a point on a constructed line, or a point that does not lie on a constructed line, or a point on a given line that does not lie on a given circle, and so on.

In principle, when you make such a free choice it is possible to make the right construction by accident.
Nevertheless, we do not accept such a coincidence as a true construction.
The result should not depend on chance.

\parit{Solution.}\label{page:solution-for-ex:circumtool}
Arguing by contradiction, 
assume we have a verifying construction. 

Apply an inversion across a circle perpendicular to $\Gamma$ to the whole construction.
According to Corollary~\ref{cor:perp-inverse-clines},
the circle
$\Gamma$ maps to itself.
Recall that the inversion sends a circline to a circline (\ref{thm:inverse-cline}) and respects inversion (\ref{cor:invese-comp}).
Therefore we get that the whole  construction is mapped to an equivalent construction; 
that is, a construction with a different choice of free points.

According to Exercise~\ref{ex:inv-center not=center-inv}, 
the inversion sends the center of $\Gamma$ to another point.
However, this construction claims that this new point is the center as well --- a contradiction.
\qeds

A similar example of impossible constructions for a ruler and a parallel tool
 is given in Exercise~\ref{ex:affine-perp}.
 
Let us discuss another example of a ruler-only construction.
Note that ruler-only constructions are invariant with respect to the projective transformations. 
In particular, to solve the following exercise, it is sufficient to construct a projective transformation that fixes two points and moves their midpoint.

\begin{thm}{Exercise}\label{ex:midpoint-proj}
Show that there is no ruler-only construction verifying that a given point is a  midpoint of a given segment.
In particular, it is impossible to construct the midpoint only with a ruler.
\end{thm}

The following theorem is a stronger version of the exercise above.

\begin{thm}{Theorem}\label{thm:circle-center-proj}
There is no ruler-only construction verifying that a given point is the center of a given circle.
In particular, it is impossible to construct the center only with a ruler.
\end{thm}

The proof uses the construction in Exercise~\ref{ex:cone}.

\parit{Sketch of the proof.}
The same argument as in the problem above shows that 
it is sufficient to construct a projective transformation 
that sends the given circle $\Gamma$ to a circle $\Gamma'$ such that the center of $\Gamma'$ is not the image of the center of~$\Gamma$.

Choose a circle $\Gamma$ that lies in the plane $\Pi$ in the Euclidean space.
By Theorem~\ref{thm:inversion-3d}, 
the inverse of a circle across a sphere is a circle or a line.
Fix a sphere $\Sigma$ with the center $O$ so that the inversion $\Gamma'$ of $\Gamma$
is a circle and the plane $\Pi'$ containing $\Gamma'$ is not parallel to $\Pi$;
any sphere $\Sigma$ in a general position will do.

Let $Z$ and $Z'$ denote the centers of $\Gamma$ and~$\Gamma'$.
Note that  $Z'\z\notin(OZ)$.
It follows that the perspective projection $\Pi\to \Pi'$ with center $O$ sends $\Gamma$ to $\Gamma'$, but $Z'$ is not the image of~$Z$.
\qeds

\section{Comparison of construction tools}

We say that one set of tools is {}\emph{stronger} than another if any configuration of points that can be constructed with the second set can also be constructed with the first set.
If in addition, there is a configuration constructible with the first set, but not constructible with the second, then we say that the first set is {}\emph{strictly stronger} than the second.
Otherwise (that is, if any configuration that can be constructed with the first set can be constructed with the second), we say that the sets of tools are {}\emph{equivalent}. 
Two sets of tools might also be {}\emph{incomparable};
that is, there are constructions possible with the first set of tools and impossible with the second, and the other way around.

As an example, consider the following classical result:

\begin{thm}{Mohr--Mascheroni theorem}\index{Mohr--Mascheroni theorem}
Compass alone is equivalent to ruler and compass.
\end{thm}

The theorem does \textit{not} state that one can construct a whole line with a compass alone!
--- since we consider only configurations of points we do not have to.
One may think that a \textit{line is constructed} if we construct two points on it.

For sure ruler and compass form a stronger set than a compass alone.
Therefore the Mohr--Mascheroni theorem will follow once we solve the following two construction problems:
\begin{enumerate}[(i)]
\item Given four points $X$, $Y$, $P$, and $Q$, construct the intersection of the lines $(XY)$ and $(PQ)$ with  a compass only.
\item Given two points $X$, $Y$, and a circle $\Gamma$, construct the intersection of the lines $(XY)$ and $\Gamma$ with a compass only.
\end{enumerate}
Indeed, once we have these two constructions, we can do every step of a ruler-and-compass construction using a compass alone.

If you wonder how such a theorem can be proved, read about the \textit{Peaucellier--Lipkin inversor};
it is a planar linkage capable of transforming rotary motion into perfect straight-line motion.
Another classical theorem that can be proved using this linkage is the so-called \index{Poncelet--Steiner theorem}\emph{Poncelet--Steiner theorem};
it states that \textit{the set of ruler and compass is equivalent to the ruler alone, provided that a single circle and its center are given}.

\begin{thm}{Exercise}\label{ex:comparison}
Compare the following sets of tools:
(a) a ruler and compass, 
(b) a set-square, 
(c) a ruler and a parallel tool,
and
(d) a circumtool and an inversor.
\end{thm}

\chapter{Area}
\label{chap:area}

We will define area as a function that satisfies certain conditions (Section~\ref{sec:def(area)}).
The so-called \emph{Lebesgue measure} provides an example of such a function.
In particular, the existence of the Lebesgue measure implies the existence of an area function.
This construction can be found in any textbook on real analysis.

Relying solely on this existence result, we develop the concept of area with no cheating.
We adopt this approach because all rigorous introductions to area tend to be tedious.
(If you haven’t encountered the Lebesgue measure yet, you will sooner or later.)

\section{Solid triangles}\label{Solid triangles}

{

\begin{wrapfigure}{r}{27 mm}
\vskip-8mm
\centering
\includegraphics{mppics/pic-292}
\end{wrapfigure}

We say that a point $X$ lies \index{inside!a triangle}\emph{inside} a nondegenerate triangle $ABC$ if the following three conditions hold:

\begin{itemize}
\item $A$ and $X$ lie on the same side of the line~$(BC)$;
\item $B$ and $X$ lie on the same side of the line~$(CA)$;
\item $C$ and $X$ lie on the same side of the line~$(AB)$.
\end{itemize}

}

The set of all points inside $\triangle ABC$ 
and on its sides $[AB]$, $[BC]$, $[CA]$
will be called \index{triangle!solid triangle}\index{solid!triangle}\emph{solid triangle} $ABC$ and denoted by \index{25@$\solidtriangle$}$\solidtriangle ABC$.

\begin{thm}{Exercise}\label{ex:triangle-convex}
Show that every solid triangle is \index{convex set}\emph{convex};
that is, for every pair of points $X,Y\z\in\solidtriangle ABC$,
the line segment $[XY]$ lies in  $\solidtriangle ABC$.
\end{thm}


The notations $\triangle ABC$ and $\solidtriangle ABC$ look similar, 
and they also have closely related but different meanings, which it is better not to confuse.
Recall that $\triangle ABC$ is an ordered triple of distinct points
(see Section~\ref{sec:cong-triangles}),
while $\solidtriangle ABC$ is an infinite set of points.

In particular, $\solidtriangle ABC=\solidtriangle BAC$ for every triangle $ABC$.
Indeed, any point that belongs to the set $\solidtriangle ABC$ 
also belongs to the set $\solidtriangle BAC$
and the other way around.
On the other hand,
$\triangle ABC\ne\triangle BAC$ simply because the ordered triple of points $(A,B,C)$ is distinct from the ordered triple $(B,A,C)$.

Note that $\solidtriangle ABC\cong\solidtriangle BAC$ even if $\triangle ABC\ncong\triangle BAC$, where congruence of the sets $\solidtriangle ABC$ and $\solidtriangle BAC$ 
is understood the following way:

\begin{thm}{Definition}\label{def:cong-sets}
Two sets $\mathcal{S}$ and $\mathcal{T}$ in the plane  
are called \index{congruent!sets}\emph{congruent} 
(briefly \index{32@$\cong$}$\mathcal{S}\cong \mathcal{T}$)
if 
$\mathcal{T}=f(\mathcal{S})$ for some motion $f$ of the plane.
\end{thm}

If $\triangle ABC$ is not degenerate
and \[\solidtriangle ABC\cong \solidtriangle A'B'C',\]
then after relabeling the vertices of $\triangle ABC$ 
we will have 
\[\triangle ABC\cong \triangle A'B'C'.\]

Indeed it is sufficient to show that 
if $f$ is a motion that maps $\solidtriangle ABC$ to $\solidtriangle A'B'C'$,
then $f$ maps each vertex of $\triangle ABC$ to a vertex $\triangle A'B'C'$.
The latter follows from the characterization of vertices of solid triangles given in the following exercise:

\begin{thm}{Exercise}\label{ex:vertex}
Let $\triangle ABC$ be nondegenerate and $X\z\in \solidtriangle ABC$.
Show that $X$ is a vertex of $\triangle ABC$
if and only if there is a line $\ell$ that intersects $\solidtriangle ABC$
at the single point~$X$.
\end{thm}

\section{Polygonal sets}

An e\index{elementary set}\emph{elementary set} on the plane 
is a set of one of the following three types:
\begin{itemize}
 \item one-point set;
 \item segment;
 \item solid triangle.
\end{itemize}

\begin{wrapfigure}{r}{24 mm}
\vskip-12mm
\centering
\includegraphics{mppics/pic-294}
\end{wrapfigure}

A set in the plane is called \index{polygonal set}\emph{polygonal} if it can be presented as a union of a finite collection of elementary sets.

According to this definition, the empty set $\emptyset$
is a polygonal set.
Indeed, $\emptyset$ is a union of an empty collection of elementary sets.

A polygonal set is called \index{polygonal set!degenerate polygonal set}\index{degenerate!polygonal set}\emph{degenerate} if it can be presented as a union of a finite collection of one-point sets and segments.

If $X$ and $Y$ lie on opposite sides of the line $(AB)$,
then the polygonal set
$\solidtriangle AXB\cup \solidtriangle BYA$
is called \index{solid!quadrangle, parallelogram,\\ rectangle, square}\index{quadrangle!solid quadrangle}\emph{solid quadrangle} $AXBY$ and denoted by 
\index{27@$\solidsquare$}$\solidsquare AXBY$.
In particular, 
we can talk about \index{parallelogram!solid parallelogram}\emph{solid parallelograms}, \index{rectangle!solid rectangle}\emph{rectangles}, and \index{square!solid square}\emph{squares}.

\begin{wrapfigure}{o}{38 mm}
\centering
\includegraphics{mppics/pic-296}
\end{wrapfigure}

Typically a polygonal set admits many 
presentations as a union of a finite collection of elementary sets.
For example, if $\square AXBY$ is a parallelogram, then
\[\solidsquare AXBY=\solidtriangle AXB\cup \solidtriangle AYB=\solidtriangle XAY\cup \solidtriangle XBY.\]

\begin{thm}{Exercise}\label{ex:solid-square}
Show that a solid square is not degenerate.
\end{thm}

\begin{thm}{Exercise}\label{ex:poly-circ}
Show that a circle is not a polygonal set.
\end{thm}

\begin{thm}{Claim}\label{clm:poly-ring}
For any every polygonal sets $\mathcal{P}$ and $\mathcal{Q}$,
the union $\mathcal{P}\cup\mathcal{Q}$, 
as well as the intersection $\mathcal{P}\cap\mathcal{Q}$, 
are also polygonal sets.
\end{thm}

A class of sets that is closed with respect to union and intersection is called a {}\emph{ring of sets}.
The claim above, therefore, states that polygonal sets in the plane form a ring of sets.

\parit{Informal proof.}
Let us present $\mathcal{P}$ and $\mathcal{Q}$
as a union of a finite collection of elementary sets $\mathcal{P}_1,\dots,\mathcal{P}_k$ 
and $\mathcal{Q}_1,\dots,\mathcal{Q}_n$ respectively.

{

\begin{wrapfigure}{o}{22 mm}
\centering
\includegraphics{mppics/pic-298}
\end{wrapfigure}

Note that
\[\mathcal{P}\cup\mathcal{Q}
=
\mathcal{P}_1
\cup
\dots
\cup
\mathcal{P}_k
\cup
\mathcal{Q}_1
\cup
\dots
\cup
\mathcal{Q}_n.\]
Therefore, $\mathcal{P}\cup\mathcal{Q}$ is polygonal.

Note that $\mathcal{P}\cap \mathcal{Q}$ is the union of $\mathcal{P}_i\cap \mathcal{Q}_j$ for all $i$ an~$j$.
Therefore, to show that $\mathcal{P}\cap \mathcal{Q}$ is polygonal,
it is sufficient to show that $\mathcal{P}_i\cap \mathcal{Q}_j$ is polygonal for all $i$ and $j$.

The picture suggests a proof of the latter statement for solid triangles $\mathcal{P}_i$ and $\mathcal{Q}_j$.
The other cases are simpler; a formal proof can be built on Exercise~\ref{ex:triangle-convex}.
\qeds

\section{The definition of area}
\label{sec:def(area)}

\index{area}\emph{Area} is defined as a function $\mathcal{P}\mapsto \area\mathcal{P}$
that returns a nonnegative real number $\area\mathcal{P}$ for every polygonal set $\mathcal{P}$ and satisfies the following conditions:\label{page:area-def}
\textit{
\begin{enumerate}[(a)]
\item
$\area\mathcal{K}_1=1$
where  $\mathcal{K}_1$ is a solid unit square;
\item\label{area-3} the conditions
\begin{align*}
\mathcal{P}\cong\mathcal{Q}
\quad 
&\Rightarrow
\quad \area\mathcal{P}=\area\mathcal{Q};
\\
\mathcal{P}\subset\mathcal{Q}
\quad
&\Rightarrow
\quad 
\area\mathcal{P}\le\area\mathcal{Q};
\\
\area\mathcal{P}+\area\mathcal{Q}
&=
\area(\mathcal{P}\cup\mathcal{Q})+\area(\mathcal{P}\cap\mathcal{Q})
\end{align*}
hold 
for any two polygonal sets $\mathcal{P}$ and $\mathcal{Q}$.
\end{enumerate}
}

The first condition is called {}\emph{normalization}; essentially, it states that a solid unit square is used as the unit of area.
The three conditions in \textit{(\ref{area-3})} are called {}\emph{invariance}, {}\emph{monotonicity}, and {}\emph{additivity}.

The Lebesgue measure is an example of an area function.
Namely, if one takes the area of $\mathcal{P}$ to be its Lebesgue measure,
then the function $\mathcal{P} \mapsto \area\mathcal{P}$ satisfies the above conditions.

The construction of the Lebesgue measure can be found in any textbook on real analysis.
We do not discuss it here.

If the reader is not familiar with the Lebesgue measure, the existence of an area function can be taken for granted;
it might be added to the axioms, altho it follows from axioms \ref{def:birkhoff-axioms:0}--\ref{def:birkhoff-axioms:4}.
In any case, we will use the following statement without proof.

\begin{thm}{Claim}
An area function exists.
\end{thm}

\section{Vanishing area and subdivisions}

\begin{thm}{Proposition}\label{prop:area-segment}
For any two points $A$ and $B$, we have
\[\area[AB]=0\quad\text{and}\quad\area\{A\}=0.\]
\end{thm}

\parit{Proof.}
Fix a line segment~$[AB]$.
Consider a sold square $\solidsquare ABCD$.

For every positive integer $n$,
there are $n$ disjoint segments $[A_1B_1],\z\dots,[A_nB_n]$ 
in $\solidsquare ABCD$,
such that each $[A_iB_i]$ is congruent to $[AB]$ in the sense of the Definition~\ref{def:cong-sets}.

\begin{wrapfigure}[10]{o}{37 mm}
\vskip-2mm
\centering
\includegraphics{mppics/pic-300}
\vskip0mm
\end{wrapfigure}

Applying the invariance, additivity, and monotonicity of the area function, 
we get
\begin{align*}
n\cdot \area[AB]
&=\area\left([A_1B_1]\cup\dots\cup[A_nB_n]\right)\le
\\
&\le \area(\solidsquare ABCD)              
\end{align*}
In other words,
\[\area[AB]\le \tfrac1n\cdot\area(\solidsquare ABCD)\] 
for every positive integer~$n$.
It follows that $\area[AB]\le 0$. 
Since area cannot be negative, we get
\[\area[AB]= 0.\]

For every one-point set $\{A\}$
we have that $\{A\}\subset [AB]$.
Therefore, 
\[0\le \area\{A\}\le \area[AB]=0.\]
Hence $\area\{A\}=0$.
\qeds

\begin{thm}{Corollary}\label{cor:degenerate}
Any degenerate polygonal set has vanishing area.
\end{thm}

\parit{Proof.}
Let $\mathcal P$ be a degenerate set,
say
\[\mathcal{P}=[A_1B_1]\cup\dots\cup[A_nB_n]\cup\{C_1,\dots,C_k\}.\]
Since area is nonnegative by definition, applying additivity several times, we get that
\begin{align*}
\area\mathcal{P}\le
& \area[A_1B_1]+\dots+\area[A_nB_n]+
\\
&+\area\{C_1\}+\dots+\area\{C_k\}.
\end{align*}
By Proposition~\ref{prop:area-segment}, the right-hand side vanishes.

On the other hand, 
$\area\mathcal{P}\ge 0$,
hence the result.
\qeds

We say that a polygonal set $\mathcal{P}$ is \index{subdivision of polygonal set}\emph{subdivided} 
into two polygonal sets $\mathcal{Q}_1,\dots,\mathcal{Q}_n$ 
if $\mathcal{P}=\mathcal{Q}_1\cup\dots\cup \mathcal{Q}_n$ 
and the intersection $\mathcal{Q}_i\cap\mathcal{Q}_j$ if $i\ne j$.
(Recall that according to Claim~\ref{clm:poly-ring},
the intersections $\mathcal{Q}_i\cap\mathcal{Q}_j$ are polygonal.)

\begin{thm}{Proposition}\label{prop:subdivision}
Assume that a polygonal set $\mathcal{P}$ is subdivided into polygonal sets $\mathcal{Q}_1,\dots,\mathcal{Q}_n$.
Then 
\[\area\mathcal{P}=\area\mathcal{Q}_1+\dots+\area\mathcal{Q}_n.\]

\end{thm}

\begin{wrapfigure}{o}{26 mm}
\vskip-4mm
\centering
\includegraphics{mppics/pic-302}
\end{wrapfigure}

\parit{Proof.}
Assume \(n = 2\).
Since \(\mathcal{Q}_1 \cap \mathcal{Q}_2\) is degenerate, Corollary~\ref{cor:degenerate} implies that \(\area(\mathcal{Q}_1 \cap \mathcal{Q}_2) = 0\).
By the additivity of area,
\begin{align*}
\area \mathcal{P} &= \area \mathcal{Q}_1 + \area \mathcal{Q}_2 - \area (\mathcal{Q}_1 \cap \mathcal{Q}_2) \\
&= \area \mathcal{Q}_1 + \area \mathcal{Q}_2.
\end{align*}

By repeatedly applying this formula, we obtain the general case.
Indeed, suppose \(\mathcal{P}\) is subdivided into \(\mathcal{Q}_1, \dots, \mathcal{Q}_n\).
Observe that
\[
\mathcal{Q}_i \cap (\mathcal{Q}_{i+1} \cup \dots \cup \mathcal{Q}_n)
= (\mathcal{Q}_i \cap \mathcal{Q}_{i+1}) \cup \dots \cup (\mathcal{Q}_i \cap \mathcal{Q}_n)
\]
for each \(i\).
Since every intersection \(\mathcal{Q}_i \cap \mathcal{Q}_j\) is degenerate,
the intersection $\mathcal{Q}_i \cap (\mathcal{Q}_{i+1} \cup \dots \cup \mathcal{Q}_n)$
is also degenerate.
Consequently,
\[
\area(\mathcal{Q}_i \cap (\mathcal{Q}_{i+1} \cup \dots \cup \mathcal{Q}_n)) = 0
\]
for each \(i\).
Therefore,
\begin{align*}
\area \mathcal{P} &= \area \mathcal{Q}_1 + \area (\mathcal{Q}_2 \cup \dots \cup \mathcal{Q}_n) \\
&= \area \mathcal{Q}_1 + \area \mathcal{Q}_2 + \area (\mathcal{Q}_3 \cup \dots \cup \mathcal{Q}_n) \\
&\ \ \ \vdots \\
&= \area \mathcal{Q}_1 + \area \mathcal{Q}_2 + \dots + \area \mathcal{Q}_n.
\end{align*}
\qedsf

\parbf{Remark.}
Two polygonal sets $\mathcal{P}$ and $\mathcal{P}'$ are called \index{equidecomposable sets}\emph{equidecomposable} if they admit subdivisions into polygonal sets $\mathcal{Q}_1,\dots,\mathcal{Q}_n$ and $\mathcal{Q}'_1,\dots,\mathcal{Q}'_n$ such that 
$\mathcal{Q}_i\cong\mathcal{Q}'_i$ for each $i$.

According to the proposition, if $\mathcal{P}$ and $\mathcal{P}$ are equidecomposable, then $\area \mathcal{P}=\area\mathcal{P}'$.
A converse to this statement also holds;
namely, \textit{if two nondegenerate polygonal sets have equal area, then they are equidecomposable.}

The last statement was proved by William Wallace, Farkas Bolyai, and Paul Gerwien.
The analogous statement in three dimensions, known as {}\emph{Hilbert's third problem}, is false; it was proved by Max Dehn.

\section{Rectangles}

\begin{thm}{Theorem}\label{thm:area-rect}
The area of a solid rectangle is the product of its adjacent sides.
\end{thm}

\parit{Proof.}
Let $\mathcal{R}_{a,b}$ denote the solid rectangle with adjacent sides $a$ and~$b$,
and let $s(a,b)\df\area \mathcal{R}_{a,b}$.
We need to show that
\[s(a,b)=a\cdot b\]
for any $a>0$ and $b>0$.

Since $\mathcal{R}_{1,1}$ is a solid unit square, we have
\[s(1,1)=1.\]

Given a positive integer $n$,
we can subdivide a solid unit square into $n^2$ squares with side length $\tfrac1n$.
By Proposition~\ref{prop:subdivision} and the invariance of area, we get $n^2\cdot s(\tfrac1n,\tfrac1n)=1$ and
\[s(\tfrac1n,\tfrac1n)=\tfrac1n\cdot\tfrac1n.\]

Further, for positive integers $k$, $m$, and $n$, the rectangle $\mathcal{R}_{\frac kn,\frac mn}$ can be subdivided into $k\cdot m$ squares with side length $\tfrac1n$.
Therefore,
\[s(\tfrac kn,\tfrac mn)=\tfrac kn\cdot \tfrac mn,\eqlbl{eq:a(a,b)=ab-}\]
which is a special case of the theorem.

Suppose $a\le a'$ and $b\le b'$.
Then we can assume that $\mathcal{R}_{a,b}\subset \mathcal{R}_{a',b'}$, and by the monotonicity of area, we have
\[s(a,b)\le s(a',b').\eqlbl{s(a,b)>s(a',b')}\]

Now let us argue by contradiction.
Assume that $s(a,b)\ne a\cdot b$ for some $a>0$ and $b>0$.
Then $s(a,b)> a\cdot b$ or $s(a,b) <a\cdot b$.

\raggedcolumns\setlength{\multicolsep}{.5mm}
\setlength{\columnseprule}{1pt}
\begin{multicols}{2}
If $s(a,b)> a\cdot b$,
we can choose a positive integer $n$ such that
\[s(a,b)> (a+\tfrac1n)\cdot (b+\tfrac1n).\eqlbl{s(a,b)>}\]
Choose positive integers $k$ and $m$ such that
\begin{align*}
a< \tfrac kn&\le a+\tfrac1n,
\\
b<\tfrac mn&\le b+\tfrac1n.
\end{align*}
By \ref{s(a,b)>s(a',b')} and \ref{eq:a(a,b)=ab-}, we get that
\begin{align*}
s(a,b)&\le s(\tfrac kn,\tfrac mn)=
\\
&=\tfrac kn\cdot\tfrac mn\le
\\
&\le (a+\tfrac1n)\cdot(b+\tfrac1n),
\end{align*}
which contradicts \ref{s(a,b)>}.

\columnbreak

The case $s(a,b)< a\cdot b$ is analogous.
Choose $n$ such that $a>\tfrac1n$, $b>\tfrac1n$, and
\[s(a,b)< (a-\tfrac1n)\cdot (b-\tfrac1n).\eqlbl{s(a,b)<}\]
Choose $k$ and $m$ such that
\begin{align*}
a> \tfrac kn&\ge a-\tfrac1n,
\\
b>\tfrac mn&\ge b-\tfrac1n.
\end{align*}
By \ref{s(a,b)>s(a',b')} and \ref{eq:a(a,b)=ab-}, we get that
\begin{align*}
s(a,b)&\ge s(\tfrac kn,\tfrac mn)=
\\
&=\tfrac kn\cdot\tfrac mn\ge
\\
&\ge (a-\tfrac1n)\cdot(b-\tfrac1n),
\end{align*}
which contradicts \ref{s(a,b)<}.
\qeds
\end{multicols}
\setlength{\columnseprule}{0pt}

\section{Parallelograms}

\begin{thm}{Proposition}\label{prop:area-parallelogram}
Let $\square ABCD$ be a parallelogram in the Euclidean plane.
Then 
\[\area(\solidsquare ABCD)=a\cdot h,\]
where $a=AB$, and $h$ is the distance between the lines $(AB)$ and~$(CD)$.
\end{thm}

{

\begin{wrapfigure}{r}{41 mm}
\vskip-2mm
\centering
\includegraphics{mppics/pic-306}
\end{wrapfigure}

\parit{Proof.}
Let $A'$ and $B'$ denote the footpoints of $A$ and $B$
on the line~$(CD)$.

Note that $ABB'A'$ is a rectangle with sides $a$ and $h$.
By Proposition~\ref{thm:area-rect},
\[\area( \solidsquare ABB'A')=h\cdot a.
\eqlbl{eq:ABBA}\]

}

Without loss of generality, we may assume that  $\solidsquare ABCA'$ 
contains $\solidsquare ABCD$ and $\solidsquare ABB'A'$.
In this case, $\solidsquare ABCA'$ admits two subdivisions: 
\[\solidsquare ABCA'=\solidsquare ABCD\cup\solidtriangle AA'D=\solidsquare ABB'A'\cup\solidtriangle BB'C.\]

By Proposition~\ref{prop:subdivision},
\[\begin{aligned}
\area( \solidsquare ABCD)&+\area(\solidtriangle AA'D)=
\\
&=
\area(\solidsquare ABB'A')+ \area (\solidtriangle BB'C).   
  \end{aligned}
\eqlbl{eq:area-sum}\]

Note that 
\[\triangle AA'D\cong \triangle BB'C.\eqlbl{eq:cong-area}\]
Indeed, since the quadrangles $ABB'A'$ and $ABCD$ are parallelograms, 
by Lemma~\ref{lem:parallelogram},
we have that $AA'=BB'$, $AD=BC$, and $DC=AB\z=A'B'$.
It follows that $A'D=B'C$.
Applying the SSS congruence condition, we get \ref{eq:cong-area}.

In particular,
\[\area(\solidtriangle BB'C)=\area (\solidtriangle AA'D).
\eqlbl{eq:area-trinagles}\]

Subtracting \ref{eq:area-trinagles} from \ref{eq:area-sum},
we get that
\[\area (\solidsquare ABCD)=\area(\solidsquare ABB'D).\]
It remains to apply \ref{eq:ABBA}.
\qeds

{

\begin{wrapfigure}{r}{26 mm}
\vskip-5mm
\centering
\includegraphics{mppics/pic-308}
\end{wrapfigure}

\begin{thm}{Exercise}\label{ex:two-parallelograms}
Assume $\square ABCD$ and $\square AB'C'D'$ are two parallelograms such that $B'\in[BC]$ and $D\z\in [C'D']$.
Show that
\[\area(\solidsquare ABCD)=\area(\solidsquare AB'C'D').\]

\end{thm}

}

\section{Triangles}

\begin{thm}{Theorem}\label{thm:area-of-triangle}
Let $h_A$ be the altitude from $A$
in  $\triangle ABC$ and $a=BC$.
Then 
\[\area(\solidtriangle ABC)=\tfrac12\cdot a\cdot h_A.\]

\end{thm}

\parit{Proof.}
Draw the line $m$ thru $A$ that is parallel to $(BC)$
and line $n$ thru $C$ parallel to~$(AB)$.
Note that the lines $m$ and $n$ are not parallel;
denote by $D$ their point of intersection.
By construction, $\square ABCD$ is a parallelogram.

Note that $\solidsquare ABCD$ admits a subdivision into $\solidtriangle ABC$ and $\solidtriangle CDA$.
Therefore, 
\[\area(\solidsquare ABCD)
=
\area(\solidtriangle ABC)
+
\area(\solidtriangle CDA).\]

\begin{wrapfigure}{o}{28 mm}
\centering
\includegraphics{mppics/pic-310}
\vskip2mm
\end{wrapfigure}

Since $\square ABCD$ is a parallelogram,  Lemma~\ref{lem:parallelogram} implies that
\[AB=CD
\quad
\text{and}
\quad
BC=DA.\]
Therefore, by the SSS congruence condition, we have
$\triangle ABC\z\cong\triangle CDA$.
In particular
\[\area(\solidtriangle ABC)
=
\area(\solidtriangle CDA).\]

From above and Proposition~\ref{prop:area-parallelogram}, we get that
\begin{align*}
\area(\solidtriangle ABC)
&=\tfrac12\cdot\area(\solidsquare ABCD)=
\\
&=\tfrac12\cdot h_A\cdot a.
\end{align*}
\qedsf

\begin{thm}{Exercise}\label{ex:three-trig}
Let $h_A$, $h_B$, and $h_C$ denote the altitudes of $\triangle ABC$ from vertices $A$, $B$, and $C$ respectively.
Note that from Theorem~\ref{thm:area-of-triangle},
it follows that
\[h_A\cdot BC=h_B\cdot CA=h_C\cdot AB.\]

Give a proof of this statement without using Theorem~\ref{thm:area-of-triangle}.
\end{thm}

\begin{thm}{Exercise}\label{ex:half-parallelogram}
Assume $M$ lies inside the parallelogram $ABCD$;
that is, $M$ belongs to the solid parallelogram $\solidsquare ABCD$ but does not lie on its sides.
Show that
\[\area(\solidtriangle ABM)+\area(\solidtriangle CDM)
=\tfrac12\cdot \area(\solidsquare ABCD).\]
\end{thm}

\begin{thm}{Exercise}\label{ex:area-diag}
Assume that diagonals 
of a nondegenerate quadrangle $ABCD$ 
intersect at point $M$.
Show that 
\[\area(\solidtriangle ABM)\cdot\area(\solidtriangle CDM)
=
\area(\solidtriangle BCM)\cdot\area(\solidtriangle DAM).\]
 
\end{thm}

\begin{thm}{Exercise}\label{ex:area-inradius}
Let $r$ be the inradius of a nondegenerate triangle $ABC$,
and let $p$ be its {}\emph{semiperimeter};
that is, $p=\tfrac12\cdot(AB+BC+CA)$.
Show that
\[\area(\solidtriangle ABC)=p\cdot r.\]

\end{thm}

\begin{thm}{Exercise}\label{ex:subdivision}
Show that every polygonal set admits a subdivision into a finite collection of solid triangles and a degenerate set.
Conclude that for every polygonal set, its area is uniquely defined.
\end{thm}

\section{The area method}

In this section, we give examples of proofs using the properties of the area function.
These proofs are not truly elementary --- the price for the existence of the area function is high.

{

\begin{wrapfigure}{r}{20 mm}
\vskip-0mm
\centering
\includegraphics{mppics/pic-312}
\end{wrapfigure}

We start with the proof of the \index{Pythagorean theorem}Pythagorean theorem.
In the Elements of Euclid, the Pythagorean theorem was formulated as equality  \ref{eq:pyth+area} below,
and the proof used a similar technique.

\parit{Proof.}
We need to show that 
\[a^2+b^2=c^2,\]
where $a$ and $b$ are legs and $c$ is the hypotenuse 
of a right triangle.

Let $\mathcal{K}_{x}$ be the solid square with side~$x$.
Denote by $\mathcal{T}$ the right solid triangle with legs $a$ and $b$.

}

Let us construct two subdivisions of $\mathcal{K}_{a+b}$:
\begin{enumerate}[1.]
\item Subdivide $\mathcal{K}_{a+b}$ into two solid squares congruent to $\mathcal{K}_a$ and $\mathcal{K}_b$
and four solid triangles congruent to $\mathcal{T}$;
see the first picture.

\item Subdivide $\mathcal{K}_{a+b}$ into one solid square congruent to $\mathcal{K}_c$
and four solid right triangles congruent to $\mathcal{T}$;
see the second picture.

\end{enumerate}

Applying Proposition~\ref{prop:subdivision} a few times,
we get that
\begin{align*}
\area\mathcal{K}_{a+b}
&=
\area\mathcal{K}_{a}+\area\mathcal{K}_{b}+4\cdot\area\mathcal{T}=
\\
&=\area\mathcal{K}_{c}+4\cdot\area\mathcal{T}.
\end{align*}
Therefore, 
\[\area\mathcal{K}_{a}+\area\mathcal{K}_{b}=\area\mathcal{K}_{c}.\eqlbl{eq:pyth+area}\]
By Theorem~\ref{thm:area-rect}, we know that 
\[\area\mathcal{K}_x=x^2,\] 
for any $x>0$. 
Hence the statement follows.\qeds

{

\begin{wrapfigure}{r}{19mm}
\vskip-4mm
\centering
\includegraphics{mppics/pic-314}
\end{wrapfigure}

\begin{thm}{Exercise}\label{ex:pyth-2}
Build another proof of the Pythagorean theorem
based on the picture.

(In the notations above it shows a subdivision of $\mathcal{K}_c$ into $\mathcal{K}_{a-b}$ and four copies of~$\mathcal{T}$ if $a>b$.)
\end{thm}

} 

\begin{thm}{Exercise}\label{ex:sum-3-dist}
Show that the sum of distances from a point to the sides of an equilateral triangle is the same for all points inside the triangle.
\end{thm}

\begin{thm}{Claim}\label{clm:area-ratio}
Assume  that two triangles $ABC$ and $A'B'C'$ in the Euclidean plane 
have equal altitudes dropped from $A$ and $A'$ respectively.
Then
\[\frac{\area(\solidtriangle A'B'C')}{\area(\solidtriangle ABC)}
=
\frac{B'C'}{BC}.\]

In particular, the same identity holds if $A=A'$ and the bases $[BC]$ and $[B'C']$ lie on one line.
\end{thm}

\parit{Proof.}
Let $h$ be the altitude.
By Theorem~\ref{thm:area-of-triangle},
\[\frac{\area(\solidtriangle A'B'C')}{\area(\solidtriangle ABC)}
=
\frac{\frac12 \cdot h\cdot B'C'}{\frac12 \cdot h\cdot BC}
=
\frac{B'C'}{BC}.\]
\qedsf

\begin{thm}{Exercise}\label{ex:area-medians}
Prove that the medians of a nondegenerate triangle divide it into six triangles of equal area.
\end{thm}

Now let us show how to use this claim to prove Lemma~\ref{lem:bisect-ratio}.
First, let us recall its statement:

\begin{thm*}{Lemma}
If $\triangle A B C$ is nondegenerate and its angle bisector at $A$ intersects $[BC]$ at point~$D$, then 
$\frac{AB}{AC}=\frac{DB}{DC}$.

\end{thm*}

\begin{wrapfigure}{r}{30 mm}
\vskip-4mm
\centering
\includegraphics{mppics/pic-316}
\end{wrapfigure}

\parit{Proof.}
Applying  Claim~\ref{clm:area-ratio}, we get that
\[\frac{\area(\solidtriangle ABD)}{\area(\solidtriangle ACD)}
=\frac{BD}{CD}.\]

By Proposition~\ref{prop:angle-bisect-dist} the triangles $ABD$ and $ACD$ have equal altitudes from~$D$.
Applying  Claim~\ref{clm:area-ratio} again, we get that
\[\frac{\area(\solidtriangle ABD)}{\area(\solidtriangle ACD)}=\frac{AB}{AC}.\]
Hence the result follows.
\qeds

The second statement in the following exercise is a partial case of \index{Ceva's theorem}Ceva's theorem; see \ref{thm:ceva-affine}.

\begin{thm}{Exercise}\label{ex:ceva}
Let $ABC$ be a nondegenerate triangle.
Assume $A'$ lies between $B$ and $C$,
point $B'$ lies between $C$ and $A$,
point $C'$ lies between $A$ and $B$.
Suppose that line segments $[AA']$, $[BB']$, and $[CC']$ meet at a point $X$.
Show that 

\vskip-3mm

{

\begin{wrapfigure}[6]{r}{40 mm}
\vskip4mm
\centering
\includegraphics{mppics/pic-318}
\end{wrapfigure}

\begin{align*}
\frac{\area(\solidtriangle ABX)}{\area(\solidtriangle BCX)}&=\frac{AB'}{B'C},
\\
\frac{\area(\solidtriangle BCX)}{\area(\solidtriangle CAX)}&=\frac{BC'}{C'A},
\\
\frac{\area(\solidtriangle CAX)}{\area(\solidtriangle ABX)}&=\frac{CA'}{A'B} .
\end{align*}

}

Conclude that 
\[\frac{AB'\cdot CA'\cdot BC'}{B'C\cdot A'B\cdot C'A}=1.\]
\end{thm}

\begin{wrapfigure}[5]{r}{40 mm}
\vskip-0mm
\centering
\includegraphics{mppics/pic-319}
\end{wrapfigure}

\begin{thm}{Exercise}\label{ex:cross-ratio-area}
Suppose that points $L_1$, $L_2$, $L_3$, $L_4$ lie on a line $\ell$ 
and points $M_1$, $M_2$, $M_3$, $M_4$ lie on a line $m$. 
Assume that the lines $(L_1M_1)$, $(L_2M_2)$, $(L_3M_3)$, and $(L_4M_4)$ pass thru point $O$ that does not lie on $\ell$ nor~$m$.

\begin{enumerate}[(a)]
 \item\label{ex:cross-ratio-area:a} Apply Claim~\ref{clm:area-ratio} to show that
 \[\frac{\area\solidtriangle OL_iL_j}{\area\solidtriangle OM_iM_j}=\frac{OL_i\cdot OL_j}{OM_i\cdot OM_j}\]
 for any $i\ne j$.
 \item\label{ex:cross-ratio-area:b} Use \textit{(\ref{ex:cross-ratio-area:a})} to prove that 
 \[\frac{L_1L_2\cdot L_3L_4}{L_2L_3\cdot L_4L_1}=\frac{M_1M_2\cdot M_3M_4}{M_2M_3\cdot M_4M_1};\]
 that is, the quadruples $(L_1, L_2, L_3, L_4)$ and $(M_1, M_2, M_3, M_4)$ have the same cross-ratio.
 
\end{enumerate}

\end{thm}

\section{Neutral planes and spheres}\label{Neutral planes and spheres}

Area can be defined in the neutral planes and spheres.
In this definition,
the solid unit square $\mathcal{K}_1$ has to be  
exchanged for a fixed nondegenerate polygonal set $\mathcal{U}$.
Such a change is necessary for a good reason --- 
hyperbolic plane and sphere have no squares.
In this case, the set $\mathcal{U}$ serves as the unit measure for area;
changing $\mathcal{U}$ would require a conversion of area units.

\begin{wrapfigure}{o}{26 mm}
\vskip-0mm
\centering
\includegraphics{mppics/pic-320}
\end{wrapfigure}

According to the standard convention, the set $\mathcal{U}$
is chosen so that on small scales the area behaves like in the Euclidean plane.
For example, if $\mathcal{K}_a$ denotes the solid quadrangle $\solidsquare ABCD$ 
with right angles at $A$, $B$, and $C$ such that  $AB=BC=a$, 
then we may assume that
\[\tfrac{1}{a^2}\cdot\area \mathcal{K}_a\to 1
\quad
\text{as}
\quad 
a\to0.\]

This convention works equally well for spheres and neutral planes, including the Euclidean plane.
In spherical geometry  equivalently we may assume that if $r$ is the radius of the sphere, 
then the area of the whole sphere is $4\cdot\pi\cdot r^2$.

Recall that the \index{defect}\emph{defect of triangle} $\triangle ABC$ is defined as 
$$\defect(\triangle ABC)
\df 
\pi-|\measuredangle ABC|-|\measuredangle BCA|-|\measuredangle CAB|.$$
It turns out that for every neutral plane or sphere,
there is a real number $k$
such that 
$$k\cdot\area(\solidtriangle ABC)+\defect(\triangle ABC)=0
\eqlbl{eq:curv-defect}$$
for every $\triangle ABC$.

This number $k$ is called \index{curvature}\emph{curvature};
$k=0$ for the Euclidean plane,
$k=-1$ for the h-plane, 
$k=1$ for the unit sphere (this case is proved in Lemma \ref{lem:area-spher-triangle}),
and $k=\tfrac1{r^2}$ for the sphere of radius~$r$.

Since the angles of an ideal triangle vanish, every ideal triangle in the h-plane has area~$\pi$.
Similarly, in the unit sphere, the area of an equilateral triangle with right angles has an area of $\tfrac\pi2$;
since the whole sphere can be subdivided into eight such triangles, we get that the area of the unit sphere is $4\cdot\pi$.

The identity \ref{eq:curv-defect} can be used as an alternative way to introduce the area function; it works on spheres and all neutral planes, except for the Euclidean plane.

\section{Quadrable sets}

A set $\mathcal{S}$ 
in the plane is called \index{quadrable set}\emph{quadrable}
if, for any $\epsilon>0$, there are two polygonal sets 
$\mathcal{P}$ and $\mathcal{Q}$
such that 
\[\mathcal{P}
\subset
\mathcal{S}\subset\mathcal{Q}
\quad
\text{and}
\quad
\area\mathcal{Q}-\area\mathcal{P}
<
\epsilon.\]

If $\mathcal{S}$ is quadrable,
its area can be defined
as the necessarily unique real number $s=\area\mathcal{S}$
such that the inequality
\[\area\mathcal{Q}\le s\le \area\mathcal{P}
\]
holds for any polygonal sets $\mathcal{P}$ and $\mathcal{Q}$ such that $\mathcal{P}\subset\mathcal{S}\subset\mathcal{Q}$.

\begin{thm}{Exercise}\label{ex:circle-is-quadrable}
Let $\mathcal{D}$ be a \index{disc}\emph{closed unit disc};
that is, $\mathcal{D}$ is a set that contains 
the unit circle $\Gamma$ and all the points inside~$\Gamma$.

Show that $\mathcal{D}$ is a quadrable set.
\end{thm}

Since $\mathcal{D}$ is quadrable, the expression $\area\mathcal{D}$ makes sense and the constant $\pi$ can be defined as $\pi\df\area\mathcal{D}$.

\medskip

It turns out that the class of quadrable sets is the largest class for which 
the area function satisfying the conditions in Section~\ref{sec:def(area)} is \textit{uniquely} defined.

If you do not require uniqueness, then there are ways to extend the area function to all bounded sets.
(A set in the plane is called \index{bounded set}\emph{bounded} if it lies inside a circle.)
On the sphere and hyperbolic plane, 
there is no similar construction.
If you wonder why,
read about the \index{doubling the ball}\emph{doubling the ball} --- a paradox of Felix Hausdorff, Stefan Banach, and Alfred Tarski.

\newgeometry{top=0.9in, bottom=0.9in,inner=0.55in, outer=0.45in}
{\footnotesize

\backmatter

\chapter{Hints}

\raggedcolumns\setlength{\multicolsep}{-7mm}
\spell{\begin{multicols}{2}}{}

\refstepcounter{chapter}
\setcounter{eqtn}{0}

\parbf{\ref{ex:dist-square}.} Check the triangle inequality for $A\z=0$, $B\z=1$, and $C\z=2$.

\parbf{\ref{ex:d_1+d_2+d_infty}.} 
Check all the conditions in \ref{def:metric-space}.
Let us discuss the triangle inequality --- the remaining conditions are self-evident.

Let $A=(x_A,y_A)$, $B=(x_B,y_B)$, and $C=(x_C,y_C)$.
Set 
\begin{align*}
x_1&=x_B-x_A, 
&
y_1&=y_B-y_A,
\\
x_2&=x_C-x_B,
&
y_2&=y_C-y_B.
\end{align*}

\parit{(a).}
The inequality
$$d_1(A,C)\le d_1(A,B)+d_1(B,C)$$
can be written as 
$$|x_1+x_2|+|y_1+y_2|
\le 
|x_1|+|y_1|+|x_2|+|y_2|.$$
The latter follows since $|x_1+x_2|\le |x_1|+|x_2|$ 
and
$|y_1+y_2|\le |y_1|+|y_2|$.

\parit{(b).}
The inequality
$$d_2(A,C)\le d_2(A,B)+d_2(B,C)\eqlbl{eq:trig-inq-d2}$$
can be written as 
\begin{align*}
&\sqrt{\bigl(x_1+x_2\bigr)^2+\bigl(y_1+y_2\bigr)^2}\le
\\
&\qquad\le
\sqrt{x_1^2+y_1^2}+\sqrt{x_2^2+y_2^2}.
\end{align*}
Take the square of the left and the right-hand sides,
simplify,
take the square again, and simplify again.
You should get the following inequality:
$$0
\le 
(x_1\cdot y_2-x_2\cdot y_1)^2,$$
which is equivalent to \ref{eq:trig-inq-d2}
and evidently true.

\parit{(c).}
The inequality
$$d_\infty(A,C)\le d_\infty(A,B)+d_\infty(B,C)$$
can be written as 
$$
\begin{aligned}
&\max\{|x_1+x_2|,|y_1+y_2|\}\le
\\
&\qquad\le 
\max\{|x_1|,|y_1|\}+\max\{|x_2|,|y_2|\}.
\end{aligned}
\eqlbl{eq:max-trig}$$
Without loss of generality, we may assume that 
$$\max\{|x_1+x_2|,|y_1+y_2|\}=|x_1+x_2|.$$
Furthermore,
\begin{align*}
|x_1+x_2|&\le |x_1|+|x_2|\le 
\\
&\le\max\{|x_1|,|y_1|\}+\max\{|x_2|,|y_2|\}.
\end{align*}
Hence \ref{eq:max-trig} follows.

\parbf{\ref{ex:4-triangle}.} Sum up four triangle inequalities.

\parbf{\ref{ex:dist-preserv=>injective}.}
If $A\ne B$, then $d_\mathcal{X}(A,B)>0$.
Since $f$ is distance-preserving,
$$d_\mathcal{Y}(f(A),f(B))=d_\mathcal{X}(A,B).$$
Therefore, $d_\mathcal{Y}(f(A),f(B))>0$; hence $f(A)\z\ne f(B)$.

\parbf{\ref{ex:motion-of-R}.}
Set $f(0)=a$ and $f(1)=b$.
Show that $b=a+1$ or $a-1$.
Moreover, $f(x)=a\pm x$, and at the same time, $f(x)=b\pm(x-1)$ for any~$x$.

Suppose $b=a+1$. 
Show that 
$f(x)=a+x$ for any~$x$.

In the same way, if $b=a-1$, 
show that 
$f(x)=a-x$ for any~$x$.

\parbf{\ref{ex:d_1=d_infty}.} 
Show that the map $(x,y)\mapsto (x+y,x-y)$ is an isometry  $(\mathbb{R}^2,d_1)\z\to (\mathbb{R}^2,d_\infty)$.
You need to check if this map is bijective and distance-preserving.

\parbf{\ref{ad-ex:motions of Manhattan plane}.} 
First prove that \textit{two points $A=(x_A,y_A)$ and $B\z=(x_B,y_B)$ on the Manhattan plane have a unique midpoint if and only if $x_A=x_B$ or $y_A=y_B$}; compare with the example in \ref{sec:cong-triangles}. 

Use the above statement to prove that
any motion of the Manhattan plane 
can be written in one of the following eight ways:
\begin{align*}
(x,y)&\mapsto (\pm x+a,\pm y+b)
\shortintertext{or} 
(x,y)&\mapsto (\pm y+b,\pm x+a),
\end{align*}
for fixed real numbers $a$ and~$b$.
In each case, we have 4 choices of signs, so for a fixed pair $(a,b)$ we have 8 distinct motions.

\parbf{\ref{ex:y=|x|}.}
Assume three points $A$, $B$, and $C$ lie on one line.
Note that in this case one of three triangle inequalities involving the points $A$, $B$, and $C$ becomes an equality.

Set $A=(-1,1)$, $B=(0,0)$, and $C=(1,1)$.
Show that for $d_1$ and $d_2$
all the triangle inequalities with the points $A$, $B$, and $C$ are strict.
It follows that the graph is not a line.

For $d_\infty$ show that $(x,|x|)\mapsto x$ is a isometry from the graph to~$\mathbb{R}$.
Conclude that the graph is a line in $(\mathbb{R}^2,d_\infty)$.

\parbf{\ref{ex:line-motion}.}
Spell the definitions of line and motion.

\parbf{\ref{ex:trig==}.}
According to the definition of a half-line,
there is an isometry
$$f\:[P Q)\to [0,\infty)$$
such that $f(P)=0$.
By the definition of an isometry, $P X=f(X)$ for every $X\z\in [P Q)$.
Thus, $P X=r$ if and only if $f(X)=r$.

Since an isometry must be bijective, the statement follows.




\parbf{\ref{ex:2a=0}.}
The equation
$2\cdot\alpha\equiv 0$
means that $2\cdot\alpha\z=2\cdot k\cdot\pi$ for an integer~$k$.
Therefore,
$\alpha=k\cdot\pi$.

Equivalently, $\alpha=2\cdot n\cdot \pi$ or $\alpha=(2\cdot n+1)\cdot \pi$ for an integer~$n$.
In these cases, we have $\alpha\equiv 0$ or $\alpha\equiv \pi$ respectively.

\parbf{\ref{ex:a+b==c}.}
Observe that $\gamma'=\alpha+\beta\in [0,2\cdot \pi]$.
Show and use that if $\gamma'\equiv\gamma$, then $\gamma'=\gamma$.

\parbf{\ref{ex:dist-cont}.} \textit{(a).}
By the triangle inequality,
$|f(A')\z-f(A)|\le d(A',A)$.
Therefore, we can take $\delta=\epsilon$.

\parit{(b).}
By the triangle inequality,
\begin{align*}
&|g(A',B')-g(A,B)|
\le 
\\
&\le|g(A',B')-g(A,B')|+
\\
&\quad+
|g(A,B')-g(A,B)|\le
\\
&\le d(A',A)+d(B',B).
\end{align*}
Therefore, we can take $\delta=\tfrac\epsilon2$.

\parbf{\ref{ex:comp+cont}.}
Fix $A\in \mathcal{X}$ and $B\in\mathcal{Y}$
such that $f(A)=B$.

Fix $\epsilon>0$.
Since $g$ is continuous at $B$, there is a positive value $\delta_1$ such that 
$$d_{\mathcal{Z}}(g(B'),g(B))<\epsilon
\quad
\text{if}
\quad
d_{\mathcal{Y}}(B',B)<\delta_1.$$ 

Since $f$ is continuous at $A$, there is $\delta_2>0$ such that 
$$d_{\mathcal{Y}}(f(A'),f(A))\z<\delta_1
\quad
\text{if}
\quad
d_{\mathcal{X}}(A',A)<\delta_2.$$ 

Since $f(A)=B$, we get that
$$d_{\mathcal{Z}}(h(A'),h(A))<\epsilon
\quad
\text{if}
\quad
d_{\mathcal{X}}(A',A)<\delta_2.$$ 
Hence the result.

\parbf{\ref{ex:ncong}.} \textit{(a).}
Show that $A\mapsto B$, $B\mapsto A$, $C\mapsto C$, $D\mapsto D$ defines a motion of the space;
conclude that $\triangle ABC\cong \triangle BAC$.

\parit{(b).} Suppose $\triangle ABC\cong \triangle BCA$; so, there is a motion $m$ that maps $A\mapsto B$, $B\mapsto C$, and $C\mapsto A$.
Show that $m\:D\mapsto D$; arrive at a contradiction.

\refstepcounter{chapter}
\setcounter{eqtn}{0}

\parbf{\ref{ex:infinite}.} By Axiom~\ref{def:birkhoff-axioms:0}, there are at least two points in the plane.
Therefore, by Axiom~\ref{def:birkhoff-axioms:1}, 
the plane has a line.
To prove \textit{(a)}, it remains to note that a line is an infinite set of points.
To prove \textit{(b)} apply Axiom~\ref{def:birkhoff-axioms:2} in addition.

\parbf{\ref{ex:[OA)=[OA')}.}
By Axiom~\ref{def:birkhoff-axioms:1},
$(OA)=(OA')$.
Therefore, the statement boils down to the following:

\textit{Assume $f\:\mathbb{R}\to \mathbb{R}$ is a motion of the line that sends $0\mapsto 0$ and one positive number to a positive number, then $f$ is an identity map.}

The latter follows from \ref{ex:motion-of-R}.

\parbf{\ref{ex:2.4}.}
By \ref{lem:AOA=0},
$\measuredangle AOA=0$.
It remains to apply Axiom~\ref{def:birkhoff-axioms:2a}.

\parbf{\ref{ex:lineAOB}.}
Apply \ref{lem:AOA=0},
\ref{thm:straight-angle},
and \ref{ex:2a=0}.

\parbf{\ref{ex:ABCO-line}.}
By Axiom~\ref{def:birkhoff-axioms:2b},
$2\cdot\measuredangle BOC
\equiv 
2\cdot\measuredangle AOC\z-2\cdot \measuredangle AOB
\equiv 0$.
By \ref{ex:2a=0}, 
this implies that 
$\measuredangle BOC$ is either $0$ or~$\pi$.
It remains to apply \ref{ex:2.4} and \ref{thm:straight-angle} respectively in these two cases.

\parbf{\ref{ex:infinite-number-of-lines}.}
Fix two points $A$ and $B$ as provided by Axiom~\ref{def:birkhoff-axioms:0}.

Choose a real number $0<\alpha<\pi$.
By Axiom~\ref{def:birkhoff-axioms:2a} there is a point $C$ such that $\measuredangle ABC\z=\alpha$.

Use \ref{cor:degenerate=pi} to show that $\triangle ABC$ is nondegenerate.


\parbf{\ref{ex:refelection-of-line}.} Apply \ref{prop:point-reflection} and \ref{ex:line-motion}.

\refstepcounter{chapter}
\setcounter{eqtn}{0}

\parbf{\ref{ex:AOB+<=>BOA-}.}
Set $\alpha=\measuredangle AOB$ 
and 
$\beta=\measuredangle BOA$.
Note that $\alpha=\pi$ if and only if $\beta=\pi$.
Otherwise, $\alpha=-\beta$.
Hence the result.

\parbf{\ref{ex:PP(PN)}.}
Set $\alpha=\measuredangle ABC$, $\beta=\measuredangle A'B'C'$.
Since $2\cdot\alpha\equiv 2\cdot \beta$, \ref{ex:2a=0} implies that
 $\alpha\equiv \beta$ or $\alpha\equiv \beta+\pi$.
In the latter case, the angles have opposite signs which is impossible.

Since $\alpha,\beta\in(-\pi,\pi]$, the equality $\alpha\equiv \beta$ implies $\alpha= \beta$.

\parbf{\ref{ex:between}.}
Apply \ref{thm:signs-of-triug} to $\triangle PAB$, $\triangle PBC$, and $\triangle PAC$.

For the second part, argue as in \ref{ex:a+b==c} to show that if three numbers $\alpha,\beta,\gamma\in (-\pi,\pi)$ have the same sign and
$\alpha+\beta\equiv\gamma$, then
$|\alpha|+|\beta|=|\gamma|$.
Apply it to $\alpha=\measuredangle APB$,
$\beta=\measuredangle BPC$,
and
$\gamma=\measuredangle APC$.

\parbf{\ref{ex:vert-intersect}.}
Note that $O$ and $A'$
lie on the same side of~$(AB)$.
Similarly, $O$ and $B'$
lie on the same side of~$(AB)$.
Hence the result.

\parbf{\ref{ex:signs-PXQ-PYQ}.}
Apply \ref{thm:signs-of-triug} for $\triangle PQX$ and $\triangle PQY$ and then 
apply \ref{cor:half-plane}\textit{\ref{cor:half-plane:angle}}.

\parbf{\ref{ex:chevinas}.} 
We can assume that $A'\z\ne B,C$ and $B'\ne A, C$;
otherwise, the statement trivially holds.

Note that $(BB')$ does not intersect $[A'C]$.
Applying Pasch's theorem (\ref{thm:pasch}) for $\triangle AA'C$ and $(BB')$, we get that 
$(BB')$ intersects $[AA']$; denote the point of intersection by $M$.

Тhe same way, we get that $(AA')$ intersects $[BB']$;
that is, $M$ lies on $[AA']$ and $[BB']$.

\parbf{\ref{ex:Z}.}
Assume that $Z$ is the point of intersection.

Note that $Z\ne P$ and $Z\ne Q$.
Therefore, $Z\notin (PQ)$.

Show that $Z$ and $X$ lie on one side of~$(PQ)$.
Repeat the argument to show that $Z$ and $Y$ lie on one side of~$(PQ)$.
It follows that $X$ and $Y$ lie on the same side of $(PQ)$ --- a contradiction.

\parbf{\ref{ex:intersecting-circles-3}.} The ``only-if'' part follows from the triangle inequality.
To prove the ``if'' part,  
observe that \ref{thm:abc} implies the existence of a triangle with sides $r_1$, $r_2$, and $d$.
Use this triangle to show that there is a point $X$ such that $O_1X=r_1$ and $O_2X=r_2$, where $O_1$ and $O_2$ are the centers of the corresponding circles.

\refstepcounter{chapter}
\setcounter{eqtn}{0}

\parbf{\ref{ex:equilateral}.}
Apply \ref{thm:isos} twice.

\parbf{\ref{ex:SMS}.} 
Let $D$ and $D'$ be reflections of
$C$ and $C'$ across $M$ and $M'$ respectively.
Show that $\triangle BCD\z\cong \triangle B'C'D'$ and use it to prove that $\triangle A' B' C'\z\cong\triangle A B C$.

\parbf{\ref{ex:isos-sides}.} \textit{(a)} Apply SAS.

\parit{(b)} Use \textit{(a)} and apply SSS.

\parbf{\ref{ex:ABC-motion}.}
Without loss of generality, we may assume that $X$ is distinct from $A$, $B$, and~$C$.
Set $f(X)=X'$; assume $X'\ne X$.

Note that $AX=AX'$, $BX=BX'$, and $CX=CX'$.
By SSS we get that $\measuredangle ABX\z=\pm\measuredangle ABX'$.
Since $X\ne X'$, we get that
$\measuredangle ABX\equiv - \measuredangle ABX'$.
In the same way, we get that 
$\measuredangle CBX\equiv - \measuredangle CBX'$.
Subtracting these two identities from each other, we get that
$\measuredangle ABC\equiv -\measuredangle ABC$.
Conclude that $\measuredangle ABC=0$ or $\pi$.
That is, $\triangle ABC$ is degenerate --- a contradiction. 

\parbf{\ref{ex:motion}.} Construct three circles with centers at $A'$, $B'$, and $C'$ and radii $AX$, $BX$, and $CX$, respectively.
The point $X'$ is their common point.


\refstepcounter{chapter}
\setcounter{eqtn}{0}

\parbf{\ref{ex:acute-obtuce}.} 
By Axiom~\ref{def:birkhoff-axioms:2b} and Theorem \ref{thm:straight-angle}, we have
$\measuredangle AOX\z+\measuredangle XOB\z\equiv\pi$; in particular, $|\measuredangle AOX\z+\measuredangle XOB|\ge\pi$.
Conclude that both values $\measuredangle AOX$ and $\measuredangle XOB$ are nonnegative or nonpositive.

In each case, argue as in \ref{ex:a+b==c} to show that $|\measuredangle AOX|\z+|\measuredangle XOB|=\pi$ and draw a conclusion.

\parbf{\ref{ex:pbisec-side}.}
Assume $X$ and $A$ lie on the same side of~$\ell$.

\begin{wrapfigure}{r}{27mm}
\vskip-4mm
\centering
\includegraphics{mppics/pic-322}
\end{wrapfigure}

Note that $A$ and $B$ lie on opposite sides of~$\ell$.
Therefore, by \ref{cor:half-plane}, 
$[AX]$ does not intersect $\ell$ 
and $[BX]$ intersects $\ell$;
let $Y$ be the intersection point.

By the triangle inequality, $BX=AY\z+YX\z\ge AX$.
Since $X\notin\ell$, by \ref{thm:perp-bisect} we have $AX\ne BX$.
Therefore $BX> AX$.

This way we have proved the ``if'' part.
To prove the ``only if'' part, you need to switch $A$ and $B$ and
repeat the above argument.

\parbf{\ref{ex:side-angle}.}
Apply \ref{ex:pbisec-side}, \ref{ex:vert-intersect}, and \ref{thm:isos}.

\parbf{\ref{ex:pbisec-motion}.} Show that $FX=FY$ and apply \ref{thm:perp-bisect}.

\begin{wrapfigure}{r}{27mm}
\vskip-4mm
\centering
\includegraphics{mppics/pic-326}
\end{wrapfigure}

\parbf{\ref{ex:construction-perpendicular}.}
See the picture; the given data is marked in bold.

\parbf{\ref{ex:2-reflections}.}
Note that
\begin{align*}
\measuredangle XBA&=\measuredangle ABP,
\\
\measuredangle PBC&=\measuredangle CBY.
\end{align*}
Therefore,
\begin{align*}
\measuredangle XBY
&\equiv
\measuredangle XBP+\measuredangle PBY\equiv
\\
&\equiv
 2\cdot(\measuredangle ABP+\measuredangle PBC)\equiv
\\
&
\equiv
 2\cdot \measuredangle ABC.
\end{align*}

\parbf{\ref{ex:3-reflections}.}
Choose an arbitrary nondegenerate triangle $ABC$.
Suppose that $\triangle \hat A \hat B\hat C$ denotes its image after the motion.

If $A\ne \hat A$, apply the reflection across the perpendicular bisector of~$[A\hat A]$.
This reflection sends $A$ to~$\hat A$.
Let $B'$ and $C'$ denote the reflections of $B$ and $C$ respectively.

If $B'\ne \hat B$, apply the reflection across the perpendicular bisector of~$[B'\hat B]$.
This reflection sends $B'$ to~$\hat B$.
Note that $\hat A\hat B=\hat AB'$;
that is, $\hat A$ lies on the perpendicular bisector. 
Therefore, $\hat A$ reflects to itself.
Suppose that $C''$ denotes the reflection of~$C'$.

Finally, if $C''\ne \hat C$, apply the reflection across $(\hat A\hat B)$.
Note that $\hat A\hat C\z=\hat AC''$ and $\hat B\hat C\z=\hat BC''$;
that is, $(\hat A\hat B)$ is the perpendicular bisector of $[C''\hat C]$.
Therefore, this reflection sends $C''$ to~$\hat C$.

Apply \ref{ex:ABC-motion} to show that the composition of the constructed reflections coincides with the given motion.

\parbf{\ref{ex:right-acute}.}
\textit{(a).} Apply \ref{lem:perp<oblique} and \ref{ex:side-angle}.

\parit{(b).} If the angle at $C$ is straight, then apply \ref{cor:degenerate=pi}.
Otherwise, show that there is a point $D\in [AB]$ such that $\angle ACD$ is right.
Applying \textit{(a)}, show that the angle at $A$ is acute. Reapeed it for~$B$.

\parbf{\ref{ex:obtuce}.}
We may assume that $P \ne X$; otherwise, the statement is trivial.

If $\angle PXV$ is right, then the statement follows from \ref{lem:perp<oblique}.
Otherwise, draw the line $(XY)$ perpendicular to~$(PX)$.
Show that $V$ and $W$ lie on opposite sides of $(XY)$.

Without loss of generality, assume that $P$ and $W$ lie on opposite sides of $(XY)$.
Let $W'$ be the point of intersection of $[PW]$ and $(XY)$.
Show and use that $PW>PW'>PX$.

\parbf{\ref{ex:PMQ}.}
Observe that $PM+MQ=P'M+MQ\z\ge P'Q$ and argue as in \ref{lem:perp<oblique}.

\parbf{\ref{ex:inside-outside}.}
Apply \ref{ex:obtuce} twice to $X$, $Y$, $P$, and the center.

{

\begin{wrapfigure}[6]{r}{27mm}
\vskip-6mm
\centering
\includegraphics{mppics/pic-328}
\end{wrapfigure}

\parbf{\ref{ex:chord-perp}.} Apply \ref{thm:perp-bisect}.

\parbf{\ref{ex:center}.}
See the picture; the given circle is marked in bold.

\parbf{\ref{ex:two-circ}.} Use \ref{ex:chord-perp} and \ref{perp:ex+un}.

}

\parbf{\ref{ex:tangent-circles}}; \textit{(\ref{ex:tangent-circles:a}).}
Assume that $Q$ is another point on both circles.
Show and use that $Q$ is the reflection of $P$ across~$(OO')$.

\parit{(\ref{ex:tangent-circles:b}).} Apply \textit{(\ref{ex:tangent-circles:a})} and \ref{lem:tangent}.

\parit{(\ref{ex:tangent-circles-2}).} Apply \textit{(\ref{ex:tangent-circles:a})} and \ref{cor:degenerate-trig}.

\parbf{\ref{ex:tangent-circles-3}.}
Let $A$ and $B$ be the points of intersection.
Note that the centers lie on the perpendicular bisector of the segment~$[AB]$.

\parbf{\ref{ex:tangent}+\ref{ex:tangent-circle}.}
The given data is marked in bold.

\begin{Figure}
\begin{minipage}{.49\textwidth}
\centering
\includegraphics{mppics/pic-330}
\end{minipage}
\hfill
\begin{minipage}{.49\textwidth}
\centering
\includegraphics{mppics/pic-332}
\end{minipage}
\end{Figure}

\refstepcounter{chapter}
\setcounter{eqtn}{0}

\parbf{\ref{ex:mid-triangle}.} Apply the SAS similarity condition to show that
$\triangle AC'B'\sim\triangle ABC$, $\triangle C'BA'\sim\triangle ABC$, and $\triangle B'A'C\sim\triangle ABC$.
Show that in each case, the similarity coefficient is $\tfrac12$.
Conclude that $A'B'=\tfrac12\cdot AB$, $B'C'=\tfrac12\cdot BC$, and $C'A'\z=\tfrac12\cdot CA$.
Apply the SSS similarity condition.

\parbf{\ref{ex:k*triangle}.} Apply the SAS similarity condition to show that
$\triangle OA'B'\sim\triangle OAB$,
$\triangle OB'C'\sim\triangle OBC$,
and $\triangle OC'A'\sim\triangle OCA$.
Show that in each case, the similarity coefficient is $k$.
Conclude that $A'B'=k\cdot AB$, $B'C'=k\cdot BC$, and $C'A'\z=k\cdot CA$.
Apply the SSS similarity condition.

\parbf{\ref{ex:angle-preserving-euclid}.}
By the AA similarity condition, the transformation multiplies the sides of any nondegenerate triangle by the same number; we need to show that this number does not depend on the triangle. 

Note that for any two nondegenerate triangles that share one side, this number is the same.
Apply this observation to a chain of triangles.

\parbf{\ref{ex:pyth}.}
Apply that $\triangle ADC\sim \triangle CDB$.

\parbf{\ref{ex:pyth-conv}.}
Apply the Pythagorean theorem (\ref{thm:pyth}) and the SSS congruence condition.

\parbf{\ref{ex:two-pairs-sim}.}
By the AA similarity condition (\ref{prop:sim}), $\triangle AYC\z\sim \triangle BXC$.
Conclude that 
$\frac{YC}{AC}=\frac{XC}{BC}$.
Apply the SAS similarity condition to show that $\triangle ABC\z\sim \triangle YXC$.

Similarly, apply AA and the equality of vertical angles to prove that $\triangle AZX\sim \triangle BZY$ and use SAS to show that $\triangle ABZ\z\sim \triangle YXZ$.

\parbf{\ref{ex:ABC+D}.}
Show and use that $\triangle ABC\sim \triangle CBD$.

\parbf{\ref{ex:right-perp-bi}.}
Show and use that $\triangle AQM\sim \triangle ABC\z\sim\triangle PBM$.


\refstepcounter{chapter}
\setcounter{eqtn}{0}

\parbf{\ref{ex:perp-perp}.}
Apply \ref{prop:perp-perp} to show that $k\parallel m$.
By \ref{cor:parallel-1}, $k\parallel n\z\Rightarrow m\parallel n$.
The latter contradicts that $m\perp n$.

\parbf{\ref{ex:construction-parallel}.}
Repeat the construction in \ref{ex:construction-perpendicular} twice.

\parbf{\ref{ex:parallel-angles}.}
Apply Axiom~\ref{def:birkhoff-axioms:2b}, \ref{thm:parallel-2}, and \ref{ex:2a=0}.

\begin{wrapfigure}{r}{23mm}
\vskip-0mm
\centering
\includegraphics{mppics/pic-333}
\end{wrapfigure}

\parbf{\ref{ex:smililar+parallel}.}
Since  $\ell\parallel (AC)$, it cannot cross $[AC]$.
By Pasch's theorem (\ref{thm:pasch}), $\ell$ has to cross another side of $\triangle ABC$.
Therefore $\ell$ crosses $[BC]$; denote the point of intersection by $Q$.

Use the transversal property (\ref{thm:parallel-2}) to show that $\measuredangle BAC= \measuredangle BPQ$.
The same argument shows that $\measuredangle ACB\z= \measuredangle PQB$; it remains to apply the AA similarity condition.

\parbf{\ref{ex:trisection}.}
Assume we need to trisect segment $[AB]$.
Construct a line $\ell\ne (AB)$ with four points $A,C_1,C_2, C_3$
such that $C_1$ and $C_2$ trisect $[AC_3]$.
Draw the line $(BC_3)$
and draw parallel lines thru $C_1$ and~$C_2$.
The points of intersections of these two lines with $(AB)$ trisect the segment $[AB]$.

\parbf{\ref{ex:|3sum|}.}
If $\triangle ABC$ is degenerate, then one of the angle measures is $\pi$, and the other two are~$0$.
Hence the result.

Assume $\triangle ABC$ is nondegenerate.
Set $\alpha\z=\measuredangle CAB$, $\beta=\measuredangle ABC$, and $\gamma=\measuredangle BCA$.

By \ref{thm:signs-of-triug},
we may assume that $0<\alpha,\beta,\gamma<\pi$.
Therefore, 
\[0<\alpha+\beta+\gamma<3\cdot\pi.
\eqlbl{eq:|3|<3pi}\]

By \ref{thm:3sum},
$$\alpha+\beta+\gamma\equiv\pi.\eqlbl{eq:|3|==pi}$$

From \ref{eq:|3|<3pi} and \ref{eq:|3|==pi} the result follows.

\parbf{\ref{ex:pent}.}
Apply \ref{thm:3sum} or \ref{ex:|3sum|} twice and \ref{thm:isos} trice.

\parbf{\ref{ex:right-isos}.}
Apply \ref{thm:isos} twice and \ref{thm:3sum} once. 





\parbf{\ref{ex:quadrangle}.}
Apply \ref{thm:3sum} to $\triangle ABC$ and $\triangle BDA$.

\parbf{\ref{ex:4parallels}.}
Denote by $M$ the center of symmetry of $\square ABCD$;
it exists by \ref{lem:parallelogram}.
Let $P'$ be the reflection of $P$ across $M$.
Show and use that $a=(AP')$, $b=(BP')$, $c=(CP')$, and $d=(DP')$.

\parbf{\ref{ex:romb}.}
Since $\triangle ABC$ is isosceles, $\measuredangle CAB\z=\measuredangle BCA$.
 
By SSS, $\triangle ABC\cong \triangle CDA$.
Therefore, 
$\pm\measuredangle DCA= \measuredangle BCA=\measuredangle CAB$.

Since $D\ne C$, we get ``$-$'' in the last formula.
Use the transversal property (\ref{thm:parallel-2}) to show that $(AB)\parallel (CD)$. Repeat the argument to show that $(AD)\z\parallel(BC)$ 

\parbf{\ref{ex:rectangle}.} 
By \ref{lem:parallelogram} and SSS, 
$AC=BD$
if and only if
$\measuredangle ABC=\pm \measuredangle BCD$.
By the transversal property~(\ref{thm:parallel-2}), 
$\measuredangle ABC+\measuredangle BCD\equiv \pi$.

Therefore, 
$AC=BD$
if and only if
$\measuredangle ABC
=\measuredangle BCD
=\pm\tfrac\pi2$.

\parbf{\ref{ex:romb2}.} 
Fix a parallelogram $ABCD$.
By \ref{lem:parallelogram},
its diagonals $[AC]$ and $[BD]$ have a common midpoint; denote it by~$M$.

Use SSS and \ref{lem:parallelogram} to show the following:
$AB=CD\ \iff\ \triangle AMB
\cong
\triangle AMD\ \z\iff\ \measuredangle AMB
=
\pm\tfrac\pi2$.

\parbf{\ref{ex:inscribed-rhombus}.}
Show and use that $\triangle CZY\sim\triangle YXB$.

\parbf{\ref{ex:coordinates}.} \textit{(a).} Use the uniqueness of the parallel line (\ref{thm:parallel}).

\parit{(b)} Use \ref{lem:parallelogram} and the Pythagorean theorem (\ref{thm:pyth}).

\parbf{\ref{ex:abc}.}
Set $A=(0,0)$, $B=(c,0)$, and $C=(x,y)$.
Clearly, $AB=c$,
$AC^2=x^2+y^2$ and $BC^2\z=(c-x)^2+y^2$.

It remains to show that there is a pair of real numbers $(x,y)$ 
that satisfy the following system of equations:
$$
\left\{
\begin{aligned}
b^2&=x^2+y^2
\\
a^2&=(c-x)^2+y^2
\end{aligned}
\right.
$$
if $0<a\le b\le c\le a+c$.

\parbf{\ref{ex:line-coord}.}
Notice that $MA=MB$ if and only if
\[(x-x_A)^2+(y-y_A)^2=(x-x_B)^2+(y-y_B)^2,\]
where $M=(x,y)$. 
To prove the first part, simplify this equation.
For the remaining parts use that any line is a perpendicular bisector of some line segment.

\parbf{\ref{ex:circle-coord}.} Rewrite it the following way and think 
\[(x+\tfrac a2)^2+(y+\tfrac b2)^2=(\tfrac a2)^2+(\tfrac b2)^2-c.\]

\parbf{\ref{ex:apolonnius}.}
We can choose the coordinates so that $B=(0,0)$ and $A=(a,0)$ for some $a>0$.
If $M=(x,y)$, then the equation $AM=k\cdot BM$ can be written in coordinates as 
\[k^2\cdot(x^2+y^2)=(x-a)^2+y^2.\]
It remains to rewrite this equation as in \ref{ex:circle-coord}.

\parbf{\ref{ex:apolonnius-construction}.}
Assume $M\notin(AB)$.
Show and use that the points $P$ and $Q$ constructed in the following picture lie on the Apollonian circle.

\begin{Figure}
\centering
\includegraphics{mppics/pic-335}
\end{Figure}

\refstepcounter{chapter}
\setcounter{eqtn}{0}

\parbf{\ref{ex:unique-cline}.}
Apply \ref{thm:circumcenter} and \ref{thm:perp-bisect}.

\parbf{\ref{ex:orthic-4}.}
Note that $(AC)\perp (BH)$ and $(BC)\z\perp (AH)$ and apply \ref{thm:orthocenter}.

(Note that each of $A,B,C,H$ is the orthocenter of the remaining three; such a quadruple of points $A,B,C,H$ is called an \index{orthocentric system}\emph{orthocentric system}.)

\parbf{\ref{ex:orthic-sim}.}
Show that $\triangle AA'C \sim \triangle BB'C$.
Use this to show that $\triangle ABC \sim \triangle A'B'C$.
Repeat the argument for $\triangle AB'C'$ and $\triangle A'BC'$.

\parbf{\ref{ex:midle}.}
Use the idea from the proof of \ref{thm:centroid}
to show that $(XY)\z\parallel (AC)\z\parallel (VW)$ and
$(XW)\z\parallel (BD)\z\parallel (YV)$.

\parbf{\ref{ex:euler-line}.}
Read about homothety.

Applying \ref{thm:centroid}, show that the homothety centered at $M$ with a coefficient of $-\tfrac{1}{2}$ maps $\triangle ABC$ onto its medial triangle $\triangle A'B'C'$.
Show that the homothety also maps the orthocenter of $\triangle ABC$ to the orthocenter of $\triangle A'B'C'$.

Arguing as in \ref{thm:orthocenter}, show that the orthocenter of $\triangle A'B'C'$ coincides with the circumcenter of $\triangle ABC$.
Make a conclusion.

\parbf{\ref{ex:perp-bisectors}.}
Let $(BX)$ and $(BY)$ be the internal and external bisectors of $\angle ABC$.
Then 
\begin{align*}
2\cdot \measuredangle XBY&\equiv2\cdot \measuredangle XBA+2\cdot \measuredangle ABY\equiv
\\
&\equiv
\measuredangle CBA+\pi+2\cdot \measuredangle ABC\equiv
\\
&\equiv\pi+\measuredangle CBC=\pi
\end{align*}
and hence the result.

\parbf{\ref{ex:bisect=altitude}.}
Apply ASA to the two triangles that the bisector cuts from the original triangle. 

\parbf{\ref{ex:ext-disect}.} 
If $E$ is the point of intersection of $(BC)$ 
with the external bisector of $\angle BAC$, then 
$\frac{AB}{AC}\z=\frac{EB}{EC}$.
It can be proved along the same lines as \ref{lem:bisect-ratio}.

\parbf{\ref{ex:bisect=median}.}
Apply \ref{lem:bisect-ratio}.
See also the solution of \ref{ex:abs-bisect=median}.

\parbf{\ref{ex:bisector-parallel}.}
Apply \ref{thm:isos}, \ref{thm:parallel-2}, and \ref{lem:parallelogram}.

\parbf{\ref{ex:2x=b+c-a}.}
Let $I$ be the incenter.
By SAS, we get that $\triangle AIZ\z\cong\triangle AIY$.
Therefore, $AY=AZ$.
In the same way, we get that $BX=BZ$ and $CX=CY$.
Hence the result.

\parbf{\ref{ex:orthic-triangle}.}
Argue as in \ref{ex:orthic-sim} to show that $\triangle A'B'C\z\sim \triangle ABC\z\sim\triangle AB'C'$.
Conclude that $(BB')$ bisects $\angle A'B'C'$.

If $\triangle ABC$ is obtuse, then its orthocenter coincides with one of the \index{excenter}\emph{excenters} of $\triangle A'B'C'$;
that is, 
the point of intersection of two external and one internal bisectors of $\triangle A'B'C'$.

\parbf{\ref{ex:bisector-incenter}.}
Apply \ref{lem:bisect-ratio} twice.
Use the obtained identity to show that the angle bisector at $C$ passes thru $I$.

\refstepcounter{chapter}
\setcounter{eqtn}{0}

\parbf{\ref{ex:inscribed-angle}.} \textit{(a).}
Apply \ref{thm:inscribed-angle} for $\angle XX'Y$ and $\angle X'YY'$
and \ref{thm:3sum} for $\triangle PYX'$.

\parit{(b)} If $P$ is inside of $\Gamma$, then $P$ lies between $X$ and $X'$ and between $Y$ and $Y'$.
In this case, $\angle XPY$ is vertical to $\angle X'PY'$.

If $P$ is outside of $\Gamma$ then $[PX)\z=[PX')$ and $[PY)=[PY')$.
In both cases we have that $\measuredangle XPY=\measuredangle X'PY'$.

Applying \ref{thm:inscribed-angle} and \ref{ex:ABCO-line}, we get
\pagebreak[0]
\begin{align*}
2\cdot \measuredangle Y'X'P
&\equiv
2\cdot \measuredangle Y'X'X\equiv 
\\
\equiv
2\cdot\measuredangle Y'YX
&\equiv
2\cdot\measuredangle PYX.
\end{align*}
According to \ref{thm:signs-of-triug}, $\angle Y'X'P$ and $\angle PYX$ have the same sign;
therefore
$\measuredangle Y'X'P\z= \measuredangle PYX$.
It remains to apply the AA similarity condition.

\parit{(c)} Apply \textit{(b)} assuming $[YY']$ is the diameter of~$\Gamma$. 

\parbf{\ref{ex:inscribed-hex}.} Apply \ref{ex:inscribed-angle}\textit{\ref{ex:inscribed-angle:b}}
thrice.

\parbf{\ref{ex:altitudes-circumcircle}.}
Let $X$ and $Y$ be the footpoints of the altitudes from $A$ and~$B$.
Suppose that $O$ denotes the circumcenter.
 
By AA, $\triangle A X C\sim \triangle B Y C$.
Thus 
\begin{align*}
\measuredangle A'OC
&\equiv 
2\cdot \measuredangle A' A C
\equiv
\\
&\equiv-2\cdot\measuredangle B' B C
\equiv
\\
&\equiv-\measuredangle B'OC.
\end{align*}

By SAS, $\triangle A'OC\cong\triangle B'OC$.
Therefore, $A'C=B'C$.

\parbf{\ref{ex:two-chords}.}
Apply the transversal property (\ref{thm:parallel-2}) and the theorem on inscribed angles (\ref{thm:inscribed-angle}).

Alternatively, apply \ref{cor:reflection+angle} to the reflection across the common perpendicular to the lines passing thru the center of the circle.

\begin{wrapfigure}[9]{r}{20mm}
\vskip-2mm
\centering
\includegraphics{mppics/pic-336}
\end{wrapfigure}

\parbf{\ref{ex:tangent-construction-inscribed}.} Apply \ref{cor:right-angle-diameter} and \ref{lem:tangent}.

\parbf{\ref{ex:perp-construction-inscribed}.} Apply \ref{cor:right-angle-diameter}.

\parbf{\ref{ex:altitude+circles}.}
Show and use that $\angle ABX$ and $\angle ABY$ are right.
(We may assume that the points $A$, $B$, $X$, and $Y$ are distinct; otherwise, there is nothing to prove.)

\parbf{\ref{ex:perpendicular-ruler}.}
Guess the construction from the picture.
To prove it,
apply \ref{thm:orthocenter} and \ref{cor:right-angle-diameter}.

\parbf{\ref{ex:tnagents+midpoint}.}
Let $O$ be the center of $\Gamma$.
Use \ref{cor:right-angle-diameter} to show that the points lie on the circle with diameter $[PO]$.

\parbf{\ref{ex:VVAA}.} 
Note that $\measuredangle AA'B\z=\pm\tfrac\pi2$ and $\measuredangle AB'B\z=\pm\tfrac\pi2$.
Then apply \ref{cor:inscribed-quadrangle}
to $\square AA'BB'$.

\parbf{\ref{ex:secant-circles}.}
Apply \ref{cor:inscribed-quadrangle} twice for $\square ABYX$ and $\square ABY'X'$ and use the transversal property (\ref{thm:parallel-2}).

\parbf{\ref{ex:perim+angle+side}.}
Construct $\triangle AXC$ such that $AC=b$, $AX=p-b$, and $\measuredangle AXC=\tfrac12\cdot \beta$.
Note that point $B$ at the intersection of $AX$ and the perpendicular bisector of $[CX]$ solves the problem.

\parbf{\ref{ex:inaccuracy}.}
One needs to show that $(A'B') \nparallel (XP)$; otherwise, the first line in the proof does not make sense.
By the transversal property, this means that
\[
2 \cdot \measuredangle XPA' + 2 \cdot \measuredangle PA'B' \not\equiv 0.
\eqlbl{eq:2<B'A'P+2<A'PX}
\]
This should be done later in the proof, right before we use point $Y$.

Instead of $2 \cdot \measuredangle AXY \equiv 2 \cdot \measuredangle AA'Y$, we should write
\[
2 \cdot \measuredangle AXP \equiv 2 \cdot \measuredangle PA'B'.
\eqlbl{eq:2<AXP=2<AA'B'}
\]
Since $\triangle XAP$ has a right angle at $A$, we have
\[
\measuredangle AXP + \measuredangle XPA \equiv \pm \tfrac{\pi}{2}.
\]
Since $2 \cdot \measuredangle XPA' \equiv 2 \cdot \measuredangle XPA$, \ref{eq:2<AXP=2<AA'B'} implies that
\[
2 \cdot \measuredangle XPA' + 2 \cdot \measuredangle PA'B' \equiv \pi,
\]
and \ref{eq:2<B'A'P+2<A'PX} follows.

Additionally, we implicitly use the following identities:
\begin{align*}
2 \cdot \measuredangle AXP &\equiv 2 \cdot \measuredangle AXY, \\
2 \cdot \measuredangle ABP &\equiv 2 \cdot \measuredangle ABB', \\
2 \cdot \measuredangle AA'B' &\equiv 2 \cdot \measuredangle AA'Y.
\end{align*}

\parbf{\ref{ex:equilateral-2}.}
By \ref{cor:right-angle-diameter},
the points $L$, $M$, and $N$ lie on the circle $\Gamma$ with diameter~$[OX]$.
It remains to apply \ref{thm:inscribed-angle} for the circle $\Gamma$ 
and two inscribed angles with a common vertex at~$O$.

\parbf{\ref{ex:median-angle}.}
Observe that the points $A$, $A'$, $B$, and $B'$ line on one circle.
Show that $(AB)\parallel (B'A')$.
Applying \ref{ex:two-chords}, show that $AB'=BA'$ and therefore $AC=BC$.

\parbf{\ref{ex:simson}.}
Let $X$, $Y$, and $Z$ denote the footpoints of $P$ on $(BC)$, $(CA)$, and $(AB)$ respectively.

\begin{Figure}
\centering
\includegraphics{mppics/pic-338}
\end{Figure}

Show that $\square AZPY$, $\square BXPZ$, $\square CYPX$, and $\square ABCP$ are inscribed.
Use this to show that
\begin{align*}
2\cdot \measuredangle CXY&\equiv 2\cdot \measuredangle CPY,
&
2\cdot \measuredangle BXZ&\equiv 2\cdot \measuredangle BPZ,
\\
2\cdot \measuredangle YAZ&\equiv 2\cdot \measuredangle YPZ,
&
2\cdot \measuredangle CAB&\equiv 2\cdot \measuredangle CPB.
\end{align*}

Conclude that 
$2\cdot \measuredangle CXY\equiv 2\cdot \measuredangle BXZ$
and hence $X$, $Y$, and $Z$ lie on one line.


\parbf{\ref{ex:a+b=c}.}
Show that $P$ lies on the arc opposite from $ACB$;
conclude that
$\measuredangle APC\z=\measuredangle CPB\z=\pm\tfrac\pi3$.

{

\begin{wrapfigure}[8]{r}{25mm}
\vskip-3mm
\centering
\includegraphics{mppics/pic-337}
\end{wrapfigure}

Choose a point $A'\z\in [PC]$ such that $PA'\z=PA$.
Note that $\triangle APA'$ is equilateral.
Prove and use that $\triangle AA'C\z\cong \triangle APB$.

\parbf{\ref{ex:tangent-arc}.}
If $C\in (AX)$, then the arc is the line segment $[AC]$ or the union of two half-lines in $(AX)$ with vertices at $A$ and~$C$.

}

Assume $C\notin (AX)$.
Let $\ell$ be the perpendicular line dropped from $A$ to $(AX)$ and $m$ be the perpendicular bisector of~$[AC]$.

Note that $\ell\nparallel m$;
suppose they meet at $O$.
Note that the circle with center $O$ passing thru $A$ is also passing thru $C$ and tangent to~$(AX)$.

Note that one of the two arcs with endpoints $A$ and $C$ is tangent to~$[AX)$.

The uniqueness follows from \ref{prop:arc(angle=tan)}.

\parbf{\ref{ex:tangent-lim}.} Use \ref{prop:arc(angle=tan)} and \ref{thm:3sum} to show that 
$\measuredangle XAY\z=\measuredangle ACY$.
By Axiom~\ref{def:birkhoff-axioms:2c}, $\measuredangle ACY\to 0$ as $AY\to 0$;
hence the result.

\begin{wrapfigure}[7]{r}{25mm}
\vskip-6mm
\centering
\includegraphics{mppics/pic-340}
\end{wrapfigure}

\parbf{\ref{ex:two-arcs}.} 
Apply \ref{prop:arc(angle=tan)} twice.

(Alternatively, apply \ref{cor:reflection+angle} for the reflection across the perpendicular bisector of $[AC]$.)

\parbf{\ref{ex:3x120}.} Guess a construction from the picture.
To show that it produces the needed point, apply~\ref{thm:inscribed-angle}.

\refstepcounter{chapter}
\setcounter{eqtn}{0}

\parbf{\ref{ex:constr-inversion}.}
By \ref{lem:tangent}, $\angle OTP'$ is right. 
Therefore, $\triangle OPT\z\sim \triangle OTP'$
and in particular
$OP\cdot OP'\z=OT^2$,
and hence the result.

\parbf{\ref{ex:appolo-circ}.}
Suppose that $O$ denotes the center of $\Gamma$.
Assume that $X,Y\z\in \Gamma$;
in particular, $OX\z=OY$.

Note that the inversion sends $X$ and $Y$ to themselves.
By \ref{lem:inversion-sim},
$$\triangle OPX\z\sim \triangle OXP'
\quad
\text{and}
\quad
\triangle OPY\sim \triangle OYP'.$$
Therefore, 
$\frac{PX}{P'X}=\frac{OP}{OX}=\frac{OP}{OY}=\frac{PY}{P'Y}$
and hence the result.

\parbf{\ref{ex:incenter+inversion=orthocenter}.}
By \ref{lem:inversion-sim},
\begin{align*}
\measuredangle IA'B'&\equiv -\measuredangle IBA,
&
\measuredangle IB'A'&\equiv -\measuredangle IAB,
\\
\measuredangle IB'C'&\equiv -\measuredangle ICB,
&
\measuredangle IC'B'&\equiv -\measuredangle IBC,
\\
\measuredangle IC'A'&\equiv -\measuredangle IAC,
&
\measuredangle IA'C'&\equiv -\measuredangle ICA.
\end{align*}

It remains to apply the theorem on the sum of angles of triangle (\ref{thm:3sum})
to show that $(A'I)\z\perp (B'C')$, 
$(B'I)\z\perp (C'A')$
and
$(C'I)\z\perp (B'A')$.

\parbf{\ref{ex:consturuction-of-inversion}.}
Guess the construction from the picture.

\begin{Figure}
\vskip-0mm
\centering
\includegraphics{mppics/pic-342}
\end{Figure}

\parbf{\ref{ex:inv-center not=center-inv}.}
Assume $r\z>0$, $x\ne 0$, and $y\ne 0$.
Show that
$$\frac{r^2}{(x+y)/2}
=
\left(\frac {r^2}x+\frac {r^2}y\right)/2\ \iff\ x=y.$$

Suppose $\ell$ denotes the line passing thru $Q$, $Q'$, and the center of the inversion $O$.
Choose an isometry $\ell\to\mathbb{R}$ that sends $O$ to $0$;
assume $x,y\in \mathbb{R}$ are the values of $\ell$ for the two points in $\ell\cap\Gamma$;
note that $x\ne y$.
Assume $r$ is the radius of the circle of inversion.
Then the left-hand side above is the coordinate of $Q'$ 
and the right-hand side is the coordinate of the center of $\Gamma'$.

\parbf{\ref{ex:circumtool}.}
A solution is given in \ref{sec:verification}.

\parbf{\ref{ex:tangent-circ->parallels}.}
Apply an inversion across a circle with the center at the only point of intersection of the circles;
then use \ref{thm:inverse}.

\parbf{\ref{ex:4-circles}.}
Label the points of tangency as in the picture.
Apply an inversion with the center at $P$. 
Observe that the two circles that tangent are at $P$ become parallel lines and 
the remaining two circles are tangent to each other and these two parallel lines.

\begin{Figure}
\centering
\includegraphics{mppics/pic-344}
\end{Figure}

Note that the points of tangency $A'$, $B'$, $X'$, and $Y'$ with the parallel lines are vertices of a square;
in particular, they lie on one circle.
These points are images of $A$, $B$, $X$, and $Y$ under the inversion.
By \ref{thm:inverse-cline}, the points $A$, $B$, $X$, and $Y$ also lie on one circline.

\parbf{\ref{ex:inverse}.} 
Apply the inversion across a circle with center~$A$. 
Point $A$ will go to infinity; the two circles tangent at $A$ will become parallel lines
while the two parallel lines will become circles tangent at~$A$.

\begin{Figure}
\centering
\includegraphics{mppics/pic-346}
\end{Figure}

It remains to show that the dashed line $(AB')$ in the picture is parallel to the other two lines.

\parbf{\ref{ex:inscribed+inv}.}
Apply \ref{lem:inverse-4-angle}\textit{\ref{lem:inverse-4-angle:angle}}, 
\ref{ex:quadrangle},
and \ref{thm:inscribed-angle}.

\parbf{\ref{ex:centers-of-perp-circles}.}
Suppose that $T\in \Omega_1\cap\Omega_2$.
Let $P$ be the footpoint of $T$ on~$(O_1O_2)$.
Show that
$\triangle O_1PT
\z\sim \triangle O_1TO_2
\z\sim \triangle TPO_2$.
Conclude that $P$ coincides with the inverses of $O_1$ across $\Omega_2$ and of $O_2$ across~$\Omega_1$.

\parbf{\ref{ex:4-th-perp-circ}.}
Since $\Gamma\perp\Omega_1$ and $\Gamma\perp\Omega_2$,
Corollary~\ref{cor:perp-inverse-clines} 
implies that
the circles $\Omega_1$ and $\Omega_2$ are inverted in $\Gamma$ 
to themselves.
Conclude that $A$ and $B$ are inverses of each other.

Since $\Omega_3\ni A,B$,
Corollary~\ref{cor:perp-inverse} implies that
$\Omega_3\perp \Gamma$.

\parbf{\ref{ex:construction-perp-clines}.}
Let $P_1$ and $P_2$ be the inverses of $P$ 
across $\Omega_1$ and~$\Omega_2$.
Apply \ref{cor:perp-inverse} and \ref{thm:perp-inverse}
to show that a circline $\Gamma$ passing thru $P$, $P_1$, and $P_3$ is a solution.

\parbf{\ref{ex:3-construction-perp-clines}.}
All circles perpendicular to $\Omega_1$ and $\Omega_2$ pass thru a fixed point~$P$.
Try to construct~$P$.

If two of the circles intersect, try to apply \ref{cor:invese-comp}.

\refstepcounter{chapter}
\setcounter{eqtn}{0}

\parbf{\ref{ex:tangent-angle-neutral}.}
One cannot apply \ref{thm:3sum} (its proof relies on \ref{thm:parallel-2}, which in turn uses \ref{thm:parallel}, and that depends on the SAS similarity condition, which is essentially Axiom~\ref{def:birkhoff-axioms:4}).

\parbf{\ref{ex:abs-bisect=median}.}
Suppose $D$ is the midpoint of~$[BC]$.
Assume that $(AD)$ is the angle bisector at~$A$.

Let $A'\in [AD)$ be the reflection of $A$ across~$D$.
Note that $\triangle CAD\z\cong\triangle BA'D$.
In particular, $\measuredangle BAA'=\measuredangle AA'B$.
Apply \ref{thm:isos} to  $\triangle ABA'$.

\parbf{\ref{ex:abs-inscibed}.}
The statement is evident if $A$, $B$, $C$, and $D$ lie on one line.

In the remaining case, suppose that $O$ denotes the circumcenter.
Apply \ref{thm:isos} to
$\triangle AOB$,
$\triangle BOC$, 
$\triangle COD$, 
and
$\triangle DOA$. 

\textit{(In the Euclidean plane the statement follows from \ref{cor:inscribed-quadrangle} and \ref{ex:quadrangle},
but these statements cannot be used in the neutral plane.)}

\parbf{\ref{ex:parallel-abs}.}
Arguing by contradiction, 
assume 
$2\cdot(\measuredangle ABC+\measuredangle BCD)\equiv0$, 
but $(AB)\z\nparallel(C D)$.
Let $Z\in (AB)\cap(CD)$.

Note that 
$
2\cdot \measuredangle ABC\equiv 2\cdot \measuredangle ZBC,
$ and
$2\cdot \measuredangle BCD\equiv 2\cdot \measuredangle BCZ$.

Apply \ref{prop:2sum} to $\triangle ZBC$ and try to reach a contradiction.

\parbf{\ref{ex:SAA}.}
Choose $C''\in [B'C')$ such that $B'C''\z=BC$.

Use SAS to show that $\triangle ABC\cong \triangle A'B'C''$.
Conclude that $\measuredangle B'C'A'\z= \measuredangle B'C''A'$.

\begin{Figure}
\vskip-0mm
\centering
\includegraphics{mppics/pic-348}
\end{Figure}

Therefore, it is sufficient to show that $C''\z=C'$.
If $C'\z\ne C''$ apply \ref{prop:2sum} to $\triangle A'C'C''$ and try to reach at a contradiction.


\parbf{\ref{ex:chev<side}.} 
Use \ref{ex:side-angle} and \ref{prop:2sum}.

(Alternatively, follow the solution of \ref{ex:inside-outside}.)

\parbf{\ref{ex:neutral-quadrangle}.}
Set $a=AB$, $b=BC$, $c=CD$, and $d=DA$; we need to show that $c\ge a$.

Mimic the proof of \ref{thm:3sum-a} for the shown fence made from copies of quadrangle $ABCD$.
You should get that $a\le c+\tfrac dn$ for any integer $n\ge 1$, and this implies the needed inequality.

\begin{Figure}
\vskip-0mm
\centering
\includegraphics{mppics/pic-349}
\end{Figure}

\parit{Alternative way.}
By \ref{lem:perp<oblique}, we can assume that $\angle BCD$ is right.
Use \ref{thm:3sum-a} to show that $|\measuredangle CBD|\ge|\measuredangle ADB|$.
Use \ref{prop:angle-side} and the construction in the neutral proof of \ref{thm:hypotenuse-leg} to show that $2\cdot c\ge 2\cdot a$.

\parbf{\ref{ex:defect}.}
Note that 
$|\measuredangle ADC|+|\measuredangle CDB|=\pi$.
Then apply the definition of the defect.

\parbf{\ref{ex:defect=}.}
Show that $\triangle AMX\cong \triangle BMC$. 
Apply \ref{ex:defect} to $\triangle ABC$ and $\triangle AXC$.

\parbf{\ref{ex:neutral-rectangle}.}
Observe that the total sum of absolute values of angle measures in $\triangle ABC$ and $\triangle CDA$ is at least $2\cdot\pi$.
Apply \ref{thm:3sum-a} to show that
\[\defect(\triangle ABC)=\defect(\triangle CDA)=0.\]
Use it to show that $\measuredangle CAB=\measuredangle ACD$ and $\measuredangle ACB=\measuredangle CAD$.
By ASA, $\triangle ABC\cong\triangle CDA$, and, in particular, $AB=CD$.

(Note that it also follows from \ref{ex:neutral-quadrangle}.)

\parbf{\ref{ex:neutral-rectangle+}.}
Choose a triangle $\triangle_0$.
Show that $\triangle_0$ can be covered by a large rectangle $\square$ that is a union of several copies of the given rectangle.
Further, subdivide $\square$ into triangles $\triangle_0,\dots, \triangle_n$.
Show that the sum of the defects of $\triangle_i$ is zero and use that defects are nonnegative.

\begin{Figure}
\vskip-0mm
\centering
\includegraphics{mppics/pic-1114}
\end{Figure}

\parit{Remark.}
This is a key step in the proofs of parts \textit{(\ref{thm:=IV:defect})} and \textit{(\ref{thm:=IV:rectangle})} in Theorem~\ref{thm:=IV}.

\refstepcounter{chapter}
\setcounter{eqtn}{0}

\parbf{\ref{ex:ideal-line-unique}.} 
Consider the ideal points $A$ and $B$ be the ideal points of the h-line~$\ell$. 
Note that the center of the Euclidean circle containing $\ell$ lies 
at the intersection of the lines tangent to the absolute at the ideal points of~$\ell$.

\parbf{\ref{ex:1ideal-line-unique}.}
Assume $A$ is an ideal point of the h-line $\ell$
and $P\in \ell$.
Suppose that $P'$ denotes the inverse of $P$ across the absolute.
By \ref{cor:perp-inverse-clines},
$\ell$ lies in the intersection of the h-plane and the (necessarily unique) circline 
passing thru $P$, $A$, and~$P'$.

\parbf{\ref{ex:line/h-line}.} 
Let $\Omega$ and $O$ denote the absolute and its center. 

Let $\Gamma$ be the circline containing~$[PQ]_h$.
Note that $[PQ]_h=[PQ]$ if and only if $\Gamma$ is a line.

Suppose that $P'$ denotes the inverse of $P$ across~$\Omega$.
Note that $O$, $P$, and $P'$ lie on one line.

By the definition of an h-line, $\Omega\perp \Gamma$.
By \ref{cor:perp-inverse-clines}, $\Gamma$ passes thru $P$ and~$P'$. 
Therefore, $\Gamma$
is a line if and only if it passes thru~$O$.

\parbf{\ref{ex:h-dist-eq}.}
Assume that the absolute is a unit circle.

Set $a\z=OX\z=XY$.
Note that $0<a<\tfrac12$,
$
OX_h=\ln \tfrac{1+a}{1-a}$,
and
$XY_h=\ln \tfrac{(1+2\cdot a)\cdot(1-a)}{(1-2\cdot a)\cdot(1+a)}$.
Verify that the inequalities 
\[
1
<
\tfrac{1+a}{1-a}
<
\tfrac{(1+2\cdot a)\cdot(1-a)}{(1-2\cdot a)\cdot(1+a)}\]
hold when $0<a<\tfrac12$.

\parbf{\ref{ex:h-perp-unique}.} 
Spell the meaning of the terms ``perpendicular'' and ``h-line'' and then apply \ref{ex:construction-perp-clines}.

\parbf{\ref{ex:h-circle=circle}.}
Choose the vertices $P$, $Q$, and $R$ on a Euclidean circle that intersects the absolute and is not orthogonal to it.
Apply~\ref{lem:h-circle=circle}.

\begin{wrapfigure}[9]{r}{34mm}
\vskip-5mm
\centering
\includegraphics{mppics/pic-350}
\end{wrapfigure}

\parbf{\ref{ex:3-h-lines}.}
Choose the required h-lines from the picture.

\parbf{\ref{ex:O-h-dist}.} Use \ref{lem:O-h-dist}.

\parbf{\ref{ex:cosh}.}
By \ref{cor:invese-comp} and \ref{lem:inverse-4-angle},
the right-hand sides of the identities 
remain unchanged under an inversion across a circle perpendicular to the absolute.

As usual, we assume that the absolute is a unit circle.
Let $O$ be the h-midpoint of $[PQ]_h$.
By the main observation (\ref{thm:main-observ})
we can assume that $O$ is the center of the absolute.
In this case, $O$ is also the Euclidean midpoint of $[PQ]$.%

Set $a=OP=OQ$; in this case, we have
\begin{align*}
PQ&=2\cdot a,
&
PP'=QQ'&=\tfrac1a-a,
\\
P'Q'&=2\cdot \tfrac1a,
&
PQ'=QP'&=\tfrac1a+a.
\end{align*}
and 
\[PQ_h=\ln \tfrac{(1+a)^2}{(1-a)^2}=2\cdot \ln \tfrac{1+a}{1-a}.\]
Therefore
\begin{align*}
\cosh[\tfrac12\cdot PQ_h]
&=\tfrac12\cdot(\tfrac{1+a}{1-a}+\tfrac{1-a}{1+a})=
\\
&=\tfrac{1+a^2}{1-a^2};
\end{align*} 
\begin{align*}
\sqrt{\frac{PQ'\cdot P'Q}{PP'\cdot QQ'}}
&=\frac{\frac1a+a}{\frac1a-a}=
\\
&=\tfrac{1+a^2}{1-a^2}.
\end{align*} 
Hence part \textit{(\ref{ex:cosh/2})} follows.
Similarly,
\begin{align*}
\sinh[\tfrac12\cdot PQ_h]
&=\tfrac12\cdot\left(\tfrac{1+a}{1-a}-\tfrac{1-a}{1+a}\right)=
\\
&=\tfrac{2\cdot a}{1-a^2};
\\
\sqrt{\frac{PQ\cdot P'Q'}{PP'\cdot QQ'}}
&=\frac{2}{\frac1a-a}=
\\
&=\tfrac{2\cdot a}{1-a^2}.
\end{align*} 
Hence part \textit{(\ref{ex:coshsinh})} follows.

Parts \textit{(\ref{ex:coshtanh})} and \textit{(\ref{ex:coshcosh})} follow from \textit{(\ref{ex:cosh/2})}, \textit{(\ref{ex:coshsinh})}, the definition of a hyperbolic tangent, and the double-argument identity for hyperbolic cosine, see \ref{double-argument}.

(We could also move $Q$ to the center of absolute.
In this case, the derivations are simpler.
But since $Q'Q=Q'P=Q'P'=\infty$, one has to justify that $\tfrac\infty\infty=1$ every time.)

\refstepcounter{chapter}
\setcounter{eqtn}{0}

\parbf{\ref{ex:lambert-parallelism}.}
Use \ref{prop:perp-perp} to show that $(AB)_h\parallel (CD)_h$.
Apply the definition of the angle of parallelism.

\parbf{\ref{ex:ultra-parallel}}; \textit{``only-if'' part.}
Suppose $\ell$ and $m$ are ultraparallel; that is, they do not intersect and have distinct ideal points.

Denote the ideal points by $A$, $B$, $C$, and $D$;
we may assume that they appear on the absolute in the same order.
In this case, the h-line with ideal points $A$ and $C$ intersects the h-line with ideal points $B$ and $D$.
Let $O$ be their point of intersection.
By \ref{thm:main-observ}, we can assume that $O$ is the center of absolute.
Note that $\ell$ is the reflection of $m$ across $O$ in the Euclidean sense.

\begin{Figure}
\vskip-0mm
\centering
\includegraphics{mppics/pic-359}
\end{Figure}

Drop an h-perpendicular $n$ from $O$ to $\ell$ and
show that $n\perp m$.

\parit{``If'' part.} 
Suppose $n$ is a common perpendicular.
Let $L$ and $M$ be its points of intersection with $\ell$ and $m$ respectively.

Let $O$ be the center of the absolute.
By \ref{thm:main-observ}, we can assume that $O$ is the h-midpoint of $L$ and $M$.
Notice that in this case $\ell$ is the reflection of $m$ across $O$ in the Euclidean sense.
It follows that the ideal points of the h-lines $\ell$ and $m$ are symmetric to each other.
Therefore, if one pair of them coincides, then the other pair does too. 
By \ref{ex:ideal-line-unique}, $\ell=m$, which contradicts the assumption $\ell\ne m$.

\parbf{\ref{ex:right-angle-parallelism}.} Show that the angle of parallelism at $C$ with respect to $(AB)_h$ is less than $\tfrac\pi4$, and apply \ref{prop:angle-parallelism}.
You can use approximations such as $\sqrt2\approx 1.414$ and $e\approx2.718$.

\parbf{\ref{ex:small-angle}.}
By the triangle inequality, the h-distance from $B$ to $(AC)_h$ is at least 50.
By \ref{prop:angle-parallelism}, $|\measuredangle_h ABC|<2\cdot\phi$, where $\phi$ is the angle of parallelism at the distance $50$.
Show that $\cos\phi\z=\tfrac{e^{100}-1}{e^{100}+1}$ and
use it to estimate~$\phi$.
You may use that $\cos\phi\z\le 1-\tfrac1{10}\cdot\phi^2$ for $|\phi|\z<\tfrac\pi2$ and $e^3>10$.

\parbf{\ref{ex:side-sup}.}
Note that the angle of parallelism at $B$ with respect to $(CD)_h$ is greater than $\tfrac\pi4$,
and it approaches to $\tfrac\pi4$ as $CD_h\to\infty$.

Applying \ref{prop:angle-parallelism},
we get that
$$BC_h<\tfrac12\cdot\ln\frac{1+\frac1{\sqrt{2}}}{1-\frac1{\sqrt{2}}}=\ln\left(1+\sqrt{2}\right).$$

The right-hand side is the limit of $BC_h$ as $CD_h\to\infty$.
Therefore, $\ln\left(1+\sqrt{2}\right)$ is the optimal upper bound.

\parbf{\ref{ex:equidistant-reflection}.}
Let $P'$ be another h-point.
Denote by $R$ and $R'$ the h-footpoint of $P$ and $P'$ on $m$.

Show and use that $P'$ is an h-reflection of $P$ across $Q\in m$ if and only if $P$ and $P'$ lie on the opposite sides of $m$, 
$PR_h=P'R'_h$, and $Q$ is the h-midpoint of $[RR']_h$. 
You may need \ref{prop:vert}, SAS, and SAA (see \ref{ex:SAA}).

\parbf{\ref{ex:right-trig-horocycle}.}
As usual, we assume that the absolute is a unit circle. 

Consider a hyperbolic triangle $PQR$
with a right angle at $Q$, where  $PQ_h\z=QR_h$
and the vertices $P$, $Q$, and $R$ 
lie on a horocycle.

\begin{wrapfigure}{r}{32mm}
\vskip-2mm
\centering
\includegraphics{mppics/pic-352}
\end{wrapfigure}

We may assume that $Q$ is the center of the absolute.
In this case, $\measuredangle_hPQR\z=\measuredangle PQR\z=\pm\tfrac\pi2$ and $PQ=QR$.

Note that the Euclidean circle passing thru $P$, $Q$, and $R$ is tangent to the absolute.
Conclude that $PQ=\tfrac1{\sqrt2}$. 
Apply \ref{lem:O-h-dist} to find $PQ_h$.

\parbf{\ref{ex:angle-preserving-hyp}.}
Apply the AAA-congruence condition (\ref{thm:AAA}).

\parbf{\ref{ex:circum}.}
Apply \ref{prop:circum}.
Use that the function $r\mapsto e^{-r}$ is decreasing and $e\z>2$.

\parbf{\ref{ex:c+1>a+b}.}
Apply the hyperbolic Pythagorean theorem and the definition of a hyperbolic cosine.
The following observations should help:
\begin{itemize}
 \item The function $x\mapsto e^x$ is an increasing positive function.
 \item By the triangle inequality,  we have
 \[-c\le a-b\quad \text{and}\quad  -c\le b-a.\]
\end{itemize}

\refstepcounter{chapter}
\setcounter{eqtn}{0}

\parbf{\ref{ex:affine-par}.}
Assume the distinct lines $\ell$ and $m$ 
are mapped to intersecting lines $\ell'$ and~$m'$.
Suppose $P'$ denotes their point of intersection.

Let $P$ be the inverse image of~$P'$.
By the definition of an affine map, $P$ must lie on both $\ell$ and $m$.
Make a concussion.

\parbf{\ref{ex:afine-linear}.}
In each case, check that the map is a bijection and apply \ref{ex:line-coord}.

\parbf{\ref{ex:collinear=affine}.}
Choose a line $(AB)$.

Assume $X'\z\in (A'B')$ for some $X\z\notin(AB)$.
Since $P\mapsto P'$ maps collinear points to collinear, 
the three lines $(AB)$, $(AX)$, and $(BX)$ are mapped to~$(A'B')$.
Furthermore, any line connecting a pair of points on these three lines is also mapped to~$(A'B')$.
Use it to show that the entire plane is mapped to $(A'B')$.
The latter contradicts that the map is a bijection.

{

\begin{wrapfigure}{r}{22mm}
\vskip-3mm
\centering
\includegraphics{mppics/pic-353}
\end{wrapfigure}

By the assumption, if $X\z\in (AB)$, then $X'\z\in (A'B')$.
From above, if $X\notin (AB)$, then $X'\z\notin (A'B')$.
Use it to prove the second statement.

}

\parbf{\ref{ex:circle=affine}.}
Observe that $\alpha$ maps noncollinear triples of points to noncollinear ones.
Therefore, $\alpha^{-1}$ maps collinear triples to collinear ones.
It remains to apply \ref{ex:collinear=affine}.

{

\begin{wrapfigure}{r}{28mm}
\vskip-6mm
\centering
\includegraphics{mppics/pic-354}
\end{wrapfigure}

\parbf{\ref{ex:midpoint-affine}.}
It is sufficient to construct the midpoint of $[AB]$
with a ruler and a parallel tool.
Guess the construction from the picture.

}

\parbf{\ref{ex:R-hom}.}
Let $O$, $E$, $A$, and $B$ be the points with coordinates $(0,0)$, $(1,0)$, $(a,0)$, and $(b,0)$ respectively.

To construct a point $W$ with coordinates $(a+b,0)$ (or $(a-b,0)$) try to construct two parallelograms $OAPQ$ and $BWPQ$ (or $WBPQ$).

To construct $Z$ with coordinates $(a\cdot b,0)$
choose a line $(OE')\ne (OE)$
and try to construct points $A'\in (OE')$
and $Z \in(OE)$
so that 
$\triangle OEE'\z\sim \triangle OAA'$ and $\triangle OE'B\z\sim \triangle OA'Z$.
Likewise, construct the point $(\tfrac ab,0)$.

\parbf{\ref{ex:center-circ-affine}.}
Draw two parallel chords $[XX']$ and~$[YY']$.
Let $Z\z\in(XY)\z\cap (X'Y')$ and $Z'\z\in (XY')\cap (X'Y)$.
Note that $(ZZ')$ passes thru the center.

Repeat the same construction for another pair of parallel chords.
The center lies at the intersection of the obtained lines.

\parbf{\ref{ex:affine-perp}.}
Assume a construction produces two perpendicular lines.
Apply a shear map that changes the angle between the lines (see \ref{ex:afine-linear}\textit{\ref{ex:afine-linear:shear}}).

Note that it transforms the construction into the same construction for other free choices of points.
Therefore, our construction does not generally produce perpendicular lines.
(It might produce a perpendicular line only by coincidence.)
 
\parbf{\ref{ex:parallelogram-rule}.}
The first part follows from \ref{ex:affine-par}.

Suppose $A$, $B$, $X$, and $Y$ are not collinear;
in this case, $\square ABYX$ is a parallelogram.
By the parallelogram rule, the only-if part follows.

Now suppose $A$, $B$, $X$, and $Y$ lie on line $\ell$.
Choose two additional points $P,Q\notin\ell$ such that 
\[\overrightarrow{XY}=\overrightarrow{PQ}
\quad\text{and therefore}\quad 
\overrightarrow{PQ}=\overrightarrow{AB}.\]
From above we get 
\[\overrightarrow{X'Y'}=\overrightarrow{P'Q'}
\quad\text{and}\quad 
\overrightarrow{P'Q'}=\overrightarrow{A'B'}.\]
Hence the only-if part follows in the general case.

The if part follows since the inverse of an affine transformation is also affine.
 
\parbf{\ref{ex:affine-continuous} and \ref{ex:affine-coordinates}.}
Choose a coordinate system and apply the fundamental theorem of affine geometry (\ref{thm:fundamental-theorem-of-affine-geometry}) for the points $O=(0,0)$, $X=(1,0)$, and $Y=(0,1)$.

\parbf{\ref{ex:preserved-circle}.} 
Use \ref{ex:center-circ-affine}, to show that $O'=O$.
Choose a coordinate system with origin at $O$.
By \ref{ex:affine-coordinates}, 
\[\alpha\:(x,y)\mapsto(a\cdot x+b\cdot y,\ c\cdot x+d\cdot y).\]
Use that $\alpha(\Gamma)=\Gamma$ to show at $(\begin{smallmatrix}a&b\\c&d\end{smallmatrix})\cdot (\begin{smallmatrix}a&c\\b&d\end{smallmatrix})\z=(\begin{smallmatrix}1&0\\0&1\end{smallmatrix})$.
Conclude that $\alpha$ is a motion.

\parbf{\ref{ex:inversions-inversive}.}
Observe that any reflection meets the condition.
By Theorem~\ref{thm:inverse-cline}, any inversion meets the condition as well.
Therefore, the same holds for any composition of inversions and reflections.

To prove the converse, choose a bijection $\alpha$ that maps circlines to circlines.
Show that we can compose $\alpha$ with several inverses and reflections to obtain a bijection $\alpha'$ such that $\alpha'(\infty)=\infty$ and $\alpha(\Gamma)=\Gamma$ for some circle $\Gamma$.

By \ref{ex:preserved-circle}, $\alpha'$ is a motion.
By \ref{ex:3-reflections}, $\alpha'$ is a composition of reflections.
Therefore so is $\alpha$.

\parbf{\ref{ex:f(1)=1}.}
Set $a=f(1)$ and $b=f(0)$.
Show and use that $a=a^2$ and $b=b^2$.

\parbf{\ref{ex:ceva-affine}.} Apply Menelaus's theorem
for $\triangle AA'B$ with $(CC')$
and
for $\triangle AA'C$ with $(BB')$.

\refstepcounter{chapter}
\setcounter{eqtn}{0}

\parbf{\ref{ex:proj-cross-ratio}.}
To prove \textit{(a)}, apply \ref{prop:affine-linear}.

To prove \textit{(b)}, suppose $P_i=(x_i,y_i)$;
show and use that 
\[\frac{P_1P_2\cdot P_3P_4}{P_2P_3\cdot P_4P_1}
=\left|\frac{(x_1-x_2)\cdot(x_3-x_4)}{(x_2-x_3)\cdot (x_4-x_1)}\right|\]
if all $P_i$ lie on a horizontal line $y=b$, and
\[\frac{P_1P_2\cdot P_3P_4}{P_2P_3\cdot P_4P_1}
=\left|\frac{(y_1-y_2)\cdot(y_3-y_4)}{(y_2-y_3)\cdot (y_4-y_1)}\right|\]
otherwise. (See \ref{ex:cross-ratio-area} for another proof.)

To prove \textit{(c)}, apply \textit{(a)}, \textit{(b)}, and \ref{thm:moving}.

\parbf{\ref{ex:proj-cross-ratio=1}.}
Observe that for the perspective projection from $(AB)$ to $(AC)$ with center at $W$ we have
$A\mapsto A$, $B\mapsto C$, $X\mapsto P$, and $Y\mapsto Q$.
Therefore \ref{ex:proj-cross-ratio} implies \textit{(\ref{ex:proj-cross-ratio=1:=})}.

Applying the same argument with center $V$, we get  
$\frac{AY\cdot BX}{AX\cdot BY}=\frac{AP\cdot CQ}{AQ\cdot CP}$.
Therefore 
$\frac{AX\cdot BY}{AY\cdot BX}=\frac{AY\cdot BX}{AX\cdot BY}$,
hence \textit{(\ref{ex:proj-cross-ratio=1:1})}.

\parbf{\ref{ex:pappus}.}
Assume that $(AB)$ meets $(A'B')$ at~$O$.
Since $(AB')\parallel (BA')$, we get that $\triangle OAB'\z\sim\triangle OBA'$
and
$\frac{OA}{OB}=\frac{OB'}{OA'}$.

Similarly, $(AC')\parallel (CA')\ \Longrightarrow\ \frac{OA}{OC}=\frac{OC'}{OA'}$.

Therefore
$\frac{OB}{OC}=\frac{OC'}{OB'}$.
Applying the SAS similarity condition, we get that
$\triangle OBC'\z\sim\triangle OCB'$.
Consequently, $(BC')\parallel (CB')$.

The case $(AB)\parallel(A'B')$ is similar.

\parbf{\ref{ex:pappus-converse}.} Observe that the statement is equivalent to Pappus' theorem.

\parbf{\ref{ex:desargues-construction};} \textit{(\ref{ex:desargues-construction:desargues}).}
Assume that the parallelogram is formed by the two pairs of parallel lines $(AB)\z\parallel (A'B')$ and $(BC)\parallel(B'C')$ and $\ell=(AC)$ in the notation of Desargues' theorem (\ref{thm:desargues}).

\parit{(\ref{ex:desargues-construction:pappus}).} Suppose that the parallelogram is formed by the two pairs of parallel lines $(AB')\z\parallel (A'B)$ and $(BC')\parallel(B'C)$ and $\ell=(AC')$ in the notation of Pappus' theorem (\ref{thm:pappus}).

\parbf{\ref{ex:dual-configurations};}
\textit{(\ref{ex:dual-configurations:infty}).} Observe and use that $\square ABCD$ and $\square ABDX$ are parallelograms.

\parit{(\ref{ex:dual-configurations:dual}).}
Draw $a=(KN)$, $b=(KL)$, $c=(LM)$, $d=(MN)$, mark $P\in b\cap d$, and continue.

\parbf{\ref{ex:dual-euclid}.}
Assume there is a duality.
Choose two distinct parallel lines $\ell$ and~$m$.
Let $L$ and $M$ be their dual points.
Set $s=(ML)$, then its dual point $S$ has to lie on both $\ell$ and $m$ --- a contradiction.

\parbf{\ref{ex:dula-coordinates}.}
Assume $M=(a,b)$ 
and the line $s$ is given by the equation $p\cdot x+q\cdot y=1$.
Then $M\in s$ is equivalent to $p\cdot a+q\cdot b=1$.

This, in turn, is equivalent to $m\ni S$
where $m$ is the line given by the equation 
$a\cdot x+b\cdot y=1$ and $S=(p,q)$.

To extend this bijection to the whole projective plane, assume that 
(1) the ideal line corresponds to the origin 
and (2) the ideal point is given by the pencil of lines $b\cdot x-a\cdot y=c$ for different values of $c$ corresponds to the line defined by $a\cdot x+b\cdot y=0$.

\parbf{\ref{ex:dual-pappus}.}
Assume one set of concurrent lines $a$, $b$, $c$, 
and another set of concurrent lines $a'$, $b'$, $c'$ are given.
Let
\begin{align*}
P&\in b\cap c',
&
Q&\in c\cap a',
&
R&\in a\cap b',\\
P'&\in b'\cap c,
&
Q'&\in c'\cap a,
&
R'&\in a'\cap b.
\end{align*}
Then the lines $(PP')$, $(QQ')$, and $(RR')$ are concurrent.

\begin{Figure}
\vskip-0mm
\centering
\includegraphics{mppics/pic-356}
\vskip-0mm
\end{Figure}

(The obtained configuration of nine points and nine lines is the same as in the original theorem and the obtained result is its reformulation.)

\parbf{\ref{ex:dual-desargues-construction},} \textit{(\ref{ex:dual-desargues-construction:desargues})}.
Assume $(AA')$ and $(BB')$ are the given lines and $C$ is the given point.
Apply the dual Desargues' theorem (\ref{thm:dual-desargues}) to construct $C'$ so that $(AA')$, $(BB')$, and $(CC')$ are concurrent. 
Since $(AA')\z\parallel (BB')$, 
we get that 
$(AA')\z\parallel (BB')\z\parallel (CC')$.

\parit{\textit{(\ref{ex:dual-desargues-construction:pappus})}.} Assume that $P$ is the given point and $(R'Q)$, $(P'R)$ are the given parallel lines.
Try to construct point $Q'$ as in the dual Pappus' theorem (see the solution of \ref{ex:dual-pappus}).

\parbf{\ref{ex:revert}.} Suppose $p=(QR)$; denote by $q$ and $r$ the dual lines produced by the construction.
Then, by \ref{clm:polar}, $P$ is the point of intersection of $q$ and $r$.

\parbf{\ref{ex:tangent ruler}.}
The line $v$ polar to $V$ is tangent to~$\Gamma$.
Since $V\in p$, by \ref{clm:polar}, we get that $P\in v$;
that is, $(PV)=v$.
Hence the statement follows.

\parbf{\ref{ex:concentric-circ}.}
Choose a point $P$ outside of the larger circle.
Construct the lines dual to $P$ for both circles.
Note that these two lines are parallel. 

Assume that the lines intersect the bigger circle at two pairs of points $X$, $X'$ and $Y$, $Y'$.
Let $Z\z\in (XY)\cap (X'Y')$.
Note that the line $(PZ)$ passes thru the common center.

\begin{Figure}
\vskip-0mm
\centering
\includegraphics{mppics/pic-370}
\vskip-6mm
\end{Figure}

The center is the intersection of $(PZ)$ and another line constructed in the same way.

\parbf{\ref{ex:proj-perp}.} 
Construct polar lines for two points on~$\ell$.
Denote their intersection by $L$.
Note that $\ell$ is polar to $L$ and therefore $(OL)\perp \ell$.

\refstepcounter{chapter}
\setcounter{eqtn}{0}

\parbf{\ref{ex:defect-sphere}.}
Apply \ref{lem:area-spher-triangle}.

\parbf{\ref{ex:s-medians}.}
\textit{(\ref{ex:s-medians:a})}.
Observe and use that 
$OA'\z=OB'\z=OC'$.

\parit{(\ref{ex:s-medians:b}).} Notice that the medians of spherical triangle $ABC$ 
map to the medians of Euclidean a triangle $A'B'C'$.
It remains to apply \ref{thm:centroid} for $\triangle A'B'C'$.

\parbf{\ref{ex:s-altitudes}.}
Apply the reflection in the plane thru $O$, $N$, and $P$.

\parbf{\ref{ex:stereographic-inversion}.}
Apply \ref{thm:inverse}\textit{\ref{thm:inverse:line}}
and \ref{thm:inversion-3d}\textit{\ref{thm:inversion-3d:angle}}.

\parbf{\ref{ex:great-circ}.}
Apply \ref{thm:inversion-3d}\textit{\ref{thm:inversion-3d:b}}.

\parbf{\ref{ex:conform-sphere}.}
Set $z=P'Q'$.
Note that $\tfrac yz\to 1$ as $x\to 0$.

It remains to show that 
$$\lim_{x\to 0} \frac{z}{x}=\frac{2}{1+OP^2}.$$

Recall that the stereographic projection is the inversion across the sphere $\Upsilon$ with the center at the south pole $S$ restricted to the plane $\Pi$.
Show that there is a plane $\Lambda$ passing thru $S$, $P$, $Q$, $P'$, and~$Q'$.
In the plane $\Lambda$, the map $Q\mapsto Q'$ is an inversion across the circle $\Upsilon\cap \Lambda$.

This reduces the problem to Euclidean plane geometry.
The remaining calculations in $\Lambda$ are similar to those in the proof of \ref{lem:conformal}.

\parbf{\ref{ex:cone}.}
Consider the inversion of the base circle across a sphere with its center at the tip of the cone and apply \ref{thm:inversion-3d}.

\parbf{\ref{ex:2(pi/4)=pi/3}.} 
Apply the spherical Pythagorean theorem to show that
$
\cos AB_s\z=\tfrac12
$,
and conclude that $AB_s\z=\tfrac\pi3$.

{

\begin{wrapfigure}{r}{28mm}
\vskip-0mm
\centering
\includegraphics{mppics/pic-358}
\end{wrapfigure}

Alternatively, 
look at the tessellation of the half-sphere in the picture;
it is made from 12 copies of $\triangle_s A B C$ and 4 equilateral spherical triangles.
Due to the symmetry of this tessellation, it follows that $[AB]_s$ occupies $\tfrac16$ of the equator, meaning $AB_s=\tfrac\pi3$.

}

\parbf{\ref{ex:taurinus}.}
Use Euler's formula (\ref{sec:Euler's formula}) to show that $\sin(i\cdot x)=i\cdot\sinh x$ and apply it.
You should get the following answers:
\[\cosh c=\cosh a \cdot \cosh b-\sinh a\cdot \sinh b\cdot \cos\gamma,\]
\[\cos \gamma=-\cos \alpha \cdot \cos \beta+\sin \alpha\cdot \sin \beta \cdot \cosh c,\]
\[\frac{\sin \alpha}{\sinh a}=\frac{\sin \beta}{\sinh b}=\frac{\sin \gamma}{\sinh c}.\]

\refstepcounter{chapter}
\setcounter{eqtn}{0}



\parbf{\ref{ex:hex}.}
Consider the bijection $P\z\leftrightarrow \hat P$ of the h-plane with the absolute~$\Omega$.
Note that $\hat P\z\in [A_iB_i]$ if and only if $P\in\Gamma_i$.

\parbf{\ref{ex:P<->hatP}.} Apply \ref{lem:P-hat-chord} and the definition of h-distance in \ref{sec:conformal-model}.

\parbf{\ref{ex:h-median}.} 
The observation follows since the reflection across the perpendicular bisector of $[PQ]$ is a motion of the Euclidean plane and a motion of the h-plane as well.

We can assume that the center of the circumcircle coincides with the center of the absolute.
In this case, the h-medians of the triangle coincide with the Euclidean medians.
It remains to apply \ref{thm:centroid}.

\begin{Figure}
\centering
\vskip-0mm
\includegraphics{mppics/pic-360}
\end{Figure}

\parbf{\ref{ex:h-altitudes}.} 
Use the projective model.
Assume that two h-altitudes intersect at a point~$H$.
Move $H$ to the center of the absolute.
By \ref{obs:h-p-perp} these h-altitudes become Euclidean altitudes of the triangle.
By \ref{thm:orthocenter}, the remaining Euclidean altitude passes thru $H$.
By \ref{obs:h-p-perp} this Euclidean altitude is also an h-altitude.

\parbf{\ref{ex:klein-perp}.} 
Let $\hat\ell$ and $\hat m$ denote the h-lines in the conformal model that correspond to $\ell$ and $m$.
We need to show that $\hat\ell\perp\hat m$ as arcs in the Euclidean plane.

The point $Z$, where $s$ meets $t$, is the center of the circle $\Gamma$ containing~$\hat\ell$.

If $\bar m$ is passing thru $Z$, then the inversion across $\Gamma$ exchanges the ideal points of~$\hat\ell$.
Consequently, $\hat\ell$ maps to itself. 
Hence the result.

\parbf{\ref{ex:klein-for-angle-parallelism}.}
Let $Q$ be the footpoint of $P$ on the line and $\phi$ be the angle of parallelism. 

{

\begin{wrapfigure}{r}{28mm}
\vskip-2mm
\centering
\includegraphics{mppics/pic-362}
\end{wrapfigure}

We can assume that $P$ is the center of the absolute.
Therefore $PQ=\cos\phi$ and
\[PQ_h=\tfrac12\cdot\ln\frac{1+\cos\phi}{1-\cos\phi}.\]

\parbf{\ref{ex:klein-inradius}.} 
Apply \ref{ex:klein-for-angle-parallelism} for $\phi=\tfrac\pi3$.

}

\parbf{\ref{ex:pyth-h-proj}.}
Note that
$
b=\tfrac12\cdot\ln\frac{1+t}{1-t}$;
therefore
\[
\cosh b
=
\tfrac12\cdot\left(\sqrt{\tfrac{1+t}{1-t}}+\sqrt{\tfrac{1-t}{1+t}}\right)
=
\frac1{\sqrt{1-t^2}}.
\eqlbl{cosh-b}
\]
The same way, we get
\[\begin{aligned}\cosh c&=\frac1{\sqrt{1-u^2}}.
\end{aligned}
\eqlbl{cosh-c}
\]

Let $X$ and $Y$ be the ideal points of~$(BC)_h$.
Applying the Pythagorean theorem (\ref{thm:pyth}) again,
we get that
$CX=CY=\sqrt{1-t^2}$.
Therefore, 
\[
a
=
\tfrac12\cdot\ln\frac{\sqrt{1-t^2}+s}{\sqrt{1-t^2}-s},\]
and
\[
\begin{aligned}
\cosh a&=\tfrac12\cdot
\sqrt{\frac{\sqrt{1-t^2}+s}{\sqrt{1-t^2}-s}}+
\\
&+
\tfrac12\cdot\sqrt{\frac{\sqrt{1-t^2}-s}{\sqrt{1-t^2}+s}}=
\\
&=\frac{\sqrt{1-t^2}}{\sqrt{1-t^2-s^2}}=
\\
&=\frac{\sqrt{1-t^2}}{\sqrt{1-u^2}}.
\end{aligned}
\eqlbl{cosh-a}
\]

Finally, note that \ref{cosh-b}, \ref{cosh-c}, and \ref{cosh-a} imply the theorem.

\parbf{\ref{ex:Boyai-in-Euclid}.}
In the Euclidean plane, the circle $\Gamma_2$ is tangent to $k$; 
that is, $T$ is the only point at intersection of $\Gamma_2$ and $k$.
It defines a unique line $(PT)$ parallel to~$\ell$.

\parbf{\ref{ex:common-perp}.}
Choose two points $P$ and $Q$ on $\ell$.
Drop perpendiculars $p$ and $q$ from $P$ and $Q$ to~$m$
(follow \ref{ex:construction-perpendicular}).
Drop perpendiculars from $P$ to $q$ and from $Q$ to~$p$;
label their point of intersection by~$H$.
(Assume $H$ exists, if not choose other $P$ and $Q$.) 
Drop perpendicular $n$ from $H$ to $m$.

To justify the construction, we need to show that $n\perp \ell$.
To do this, move $H$ to the center of the absolute and apply \ref{obs:h-p-perp}, \ref{ex:klein-perp}, and \ref{thm:orthocenter}.

\refstepcounter{chapter}
\setcounter{eqtn}{0}

\parbf{\ref{ex:|zw|}.} Use that $|z|^2=z\cdot \bar z$ and $\bar z\cdot \bar u=\overline{z\cdot u}$.

\parbf{\ref{ex:ptolemy}.}
Choose a quadrangle $ABCD$.
Assume that $0$, $u$, $v$, and $w$ are complex coordinates of $A$, $B$, $C$, and $D$ respectively.
Rewrite Ptolemy's inequality using $u$, $v$, and $w$.
Deduce this inequality from the provided identity and the triangle inequality.

\parbf{\ref{ex:3-sum-C}.} 
Let $z$, $v$, and $w$ denote the complex coordinates of $Z$, $V$, and $W$ respectively.
Then 
\begin{align*}
&\qquad \measuredangle ZVW+\measuredangle VWZ+\measuredangle WZV\equiv
\\
&\equiv
\arg \tfrac{w-v}{z-v}+\arg \tfrac{z-w}{v-w}+\arg \tfrac{v-z}{w-z}\equiv
\\
&\equiv
\arg \tfrac{(w-v)\cdot(z-w)\cdot(v-z)}{(z-v)\cdot(v-w)\cdot(w-z)}\equiv
\\
&\equiv\arg (-1)\equiv
\pi.
\end{align*}

\parbf{\ref{ex:C-sim}.}
Note and use that 
$
\measuredangle EOV=\measuredangle WOZ\z=\arg v$
and
$\frac{OW}{OZ}=\frac{OZ}{OW}=|v|$.

\parbf{\ref{ex:real-cross-ratio}.}
Note that 
\begin{align*}
&\qquad\arg\frac{(v-u)\cdot(z-w)}{(v-w)\cdot(z-u)}\equiv
\\
&\equiv
\arg\frac{v-u}{z-u}
+
\arg\frac{z-w}{v-w}=
\\
&= \measuredangle ZUV+\measuredangle VWZ.
\end{align*}

The statement follows since the value $\tfrac{(v-u)\cdot(z-w)}{(v-w)\cdot(z-u)}$ is real if and only if 
\[2\cdot\arg\frac{(v-u)\cdot(z-w)}{(v-w)\cdot(z-u)}\equiv0.\]

\parbf{\ref{ex:3-squares}.}
We can choose the complex coordinates so that the points $O$, $E$, $A$, $B$, and $C$ have coordinates
$0$, $1$, $1+i$, $2+i$, and $3+i$ respectively.
Set $\measuredangle EOA=\alpha$, $\measuredangle EOB=\beta$, and $\measuredangle EOC\z=\gamma$.
Note that
\begin{align*}
&\ \ \ \ \alpha+\beta+\gamma\equiv
\\
&\equiv\arg(1+i)+\arg(2+i)+\arg(3+i)\equiv
\\
&\equiv\arg[(1+i)\cdot(2+i)\cdot(3+i)]\equiv
\\
&\equiv\arg [10\cdot i]=
\\
&=\tfrac\pi2.
\end{align*}
Note that these three angles are acute and conclude that $\alpha+\beta+\gamma=\tfrac\pi2$.

\parbf{\ref{ex:6-circles}.}
The identity can be verified by straightforward computations.

By \ref{thm:inscribed-quadrangle-C}, five out of the six cross-ratios in this identity are real.
Consequently, the sixth cross-ratio must also be real.
It remains to apply the theorem again.

\parbf{\ref{ex:4-sim}.}
Use \ref{thm:signs-of-triug} and \ref{cor:half-plane} to show that $\angle UAB$, $\angle BVA$, and $\angle ABW$ have the same sign.
By SAS, we have
\[\frac{AU}{AB}=\frac{VB}{VA}=\frac{BA}{BW},\]
and
\[\measuredangle UAB=\measuredangle BVA=\measuredangle ABW.\]
The latter means that 
\[|\frac{u-a}{b-a}|=|\frac{b-v}{a-v}|=|\frac{a-b}{w-b}|,\]
and
\[\arg\frac{b-a}{u-a}=\arg\frac{a-v}{b-v}=\arg\frac{a-b}{w-b}.\]
It implies the first two equalities in 
\[\frac{b-a}{u-a}=\frac{a-v}{b-v}=\frac{w-b}{a-b}=\frac{w-v}{u-v};\eqlbl{eq:4fractions}\]
the last equality holds since 
\[\frac{(b-a)+(a-v)+(w-b)}{(u-a)+(b-v)+(a-b)}=\frac{w-v}{u-v}.\]

To prove \textit{(b)}, rewrite \ref{eq:4fractions} using angles and distances between the points and apply SAS.

\parbf{\ref{ex:inverse-Mob}.}
Show that the inverse of each elementary transformation is elementary
and use \ref{prop:mob-comp}.

\parbf{\ref{ex:3-point-Mob}.}
The fractional linear transformation
\[f(z)=\frac{(z_1-z_\infty)\cdot(z-z_0)}{(z_1-z_0)\cdot(z-z_\infty)}\]
meets the conditions.

To show the uniqueness, assume there is another fractional linear transformation
$g(z)$ that meets the conditions.
Then the composition
$h=g\circ f^{-1}$ 
is a fractional linear transformation; set
$h(z)=\tfrac{a\cdot z+b}{c\cdot z+d}$.

Note that $h(\infty)=\infty$;
therefore, $c=0$.
Furthermore, $h(0)=0$ implies $b=0$.
Finally, since $h(1)=1$, we get that $\tfrac ad=1$.
Therefore, $h$ is the \index{identity map}\emph{identity};
that is, $h(z)=z$ for any~$z$.
It follows that $g=f$.

\parbf{\ref{ex:inversion-Mob}.}
Let $Z'$ be the inverse of the point $Z$.
Assume that the circle of the inversion has center $W$ and radius~$r$.
Let $z$, $z'$, and $w$ denote the complex coordinates of the points $Z$, $Z'$, and $W$ respectively.

By the definition of an inversion, $\arg (z\z-w)\z=\arg (z'-w)$ and
$|z-w|\cdot|z'-w|=r^2$.
It follows that $(\bar z'-\bar w)\cdot ( z- w)= r^2$.
Equivalently,
\[z'=\overline{\left(\frac{\bar w\cdot z+[r^2-|w|^2]}{z- w}\right)}.\]

\parbf{\ref{ex:C-cross-ratio}.}
Check the statement for each elementary transformation.
Then apply \ref{prop:mob-comp}.

\parbf{\ref{ex:schwarz-moebius}.}
Note that $f\:z\mapsto\tfrac{a\cdot z+b}{c\cdot z+d}$ preserves the unit circle $|z|=1$.
Use \ref{cor:invese-comp} and \ref{prop:mob-comp} to show that $f$ commutes with the inversion $z\mapsto 1/\bar z$.
In other words, $1/\overline{f(z)}=f(1/\bar z)$ or
\[\frac{\bar c\cdot \bar z+\bar d}{\bar a\cdot \bar z+\bar b}
=\frac{a/\bar z+b}{c/\bar z+d}\]
for any $z\in\hat{\mathbb{C}}$.
The latter identity leads to the required statement. 
The condition $|w|<|v|$ follows since $f(0)\in\mathbb{D}$.

\parbf{\ref{ex:schwarz-tanh}.} 
Note that the inverses of the points $z$ and $w$ have complex coordinates $1/\bar z$ and $1/\bar w$.
Apply \ref{ex:cosh} and simplify.

The second part follows since the function $x\mapsto \tanh(\tfrac12\cdot x)$ is increasing.

\parbf{\ref{ex:schwarz}.}
Apply Schwarz--Pick theorem for a function $f$ such that $f(0)\z=0$ and then apply \ref{lem:O-h-dist}.

\refstepcounter{chapter}
\setcounter{eqtn}{0}

%
%
%
%

\parbf{\ref{ex:a2/b}.}
To construct  $\sqrt{a\cdot b}$:
(1) construct points $A$, $B$, and $D\z\in [AB]$
such that $AD=a$ and $BD=b$;
(2) construct the circle $\Gamma$ on the diameter $[AB]$;
(3) draw the line $\ell$ thru $D$ perpendicular to $(AB)$; 
(4) let $C$ be an intersection of $\Gamma$ and~$\ell$.
Then $DC= \sqrt{a\cdot b}$.

The construction of $\tfrac{a^2}b$ is analogous.

\parbf{\ref{ex:5-gon},} \textit{(\ref{ex:5-gon:a})}.
Look at the picture;
show that the angles marked the same way have the same angle measure.

Conclude that $XC=BC$ and $\triangle ABC\z\sim \triangle AXB$.
Therefore 
\[\frac{AB}{AC}=\frac{AX}{AB}=\frac{AC-AB}{AB}=\frac{AC}{AB}-1.\]
It remains to solve for $\frac{AC}{AB}$.

{

\begin{wrapfigure}[7]{r}{27mm}
\vskip-8mm
\centering
\includegraphics{mppics/pic-366}
\end{wrapfigure}

\parit{(\ref{ex:5-gon:b}).}
Choose two points $P$ and $Q$ and use the ruler-and-compass calculator to construct two points $V$ and $W$ such that $VW\z=\tfrac{1+\sqrt5}2\cdot PQ$.
Then construct a pentagon with the sides $PQ$ and diagonals $VW$.

}

\parbf{\ref{ex:trisect-set-square}.} 
Note that with a set-square we can construct a line parallel to a given line thru the given point.
It remains to modify the construction in \ref{ex:midpoint-affine}.

\parbf{\ref{ex:equilateral triangle}.}
Choose a coordinate system so that the given vertices are $(0,0)$ and $(1,0)$.
Show that the remaining vertex is $(\tfrac12,\pm\tfrac{\sqrt{3}}2)$.
Observe that it is an irrational point; apply \ref{thm:set-square-constructible-numbers}.%
\footnote{It is okay to use that $\sqrt{3}$ is irrational without proving it.
But let us explain why it is true. 
Assume the contrary; that is, $\tfrac mn=\sqrt{3}$ for integers $m$ and $n$.
We can assume that $m$ and $n$ do not share a prime factor; in particular, if $m$ is divisible by $3$, then $n$ is not.
Observe that $m^2=3\cdot n^2$.
It follows that $m$ is divisible by 3; that is, $m=3\cdot k$ for an integer $k$.
It follows that $3\cdot k^2=n^2$.
Therefore, $n$ is divisible by 3 --- a contradiction.} 

\parbf{\ref{ex:set-square-bisect}.}
Assume that one can construct a bisector of $\angle AOB$, where $A=(1,0)$, $O=(0,0)$, and $C=(1,1)$.
Let $D$ be the point of intersection of the bisector with the line $(AB)$.
Use \ref{lem:bisect-ratio} to show that $D$ is irrational.
Apply \ref{thm:set-square-constructible-numbers} and arrive at a contradiction.

\parbf{\ref{ex:90-60-30}.}
Suppose that every initial point has coordinates $(a,b\cdot\sqrt{3})$ for rational values $a$ and $b$.
Show and use that any point that can be constructed with the 30°-set-square has coordinates of the same type.

\begin{wrapfigure}[8]{r}{21mm}
\vskip-0mm
\centering
\includegraphics{mppics/pic-368}
\end{wrapfigure}

\parbf{\ref{ex:equilateral triangle-verify},} \textit{(\ref{ex:verify:triangle})}.
Observe that three perpendiculars in the picture meet at one point if and only if the triangle is isosceles.

Use this observation a couple of times to verify that the given triangle is equilateral.

\parit{(\ref{ex:verify:bisector}).}
Suppose that a line $\ell$ passes thru the vertex of the given angle.
Choose a point $P\in \ell$.
Suppose $X$ and $Y$ are the footpoints of $P$ on the sides of the angle.
Show and use that $(XY)\perp \ell$ if and only if $\ell$ bisects the angle.

\parbf{\ref{ex:midpoint-proj}.}
Consider the perspective projection 
$(x,y)\mapsto (\tfrac 1x,\tfrac yx)$ (see \ref{sec:perspective-projection}).
Let $A\z=(1,1)$, $B\z=(3,1)$, and $M\z=(2,1)$.
Note that $M$ is the midpoint of $[AB]$.

Their images are $A'\z=(1,1)$, $B'\z=(\tfrac13,\tfrac13)$, and $M'\z=(\tfrac12,\tfrac12)$.
Clearly, $M'$ is not the midpoint of~$[A'B']$.

\parbf{\ref{ex:comparison}.}
$(a)$ is strictly stronger than $(b)$,
$(b)$ is strictly stronger than $(c)$,
$(a)$ is strictly stronger than $(d)$,
and $(d)$ is not comparable with $(b)$ and $(c)$.
Most of these statements follow from \ref{ex:construction-perpendicular},
\ref{prop:perp-perp},
\ref{ex:consturuction-of-inversion},
\ref{ex:affine-perp},
\ref{ex:equilateral triangle}, 
\ref{prob:center-inversor+circumtool}.

To show that $(d)$ is not stronger than $(c)$, show that one cannot construct a midpoint using the set $(d)$ and use the solution of \ref{ex:midpoint-affine}.
To show that $(b)$ is not stronger than $(d)$, show that given the initial configuration of 6 points 
$(0,0)$, 
$(1,0)$,
$(2,0)$,
$(0,1)$, 
$(1,1)$,
$(2,1)$,
one can construct an equilateral triangle using the set $(d)$ and apply \ref{ex:equilateral triangle}.

\refstepcounter{chapter}
\setcounter{eqtn}{0} 

\parbf{\ref{ex:triangle-convex}.}
Assume the contrary; 
that is, there is a point $W\in [XY]$ such that $W\notin\solidtriangle ABC$.

Without  loss of generality, we may assume that $W$ and $A$ lie on opposite sides of the line~$(BC)$.

It implies that both segments $[WX]$ and $[WY]$ intersect $(BC)$.
By Axiom~\ref{def:birkhoff-axioms:1}, $W\in (BC)$ --- a contradiction.


\parbf{\ref{ex:vertex}.} 
To prove the ``only if'' part, consider the line passing thru the vertex that is parallel to the opposite side.

To prove the ``if'' part, use Pasch's theorem (\ref{thm:pasch}).

\parbf{\ref{ex:solid-square}.}
Assume the contrary; that is, suppose a solid square \( \mathcal{K} \) can be represented as the union of a finite collection of segments \( [A_1B_1],\dots,[A_nB_n] \)
and single-point sets \( \{C_1\},\dots,\{C_k\} \).

Note that \( \mathcal{K} \) contains an infinite number of mutually nonparallel segments.
Therefore, we can choose a segment \( [XY] \) in \( \mathcal{K} \) that is not parallel to any of the segments \( [A_1B_1],\dots,[A_nB_n] \).

It follows that \( [XY] \) has at most one common point with each of the sets \( [A_iB_i] \) and~\( \{C_i\} \).
Since \( [XY] \) contains infinitely many points, we reach a contradiction.

\smallskip

Alternatively, choose a circle \( \Gamma \) inside \( \mathcal{K} \).
Note that \( \Gamma \) contains infinitely many points, but by \ref{lem:line-circle}, it has only finitely many intersection points with any degenerate polygonal set.

Let us also note that the statement will follow from \ref{cor:degenerate} and \ref{thm:area-rect}.

\parbf{\ref{ex:poly-circ}.} 
Show that among elementary sets
only one-point sets can be subsets of a circle.
It remains to note that any circle contains an infinite number of points.

\parbf{\ref{ex:two-parallelograms}.}
Suppose that $E$ denotes the point of intersection of the lines $(BC)$ and~$(C'D')$.

\begin{wrapfigure}[7]{r}{28mm}
\vskip-2mm
\centering
\includegraphics{mppics/pic-372}
\end{wrapfigure}

Use \ref{prop:area-parallelogram} to prove that the following solid parallelograms have the same area:
$\solidsquare ABCD$, $\solidsquare AB'ED$, and $\solidsquare AB'C'D'$.

\parbf{\ref{ex:three-trig}.}
Without loss of generality, we may assume that the angles $ABC$ and $BCA$ are acute.

Let $A'$ and $B'$ denote the footpoints of $A$ and $B$ on $(BC)$ and $(AC)$ respectively.
Note that $h_A=AA'$ and $h_B=BB'$.

Note that $\triangle AA'C\sim \triangle BB'C$;
indeed the angle at $C$ is shared and the angles at $A'$ and $B'$ are right.
In particular
$\frac{AA'}{BB'}=\frac{AC}{BC}$,
or equivalently, $h_A\cdot BC=h_B\cdot AC$.

\begin{Figure}
\vskip-0mm
\begin{minipage}{.49\textwidth}
\centering
\includegraphics{mppics/pic-374}
\end{minipage}
\hfill
\begin{minipage}{.49\textwidth}
\centering
\includegraphics{mppics/pic-376}
\end{minipage}
\end{Figure}

\parbf{\ref{ex:half-parallelogram}.}
Draw the line $\ell$ 
thru $M$ parallel to $[AB]$ and $[CD]$;
it subdivides $\solidsquare ABCD$ into two solid parallelograms
which we will denote as
$\solidsquare ABEF$ and
$\solidsquare CDFE$.
In particular,
\begin{align*}
&\area(\solidsquare ABCD)=
\\
&\qquad=
\area(\solidsquare ABEF)+\area(\solidsquare CDFE).
\end{align*}

By \ref{prop:area-parallelogram} and \ref{thm:area-of-triangle} we get that 
\begin{align*}
\area(\solidtriangle ABM)&=\tfrac12\cdot\area(\solidsquare ABEF),
\\
\area(\solidtriangle CDM)&=\tfrac12\cdot\area(\solidsquare CDFE)
\end{align*}
and hence the result.

\parbf{\ref{ex:area-diag}.}
Let $h_A$ and $h_C$ denote the distances from $A$ and $C$ to the line~$(BD)$ respectively.
According to \ref{thm:area-of-triangle},
\begin{align*}
\area(\solidtriangle ABM)&=\tfrac12\cdot h_A\cdot BM;
\\
\area(\solidtriangle BCM)&=\tfrac12\cdot h_C\cdot BM;
\\
\area(\solidtriangle CDM)&=\tfrac12\cdot h_C\cdot DM;
\\
\area(\solidtriangle ABM)&=\tfrac12\cdot h_A\cdot DM.
\end{align*}
Therefore
\begin{align*}
&\area(\solidtriangle ABM)\cdot \area(\solidtriangle CDM)=
\\
&\qquad=\tfrac14 \cdot h_A\cdot h_C\cdot DM\cdot BM=
\\
&\qquad=\area(\solidtriangle BCM)\cdot \area(\solidtriangle DAM).
\end{align*}

\parbf{\ref{ex:area-inradius}.}
Let $I$ be the incenter of $\triangle ABC$.
Note that $\solidtriangle ABC$
can be subdivided into 
$\solidtriangle IAB$, 
$\solidtriangle IBC$,
and $\solidtriangle ICA$.

It remains to apply \ref{thm:area-of-triangle} 
to each of these triangles and sum up the results.

\parbf{\ref{ex:subdivision}.} Fix a polygonal set $\mathcal{P}$.
Without loss of generality, we may assume that $\mathcal{P}$ is a union of a finite collection of solid triangles.
Cut $\mathcal{P}$ along the extensions of the sides of all the triangles,
it subdivides $\mathcal{P}$ into convex polygons.
Cutting each polygon by diagonals from one vertex produces a subdivision into solid triangles.

\parbf{\ref{ex:pyth-2}.}
Assuming $a>b$,
we have subdivided $\mathcal{K}_c$ into $\mathcal{K}_{a-b}$ and four triangles congruent to~$\mathcal{T}$.
Therefore
\[\area\mathcal{K}_c=\area\mathcal{K}_{a-b}+4\cdot\area\mathcal{T}.
\eqlbl{eq:pyth-2}\]

According to \ref{thm:area-of-triangle},
$\area\mathcal{T}=\tfrac12\cdot a\cdot b$. 
Therefore, the identity \ref{eq:pyth-2} can be written as 
\[c^2=(a-b)^2+2\cdot a\cdot b.\]
Simplifying, we get the Pythagorean theorem.

Case $a=b$ is simpler.
Case $b>a$ can be done in the same way.

\parbf{\ref{ex:sum-3-dist}.} 
If $X$ is a point inside of $\triangle ABC$, then $\solidtriangle ABC$ is subdivided into $\solidtriangle ABX$, $\solidtriangle BCX$, and $\solidtriangle CAX$.
Therefore
\begin{align*}
&\area(\solidtriangle ABX)
+\area(\solidtriangle BCX)+
\\
&\qquad+\area(\solidtriangle CAX)
=\area(\solidtriangle ABC).
\end{align*}

Set $a=AB=BC=CA$.
Let $h_1$, $h_2$, and $h_3$ denote the distances from $X$ to the sides $[AB]$, $[BC]$, and~$[CA]$. 
Then by \ref{thm:area-of-triangle},
\begin{align*}
\area(\solidtriangle ABX)&=\tfrac12\cdot h_1\cdot a,
\\
\area(\solidtriangle BCX)&=\tfrac12\cdot h_2\cdot a,
\\
\area(\solidtriangle CAX)&=\tfrac12\cdot h_3\cdot a.
\end{align*}
Therefore, 
\[h_1+h_2+h_3=\tfrac2a\cdot\area(\solidtriangle ABC).\]

\parbf{\ref{ex:area-medians}.}
Apply \ref{thm:centroid} and \ref{clm:area-ratio} several times.

\parbf{\ref{ex:ceva}.}
Apply \ref{clm:area-ratio} to show that
\[\frac{\area(\solidtriangle ABX)}{\area(\solidtriangle ABB')}=\frac{BX}{BB'}=\frac{\area(\solidtriangle BCX)}{\area(\solidtriangle CBB')}\]
and, therefore,
\[\frac{\area(\solidtriangle ABX)}{\area(\solidtriangle BCX)}=\frac{\area(\solidtriangle ABB')}{\area(\solidtriangle CBB')}=\frac{AB'}{B'C}.\]

It implies the first identity; the rest is analogous.
Multipying the identities, we get the last statement.

\parbf{\ref{ex:cross-ratio-area}.}
To prove \textit{(\ref{ex:cross-ratio-area:a})}, apply \ref{clm:area-ratio} twice to $\solidtriangle OL_iL_j$, $\solidtriangle OL_jM_i$, and $\solidtriangle OM_iM_j$.

To prove part~\textit{(\ref{ex:cross-ratio-area:b})}, use \ref{clm:area-ratio} to rewrite the left-hand side using the areas of $\solidtriangle OL_1L_2$, $\solidtriangle OL_2L_3$, $\solidtriangle OL_3L_4$, and $\solidtriangle OL_4L_1$.
Furthermore, use part \textit{(\ref{ex:cross-ratio-area:a})} to rewrite it using areas of $\solidtriangle OM_1M_2$, $\solidtriangle OM_2M_3$, $\solidtriangle OM_3M_4$, and $\solidtriangle OM_4M_1$ and apply \ref{clm:area-ratio} again to get the right-hand side.

\parbf{\ref{ex:circle-is-quadrable}.}
Let $\mathcal{P}_n$ and $\mathcal{Q}_n$ be the solid regular $n$-gons
so that $\Gamma$ is inscribed in $\mathcal{Q}_n$ and circumscribed around~$\mathcal{P}_n$.
Clearly,
$\mathcal{P}_n\subset\mathcal{D}\subset\mathcal{Q}_n$.

Show that 
$\tfrac{\area\mathcal{P}_n}{\area\mathcal{Q}_n}=(\cos\tfrac\pi n)^2$;
in particular, 
$$\frac{\area\mathcal{P}_n}{\area\mathcal{Q}_n}\to 1
\quad
\text{as}
\quad
n\to\infty.$$
Next show that $\area\mathcal{Q}_n<100$ foe any large $n$.

These two statements imply that
\[(\area\mathcal{Q}_n\z-\area\mathcal{P}_n)\to 0.\]
Hence the result.

\spell{\end{multicols}}{}

\newpage

}
\newpage
\phantomsection
{\scriptsize
\input{all-lectures.ind}
}
\renewcommand{\bibname}{Used resources}
{

\def\emph{\textit}

\printbibliography[heading=bibintoc]

@book{akopyan,
  title={Geometry in figures},
  author={Akopyan, A.},
  year={2017},
  addendum={[Translated to Bulgarian, Chinese, French, Hebrew, Polish, Russian, and Spanish.]}
}

@incollection{engeler,
  title={Remarks on the theory of geometrical constructions},
  author={Engeler, E.},
  booktitle={The syntax and semantics of infinitary languages},
  pages={64--76},
  year={1968}
}

@article {alexandrov,
    AUTHOR = {Aleksandrov, A. D.},
     TITLE = {Minimal foundations of geometry},
   JOURNAL = {Siberian Math. J.},
    VOLUME = {35},
      YEAR = {1994},
    NUMBER = {6},
     PAGES = {1057--1069},
}

@book{bachmann,
  title={Aufbau der Geometrie aus dem Spiegelungsbegriff},
  author={Bachmann, F.},
  year={1959},
  addendum={[Translated to Russian.]}
}

@online{euclidea,
title={Euclidea},
  url={https://www.euclidea.xyz}
}

@article{beltrami,
  title={Teoria fondamentale degli spazii di curvatura costante},
  author={Beltrami, E.},
  journal={Annali. di Mat., ser II},
  volume={2},
  pages={232--255},
  year={1868},
  addendum={[Translated by J. Stillwell in \textit{Sources of Hyperbolic Geometry}, pp. 41--62 (1996).]}
}

@article {birkhoff,
    AUTHOR = {Birkhoff, G. D.},
     TITLE = {A set of postulates for plane geometry, based on scale and
              protractor},
   JOURNAL = {Ann. of Math. (2)},
  FJOURNAL = {Annals of Mathematics. Second Series},
    VOLUME = {33},
      YEAR = {1932},
    NUMBER = {2},
     PAGES = {329--345},
}

@book {bolyai,
    AUTHOR = {Bolyai, J.},
     TITLE = {Appendix},
     YEAR ={1832},
    addendum={[Translated by Ferenc Kárteszi in \textit{Appendix.
The theory of space}, (1987).]}   
}

@book {euclid,
 TITLE = {Euclid's Elements}
 }

@book{byrne,
  title={The first six books of the Elements of Euclid: in which coloured diagrams and symbols are used instead of letters for the greater ease of learners},
  author={Byrne, O.},
  year={1847},
  url={https://github.com/jemmybutton/byrne-euclid/releases}
}

@book{hadamard,
  title={Le{\c{c}}ons de g{\'e}om{\'e}trie {\'e}l{\'e}mentaire: G{\'e}om{\'e}trie plane},
  author={Hadamard, J.},
  year={1906},
addendum={[Translated  by M. Saul in \textit{Lessons in geometry. I.
Plane geometry}, (2008).]}
}

@book{kiselev,
  title={Элементарная геометрия},
  author={Киселёв, А. П.},
  addendum={[Translated by A. Givental in \textit{Kiselev's geometry}, (2006).]}
}

@article{lambert,
  title={Theorie der parallellinien},
  author={Lambert, J. H.},
  journal={Leipziger Magazin f{\"u}r reine und angewandte Mathematik},
  volume={1},
  number={2},
  pages={137--164},
  year={1786}
}

@article{legendre,
  title={El{\'e}ments de g{\'e}om{\'e}trie},
  author={Legendre, A.-M.},
  year={1794}
}

@article{lobachevsky,
author={Лобачевский, Н. И.},
title={О началах геометрии}, 
journal={Казанский вестник},
number={25---28},
year={1829---1830}
}

@incollection{lobachevsky-1840,
  title={Geometrische Untersuchungen zur Theorie der Parallellinien},
  author={Lobatschewsky, N. I.},
  year={1840},
addendum={[Translated  by G.~B. Halsted in \textit{The Theory of Parallels}, (2015).]}
}

@book{prasolov,
  title={Задачи по планиметрии},
  author={Прасолов, В. В.},
  year={1986},
  addendum={[Translated by D. Leites in \textit{Problems in plane and solid geometry}, (2006).]}
}

@book{saccheri,
  title={Euclides ab omni n\ae vo vindicatus},
  author={Saccheri, G. G.},
  year={1733},
addendum={[Translated by G. B. Halsted in \textit{Euclides vindicatus}, (1986).]},
}

@book{sharygin,
  title={Геометрия 7--9},
  author={Шарыгин, И. Ф.},
  year={1997}
}

@inproceedings {tarski,
    AUTHOR = {Tarski, A.},
     TITLE = {What is elementary geometry?},
 BOOKTITLE = {The axiomatic method (edited by {L}. {H}enkin, {P}. {S}uppes and {A}.
              {T}arski)},
     PAGES = {16--29},
      YEAR = {1959},
}

@book{taurinus,
  title={Geometriae: Prima elementa. Recensuit et novas observationes adjecit},
  author={Taurinus, F. A.},
  year={1826},
  publisher={Bachem}
}

@online{zadachi,
title={Задачи},
url={https://www.problems.ru/}
}
\fussy
}
\end{document}